# Representation of Universal Algebra

Aleks Kleyn

Aleks_Kleyn@MailAPS.org
http://AleksKleyn.dyndns-home.com:4080/
http://sites.google.com/site/AleksKleyn/
http://arxiv.org/a/kleyn_a_1
http://AleksKleyn.blogspot.com/

Abstract. Theory of representations of universal algebra is a natural development of the theory of universal algebra. Morphism of the representation is the map that conserve the structure of the representation. Exploring of morphisms of the representation leads to the concepts of generating set and basis of representation. In the book I considered the notion of tower of left-side representations of $\Omega_i$-algebras, $i = 1, ..., n$, as the set of coordinated left-side representations of $\Omega_i$-algebras.

# Contents

CHAPTER 1

# Preface

## 1.1. Preface

In my papers, I often explore problems relating to the representation of universal algebra. Initially it was small sketches which I repeatedly corrected and rewrote. However gradually there were new observations. As a result, auxiliary tool became a consistent theory.

I realized this when I was writing book [9], and I decided to dedicate a separate book to the questions related with representation of universal algebra. Exploring of the theory of representations of universal algebra shows that this theory has a lot of common with theory of universal algebra.

The definition of vector space as representation of field in the Abelian group was the main impetus of deeper exploring of representations of universal algebra. I put attention that this definition changes role of linear map. It was found that linear map is the map that preserves the structure of the representation. It is easy to generalize this structure for an arbitrary representation of universal algebra. Thus I came to the notion of morphism of representation.

The set of regular automorphisms of vector space forms a group. This group is single transitive on the set of basises of vector space. This statement is the foundation of the theory of invariants of vector space.

The natural question arises. Can we generalize this structure to arbitrary representation? The basis is not the only set that forms the vector space. If we add an arbitrary vector to the set of vectors of basis, then a new set also generates the same vector space, however this set is not basis. This statement is initial point where I started exploring of generating set of representation. Generating set of representation is one more interesting parallel between theory of representations and theory of universal algebra.

The set of automorphisms of representations is loop. Nonassociativity of the product is the source of numerous questions which require additional research. All these questions lead to the need to understand the theory of invariants of a given representation.

If we consider the theory of representations of universal algebra as an extension of the theory of universal algebra, then why not consider the representation of one representation in another representation. Thus the concept of the tower representations appeared. The most amazing fact is the statement that all maps in the tower of representations are coordinated.

## 1.2. Conventions

CONVENTION 1.2.1. *Function and map are synonyms. However according to tradition, correspondence between either rings or vector spaces is called map and a map of either real field or quaternion algebra is called function.* □

CONVENTION 1.2.2. *Since the number of universal algebras in the tower of representations is varying, then we use vector notation for a tower of representations. We denote the set* $(A_1, ..., A_n)$ *of* $\Omega_i$*-algebras* $A_i$, $i = 1, ..., n$ *as* $\overline{A}$. *We denote the set of representations* $(f_{1,2}, ..., f_{n-1,n})$ *of these algebras as* $\overline{f}$. *Since different algebras have different type, we also talk about the set of* $\overline{\Omega}$*-algebras. In relation to the set* $\overline{A}$*, we also use matrix notations that we discussed in section [7]-2.1. For instance, we use the symbol* $\overline{A}_{[1]}$ *to denote the set of* $\overline{\Omega}$*-algebras* $(A_2, ..., A_n)$*. In the corresponding notation* $(\overline{A}_{[1]}, \overline{f})$ *of tower of representation, we assume that* $\overline{f} = (f_{2,3}, ..., f_{n-1,n})$. □

CONVENTION 1.2.3. *Since we use vector notation for elements of the tower of representations, we need convention about notation of operation. We assume that we get result of operation componentwise. For instance,*

$$\overline{r}(\overline{a}) = (r_1(a_1), ..., r_n(a_n))$$

□

REMARK 1.2.4. I believe that diagrams of maps are an important tool. However, sometimes I want to see the diagram as three dimensional figure and I expect that this would increase its expressive power. Who knows what surprises the future holds. In 1992, at a conference in Kazan, I have described to my colleagues what advantages the computer preparation of papers has. 8 years later I learned from the letter from Kazan that now we can prepare paper using LaTeX. □

CHAPTER 2

# Representation of Universal Algebra

## 2.1. Representation of Universal Algebra

Definition 2.1.1. Let the structure of $\Omega_2$-algebra be defined on the set $M$ ([2, 13]). Endomorphism of $\Omega_2$-algebra

$$t : M \to M$$

is called **transformation of universal algebra** $M$.[2.1] □

We denote $\delta$ identical transformation.

Definition 2.1.2. Let ${}^*M$ be the set of **left-side transformations**

$$u' = tu$$

of $\Omega_2$-algebra $M$. Let ${}^*M$ be $\Omega_1$-algebra. The homomorphism

$$f : A \to {}^*M \tag{2.1.1}$$

of $\Omega_1$-algebra $A$ into $\Omega_1$-algebra ${}^\star M$ is called **left-side representation of $\Omega_1$-algebra** $A$ or $A*$**-representation** in $\Omega_2$-algebra $M$. □

Definition 2.1.3. Let $M^*$ be the set of **right-side transformations**

$$u' = ut$$

of $\Omega_2$-algebra $M$. Let $M^*$ be $\Omega_1$-algebra. The homomorphism

$$f : A \to M^*$$

of $\Omega_1$-algebra $A$ into $\Omega_1$-algebra $M^*$ is called **right-side representation of $\Omega_1$-algebra** $A$ or $*A$**-representation** in $\Omega_2$-algebra $M$. □

We extend the convention described in remarks [7]-2.2.15, [11]-2.2.15 to representation theory. We can write duality principle in the following form

Theorem 2.1.4 (duality principle). *Any statement which holds for left-side representation of $\Omega_1$-algebra $A$ holds also for right-side representation of $\Omega_1$-algebra $A$.*

Remark 2.1.5. There exist two forms of notation for transformation of $\Omega_2$-algebra $M$. In operational notation, we write the transformation $A$ as either $Aa$ which corresponds to the left-side transformation or $aA$ which corresponds to the right-side transformation. In functional notation, we write the transformation $A$

[2.1] If the set of operations of $\Omega_2$-algebra is empty, then

$$t : M \to M$$

is a map.

as $A(a)$ regardless of the fact whether this is left-side or right-side transformation. This notation is in agreement with duality principle.

This remark serves as a basis for the following convention. When we use functional notation we do not make a distinction whether this is left-side or right-side transformation. We denote $^*M$ the set of transformations of $\Omega_2$-algebra $M$. Suppose we defined the structure of $\Omega_1$-algebra on the set $^*M$. Let $A$ be $\Omega_1$-algebra. We call homomorphism

$$f : A \to {}^*M \tag{2.1.2}$$

**representation of $\Omega_1$-algebra $A$ in $\Omega_2$-algebra $M$.** We also use record

$$f : A \xrightarrow{*} M$$

to denote the representation of $\Omega_1$-algebra $A$ in $\Omega_2$-algebra $M$.

Correspondence between operational notation and functional notation is unambiguous. We can select any form of notation which is convenient for presentation of particular subject. □

There are several ways to describe the representation. We can define the map $f$ keeping in mind that the domain is $\Omega_1$-algebra $A$ and range is $\Omega_1$-algebra $^*M$. Either we can specify $\Omega_1$-algebra $A$ and $\Omega_2$-algebra $M$ keeping in mind that we know the structure of the map $f$.[2.2]

Diagram

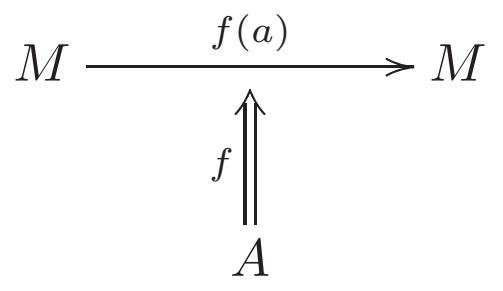

means that we consider the representation of $\Omega_1$-algebra $A$. The map $f(a)$ is image of $a \in A$.

Definition 2.1.6. Let the map (2.1.2) be an isomorphism of the $\Omega_1$-algebra $A$ into $^*M$. Then the representation of the $\Omega_1$-algebra $A$ is called **effective**. □

Remark 2.1.7. If the left-side representation of $\Omega_1$-algebra is effective, then we identify an element of $\Omega_1$-algebra and its image and write left-side transformation caused by element $a \in A$ as

$$v' = av$$

If the right-side representation of $\Omega_1$-algebra is effective, then we identify an element of $\Omega_1$-algebra and its image and write right-side transformation caused by element $a \in A$ as

$$v' = va$$

□

Definition 2.1.8. We call a representation of $\Omega_1$-algebra **transitive** if for any $a, b \in V$ exists such $g$ that

$$a = f(g)(b)$$

We call a representation of $\Omega_1$-algebra **single transitive** if it is transitive and effective. □

[2.2] For instance, we consider vector space $\overline{V}$ over field $D$ ( definition [7]-5.1.4 ).

Theorem 2.1.9. *Representation is single transitive iff for any $a, b \in M$ exists one and only one $g \in A$ such that $a = f(g)(b)$*

Proof. Corollary of definitions 2.1.6 and 2.1.8. □

## 2.2. Morphism of Representations of Universal Algebra

Theorem 2.2.1. *Let $A$ and $B$ be $\Omega_1$-algebras. Representation of $\Omega_1$-algebra $B$*

$$g : B \to {}^\star M$$

*and homomorphism of $\Omega_1$-algebra*

$$h : A \to B \tag{2.2.1}$$

*define representation $f$ of $\Omega_1$-algebra $A$*

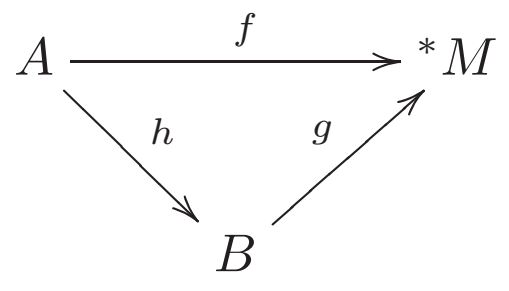

Proof. Since map $g$ is homomorphism of $\Omega_1$-algebra $B$ into $\Omega_1$-algebra ${}^*M$, the map $f$ is homomorphism of $\Omega_1$-algebra $A$ into $\Omega_1$-algebra ${}^*M$. □

Considering representations of $\Omega_1$-algebra in $\Omega_2$-algebras $M$ and $N$, we are interested in a map that preserves the structure of representation.

Definition 2.2.2. Let

$$f : A \xrightarrow{*} M$$

be representation of $\Omega_1$-algebra $A$ in $\Omega_2$-algebra $M$ and

$$g : B \xrightarrow{*} N$$

be representation of $\Omega_1$-algebra $B$ in $\Omega_2$-algebra $N$. Tuple of maps

$$(r : A \to B, R : M \to N) \tag{2.2.2}$$

such, that

- $r$ is homomorphism of $\Omega_1$-algebra
- $R$ is homomorphism of $\Omega_2$-algebra
- 

$$R \circ f(a) = g(r(a)) \circ R \tag{2.2.3}$$

is called **morphism of representations from $f$ into $g$**. We also say that **morphism of representations of $\Omega_1$-algebra in $\Omega_2$-algebra** is defined. □

Remark 2.2.3. We may consider a pair of maps $r$, $R$ as map

$$F : A \cup M \to B \cup N$$

such that

$$F(A) = B \quad F(M) = N$$

Therefore, hereinafter we will say that we have the map $(r, R)$. □

Definition 2.2.4. If representation $f$ and $g$ coincide, then morphism of representations $(r, R)$ is called **morphism of representation $f$**. □

For any $m \in M$ equation (2.2.3) has form

$$R(f(a)(m)) = g(r(a))(R(m)) \tag{2.2.4}$$

Remark 2.2.5. Consider morphism of representations (2.2.2). We denote elements of the set $B$ by letter using pattern $b \in B$. However if we want to show that $b$ is image of element $a \in A$, we use notation $r(a)$. Thus equation

$$r(a) = r(a)$$

means that $r(a)$ (in left part of equation) is image $a \in A$ (in right part of equation). Using such considerations, we denote element of set $N$ as $R(m)$. We will follow this convention when we consider correspondences between homomorphisms of $\Omega_1$-algebra and maps between sets where we defined corresponding representations. □

Remark 2.2.6. There are two ways to interpret (2.2.4)

- Let transformation $f(a)$ map $m \in M$ into $f(a)(m)$. Then transformation $g(r(a))$ maps $R(m) \in N$ into $R(f(a)(m))$.
- We represent morphism of representations from $f$ into $g$ using diagram

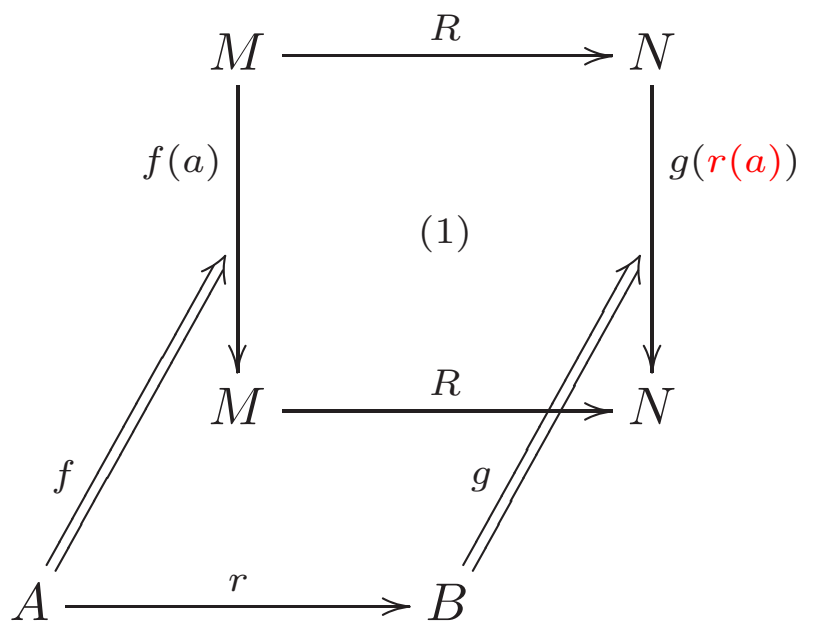

From (2.2.3), it follows that diagram (1) is commutative.

□

Theorem 2.2.7. *Consider representation*

$$f : A \xrightarrow{*} M$$

*of $\Omega_1$-algebra $A$ and representation*

$$g : B \xrightarrow{*} N$$

*of $\Omega_1$-algebra $B$. Morphism*

$$h : A \longrightarrow B \qquad H : M \longrightarrow N$$

*of representations from $f$ into $g$ satisfies equation*

$$H \circ \omega(f(a_1), ..., f(a_n)) = \omega(g(h(a_1)), ..., g(h(a_n))) \circ H \tag{2.2.5}$$

*for any $n$-ary operation $\omega$ of $\Omega_1$-algebra.*

Proof. Since $f$ is homomorphism, we have

$$H \circ \omega(f(a_1), ..., f(a_n)) = H \circ f(\omega(a_1, ..., a_n)) \tag{2.2.6}$$

From (2.2.3) and (2.2.6) it follows that

$$H \circ \omega(f(a_1), ..., f(a_n)) = g(h(\omega(a_1, ..., a_n))) \circ H \tag{2.2.7}$$

Since $h$ is homomorphism, from (2.2.7) it follows that

$$H\circ\omega(f(a_1),...,f(a_n))=g(\omega(h(a_1),...,h(a_n)))\circ H \tag{2.2.8}$$

Since $g$ is homomorphism, (2.2.5) follows from (2.2.8). □

Theorem 2.2.8. *Let the map*

$$h:A\rightarrow B\quad H:M\rightarrow N$$

*be morphism from representation*

$$f:A\longrightarrow\!\!\!\!\!*\;\; M$$

*of* $\Omega_1$*-algebra* $A$ *into representation*

$$g:B\longrightarrow\!\!\!\!\!*\;\; N$$

*of* $\Omega_1$*-algebra* $B$*. If representation* $f$ *is effective, then the map*

$${}^*H:{}^*M\rightarrow{}^*N$$

*defined by equation*

$${}^*H(f(a))=g(h(a)) \tag{2.2.9}$$

*is homomorphism of* $\Omega_1$*-algebra.*

Proof. Because representation $f$ is effective, then for given transformation $f(a)$ element $a$ is determined uniquely. Therefore, transformation $g(h(a))$ is properly defined in equation (2.2.9).

Since $f$ is homomorphism, we have

$${}^*H(\omega(f(a_1),...,f(a_n)))={}^*H(f(\omega(a_1,...,a_n))) \tag{2.2.10}$$

From (2.2.9) and (2.2.10) it follows that

$${}^*H(\omega(f(a_1),...,f(a_n)))=g(h(\omega(a_1,...,a_n))) \tag{2.2.11}$$

Since $h$ is homomorphism, from (2.2.11) it follows that

$${}^*H(\omega(f(a_1),...,f(a_n)))=g(\omega(h(a_1),...,h(a_n))) \tag{2.2.12}$$

Since $g$ is homomorphism,

$${}^*H(\omega(f(a_1),...,f(a_n)))=\omega(g(h(a_1)),...,g(h(a_n)))=\omega({}^*H(f(a_1)),...,{}^*H(f(a_n)))$$

follows from (2.2.12). Therefore, the map ${}^*H$ is homomorphism of $\Omega_1$-algebra. □

Theorem 2.2.9. *Let*

$$f:A\longrightarrow\!\!\!\!\!*\;\; M$$

*be single transitive representation of* $\Omega_1$*-algebra* $A$ *and*

$$g:B\longrightarrow\!\!\!\!\!*\;\; N$$

*be single transitive representation of* $\Omega_1$*-algebra* $B$*. Given homomorphism of* $\Omega_1$ *algebra*

$$h:A\rightarrow B$$

*there exists morphism of representations from* $f$ *into* $g$

$$h:A\rightarrow B\quad H:M\rightarrow N$$

Proof. Let us choose homomorphism $h$. Let us choose element $m\in M$ and element $n\in N$. To define map $H$, consider following diagram

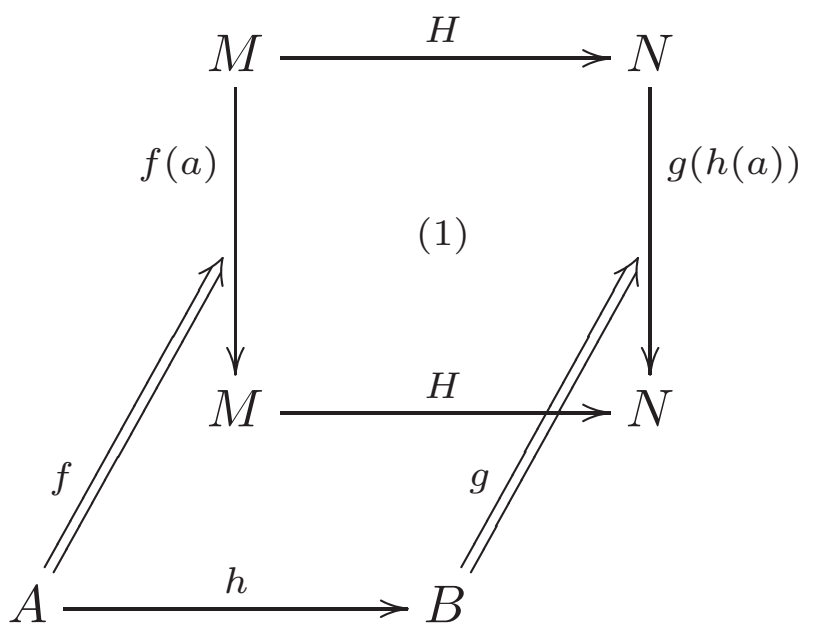

From commutativity of diagram (1), it follows that

$$H(f(a)m)=h(g(a))H(m)$$

For arbitrary $m'\in M$, we defined unambiguously $a\in A$ such that $m'=f(a)m$. Therefore, we defined map $H$ which satisfies to equation (2.2.3). □

Theorem 2.2.10. *Let the representation*

$$f:A\xrightarrow{*}M$$

*of $\Omega_1$-algebra $A$ be single transitive representation and the representation*

$$g:B\xrightarrow{*}N$$

*of $\Omega_1$-algebra $B$ be single transitive representation. Given homomorphism of $\Omega_1$-algebra*

$$h:A\longrightarrow B$$

*consider a homomorphism of $\Omega_2$-algebra*

$$H:M\longrightarrow N$$

*such that $(h,H)$ is morphism of representations from $f$ into $g$. This map is unique up to choice of image $n=H(m)\in N$ of given element $m\in M$.*

Proof. From proof of theorem 2.2.9, it follows that choice of homomorphism $h$ and elements $m\in M$, $n\in N$ uniquely defines the map $H$. □

Theorem 2.2.11. *Given single transitive representation*

$$f:A\xrightarrow{*}M$$

*of $\Omega_1$-algebra $A$, for any endomorphism $h$ of $\Omega_1$-algebra $A$ there exists morphism of representation $f$*

$$h:A\rightarrow A\quad H:M\rightarrow M$$

Proof. Consider following diagram

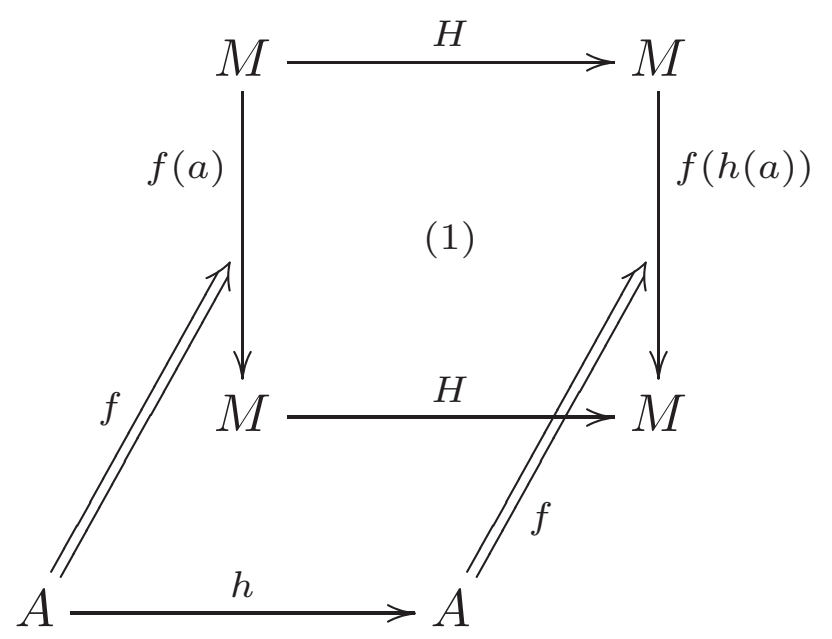

Statement of theorem is corollary of the theorem 2.2.9. □

Theorem 2.2.12. *Let*

$$f : A \longrightarrow\!\!\!\!\!\!\!\!*\ \ \ \ M$$

*be representation of $\Omega_1$-algebra $A$,*

$$g : B \longrightarrow\!\!\!\!\!\!\!\!*\ \ \ \ N$$

*be representation of $\Omega_1$-algebra $B$,*

$$h : C \to {}^*L$$

*be representation of $\Omega_1$-algebra $C$. Given morphisms of representations of $\Omega_1$-algebra*

$$p : A \longrightarrow B \qquad P : M \longrightarrow N$$

$$q : B \longrightarrow C \qquad Q : N \longrightarrow L$$

*There exists morphism of representations of $\Omega_1$-algebra*

$$r : A \longrightarrow C \qquad R : M \longrightarrow L$$

*where $r = qp$, $R = QP$. We call morphism $(r, R)$ of representations from $f$ into $h$* **product of morphisms** $(p, P)$ **and** $(q, Q)$ **of representations of universal algebra**.

Proof. We represent statement of theorem using diagram

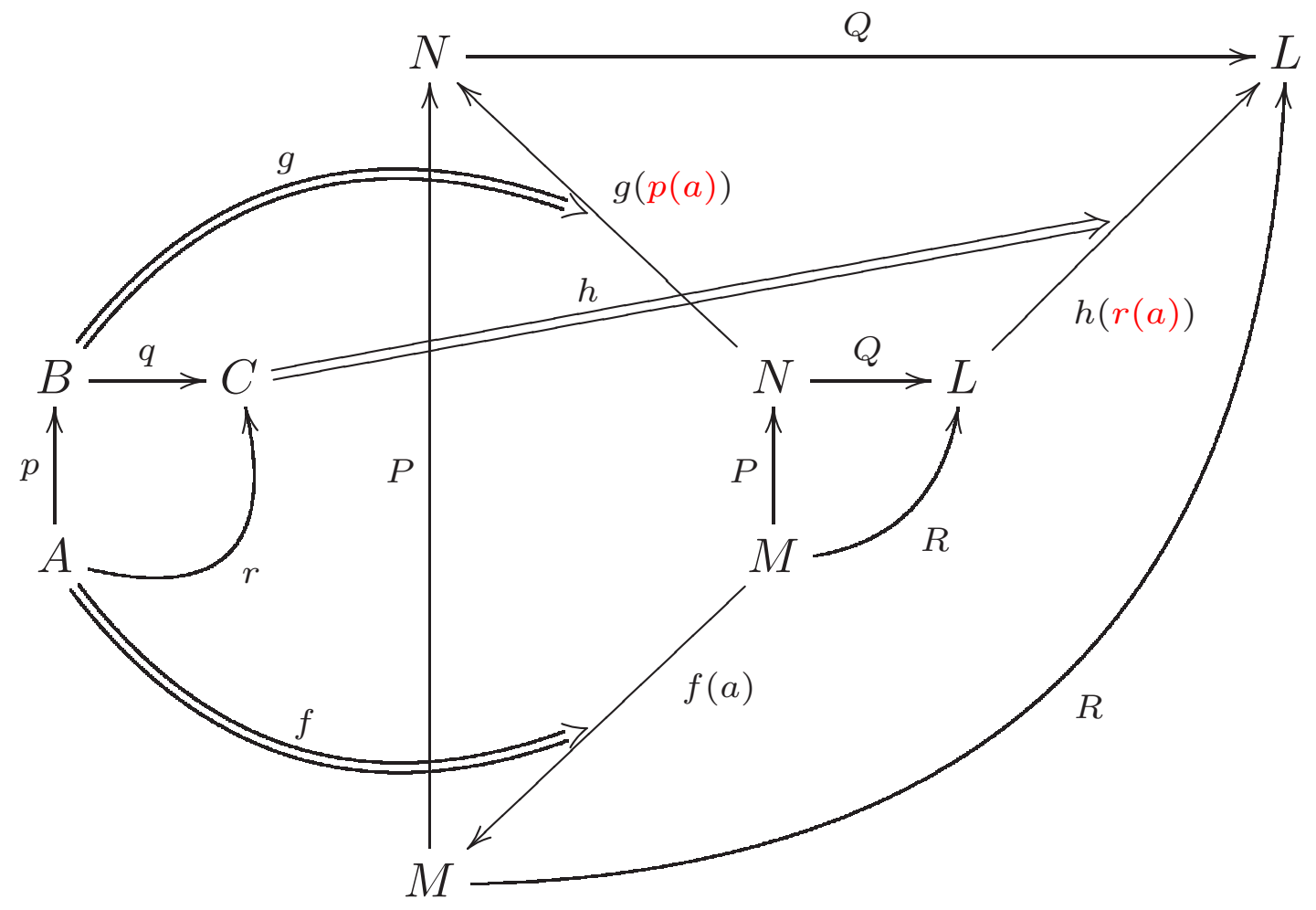

Map $r$ is homomorphism of $\Omega_1$-algebra $A$ into $\Omega_1$-algebra $C$. We need to show that tuple of maps $(r,R)$ satisfies to (2.2.3):

$$\begin{aligned}R(f(a)m)&=QP(f(a)m)\\&=Q(g(p(a))P(m))\\&=h(qp(a))QP(m))\\&=h(r(a))R(m)\end{aligned}$$

□

Definition 2.2.13. Let $\mathcal{A}$ be category of $\Omega_1$-algebras. We define **category $\mathcal{A}*$ of left-side representations** of $\Omega_1$-algebra from category $\mathcal{A}$. Left-side representations of $\Omega_1$-algebra are objects of this category. Morphisms of left-side representations of $\Omega_1$-algebra are morphisms of this category. □

Definition 2.2.14. Let us define equivalence $S$ on the set $M$. Transformation $f$ is called **coordinated with equivalence** $S$, when $f(m_1)\equiv f(m_2)(\mathrm{mod}\ S)$ follows from condition $m_1\equiv m_2(\mathrm{mod}S)$. □

Theorem 2.2.15. *Consider equivalence $S$ on set $M$. Consider $\Omega_1$-algebra on set ${}^*M$. If any transformation $f\in {}^*M$ is coordinated with equivalence $S$, then we can define the structure of $\Omega_1$-algebra on the set ${}^*(M/S)$.*

Proof. Let $h=\mathrm{nat}\ S$. If $m_1\equiv m_2(\mathrm{mod}S)$, then $h(m_1)=h(m_2)$. Since $f\in {}^*M$ is coordinated with equivalence $S$, then $h(f(m_1))=h(f(m_2))$. This allows us to define transformation $F$ according to rule

$$F([m])=h(f(m)) \tag{2.2.13}$$

Let $\omega$ be n-ary operation of $\Omega_1$-algebra. Suppose $f_1$, ..., $f_n\in {}^\star M$ and

$$F_1([m])=h(f_1(m))\quad ...\quad F_n([m])=h(f_n(m))$$

According to condition of theorem, the transformation

$$f=\omega(f_1,...,f_n)\in {}^*M$$

is coordinated with equivalence $S$. Therefore,

$$\begin{gathered}f(m_1)\equiv f(m_2)(\mathrm{mod}S)\\ \omega(f_1,...,f_n)(m_1)\equiv\omega(f_1,...,f_n)(m_2)(\mathrm{mod}S)\end{gathered}\tag{2.2.14}$$

follows from condition $m_1\equiv m_2(\mathrm{mod}S)$ and the definition 2.2.14. Therefore, we can define operation $\omega$ on the set ${}^\star(M/S)$ according to rule

$$\omega(F_1,...,F_n)[m]=h(\omega(f_1,...,f_n)(m)) \tag{2.2.15}$$

From the definition (2.2.13) and equation (2.2.14), it follows that we properly defined operation $\omega$ on the set ${}^\star(M/S)$. □

Theorem 2.2.16. *Let*

$$f:A\xrightarrow{*} M$$

*be representation of $\Omega_1$-algebra $A$,*

$$g:B\xrightarrow{*} N$$

*be representation of $\Omega_1$-algebra $B$. Let*

$$r : A \longrightarrow B \qquad R : M \longrightarrow N$$

*be morphism of representations from $f$ into $g$. Suppose*

$$s = rr^{-1} \qquad S = RR^{-1}$$

*Then there exist decompositions of $r$ and $R$, which we describe using diagram*

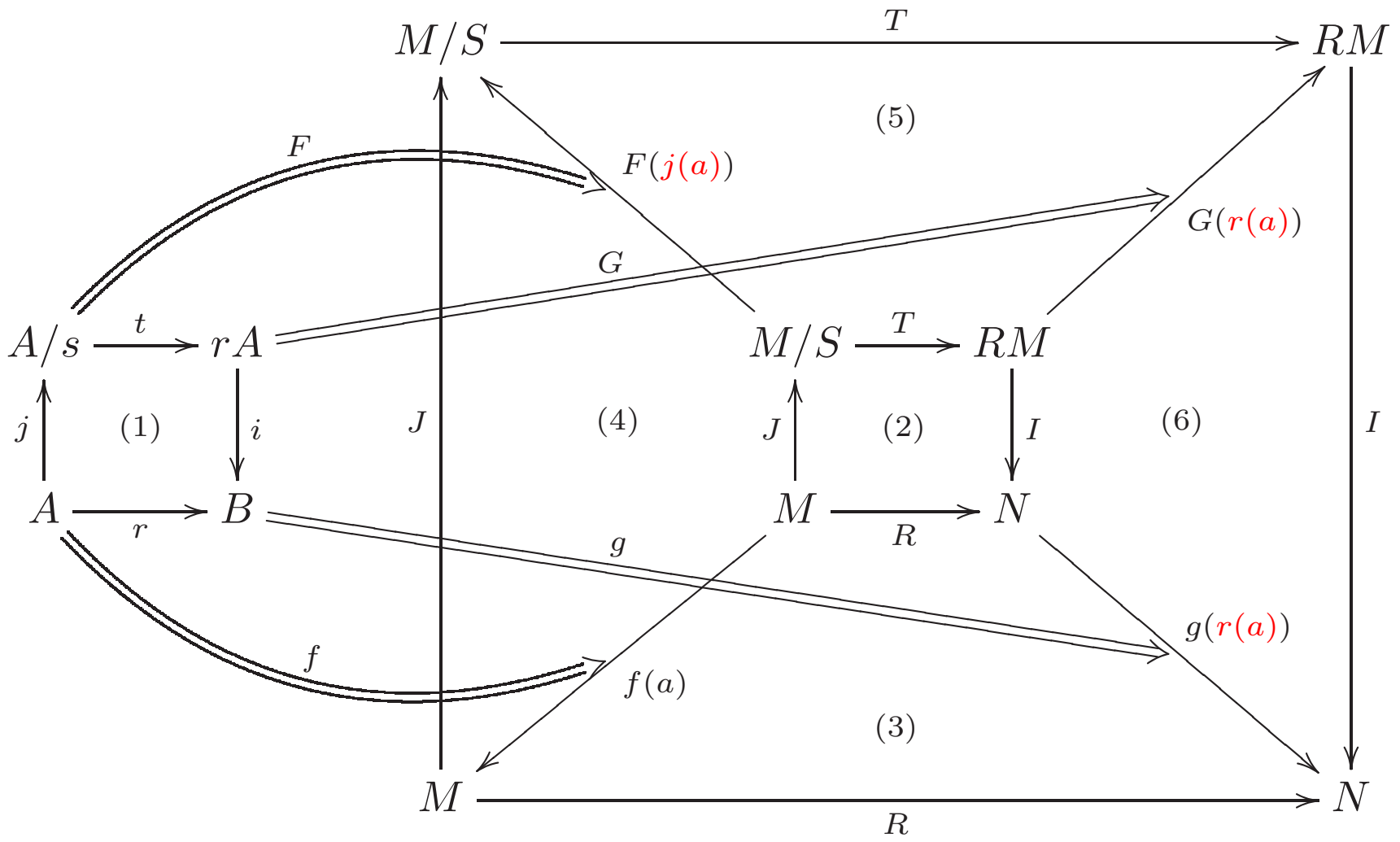

(1) *$s = \ker\ r$ is a congruence on $A$. There exists decompositions of homomorphism $r$*

(2.2.16) $$r = itj$$

*$j = \text{nat}\ s$ is the natural homomorphism*

$$j(a) = j(a)$$

*$t$ is isomorphism*

(2.2.17) $$r(a) = t(j(a))$$

*$i$ is the inclusion map*

(2.2.18) $$r(a) = i(r(a))$$

(2) *$S = \ker\ R$ is an equivalence on $M$. There exists decompositions of homomorphism $R$*

(2.2.19) $$R = ITJ$$

*$J = \text{nat}\ S$ is surjection*

$$J(m) = J(m)$$

*$T$ is bijection*

(2.2.20) $$R(m) = T(J(m))$$

*$I$ is the inclusion map*

(2.2.21) $$R(m) = I(R(m))$$

(3) *$F$ is left-side representation of $\Omega_1$-algebra $A/s$ in $M/S$*

(4) *$G$ is left-side representation of $\Omega_1$-algebra $rA$ in $RM$*
(5) *$(j, J)$ is morphism of representations $f$ and $F$*
(6) *$(t, T)$ is morphism of representations $F$ and $G$*
(7) *$(t^{-1}, T^{-1})$ is morphism of representations $G$ and $F$*
(8) *$(i, I)$ is morphism of representations $G$ and $g$*
(9) *There exists decompositions of morphism of representations*

$$(r, R) = (i, I)(t, T)(j, J) \tag{2.2.22}$$

Proof. Existence of diagrams (1) and (2) follows from theorem II.3.7 ([13], p. 60).

We start from diagram (4).

Let $m_1 \equiv m_2 (\text{mod } S)$. Then

$$R(m_1) = R(m_2) \tag{2.2.23}$$

Since $a_1 \equiv a_2 (\text{mod } s)$, then

$$r(a_1) = r(a_2) \tag{2.2.24}$$

Therefore, $j(a_1) = j(a_2)$. Since $(r, R)$ is morphism of representations, then

$$R(f(a_1)(m_1)) = g(r(a_1))(R(m_1)) \tag{2.2.25}$$

$$R(f(a_2)(m_2)) = g(r(a_2))(R(m_2)) \tag{2.2.26}$$

From (2.2.23), (2.2.24), (2.2.25), (2.2.26), it follows that

$$R(f(a_1)(m_1)) = R(f(a_2)(m_2)) \tag{2.2.27}$$

From (2.2.27) it follows

$$f(a_1)(m_1) \equiv f(a_2)(m_2)(\text{mod } S) \tag{2.2.28}$$

and, therefore,

$$J(f(a_1)(m_1)) = J(f(a_2)(m_2)) \tag{2.2.29}$$

From (2.2.29) it follows that map

$$F(j(a))(J(m)) = J(f(a)(m))) \tag{2.2.30}$$

is well defined and this map is transformation of set $M/S$.

From equation (2.2.28) (in case $a_1 = a_2$) it follows that for any $a$ transformation is coordinated with equivalence $S$. From theorem 2.2.15 it follows that we defined structure of $\Omega_1$-algebra on the set ${}^*(M/S)$. Consider $n$-ary operation $\omega$ and $n$ transformations

$$F(j(a_i))(J(m)) = J(f(a_i)(m))) \quad i = 1, ..., n$$

of the set $M/S$. We assume

$$\omega(F(j(a_1)), ..., F(j(a_n)))(J(m)) = J(\omega(f(a_1), ..., f(a_n)))(m))$$

Therefore, map $F$ is representations of $\Omega_1$-algebra $A/s$.

From (2.2.30) it follows that $(j, J)$ is morphism of representations $f$ and $F$ (the statement (5) of the theorem).

Consider diagram (5).

Since $T$ is bijection, then we identify elements of the set $M/S$ and the set $MR$, and this identification has form

$$T(J(m))=R(m) \tag{2.2.31}$$

We can write transformation $F(j(a))$ of the set $M/S$ as

$$F(j(a)):J(m)\rightarrow F(j(a))(J(m)) \tag{2.2.32}$$

Since $T$ is bijection, we define transformation

$$T(J(m))\rightarrow T(F(j(a))(J(m))) \tag{2.2.33}$$

of the set $RM$. Transformation (2.2.33) depends on $j(a)\in A/s$. Since $t$ is bijection, we identify elements of the set $A/s$ and the set $rA$, and this identification has form

$$t(j(a))=r(a) \tag{2.2.34}$$

Therefore, we defined map

$$G:rA\rightarrow {}^\star RM$$

according to equation

$$G(t(j(a)))(T(J(m)))=T(F(j(a))(J(m))) \tag{2.2.35}$$

Consider $n$-ary operation $\omega$ and $n$ transformations

$$G(r(a_i))(R(m))=T(F(j(a_i))(J(m)))\quad i=1,...,n$$

of space $RM$. We assume

$$\omega(G(r(a_1)),...,G(r(a_n)))(R(m))=T(\omega(F(j(a_1),...,F(j(a_n)))(J(m))) \tag{2.2.36}$$

According to (2.2.35) operation $\omega$ is well defined on the set ${}^\star RM$. Therefore, the map $G$ is representations of $\Omega_1$-algebra.

From (2.2.35), it follows that $(t,T)$ is morphism of representations $F$ and $G$ (the statement (6) of the theorem).

Since $T$ is bijection, then from equation (2.2.31), it follows that

$$J(m)=T^{-1}(R(m)) \tag{2.2.37}$$

We can write transformation $G(r(a))$ of the set $RM$ as

$$G(r(a)):R(m)\rightarrow G(r(a))(R(m)) \tag{2.2.38}$$

Since $T$ is bijection, we define transformation

$$T^{-1}(R(m))\rightarrow T^{-1}(G(r(a))(R(m))) \tag{2.2.39}$$

of the set $M/S$. Transformation (2.2.39) depends on $r(a)\in rA$. Since $t$ is bijection, then from equation (2.2.34) it follows that

$$j(a)=t^{-1}(r(a)) \tag{2.2.40}$$

Since, by construction, diagram (5) is commutative, then transformation (2.2.39) coincides with transformation (2.2.32). We can write the equation (2.2.36) as

$$T^{-1}(\omega(G(r(a_1)),...,G(r(a_n)))(R(m)))=\omega(F(j(a_1),...,F(j(a_n)))(J(m)) \tag{2.2.41}$$

Therefore $(t^{-1},T^{-1})$ is morphism of representations $G$ and $F$ (the statement (7) of the theorem).

Diagram (6) is the simplest case in our proof. Since map $I$ is immersion and diagram (2) is commutative, we identify $n \in N$ and $R(m)$ when $n \in \mathrm{Im}R$. Similarly, we identify corresponding transformations.

$$g'(i(r(a)))(I(R(m))) = I(G(r(a))(R(m))) \tag{2.2.42}$$

$$\omega(g'(r(a_1)), ..., g'(r(a_n)))(R(m)) = I(\omega(G(r(a_1), ..., G(r(a_n)))(R(m)))$$

Therefore, $(i, I)$ is morphism of representations $G$ and $g$ (the statement (8) of the theorem).

To prove the statement (9) of the theorem we need to show that defined in the proof representation $g'$ is congruent with representation $g$, and operations over transformations are congruent with corresponding operations over ${}^*N$.

$$\begin{aligned}
g'(i(r(a)))(I(R(m))) &= I(G(r(a))(R(m))) && \text{by (2.2.42)}\\
&= I(G(t(j(a)))(T(J(m)))) && \text{by (2.2.17), (2.2.20)}\\
&= IT(F(j(a))(J(m))) && \text{by (2.2.35)}\\
&= ITJ(f(a)(m)) && \text{by (2.2.30)}\\
&= R(f(a)(m)) && \text{by (2.2.19)}\\
&= g(r(a))(R(m)) && \text{by (2.2.3)}
\end{aligned}$$

$$\begin{aligned}
\omega(G(r(a_1)), ..., G(r(a_n)))(R(m)) &= T(\omega(F(j(a_1), ..., F(j(a_n)))(J(m)))\\
&= T(F(\omega(j(a_1), ..., j(a_n)))(J(m)))\\
&= T(F(j(\omega(a_1, ..., a_n)))(J(m)))\\
&= T(J(f(\omega(a_1, ..., a_n))(m)))
\end{aligned}$$

□

Definition 2.2.17. Let

$$f : A \xrightarrow{*} M$$

be representation of $\Omega_1$-algebra $A$,

$$g : B \xrightarrow{*} N$$

be representation of $\Omega_1$-algebra $B$. Let

$$r : A \longrightarrow B \qquad R : M \longrightarrow N$$

be morphism of representations from $f$ into $g$ such that $f$ is isomorphism of $\Omega_1$-algebra and $g$ is isomorphism of $\Omega_2$-algebra. Then map $(r, R)$ is called **isomorphism of repesentations**. □

Theorem 2.2.18. *In the decomposition (2.2.22), the map $(t, T)$ is isomorphism of representations $F$ and $G$.*

Proof. The statement of the theorem is corollary of definition 2.2.17 and statements (6) and (7) of the theorem 2.2.16. □

## 2.3. Reduced Morphism of Representations

From theorem 2.2.16 it follows that we can reduce the problem of studying of morphism of representations of $\Omega_1$-algebra to the case described by diagram

(2.3.1)

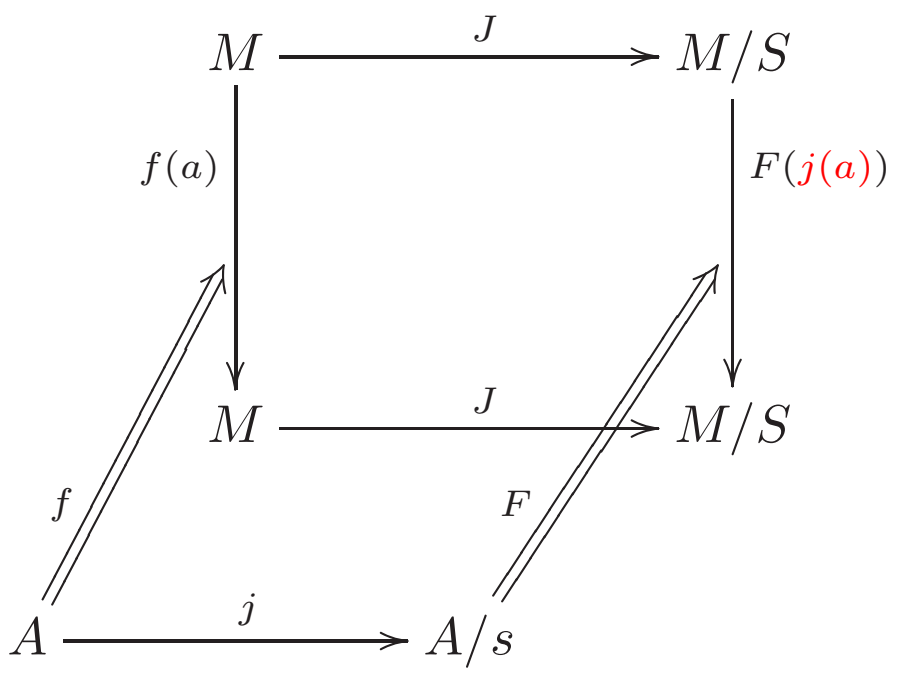

THEOREM 2.3.1. *We can supplement diagram (2.3.1) with representation $F_1$ of $\Omega_1$-algebra $A$ into set $M/S$ such that diagram*

(2.3.2)

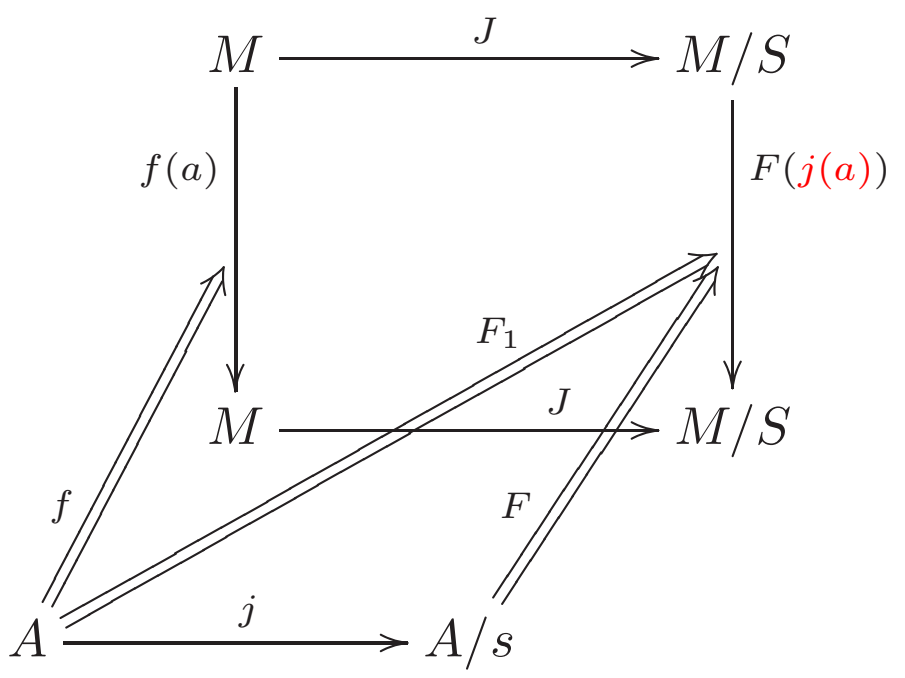

*is commutative. The set of transformations of representation $F$ and the set of transformations of representation $F_1$ coincide.*

PROOF. To prove theorem it is enough to assume

$$F_1(a) = F(j(a))$$

Since map $j$ is surjection, then $\mathrm{Im}F_1 = \mathrm{Im}F$. Since $j$ and $F$ are homomorphisms of $\Omega_1$-algebra, then $F_1$ is also homomorphism of $\Omega_1$-algebra. □

REMARK 2.3.2. Theorem 2.3.1 completes the series of theorems dedicated to the structure of morphism of representations $\Omega_1$-algebra. From these theorems it follows that we can simplify task of studying of morphism of representations $\Omega_1$-algebra and not go beyond morphism of representations of form

$$id : A \longrightarrow A \qquad R : M \longrightarrow N$$

In this case we identify morphism of $(\mathrm{id}, R)$ representations of $\Omega_1$-algebra and corresponding homomorphism $R$ of $\Omega_2$-algebra and the homomorphism $R$ is called **reduced morphism of representations**. We will use diagram

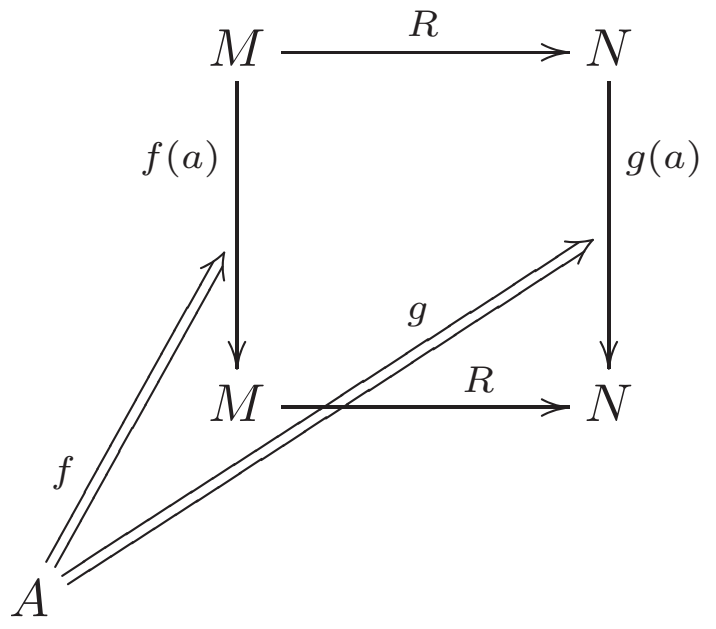

to represent reduced morphism $R$ of representations of $\Omega_1$-algebra. From diagram it follows

$$R \circ f(a) = g(a) \circ R \tag{2.3.3}$$

□

THEOREM 2.3.3. *Let the map*

$$H : M \to N$$

*be reduced morphism from representation*

$$f : A \longrightarrow\!\!\!\!\!\!*\;\; M$$

*of $\Omega_1$-algebra $A$ into representation*

$$g : A \longrightarrow\!\!\!\!\!\!*\;\; N$$

*of $\Omega_1$-algebra $A$. If representation $f$ is effective, then the map*

$$^*H : {}^*M \to {}^*N$$

*defined by equation*

$$^*H(f(a)) = g(a) \tag{2.3.4}$$

*is homomorphism of $\Omega_1$-algebra.*

PROOF. The theorem follows from the theorem 2.2.8, if we assume $h = id$. □

THEOREM 2.3.4. *Let representations*

$$f : A \longrightarrow\!\!\!\!\!\!*\;\; M$$

$$g : A \longrightarrow\!\!\!\!\!\!*\;\; N$$

*of $\Omega_1$-algebra $A$ be single transitive representations. There exists reduced morphism of representations from $f$ into $g$*

$$H : M \to N$$

Proof. Let us choose homomorphism $h$. Let us choose element $m \in M$ and element $n \in N$. To define map $H$, consider following diagram

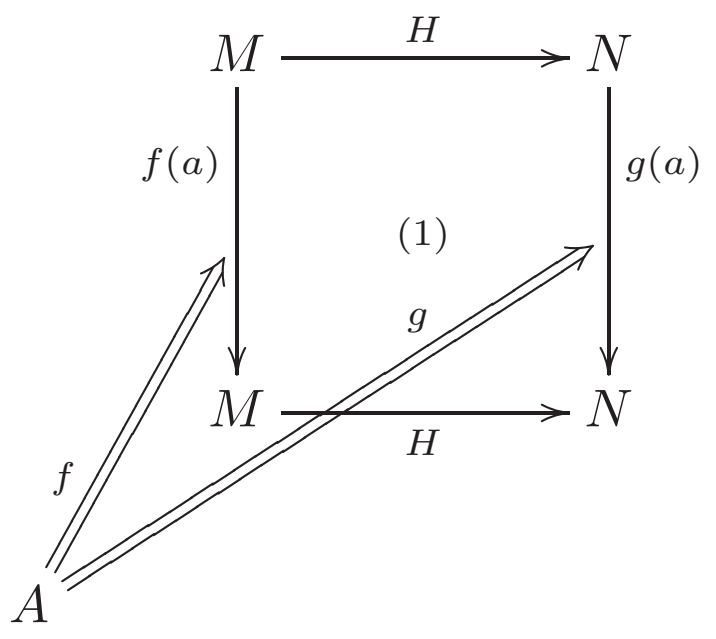

From commutativity of diagram (1), it follows that

$$H(f(a)m) = g(a)H(m)$$

For arbitrary $m' \in M$, we defined unambiguously $a \in A$ such that $m' = f(a)m$. Therefore, we defined map $H$ which satisfies to equation (2.3.3). □

Theorem 2.3.5. *Let representations*

$$f : A \twoheadrightarrow M$$

$$g : A \twoheadrightarrow N$$

*of $\Omega_1$-algebra $A$ be single transitive representations. Reduced morphism of representations from $f$ into $g$*

$$H : M \longrightarrow N$$

*is unique up to choice of image $n = H(m) \in N$ of given element $m \in M$.*

Proof. From proof of theorem 2.3.4, it follows that choice of elements $m \in M$, $n \in N$ uniquely defines the map $H$. □

By analogy with definition 2.2.13. we give following definition.

Definition 2.3.6. We define **category $A*$ of left-side representations of $\Omega_1$-algebra** $A$. Left-side representations of $\Omega_1$-algebra $A$ are objects of this category. Morphisms $(id, R)$ of left-side representations of $\Omega_1$-algebra $A$ are morphisms of this category. □

## 2.4. Automorphism of Representation of Universal Algebra

Definition 2.4.1. Let

$$f : A \twoheadrightarrow M$$

be representation of $\Omega_1$-algebra $A$ in $\Omega_2$-algebra $M$. The reduced morphism of representations of $\Omega_1$-algebra

$$R : M \to M$$

such, that $R$ is endomorphism of $\Omega_2$-algebra is called **endomorphism of representation** $f$. □

Theorem 2.4.2. *Given single transitive representation*

$$f:A\xrightarrow{*} M$$

*of $\Omega_1$-algebra $A$, for any $p$, $q\in M$ there exists unique endomorphism*

$$H:M\to M$$

*of representation $f$ such that $H(p)=q$.*

Proof. Consider following diagram

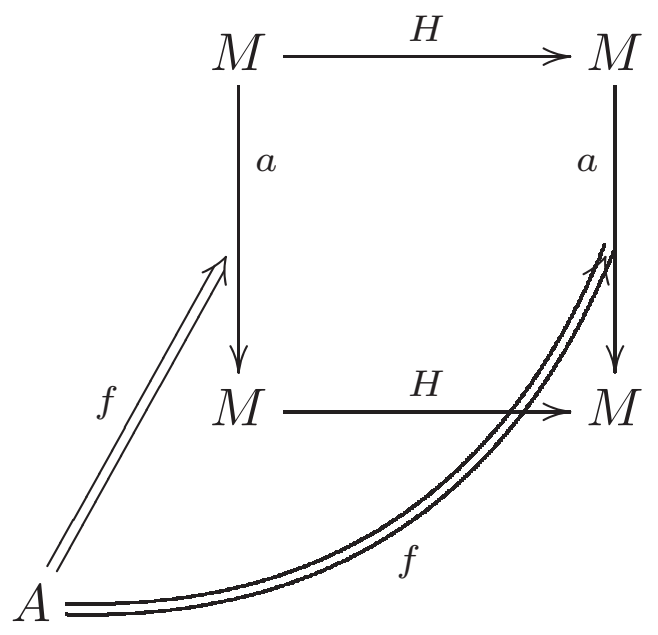

Existence of endomorphism is corollary of the theorem 2.2.9. For given $p$, $q\in M$, uniqueness of endomorphism follows from the theorem 2.2.10 when $h=\mathrm{id}$. □

Theorem 2.4.3. *Endomorphisms of representation $f$ form semigroup.*

Proof. From theorem 2.2.12, it follows that the product of endomorphisms $(p,P)$, $(r,R)$ of the representation $f$ is endomorphism $(pr,PR)$ of the representation $f$. □

Definition 2.4.4. Let

$$f:A\xrightarrow{*} M$$

be representation of $\Omega_1$-algebra $A$ in $\Omega_2$-algebra $M$. The morphism of representations of $\Omega_1$-algebra

$$R:M\to M$$

such, that $R$ is automorphism of $\Omega_2$-algebra is called **automorphism of representation** $f$. □

Theorem 2.4.5. *Let*

$$f:A\xrightarrow{*} M$$

*be representation of $\Omega_1$-algebra $A$ in $\Omega_2$-algebra $M$. The set of automorphisms of the representation $f$ forms* **group** $GA(f)$ .

Proof. Let $R$, $P$ be automorphisms of the representation $f$. According to definition 2.4.4, maps $R$, $P$ are automorphisms of $\Omega_2$-algebra $M$. According to theorem II.3.2, ([13], p. 57), the map $R\circ P$ is automorphism of $\Omega_2$-algebra $M$. From the theorem 2.2.12 and the definition 2.4.4, it follows that product of automorphisms $R\circ P$ of the representation $f$ is automorphism of the representation $f$.

Let $R$, $P$, $Q$ be automorphisms of the representation $f$. The associativity of product of maps $R$, $P$, $Q$ follows from the chain of equations[2.3]

$$\begin{aligned}((R\circ P)\circ Q)(a)&=(R\circ P)(Q(a))=R(P(Q(a)))\\&=R((P\circ Q)(a))=(R\circ(P\circ Q))(a)\end{aligned}$$

Let $R$ be an automorphism of the representation $f$. According to definition 2.4.4 the map $R$ is automorphism of $\Omega_2$-algebra $M$. Therefore, the map $R^{-1}$ is automorphism of $\Omega_2$-algebra $M$. The equation (2.2.4) is true for automorphism $R$ of representation. Assume $m'=R(m)$. Since $R$ is automorphism of $\Omega_2$-algebra, then $m=R^{-1}(m')$ and we can write (2.2.4) in the form

$$R(f(a')(R^{-1}(m')))=f(a')(m') \tag{2.4.1}$$

Since the map $R$ is automorphism of $\Omega_2$-algebra $M$, then from the equation (2.4.1) it follows that

$$f(a')(R^{-1}(m'))=R^{-1}(f(a')(m')) \tag{2.4.2}$$

The equation (2.4.2) corresponds to the equation (2.2.4) for the map $R^{-1}$. Therefore, map $R^{-1}$ of the representation $f$. □

## 2.5. Examples of Representation of Universal Algebra

**2.5.1. Vector Space.** According to definitions [7]-5.1.4, [11]-4.1.4, an effective left-side representation of division ring $D$ in Abelian group is called **left $D$-vector space**. Similarly an effective representation of associative ring $R$ in Abelian group is called **module over ring** $R$ or $R$**-module**.

**2.5.2. Representation of Group on the Set.** Let $G$ be Abelian group, and $M$ be a set. Consider effective representation of group $G$ on the set $M$. For given $a\in G$, $A\in M$ we assume $A\to A+a$. We also use notation $a=\overrightarrow{AB}$ if

$$B=A+a \tag{2.5.1}$$

Then we can represent action of group as

$$B=A+\overrightarrow{AB} \tag{2.5.2}$$

Since the representation is effective, then from equations (2.5.1), (2.5.2) and equation

$$D=C+a$$

it folows that

$$\overrightarrow{AB}=\overrightarrow{CD}$$

---

[2.3] To prove the associativity of product I follow to the example of the semigroup from [4], p. 20, 21.

## 2.6. Generating Set of Representation

Definition 2.6.1. Let

$$f : A \twoheadrightarrow^{*} M$$

be representation of $\Omega_1$-algebra $A$ in $\Omega_2$-algebra $M$. The set $N \subset M$ is called **stable set of representation** $f$, if $f(a)(m) \in N$ for each $a \in A$, $m \in N$. □

We also say that the set $M$ is stable with respect to the representation $f$.

Theorem 2.6.2. *Let*

$$f : A \twoheadrightarrow^{*} M$$

*be representation of $\Omega_1$-algebra $A$ in $\Omega_2$-algebra $M$. Let set $N \subset M$ be subalgebra of $\Omega_2$-algebra $M$ and stable set of representation $f$. Then there exists representation*

$$f_N : A \to {}^*N$$

*such that $f_N(a) = f(a)|_N$. Representation $f_N$ is called* **subrepresentation of representation** *$f$.*

Proof. Let $\omega_1$ be $n$-ary operation of $\Omega_1$-algebra $A$. Then for each $a_1$, ..., $a_n \in A$ and each $b \in N$

$$\begin{aligned}(f_N(a_1)...f_N(a_n)\omega_1)(b) &= (f(a_1)...f(a_n)\omega_1)(b)\\ &= f(a_1...a_n\omega_1)(b)\\ &= f_N(a_1...a_n\omega_1)(b)\end{aligned}$$

Let $\omega_2$ be $n$-ary operation of $\Omega_2$-algebra $M$. Then for each $b_1$, ..., $b_n \in N$ and each $a \in A$

$$\begin{aligned}f_N(a)(b_1)...f_N(a)(b_n)\omega_2 &= f(a)(b_1)...f(a)(b_n)\omega_2\\ &= f(a)(b_1...b_n\omega_2)\\ &= f_N(a)(b_1...b_n\omega_2)\end{aligned}$$

We proved the statement of theorem. □

From the theorem 2.6.2, it follows that if $f_N$ is subrepresentation of representation $f$, then the mapping $(id : A \to A, id_N : N \to M)$ is morphism of representations.

Theorem 2.6.3. *The set*[2.4] *$\mathcal{B}_f$ of all subrepresentations of representation $f$ generates a closure system on $\Omega_2$-algebra $M$ and therefore is a complete lattice.*

Proof. Let $(K_\lambda)_{\lambda\in\Lambda}$ be the set off subalgebras of $\Omega_2$-algebra $M$ that are stable with respect to representation $f$. We define the operation of intersection on the set $\mathcal{B}_f$ according to rule

$$\bigcap f_{K_\lambda} = f_{\cap K_\lambda}$$

We defined the operation of intersection of subrepresentations properly. $\cap K_\lambda$ is subalgebra of $\Omega_2$-algebra $M$. Let $m \in \cap K_\lambda$. For each $\lambda \in \Lambda$ and for each $a \in A$, $f(a)(m) \in K_\lambda$. Therefore, $f(a)(m) \in \cap K_\lambda$. Therefore, $\cap K_\lambda$ is the stable set of representation $f$. □

[2.4] This definition is similar to definition of the lattice of subalgebras ([13], p. 79, 80)

We denote the corresponding closure operator by $J(f)$. Thus $J(f,X)$ is the intersection of all subalgebras of $\Omega_2$-algebra $M$ containing $X$ and stable with respect to representation $f$.

THEOREM 2.6.4. *Let*[2.5]

$$f : A \longrightarrow\!\!\!\!\!* \;\; M$$

*be representation of $\Omega_1$-algebra $A$ in $\Omega_2$-algebra $M$. Let $X \subset M$. Define a subset $X_k \subset M$ by induction on $k$.*

2.6.4.1: $X_0 = X$
2.6.4.2: $x \in X_k \Rightarrow x \in X_{k+1}$
2.6.4.3: $x_1 \in X_k, ..., x_n \in X_k, \omega \in \Omega_2(n) \Rightarrow x_1...x_n\omega \in X_{k+1}$
2.6.4.4: $x \in X_k, a \in A \Rightarrow f(a)(x) \in X_{k+1}$

*Then*

$$\bigcup_{k=0}^{\infty} X_k = J(f,X) \tag{2.6.1}$$

PROOF. If we put $U = \cup X_k$, then by definition of $X_k$, we have $X_0 \subset J(f,X)$, and if $X_k \subset J(f,X)$, then $X_{k+1} \subset J(f,X)$. By induction it follows that $X_k \subset J(f,X)$ for all $k$. Therefore,

$$U \subset J(f,X) \tag{2.6.2}$$

If $a \in U^n$, $a = (a_1,...,a_n)$, where $a_i \in X_{k_i}$, and if $k = \max\{k_1,...,k_n\}$, then $a_1...a_n\omega \in X_{k+1} \subset U$. Therefore, $U$ is subalgebra of $\Omega_2$-algebra $M$.

If $m \in U$, then there exists such $k$ that $m \in X_k$. Therefore, $f(a)(m) \in X_{k+1} \subset U$ for any $a \in A$. Therefore, $U$ is stable set of the representation $f$.

Since $U$ is subalgebra of $\Omega_2$-algebra $M$ and is a stable set of the representation $f$, then subrepresentation $f_U$ is defined. Therefore,

$$J(f,X) \subset U \tag{2.6.3}$$

From (2.6.2), (2.6.3), it follows that $J(f,X) = U$. □

DEFINITION 2.6.5. $J(f,X)$ is called **subrepresentation generated by set** $X$, and $X$ is a **generating set of subrepresentation** $J(f,X)$. In particular, a **generating set of representation** $f$ is a subset $X \subset M$ such that $J(f,X) = M$. □

It is easy to see that the definition of generating set of representation does not depend on whether representation is effective or not. For this reason hereinafter we will assume that the representation is effective and we will use convention for effective left-side representation in remark 2.1.7. We also will use notation

$$R \circ m = R(m)$$

for image of $m \in M$ under the endomorphism of effective representation. According to the definition of product of mappings, for any endomorphisms $R$, $S$ the following equation is true

$$(R \circ S) \circ m = R \circ (S \circ m) \tag{2.6.4}$$

The equation (2.6.4) is associative law for $\circ$ and allows us to write expression

$$R \circ S \circ m$$

[2.5]The statement of theorem is similar to the statement of theorem 5.1, [13], p. 79.

without brackets.

From theorem 2.6.4, it follows next definition.

Definition 2.6.6. Let $X \subset M$. For each $x \in J(f,X)$ there exists $\Omega_2$-word defined according to following rules.

(1) If $m \in X$, then $m$ is $\Omega_2$-word.
(2) If $m_1$, ..., $m_n$ are $\Omega_2$-words and $\omega \in \Omega_2(n)$, then $m_1...m_n\omega$ is $\Omega_2$-word.
(3) If $m$ is $\Omega_2$-word and $a \in A$, then $am$ is $\Omega_2$-word.

**$\Omega_2$-word** $w(f,X,m)$ represents given element $m \in J(f,X)$. We will identify an element $m \in J(f,X)$ and corresponding it $\Omega_2$-word using equation

$$m = w(f,X,m)$$

Similarly, for an arbitrary set $B \subset J(f,X)$ we consider the set of $\Omega_2$-words[2.6]

$$w(f,X,B) = \{w(f,X,m) : m \in B\}$$

We also use notation

$$w(f,X,B) = (w(f,X,m), m \in B)$$

Denote $w(f,X)$ **the set of $\Omega_2$-words of representation** $J(f,X)$. □

Theorem 2.6.7. *Endomorphism $R$ of representation*

$$f : A \longrightarrow\!\!\!\!\!*\;\; M$$

*generates the mapping of $\Omega_2$-words*

$$w[f,X,R] : w(f,X) \to w(f,X') \quad X \subset M \quad X' = R \circ X$$

*such that*

2.6.7.1: *If* $m \in X$, $m' = R \circ m$, *then*

$$w[f,X,R](m) = m'$$

2.6.7.2: *If*

$$m_1, ..., m_n \in w(f,X)$$

$$m'_1 = w[f,X,R](m_1) \quad ... \quad m'_n = w[f,X,R](m_n)$$

*then for operation* $\omega \in \Omega_2(n)$ *holds*

$$w[f,X,R](m_1...m_n\omega) = m'_1...m'_n\omega$$

2.6.7.3: *If*

$$m \in w(f,X) \quad m' = w[f,X,R](m) \quad a \in A$$

*then*

$$w[f,X,R](am) = am'$$

Proof. Statements 2.6.7.1, 2.6.7.2 are true by definition of the endomorpism $R$. Because $r = \mathrm{id}$, the statement 2.6.7.3 follows from the equation (2.2.4). □

[2.6] The expression $w(f,X,m)$ is a special case of the expression $w(f,X,B)$, namely

$$w(f,X,\{m\}) = \{w(f,X,m)\}$$

Remark 2.6.8. Let $R$ be endomorphism of representation $f$. Let

$$m \in J(f,X) \quad m' = R \circ m \quad X' = R \circ X$$

The theorem 2.6.7 states that $m' \in J_f(X')$. The theorem 2.6.7 also states that $\Omega_2$-word representing $m$ relative $X$ and $\Omega_2$-word representing $m'$ relative $X'$ are generated according to the same algorithm. This allows considering of the set of $\Omega_2$-words $w(f,X',m')$ as mapping

$$W(f,X,m) : X' \to w(f,X',m')$$
$$W(f,X,m)(X') = W(f,X,m) \circ X'$$

such that, if for certain endomorphism $R$

$$X' = R \circ X \quad m' = R \circ m$$

then

$$W(f,X,m) \circ X' = w(f,X',m') = m'$$

The mapping $W(f,X,m)$ is called **coordinates of element $m$ relative to set $X$**. Similarly, we consider coordinates of a set $B \subset J(f,X)$ relative to the set $X$

$$W(f,X,B) = \{W(f,X,m) : m \in B\} = (W(f,X,m), m \in B)$$

Denote $W(f,X)$ **the set of coordinates of representation** $J(f,X)$. □

Theorem 2.6.9. *There is a structure of $\Omega_2$-algebra on the set of coordinates $W(f,X)$.*

Proof. Let $\omega \in \Omega_2(n)$. Then for any $m_1$, ..., $m_n \in J(f,X)$, we assume

$$W(f,X,m_1)...W(f,X,m_n)\omega = W(f,X,m_1...m_n\omega) \tag{2.6.5}$$

According to the remark 2.6.8,

$$\begin{aligned}(W(f,X,m_1)...W(f,X,m_n)\omega) \circ X &= W(f,X,m_1...m_n\omega) \circ X \\ &= w(f,X,m_1...m_n\omega)\end{aligned} \tag{2.6.6}$$

follows from the equation (2.6.5). According to rule (2) of the definition 2.6.6, from the equation (2.6.6), it follows that

$$\begin{aligned}(W(f,X,m_1)...W(f,X,m_n)\omega) \circ X &= w(f,X,m_1)...w(f,X,m_n)\omega \\ &= (W(f,X,m_1) \circ X)...(W(f,X,m_n) \circ X)\omega\end{aligned} \tag{2.6.7}$$

From the equation (2.6.7), it follows that the operation $\omega$ defined by the equation (2.6.5) on the set of coordinates is defined properly. □

Theorem 2.6.10. *There exists the representation of $\Omega_1$-algebra $A$ in $\Omega_2$-algebra $W(f,X)$.*

Proof. Let $a \in A$. Then for any $m \in J(f,X)$ we assume

$$aW(f,X,m) = W(f,X,am) \tag{2.6.8}$$

According to the remark 2.6.8,

$$(aW(f,X,m)) \circ X = W(f,X,am) \circ X = w(f,X,am) \tag{2.6.9}$$

follows from the equation (2.6.8). According to rule (3) of the definition 2.6.6, from the equation (2.6.9), it follows that

$$(aW(f,X,m))\circ X = aw(f,X,m) = a(W(f,X,m)\circ X) \tag{2.6.10}$$

From the equation (2.6.10), it follows that the representation (2.6.8) of $\Omega_1$-algebra $A$ in $\Omega_2$-algebra $W(f,X)$ is defined properly. □

Theorem 2.6.11. *Let*

$$f : A \longrightarrow\!\!\!*\!\!\!\longrightarrow M$$

*be representation of $\Omega_1$-algebra $A$ in $\Omega_2$-algebra $M$. For given sets $X\subset M$, $X'\subset M$, let map*

$$R_1 : X \to X'$$

*agree with the structure of representation $f$, i. e.*

$$\begin{aligned}
&\omega\in\Omega_2(n)\ x_1,\ ...,\ x_n,\ x_1...x_n\omega\in X,\ R_1(x_1...x_n\omega)\in X'\\
=>\ &R_1(x_1...x_n\omega)=R_1(x_1)...R_1(x_n)\omega\\
&x\in X,\ a\in A,\ R_1(ax)\in X'\\
=>\ &R_1(ax)=aR_1(x)
\end{aligned}$$

*Consider the mapping of $\Omega_2$-words*

$$w[f,X,R_1] : w(f,X)\to w(f,X')$$

*that satisfies conditions 2.6.7.1, 2.6.7.2, 2.6.7.3 and such that*

$$x\in X => w[f,X,R_1](x)=R_1(x)$$

*There exists unique endomorphism of $\Omega_2$-algebra $M$*

$$R : M\to M$$

*defined by rule*

$$R\circ m = w[f,X,R_1](w(f,X,m))$$

*which is the morphism of representations $J(f,X)$ and $J(f,X')$.*

Proof. We prove the theorem by induction over complexity of $\Omega_2$-word.

If $w(f,X,m)=m$, then $m\in X$. According to condition 2.6.7.1,

$$R\circ m = w[f,X,R_1](w(f,X,m)) = w[f,X,R_1](m) = R_1(m)$$

Therefore, mappings $R$ and $R_1$ coinside on the set $X$, and the mapping $R$ agrees with structure of representation $f$.

Let $\omega\in\Omega_2(n)$. Let the mapping $R$ be defined for $m_1$, ..., $m_n\in J(f,X)$. Let

$$w_1 = w(f,X,m_1) \quad ... \quad w_n = w(f,X,m_n)$$

If $m=m_1...m_n\omega$, then according to rule (2) of definition 2.6.6,

$$w(f,X,m)=w_1...w_n\omega$$

According to condition 2.6.7.2,

$$\begin{aligned}
R\circ m &= w[f,X,R_1](w(f,X,m)) = w[f,X,R_1](w_1...w_n\omega)\\
&= w[f,X,R_1](w_1)...w[f,X,R_1](w_n)\omega\\
&= (R\circ m_1)...(R\circ m_n)\omega
\end{aligned}$$

Therefore, the mapping $R$ is endomorphism of $\Omega_2$-algebra $M$.

Let the mapping $R$ be defined for $m_1 \in J(f,X)$, $w_1 = w(f,X,m_1)$. Let $a \in A$. If $m = am_1$, then according to rule (3) of definition 2.6.6,

$$w(f,X,am_1) = aw_1$$

According to condition 2.6.7.3,

$$\begin{aligned}R \circ m &= w[f,X,R_1](w(f,X,m)) = w[f,X,R_1](aw_1)\\ &= aw[f,X,R_1](w_1) = aR \circ m_1\end{aligned}$$

From equation (2.2.4), it follows that the mapping $R$ is morphism of the representation $f$.

The statement that the endomorphism $R$ is unique and therefore this endomorphism is defined properly follows from the following argument. Let $m \in M$ have different $\Omega_2$-words relative the set $X$, for instance

$$m = x_1 ... x_n \omega = ax \tag{2.6.11}$$

Because $R$ is endomorphism of representation, then, from the equation (2.6.11), it follows that

$$R \circ m = R \circ (x_1 ... x_n \omega) = (R \circ x_1) ... (R \circ x_n)\omega = R \circ (ax) = aR \circ x \tag{2.6.12}$$

From the equation (2.6.12), it follows that

$$R \circ m = (R \circ x_1) ... (R \circ x_n)\omega = aR \circ x \tag{2.6.13}$$

From equations (2.6.11), (2.6.13), it follows that the equation (2.6.11) is preserved under the mapping. Therefore, the image of $m$ does not depend on the choice of coordinates. □

Remark 2.6.12. The theorem 2.6.11 is the theorem of extension of mapping. The only statement we know about the set $X$ is the statement that $X$ is generating set of the representation $f$. However, between the elements of the set $X$ there may be relationships generated by either operations of $\Omega_2$-algebra $M$, or by transformation of representation $f$. Therefore, any mapping of set $X$, in general, cannot be extended to an endomorphism of representation $f$.[2.7] However, if the mapping $R_1$ is coordinated with the structure of representation on the set $X$, then we can construct an extension of this mapping and this extension is endomorphism of representation $f$. □

Definition 2.6.13. Let $X$ be the generating set of the representation $f$. Let $R$ be the endomorphism of the representation $f$. The set of coordinates $W(f,X,R \circ X)$ is called **coordinates of endomorphism of representation**. □

Definition 2.6.14. Let $X$ be the generating set of the representation $f$. Let $R$ be the endomorphism of the representation $f$. Let $m \in M$. We define **superposition of coordinates** of the representation $f$ and the element $m$ as coordinates defined according to rule

$$W(f,X,m) \circ W(f,X,R \circ X) = W(f,X,R \circ m) \tag{2.6.14}$$

[2.7] In the theorem 2.7.7, requirements to generating set are more stringent. Therefore, the theorem 2.7.7 says about extension of arbitrary mapping. A more detailed analysis is given in the remark 2.7.9.

Let $Y \subset M$. We define superposition of coordinates of the representation $f$ and the set $Y$ according to rule

$$W(f,X,Y)\circ W(f,X,R\circ X)=(W(f,X,m)\circ W(f,X,R\circ X),m\in Y) \tag{2.6.15}$$

□

THEOREM 2.6.15. *Endomorphism $R$ of representation*

$$f:A \longrightarrow\!\!\!\!\!*\;\; M$$

*generates the mapping of coordinates of representation*

$$W[f,X,R]:W(f,X)\to W(f,X) \tag{2.6.16}$$

*such that*

$$\begin{aligned} W(f,X,m)\to W[f,X,R]\star W(f,X,m)=W(f,X,R\circ m)\\ W[f,X,R]\star W(f,X,m)=W(f,X,m)\circ W(f,X,R\circ X)\end{aligned} \tag{2.6.17}$$

PROOF. According to the remark 2.6.8, we consider equations (2.6.14), (2.6.16) relative to given generating set $X$. The word

$$W(f,X,m)\circ X=w(f,X,m) \tag{2.6.18}$$

corresponds to coordinates $W(f,X,m)$; the word

$$W(f,X,R\circ m)\circ X=w(f,X,R\circ m) \tag{2.6.19}$$

corresponds to coordinates $W(f,X,R\circ m)$. Therefore, in order to prove the theorem, it is sufficient to show that the mapping $W[f,X,R]$ corresponds to mapping $w[f,X,R]$. We prove this statement by induction over complexity of $\Omega_2$-word.

If $m\in X$, $m'=R\circ m$, then, according to equations (2.6.18), (2.6.19), mappings $W[f,X,R]$ and $w[f,X,R]$ are coordinated.

Let for $m_1$, ..., $m_n\in X$ mappings $W[f,X,R]$ and $w[f,X,R]$ be coordinated. Let $\omega\in\Omega_2(n)$. According to the theorem 2.6.9

$$W(f,X,m_1...m_n\omega)=W(f,X,m_1)...W(f,X,m_n)\omega \tag{2.6.20}$$

Because $R$ is endomorphism of $\Omega_2$-algebra $M$, then from the equation (2.6.20), it follows that

$$\begin{aligned} W(f,X,R\circ(m_1...m_n\omega))&=W(f,X,(R\circ m_1)...(R\circ m_n)\omega)\\ &=W(f,X,R\circ m_1)...W(f,X,R\circ m_n)\omega\end{aligned} \tag{2.6.21}$$

From equations (2.6.20), (2.6.21) and the statement of induction, it follows that the mappings $W[f,X,R]$ and $w[f,X,R]$ are coordinated for $m=m_1...m_n\omega$.

Let for $m_1\in M$ mappings $W[f,X,R]$ and $w[f,X,R]$ are coordinated. Let $a\in A$. According to the theorem 2.6.10

$$W(f,X,am_1)=aW(f,X,m_1) \tag{2.6.22}$$

Because $R$ is endomorphism of representation $f$, then, from the equation (2.6.22), it follows that

$$W(f,X,R\circ(am_1))=W(f,X,a\,R\circ m_1)=aW(f,X,R\circ m_1) \tag{2.6.23}$$

From equations (2.6.22), (2.6.23) and the statement of induction, it follows that mappings $W[f,X,R]$ and $w[f,X,R]$ are coordinated for $m=am_1$. □

Corollary 2.6.16. *Let $X$ be the generating set of the representation $f$. Let $R$ be the endomorphism of the representation $f$. The mapping $W[f,X,R]$ is endomorphism of representation of $\Omega_1$-algebra $A$ in $\Omega_2$-algebra $W(f,X)$.* □

Hereinafter we will identify mapping $W[f,X,R]$ and the set of coordinates $W(f,X,R\circ X)$.

Theorem 2.6.17. *Let $X$ be the generating set of the representation $f$. Let $R$ be the endomorphism of the representation $f$. Let $Y\subset M$. Then*

$$W(f,X,Y)\circ W(f,X,R\circ X)=W(f,X,R\circ Y) \tag{2.6.24}$$

$$W[f,X,R]\star W(f,X,Y)=W(f,X,R\circ Y) \tag{2.6.25}$$

Proof. The equation (2.6.24) follows from the equation

$$R\circ Y=(R\circ m,m\in Y)$$

as well from equations (2.6.14), (2.6.15). The equation (2.6.25) is corollary of equations (2.6.24), (2.6.17). □

Theorem 2.6.18. *Let $X$ be the generating set of the representation $f$. Let $R$, $S$ be the endomorphisms of the representation $f$. Then*

$$W(f,X,S\circ X)\circ W(f,X,R\circ X)=W(f,X,R\circ S\circ X) \tag{2.6.26}$$

$$W[f,X,R]\star W[f,X,S]=W[f,X,R\circ S] \tag{2.6.27}$$

Proof. The equation (2.6.26) follows from the equation (2.6.24), if we assume $Y=S\circ X$. The equation (2.6.27) follows from the equation (2.6.26) and chain of equations

$$\begin{aligned}&(W[f,X,R]\star W[f,X,S])\star W(f,X,Y)\\ =&W[f,X,R]\star(W[f,X,S]\star W(f,X,Y))\\ =&(W(f,X,Y)\circ W(f,X,S\circ X))\circ W(f,X,R\circ X)\\ =&W(f,X,Y)\circ(W(f,X,S\circ X)\circ W(f,X,R\circ X))\\ =&W(f,X,Y)\circ W(f,X,R\circ S\circ X)\\ =&W(f,X,R\circ S)\star W(f,X,Y)\end{aligned}$$

□

The concept of superposition of coordinates is very simple and resembles a kind of Turing machine. If element $m\in M$ has form either

$$m=m_1...m_n\omega$$

or

$$m=am_1$$

then we are looking for coordinates of elements $m_i$ to substitute them in an appropriate expression. As soon as an element $m\in M$ belongs to the generating set of $\Omega_2$-algebra $M$, we choose coordinates of the corresponding element of the second factor. Therefore, we require that the second factor in the superposition has been the set of coordinates of the image of the generating set $X$.

We can generalize the definition of the superposition of coordinates and assume that one of the factors is a set of $\Omega_2$-words. Accordingly, the definition of the superposition of coordinates has the form

$$w(f,X,Y)\circ W(f,X,R\circ X)=W(f,X,Y)\circ w(f,X,R\circ X)=w(f,X,R\circ Y)$$

The following forms of writing an image of the set $Y$ under endomorphism $R$ are equivalent.

$$R\circ Y=W(f,X,Y)\circ(R\circ X)=W(f,X,Y)\circ(W(f,X,R\circ X)\circ X) \tag{2.6.28}$$

From equations (2.6.24), (2.6.28), it follows that

$$(W(f,X,Y)\circ W(f,X,R\circ X))\circ X=W(f,X,Y)\circ(W(f,X,R\circ X)\circ X) \tag{2.6.29}$$

The equation (2.6.29) is associative law for composition and allows us to write expression

$$W(f,X,Y)\circ W(f,X,R\circ X)\circ X$$

without brackets.

Consider equation (2.6.26), where we see change in the order of endomorphisms in a superposition of coordinates. This equation also follows from the chain of equations, where we can immediately see when order of endomorphisms changes

$$\begin{aligned}&W(f,X,m)\circ W(f,X,R\circ S\circ X)\circ X\\=&(R\circ S)\circ m=R\circ(S\circ m)\\=&R\circ((W(f,X,m)\circ W(f,X,S\circ X))\circ X)\\=&W(f,X,m)\circ W(f,X,S\circ X)\circ W(f,X,R\circ X)\circ X\end{aligned} \tag{2.6.30}$$

From the equation (2.6.30), it follows that coordinates of endomorphism act over coordinates of element of $\Omega_2$-algebra $M$ from the right.

Definition 2.6.19. Let $X\subset M$ be generating set of representation

$$f:A\longrightarrow\!\!\!\!\!*\ \ M$$

Let the mapping

$$H:M\to M$$

be endomorphism of the representation $f$. Let the set $X'=H\circ X$ be the image of the set $X$ under the mapping $H$. Endomorphism $H$ of representation $f$ is called **regular on generating set** $X$, if the set $X'$ is the generating set of representation $f$. Otherwise, endomorphism $H$ of representation $f$ is called **singular on generating set** $X$, □

Definition 2.6.20. Endomorphism of representation $f$ is called **regular**, if it is regular on every generating set. □

Theorem 2.6.21. *Automorphism $R$ of representation*

$$f:A\longrightarrow\!\!\!\!\!*\ \ M$$

*is regular endomorphism.*

Proof. Let $X$ be generating set of representation $f$. Let $X' = R \circ X$.

According to theorem 2.6.7 endomorphism $R$ forms the map of $\Omega_2$-words $w[f, X, R]$.

Let $m' \in M$. Since $R$ is automorphism, then there exists $m \in M$, $R \circ m = m'$. According to definition 2.6.6, $w(f, X, m)$ is $\Omega_2$-word, representing $m$ relative to generating set $X$. According to theorem 2.6.7, $w(f, X', m')$ is $\Omega_2$-word, representing of $m'$ relative to generating set $X'$

$$w(f, X', m') = w[f, X, R](w(f, X, m))$$

Therefore, $X'$ is generating set of representation $f$. According to definition 2.6.20, automorphism $R$ is regular. □

## 2.7. Basis of representation

Definition 2.7.1. If the set $X \subset M$ is generating set of representation $f$, then any set $Y$, $X \subset Y \subset M$ also is generating set of representation $f$. If there exists minimal set $X$ generating the representation $f$, then the set $X$ is called **basis of representation** $f$. □

Theorem 2.7.2. *The generating set $X$ of representation $f$ is basis iff for any $m \in X$ the set $X \setminus \{m\}$ is not generating set of representation $f$.*

Proof. Let $X$ be generating set of representation $f$. Assume that for some $m \in X$ there exist $\Omega_2$-word

$$w = w(f, X \setminus \{m\}, m) \tag{2.7.1}$$

Consider element $m'$ such that it has $\Omega_2$-word

$$w' = w(f, X, m') \tag{2.7.2}$$

that depends on $m$. According to the definition 2.6.6, any occurrence of $m$ into $\Omega_2$-word $w'$ can be substituted by the $\Omega_2$-word $w$. Therefore, the $\Omega_2$-word $w'$ does not depend on $m$, and the set $X \setminus \{m\}$ is generating set of representation $f$. Therefore, $X$ is not basis of representation $f$. □

Remark 2.7.3. The proof of the theorem 2.7.2 gives us effective method for constructing the basis of the representation $f$. Choosing an arbitrary generating set, step by step, we remove from set those elements which have coordinates relative to other elements of the set. If the generating set of the representation is infinite, then this construction may not have the last step. If the representation has finite generating set, then we need a finite number of steps to construct a basis of this representation.

As noted by Paul Cohn in [13], p. 82, 83, the representation may have inequivalent bases. For instance, the cyclic group of order six has bases $\{a\}$ and $\{a^2, a^3\}$ which we cannot map one into another by endomorphism of the representation. □

Remark 2.7.4. We write a basis also in following form

$$X = (x, x \in X)$$

If basis is finite, then we also use notation

$$X = (x_i, i \in I) = (x_1, ..., x_n)$$

□

Remark 2.7.5. We introduced $\Omega_2$-word of $x \in M$ relative generating set $X$ in the definition 2.6.6. From the theorem 2.7.2, it follows that if the generating set $X$ is not a basis, then a choice of $\Omega_2$-word relative generating set $X$ is ambiguous. However, even if the generating set $X$ is a basis, then a representation of $m \in M$ in form of $\Omega_2$-word is ambiguous. If $m_1$, ..., $m_n$ are $\Omega_2$-words, $\omega \in \Omega_2(n)$ and $a \in A$, then[2.8]

$$a(m_1...m_n\omega) = (am_1)...(am_n)\omega \tag{2.7.3}$$

It is possible that there exist equations related to specific character of representation. For instance, if $\omega$ is operation of $\Omega_1$-algebra $A$ and operation of $\Omega_2$-algebra $M$, then we require that $\Omega_2$-words $a_1...a_n\omega x$ and $a_1x...a_nx\omega$ describe the same element of $\Omega_2$-algebra $M$.[2.9]

In addition to the above equations in $\Omega_2$-algebra there may be relations of the form

$$w_1(f, X, m) = w_2(f, X, m) \quad m \notin X \tag{2.7.4}$$

The feature of the equation (2.7.4) is that this equation cannot be reduced.[2.10]

On the set of $\Omega_2$-words $w(f, X)$, above equations determine **equivalence** $\rho(f)$ **generated by representation** $f$. It is evident that for any $m \in M$ the choice of appropriate $\Omega_2$-word is unique up to equivalence relations $\rho(f)$. However, if during the construction, we obtain the equality of two $\Omega_2$-word relative to given basis, then we can say without worrying about the equivalence $\rho(f)$ that these $\Omega_2$-words are equal.

A similar remark concerns the mapping $W(f, X, m)$ defined in the remark 2.6.8.[2.11] □

---

[2.8] For instance, let $\{\overline{e}_1, \overline{e}_2\}$ be the basis of vector space over field $k$. The equation (2.7.3) has the form of distributive law

$$a(b^1\overline{e}_1 + b^2\overline{e}_2) = (ab^1)\overline{e}_1 + (ab^2)\overline{e}_2$$

[2.9] For vector space, this requirement has the form of distributive law

$$(a + b)\overline{e}_1 = a\overline{e}_1 + b\overline{e}_1$$

[2.10] See for instance sections 4.7.2, 4.7.3.

[2.11] If vector space has finite basis, then we represent the basis as matrix

$$\overline{\overline{e}} = \begin{pmatrix} \overline{e}_1 & ... & \overline{e}_2 \end{pmatrix}$$

We present the mapping $W(f, \overline{\overline{e}}, \overline{v})$ as matrix

$$W(f, \overline{\overline{e}}, \overline{v}) = \begin{pmatrix} v^1 \\ ... \\ v^n \end{pmatrix}$$

Then

$$W(f, \overline{\overline{e}}, \overline{v})(\overline{\overline{e}}') = W(f, \overline{\overline{e}}, \overline{v}) \begin{pmatrix} \overline{e}'_1 & ... & \overline{e}'_n \end{pmatrix} = \begin{pmatrix} v^1 \\ ... \\ v^n \end{pmatrix} \begin{pmatrix} \overline{e}'_1 & ... & \overline{e}'_n \end{pmatrix}$$

has form of matrix product.

Theorem 2.7.6. *Automorphism of the representation $f$ maps a basis of the representation $f$ into basis.*

Proof. Let the mapping $R$ be automorphism of the representation $f$. Let the set $X$ be a basis of the representation $f$. Let $X'=R\circ X$.

Assume that the set $X'$ is not basis. According to the theorem 2.7.2 there exists such $m'\in X'$ that $X'\setminus\{x'\}$ is generating set of the representation $f$. According to the theorem 2.4.5, the mapping $R^{-1}$ is automorphism of the representation $f$. According to the theorem 2.6.21 and definition 2.6.20, the set $X\setminus\{m\}$ is generating set of the representation $f$. The contradiction completes the proof of the theorem. □

Theorem 2.7.7. *Let $X$ be the basis of the representation $f$. Let*

$$R_1:X\to X'$$

*be arbitrary mapping of the set $X$. Consider the mapping of $\Omega_2$-words*

$$w[f,X,R_1]:w(f,X)\to w(f,X')$$

*that satisfies conditions 2.6.7.1, 2.6.7.2, 2.6.7.3 and such that*

$$x\in X => w[f,X,R_1](x)=R_1(x)$$

*There exists unique endomorphism of representation $f$*[2.12]

$$R:M\to M$$

*defined by rule*

$$R\circ m=w[f,X,R_1](w(f,X,m))$$

Proof. The statement of theorem is corollary theorems 2.6.7, 2.6.11. □

Corollary 2.7.8. *Let $X$, $X'$ be the bases of the representation $f$. Let $R$ be the automorphism of the representation $f$ such that $X'=R\circ X$. Automorphism $R$ is uniquely defined.* □

Remark 2.7.9. The theorem 2.7.7, as well as the theorem 2.6.11, is the theorem of extension of mapping. However in this theorem, $X$ is not arbitrary generating set of the representation, but basis. According to remark 2.7.3, we cannot determine coordinates of any element of basis through the remaining elements of the same basis. Therefore, we do not need to coordinate the mapping of the basis with representation. □

Theorem 2.7.10. *The set of coordinates $W(f,X,X)$ corresponds to identity transformation*

$$W[f,X,E]=W(f,X,X)$$

Proof. The statement of the theorem follows from the equation

$$m=W(f,X,m)\circ X=W(f,X,m)\circ W(f,X,X)\circ X$$

□

[2.12] This statement is similar to the theorem [1]-4.1, p. 135.

Theorem 2.7.11. *Let $W(f,X,R\circ X)$ be the set of coordinates of automorphism $R$. There exists set of coordinates $W(f,R\circ X,X)$, corresponding to automorphism $R^{-1}$. The set of coordinates $W(f,R\circ X,X)$ satisfy to equation*

$$W(f,X,R\circ X)\circ W(f,R\circ X,X)=W(f,X,X) \tag{2.7.5}$$
$$W[f,X,R^{-1}]=W[f,X,R]^{-1}=W(f,R\circ X,X)$$

Proof. Since $R$ is automorphism of the representation $f$, then, according to the theorem 2.7.6, the set $R\circ X$ is a basis of the representation $f$. Therefore, there exists the set of coordinates $W(f,R\circ X,X)$. The equation (2.7.5) follows from the chain of equations

$$\begin{aligned}W(f,X,R\circ X)\circ W(f,R\circ X,X)&=W(f,X,R\circ X)\circ W(f,X,R^{-1}\circ X)\\&=W(f,X,R^{-1}\circ R\circ X)=W(f,X,X)\end{aligned}$$

□

Remark 2.7.12. In $\Omega_2$-algebra $M$ there is no universal algorithm for determining the set of coordinates $W(f,R\circ X,X)$ for given set $W(f,X,R\circ X)$.[2.13] We assume that in the theorem 2.7.11 this algorithm is given implicitly. It is evident also that the set of $\Omega_2$-words

$$W(f,X,R\circ X)\circ W(f,R\circ X,X)\circ X \tag{2.7.6}$$

in general, does not coincide with the set of $\Omega_2$-words

$$W(f,X,X)\circ X \tag{2.7.7}$$

The theorem 2.7.11 states that sets of $\Omega_2$-words (2.7.6) and (2.7.7) coincide up to equivalence generated by the representation $f$. □

Theorem 2.7.13. *The group of automorphisms $G(f)$ of effective representation $f$ in $\Omega_2$-algebra $M$ generates effective representation in $\Omega_2$-algebra $M$.*

Proof. From the corollary 2.7.8, it follows that if automorphism $R$ maps a basis $X$ into a basis $X'$, then the set of coordinates $W(f,X,X')$ uniquely determines an automorphism $R$. From the theorem 2.6.15, it follows that the set of coordinates $W(f,X,X')$ determines the mapping of coordinates relative to the basis $X$ under automorphism of the representation $f$. From the equation (2.6.28), it follows that automorphism $R$ acts from the right on elements of $\Omega_2$-algebra $M$. From the equation (2.6.26), it follows that the representation of group is right-side representation. According to the theorem 2.7.10, the set of coordinates $W(f,X,X)$ corresponds to identity transformation. From the theorem 2.7.11, it follows that the set of coordinates $W(f,R\circ X,X)$ corresponds to transformation, inverse to transformation $W(f,X,R\circ X)$. □

[2.13] In vector space, the matrix of numbers corresponds to linear transformation. Accordingly, the inverse matrix corresponds to inverse transformation.

## 2.8. Examples of Basis of Representation of Universal Algebra

**2.8.1. Vector Space.** Consider the vector space $\overline{V}$ over the field $F$. Given the set of vectors $\overline{e}_1, ..., \overline{e}_n,$ according to algorithm of construction of coordinates over vector space, coordinates include such elements as $\overline{e}_1+\overline{e}_2$ and $a\overline{e}_1$. Recursively using rules, contained in the definition 2.6.6, we conclude that the set of vectors $\overline{e}_1, ..., \overline{e}_n,$ generates the set of linear combinations

$$a^1\overline{e}_1 + ... + a^n\overline{e}_n$$

According to the theorem 2.7.2, the set of vectors $\overline{e}_1, ..., \overline{e}_n,$ is a basis if for any $i$, $i = 1, ..., n$, vector $\overline{e}_i$ is not linear combination of other vectors. This requirement is equivalent to the requirement of linear independence of vectors.

**2.8.2. Representation of Group on the Set.** Consider the representation from the example 2.5.2. We can consider the set $M$ as union of orbits of the representation of the group $G$. We can select for basis of the representation the set of points such that one and only one point belongs to each orbit. If $X$ is the basis of representation, $A \in X$, $g \in G$, then $\Omega_2$-word has form $A + g$. Since there is no operations on the set $M$, then there is no $\Omega_2$-word containing different elements of the basis. If representation of group $G$ is single transitive, then basis of representation consists of one point. Any point of the set $M$ can be such point.

CHAPTER 3

# Representation of Group

## 3.1. Representation of Group

Group is among few algebras that allow somebody to consider the product of transformations of the $\Omega$-algebra $M$ in such a way that, if transformations belong to the representation, then their product also belongs to the representation. We should remember that order of maps in product depends on order of maps on diagram and how these maps act over elements of the set (from left or from right).

DEFINITION 3.1.1. Let ${}^*M$ be a group with product

$$(f\circ g)x=f(gx)$$

and $\delta$ be unit of group ${}^\star M$. Let $G$ be group. We call a homomorphism of group

$$f:G\to{}^\star M \tag{3.1.1}$$

**left-side representation of group** $G$ or $G*$**-representation** in $\Omega$-algebra $M$ if the map $f$ holds

$$f(ab)u=f(a)(f(b)u) \tag{3.1.2}$$

□

REMARK 3.1.2. Since the map (3.1.1) is homomorhism, then

$$f(ab)u=(f(a)f(b))u \tag{3.1.3}$$

We use here convention

$$f(a)f(b)=f(a)\circ f(b)$$

Thus, the idea of representation of group is that we multiply elements of group in the same order as we multiply transformations of representation. From equations (3.1.2) and (3.1.3) it follows

$$(f(a)f(b))u=f(a)(f(b)u) \tag{3.1.4}$$

Equation (3.1.4) together with associativity of product of transformations expresses **associative law** for $G*$-representation. This allows writing of equation (3.1.4) without using of brackets

$$f(ab)u=f(a)f(b)u$$

□

DEFINITION 3.1.3. Let $M^\star$ be a group with product

$$x(f\circ g)=(xf)g$$

and $\delta$ be unit of group $M^\star$. Let $G$ be group. We call a homomorphism of group

$$f:G\to M^* \tag{3.1.5}$$

**right-side representation of group** $G$ or $*G$**-representation** in $\Omega$-algebra $M$ if the map $f$ holds

$$uf(ab)=(uf(a))f(b) \tag{3.1.6}$$

□

REMARK 3.1.4. Since the map (3.1.5) is homomorhism, then

$$uf(ab)=u(f(a)f(b)) \tag{3.1.7}$$

From equations (3.1.6) and (3.1.7) it follows

$$u(f(a)f(b))=(uf(a))f(b) \tag{3.1.8}$$

Equation (3.1.8) together with associativity of product of transformations expresses **associative law** for $*G$-representation. This allows writing of equation (3.1.8) without using of brackets

$$uf(ab)=uf(a)f(b)$$

□

DEFINITION 3.1.5. We call the transformation

$$t:M\to M$$

**nonsingular transformation**, if there exists inverse map. □

THEOREM 3.1.6. *For any* $g\in G$ *transformation is nonsingular and satisfies equation*

$$f(g^{-1})=f(g)^{-1} \tag{3.1.9}$$

PROOF. Since (3.1.2) and

$$f(e)=\delta$$

we have

$$u=\delta(u)=f(gg^{-1})(u)=f(g)(f(g^{-1})(u))$$

This completes the proof. □

THEOREM 3.1.7. *The group operation determines two different representations on the group:*

- *The* **left shift** $t_\star$

$$\begin{aligned}b'&=t_\star(a)b=ab\\ b'&=t_\star(a)(b)=ab\end{aligned} \tag{3.1.10}$$

*is* $G*$*-representation on the set*[3.1] $G$

$$t_\star(ab)=t_\star(a)\circ t_\star(b) \tag{3.1.11}$$

- *The* **right shift** ${}_\star t$

$$\begin{aligned}b'&=b\ {}_\star t(a)=ba\\ b'&={}_\star t(a)(b)=ba\end{aligned} \tag{3.1.12}$$

*is* $*G$*-representation on the set* $G$

$${}_\star t(ab)={}_\star t(a)\circ\ {}_\star t(b) \tag{3.1.13}$$

[3.1]Left shift is not a representation of group in a group, because the transformation $t_\star$ is not a homomorphism of group. Similar remark is true for right shift.

PROOF. Equation (3.1.11) follows from associativity of product

$$t_\star(ab)c = (ab)c = a(bc) = t_\star(a)(t_\star(b)c) = (t_\star(a) \circ t_\star(b))c$$

In a similar manner we prove the equation (3.1.13). □

DEFINITION 3.1.8. Let $G$ be group. Let $f$ be $G*$-representation in $\Omega$-algebra $M$. For any $v \in M$ we define **orbit of representation of the group** $G$ as set

$$f(G)v = \{w = f(g)v : g \in G\}$$

□

Since $f(e) = \delta$ we have $v \in f(G)v$.

THEOREM 3.1.9. *Suppose*

$$v \in f(G)u \tag{3.1.14}$$

*Then*

$$f(G)u = f(G)v$$

PROOF. From (3.1.14) it follows that there exists $a \in G$ such that

$$v = f(a)u \tag{3.1.15}$$

Suppose $w \in f(G)v$. Then there exists $b \in G$ such that

$$w = f(b)v \tag{3.1.16}$$

If we substitute (3.1.15) into (3.1.16) we get

$$w = f(b)(f(a)u) \tag{3.1.17}$$

Since (3.1.2), we see that from (3.1.17) it follows that $w \in f(G)u$. Thus

$$f(G)v \subseteq f(G)u$$

Since (3.1.9), we see that from (3.1.15) it follows that

$$u = f(a)^{-1}v = f(a^{-1})v \tag{3.1.18}$$

From (3.1.18) it follows that $u \in f(G)v$ and therefore

$$f(G)u \subseteq f(G)v$$

This completes the proof. □

Thus, $G*$-representation $f$ in $\Omega$-algebra $M$ forms equivalence $S$ and the orbit $f(G)u$ is equivalence class. We will use notation $M/f(G)$ for quotient set $M/S$ and this set is called **space of orbits of $G*$-representation** $f$.

THEOREM 3.1.10. *Suppose $f_1$ is $G*$-representation in $\Omega$-algebra $M_1$ and $f_2$ is $G*$-representation in $\Omega$-algebra $M_2$. Then we introduce* **direct product of $G*$-representations** $f_1$ *and* $f_2$

$$f = f_1 \times f_2 : G \to M_1 \otimes M_2$$

$$f(g) = (f_1(g), f_2(g))$$

Proof. To show that $f$ is a representation, it is enough to prove that $f$ satisfies the definition 3.1.1.

$$f(e) = (f_1(e), f_2(e)) = (\delta_1, \delta_2) = \delta$$

$$\begin{aligned} f(ab)u &= (f_1(ab)u_1, f_2(ab)u_2) \\ &= (f_1(a)(f_1(b)u_1), f_2(a)(f_2(b)u_2)) \\ &= f(a)(f_1(b)u_1, f_2(b)u_2) \\ &= f(a)(f(b)u) \end{aligned}$$

□

## 3.2. Single Transitive Right-Side Representation of Group

Definition 3.2.1. We call **kernel of inefficiency of $G*$-representation** a set

$$K_f = \{g \in G : f(g) = \delta\}$$

□

Theorem 3.2.2. *A kernel of inefficiency of $G*$-representation is a subgroup of the group $G$.*

Proof. Assume $f(a_1) = \delta$ and $f(a_2) = \delta$. Then

$$f(a_1a_2)u = f(a_1)(f(a_2)u) = u$$
$$f(a^{-1}) = f^{-1}(a) = \delta$$

□

Theorem 3.2.3. *$G*$-representation is* **effective** *iff kernel of inefficiency $K_f = \{e\}$.*

Proof. Statement is corollary of definitions 2.1.6 and 3.2.1 and of the theorem 3.2.2. □

If an action is not effective we can switch to an effective one by changing group $G_1 = G|K_f$ using factorization by the kernel of inefficiency. This means that we can study only an effective action.

Definition 3.2.4. Consider $G*$-representation $f$ in $\Omega$-algebra $M$. A **little group** or **stability group** of $x \in M$ is the set

$$= \{g \in G : f(g)x = x\}$$

$G*$-representation $f$ is said to be **free**, if for any $x \in M$ stability group $G_x = \{e\}$. □

Theorem 3.2.5. *Given free $G*$-representation $f$ in the $\Omega$-algebra $A$, there exist $1-1$ correspondence between orbits of representation, as well between orbit of representation and group $G$.*

Proof. Given $a \in A$ there exist $g_1$, $g_2 \in G$

$$f(g_1)a = f(g_2)a \tag{3.2.1}$$

We multiply both parts of equation (3.2.1) by $f(g_1^{-1})$

$$a = f(g_1^{-1})f(g_2)a$$

Since the representation is free, $g_1 = g_2$. Since we established $1-1$ correspondence between orbit and group $G$, we proved the statement of the theorem. □

DEFINITION 3.2.6. We call a space $V$ **homogeneous space of group** $G$ if we have single transitive $G*$-representation on $V$. □

THEOREM 3.2.7. *If we define a single transitive representation $f$ of the group $G$ on the $\Omega$-algebra $A$ then we can uniquely define coordinates on $A$ using coordinates on the group $G$.*

*If $f$ is left-side representation then $f(a)$ is equivalent to the left shift $t_\star(a)$ on the group $G$. If $f$ is right-side representation then $f(a)$ is equivalent to the right shift ${}_\star t(a)$ on the group $G$.*

PROOF. We select a point $v \in A$ and define coordinates of a point $w \in A$ as coordinates of $a \in G$ such that $w = f(a)v$. Coordinates defined this way are unique up to choice of an initial point $v \in A$ because the action is effective.

If $f$ is left-side representation, we will use the notation

$$f(a)v = av$$

Because the notation

$$f(a)(f(b)v) = a(bv) = (ab)v = f(ab)v$$

is compatible with the group structure we see that left-side representation $f$ is equivalent to the left shift.

If $f$ is right-side representation, we will use the notation

$$vf(a) = va$$

Because the notation

$$(vf(b))f(a) = (vb)a = v(ba) = vf(ba)$$

is compatible with the group structure we see that right-side representation $f$ is equivalent to the right shift. □

REMARK 3.2.8. We will write effective $G*$-representation as

$$v' = t_\star(a)v = av$$

Orbit of this representation is

$$Gv = t_\star(G)v$$

We will use notation $M/t_\star(G)$ for the space of orbits of effective $G*$-representation. □

REMARK 3.2.9. We will write effective $*G$-representation as

$$v' = v\ {}_\star t(a) = va$$

Orbit of this representation is

$$vG = v\ {}_\star t(G)$$

We will use notation $M/{}_\star t(G)$ for the space of orbits of effective $*G$-representation. □

THEOREM 3.2.10. *Free $G**$-representation is effective. Free $G**$-representation $f$ in $\Omega$-algebra $M$ is single transitive representation on orbit.*

PROOF. The statement of theorem is the corollary of definition 3.2.4. □

THEOREM 3.2.11. *Left and right shifts on group $G$ are commuting.*

PROOF. This is the consequence of the associativity on the group $G$

$$(t_\star(a)\circ {}_\star t(b))c = a(cb) = (ac)b = ({}_\star t(b)\circ t_\star(a))c$$

□

Theorem 3.2.11 can be phrased n the following way.

THEOREM 3.2.12. *Let $G$ be group. For any $a\in G$, the map $t_\star(a)$ is automorphism of representation ${}_\star t$.*

PROOF. According to theorem 3.2.11

$$t_\star(a)\circ {}_\star t(b) = {}_\star t(b)\circ t_\star(a) \tag{3.2.2}$$

Equation (3.2.2) coincides with equation (2.2.3) from definition 2.2.2 when $r = id$, $R = t_\star(a)$. □

THEOREM 3.2.13. *Let $G*$-representation $f$ on $\Omega$-algebra $M$ be single transitive. Then we can uniquely define a single transitive $*G$-representation $h$ on $\Omega$-algebra $M$ such that diagram*

$$\begin{array}{ccc}
M & \xrightarrow{\ h(a)\ } & M \\
{\scriptstyle f(b)}\downarrow & & \downarrow{\scriptstyle f(b)} \\
M & \xrightarrow[\ h(a)\ ]{} & M
\end{array}$$

*is commutative for any $a, b\in G$.*[3.2]

PROOF. We use group coordinates for points $v\in M$. Then according to theorem 3.2.7 we can write the left shift $t_\star(a)$ instead of the transformation $f(a)$.

Let $v_0, v\in M$. Then we can find one and only one $a\in G$ such that

$$v = v_0a = v_0\ {}_\star t(a)$$

We assume

$$h(a) = {}_\star t(a)$$

For some $b\in G$ we have

$$w_0 = f(b)v_0 = t_\star(b)v_0 \qquad w = f(b)v = t_\star(b)v$$

According to the theorem 3.2.11, the diagram

$$\begin{array}{ccc}
v_0 & \xrightarrow{\ h(a)={}_\star t(a)\ } & v \\
{\scriptstyle f(b)=t_\star(b)}\downarrow & & \downarrow{\scriptstyle f(b)=t_\star(b)} \\
w_0 & \xrightarrow[\ h(a)={}_\star t(a)\ ]{} & w
\end{array} \tag{3.2.3}$$

is commutative.

Changing $b$ we get that $w_0$ is an arbitrary point of $M$.

We see from the diagram that if $v_0 = v$ then $w_0 = w$ and therefore $h(e) = \delta$. On other hand if $v_0 \ne v$ then $w_0\ne w$ because the $G*$-representation $f$ is single transitive. Therefore the $*G$-representation $h$ is effective.

[3.2] You can see this statement in [3].

In the same way we can show that for given $w_0$ we can find $a$ such that $w = h(a)w_0$. Therefore the $*G$-representation $h$ is single transitive.

In general the product of transformations of the $G*$-representation $f$ is not commutative and therefore the $*G$-representation $h$ is different from the $G*$-representation $f$. In the same way we can create a $G*$-representation $f$ using the $*G$-representation $h$. □

Representations $f$ and $h$ are called **twin representations of the group** $G$.

Remark 3.2.14. It is clear that transformations $t_\star(a)$ and ${}_\star t(a)$ are different until the group $G$ is nonabelian. However they both are maps onto. Theorem 3.2.13 states that if both right and left shift presentations exist on the set $M$, then we can define two commuting representations on the set $M$. The right shift or the left shift only cannot represent both types of representation. To understand why it is so let us change diagram (3.2.3) and assume $h(a)v_0 = t_\star(a)v_0 = v$ instead of $h(a)v_0 = v_0{}_\star t(a) = v$ and let us see what expression $h(a)$ has at the point $w_0$. The diagram

$$\begin{array}{ccc}
v_0 & \xrightarrow{h(a)=t_\star(a)} & v \\
\Big\downarrow {\scriptstyle f(b)=t_\star(b)} & & \Big\downarrow {\scriptstyle f(b)=t_\star(b)} \\
w_0 & \xrightarrow[h(a)]{} & w
\end{array}$$

is equivalent to the diagram

$$\begin{array}{ccc}
v_0 & \xrightarrow{h(a)=t_\star(a)} & v \\
\Big\uparrow {\scriptstyle f^{-1}(b)=t_\star(b^{-1})} & & \Big\downarrow {\scriptstyle f(b)=t_\star(b)} \\
w_0 & \xrightarrow[h(a)]{} & w
\end{array}$$

and we have $w = bv = bav_0 = bab^{-1}w_0$. Therefore

$$h(a)w_0 = (bab^{-1})w_0$$

We see that the representation of $h$ depends on its argument. □

Theorem 3.2.15. *Let $f$ and $h$ be twin representations of the group $G$. For any $a \in G$ the map $h(a)$ is automorphism of representation $f$.*

Proof. The statement of theorem is corollary of theorems 3.2.12 and 3.2.13. □

Remark 3.2.16. Is there a morphism of representations from $t_\star$ to $t_\star$ different from automorphism $(\mathrm{id}, {}_\star t(a))$? If we assume

$$r(g) = cgc^{-1}$$
$$R(a)(m) = cmac^{-1}$$

then it is easy to see that the map $(r, R(a))$ is morphism of the representations from $t_\star$ to $t_\star$. However this map is not automorphism of the representation $t_\star$, because $r \neq \mathrm{id}$. □

CHAPTER 4

# Tower of Representations of Universal Algebras

## 4.1. Tower of Representations of Universal Algebras

Definition 4.1.1. Consider set of $\Omega_k$-algebras $A_k$, $k = 1, ..., n$. Assume $\overline{A} = (A_1, ..., A_n)$. Assume $\overline{f} = (f_{1,2}, ..., f_{n-1,n})$. Set of representations $f_{k,k+1}$, $k = 1, ..., n$, of $\Omega_k$-algebra $A_k$ in $\Omega_{k+1}$-algebra $A_{k+1}$ is called **tower $(\overline{A}, \overline{f})$ of representations of $\overline{\Omega}$-algebras**. □

Consider following diagram for the purposes of illustration of definition 4.1.1

(4.1.1)

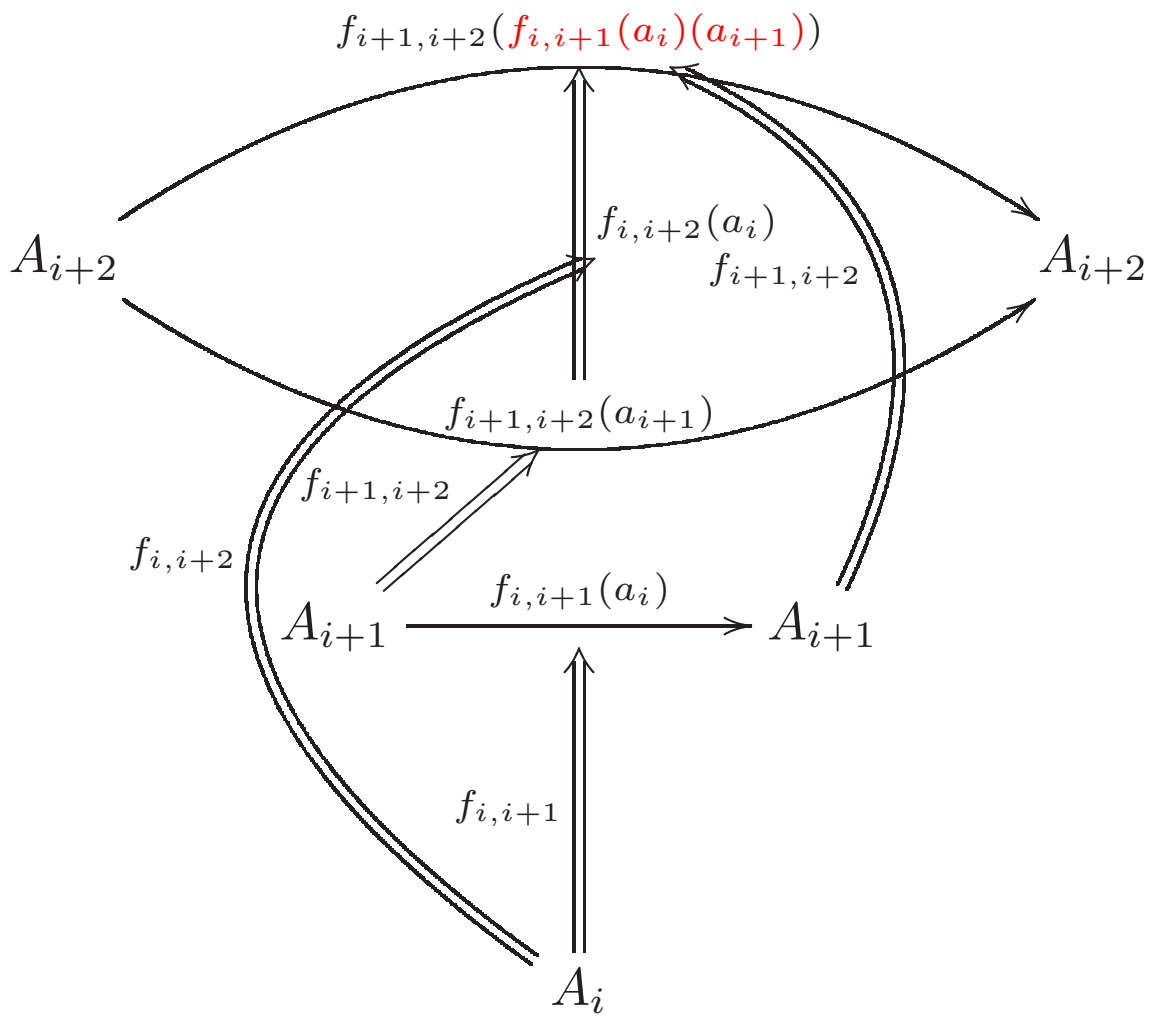

$f_{i,i+1}$ is representation of $\Omega_i$-algebra $A_i$ in $\Omega_{i+1}$-algebra $A_{i+1}$. $f_{i+1,i+2}$ is representation of $\Omega_{i+1}$-algebra $A_{i+1}$ in $\Omega_{i+2}$-algebra $A_{i+2}$.

Theorem 4.1.2. *The map*

$$f_{i,i+2} : A_i \to {}^{**}A_{i+2}$$

*is representation of $\Omega_i$-algebra $A_i$ in $\Omega_{i+1}$-algebra ${}^*A_{i+2}$.*

Proof. Automorphism $f_{i+1,i+2}(a_{i+1}) \in {}^*A_{i+2}$ corresponds to arbitrary $a_{i+1} \in A_{i+1}$. Automorphism $f_{i,i+1}(a_i) \in {}^*A_{i+1}$ corresponds to arbitrary $a_i \in A_i$. Since $f_{i,i+1}(a_i)(a_{i+1}) \in A_{i+1}$ is image of $a_{i+1} \in A_{i+1}$, then element $a_i \in A_i$ generates transformation of $\Omega_{i+1}$-algebra ${}^*A_{i+2}$ which is defined by equation

$$f_{i,i+2}(a_i)(f_{i+1,i+2}(a_{i+1})) = f_{i+1,i+2}(f_{i,i+1}(a_i)(a_{i+1})) \tag{4.1.2}$$

Let $\omega$ be $n$-ary operation of $\Omega_i$-algebra. Because $f_{i,i+1}$ is homomorphism of $\Omega_i$-algebra, then

$$f_{i,i+1}(a_{i,1}...a_{i,n}\omega) = f_{i,i+1}(a_{i,1})...f_{i,i+1}(a_{i,n})\omega \tag{4.1.3}$$

For $a_i$, ..., $a_n \in A_i$ we define operation $\omega$ on the set ${}^{**}A_{i+2}$ using equations

$$\begin{aligned}
&f_{i,i+2}(a_{i,1})(f_{i+1,i+2}(a_{i+1}))...f_{i,i+2}(a_{i,n})(f_{i+1,i+2}(a_{i+1}))\omega\\
&=f_{i+1,i+2}(f_{i,i+1}(a_{i,1})(a_{i+1}))...f_{i+1,i+2}(f_{i,i+1}(a_{i,n})(a_{i+1}))\omega\\
&=f_{i+1,i+2}(f_{i,i+1}(a_{i,1})...f_{i,i+1}(a_{i,n})\omega)(a_{i+1})
\end{aligned}\tag{4.1.4}$$

First equation follows from equation (4.1.2). We wrote second equation on the base of demand that map $f_{i+1,i+2}$ is homomorphism of $\Omega_i$-algebra. Therefore, we defined the structure of $\Omega_i$-algebra on the set ${}^{**}A_{i+2}$.

From (4.1.4) and (4.1.3) it follows that

$$\begin{aligned}
&f_{i,i+2}(a_{i,1})(f_{i+1,i+2}(a_{i+1}))...f_{i,i+2}(a_{i,n})(f_{i+1,i+2}(a_{i+1}))\omega\\
&=f_{i+1,i+2}(f_{i,i+1}(a_{i,1}...a_{i,n}\omega)(a_{i+1}))\\
&=f_{i,i+2}(a_{i,1}...a_{i,n}\omega)(f_{i+1,i+2}(a_{i+1}))
\end{aligned}\tag{4.1.5}$$

From equation (4.1.5) it follows that

$$f_{i,i+2}(a_{i,1})...f_{i,i+2}(a_{i,n})\omega=f_{i,i+2}(a_{i,1}...a_{i,n}\omega)$$

Therefore, the map $f_{i,i+2}$ is homomorphism of $\Omega_i$-algebra. Therefore, the map $f_{i,i+2}$ is representation of $\Omega_i$-algebra $A_i$ in $\Omega_{i+1}$-algebra ${}^*A_{i+2}$. □

THEOREM 4.1.3. *$(id, f_{i+1,i+2})$ is morphism of representations of $\Omega_i$-algebra from $f_{i,i+1}$ into $f_{i,i+2}$.*

PROOF. Consider diagram (4.1.1) in more detail.

(4.1.6)

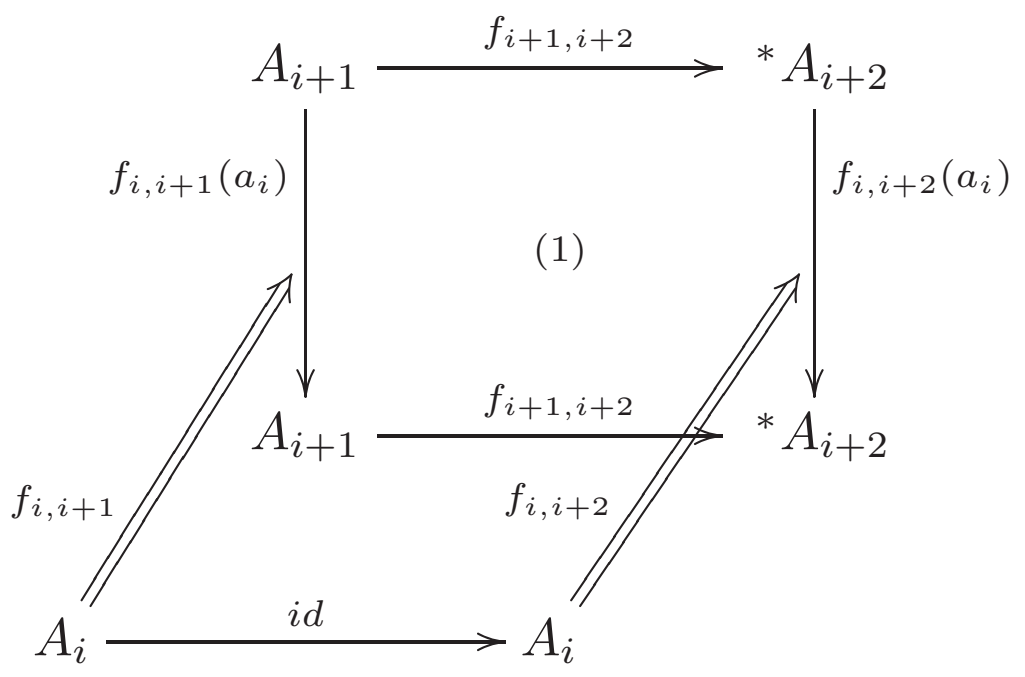

The statement of theorem follows from equation (4.1.2) and definition 2.2.2. □

THEOREM 4.1.4. *Consider tower $(\overline{A},\overline{f})$ of representations of $\overline{\Omega}$-algebras. Since identity transformation*

$$\delta_{i+2}:A_{i+2}\to A_{i+2}$$

*of $\Omega_{i+2}$-algebra $A_{i+2}$ belongs to representation $f_{i+1,i+2}$, then the representation $f_{i,i+2}$ of $\Omega_i$-algebra $A_i$ in $\Omega_{i+1}$-algebra ${}^*A_{i+2}$ can be extended to representation*

$$f'_{i,i+2}:A_i\to {}^*A_{i+2}$$

*of $\Omega_i$-algebra $A_i$ in $\Omega_{i+2}$-algebra $A_{i+2}$.*

Proof. According to the equation (4.1.2), for given $a_{i+1}\in A_{i+1}$ we can define a map

$$\begin{aligned}&f'_{i,i+2}:A_i\rightarrow {}^*A_{i+2}\\&f'_{i,i+2}(a_i)(f_{i+1,i+2}(a_{i+1})(a_{i+2}))=f_{i+1,i+2}(f_{i,i+1}(a_i)(a_{i+1}))(a_{i+2})\end{aligned}\tag{4.1.7}$$

Let there exist $a_{i+1}\in A_{i+1}$ such that

$$f_{i+1,i+2}(a_{i+1})=\delta_{i+2}\tag{4.1.8}$$

Then from (4.1.7),(4.1.8), it follows that

$$f'_{i,i+2}(a_i)(a_{i+2})=f_{i+1,i+2}(f_{i,i+1}(a_i)(a_{i+1}))(a_{i+2})\tag{4.1.9}$$

Let $\omega$ be $n$-ari operation of $\Omega_i$-algebra $A_i$. Since $f_{i,i+2}$ is homomorphism of $\Omega_i$-algebra, then, from the equation (4.1.2), it follows that

$$\begin{aligned}&f_{i,i+2}(a_{i\cdot 1}...a_{i\cdot n}\omega)(f_{i+1,i+2}(a_{i+1}))\\&=f_{i,i+2}(a_{i\cdot 1})(f_{i+1,i+2}(a_{i+1}))...f_{i,i+2}(a_{i\cdot n})(f_{i+1,i+2}(a_{i+1}))\omega\\&=f_{i+1,i+2}(f_{i,i+1}(a_{i\cdot 1})(a_{i+1})...f_{i,i+1}(a_{i\cdot n})(a_{i+1})\omega)\end{aligned}\tag{4.1.10}$$

From the equation (4.1.10), it follows that

(4.1.11)

$$\begin{aligned}&f_{i,i+2}(a_{i\cdot 1}...a_{i\cdot n}\omega)(f_{i+1,i+2}(a_{i+1})(a_{i+2}))\\&=f_{i,i+2}(a_{i\cdot 1})(f_{i+1,i+2}(a_{i+1})(a_{i+2}))...f_{i,i+2}(a_{i\cdot n})(f_{i+1,i+2}(a_{i+1})(a_{i+2}))\omega\\&=f_{i+1,i+2}(f_{i,i+1}(a_{i\cdot 1})(a_{i+1})...f_{i,i+1}(a_{i\cdot n})(a_{i+1})\omega)(a_{i+2})\end{aligned}$$

From equations (4.1.8), (4.1.9), (4.1.11), it follows that

$$\begin{aligned}&f'_{i,i+2}(a_{i\cdot 1}...a_{i\cdot n}\omega)(a_{i+2})\\&=f'_{i,i+2}(a_{i\cdot 1})(a_{i+2})...f'_{i,i+2}(a_{i\cdot n})(a_{i+2})\omega\\&=f_{i+1,i+2}(f_{i,i+1}(a_{i\cdot 1})(a_{i+1})...f_{i,i+1}(a_{i\cdot n})(a_{i+1})\omega)(a_{i+2})\end{aligned}\tag{4.1.12}$$

The equation (4.1.12) defines an operation $\omega$ on the set ${}^*A_{i+2}$

$$f'_{i,i+2}(a_{i\cdot 1})...f'_{i,i+2}(a_{i\cdot n})\omega=f'_{i,i+2}(a_{i\cdot 1}...a_{i\cdot n}\omega)\tag{4.1.13}$$

and the map $f'_{i,i+2}$ is homomorphism of $\Omega_i$-algebra. Therefore, we have a representation of $\Omega_i$-algebra $A_i$ in $\Omega_{i+2}$-algebra $A_{i+2}$. □

Definition 4.1.5. Consider the tower of representations

$$((A_1,A_2,A_3),(f_{1,2},f_{2,3}))$$

The map $f_*=(f_{1,2},f_{1,3})$ is called **representation of $\Omega_1$-algebra $A_1$ in representation $f_{2,3}$**. □

Definition 4.1.6. Consider the tower of representations $(\overline{A},\overline{f})$. The map

$$f_*=(f_{1,2},...,f_{1,n})$$

is called **representation of $\Omega_1$-algebra $A_1$ in tower of representations**

$$(\overline{A}_{[1]},\overline{f})=((A_2,...,A_n),(f_{2,3},...,f_{n-1,n}))$$

□

## 4.2. Morphism of Tower of $T*$-Representations

Definition 4.2.1. Consider the set of $\Omega_k$-algebr $A_k$, $B_k$, $k = 1, ..., n$. The set of maps $(h_1, ..., h_n)$ is called **morphism from tower of representations** $(\overline{A}, \overline{f})$ **into tower of representations** $(\overline{B}, \overline{g})$, if for any $k$, $k = 1, ..., n-1$, the tuple of maps $(h_k, h_{k+1})$ is morphism of representations from $f_{k,k+1}$ into $g_{k,k+1}$. □

For any $k$, $k = 1, ..., n-1$, we have diagram

(4.2.1)

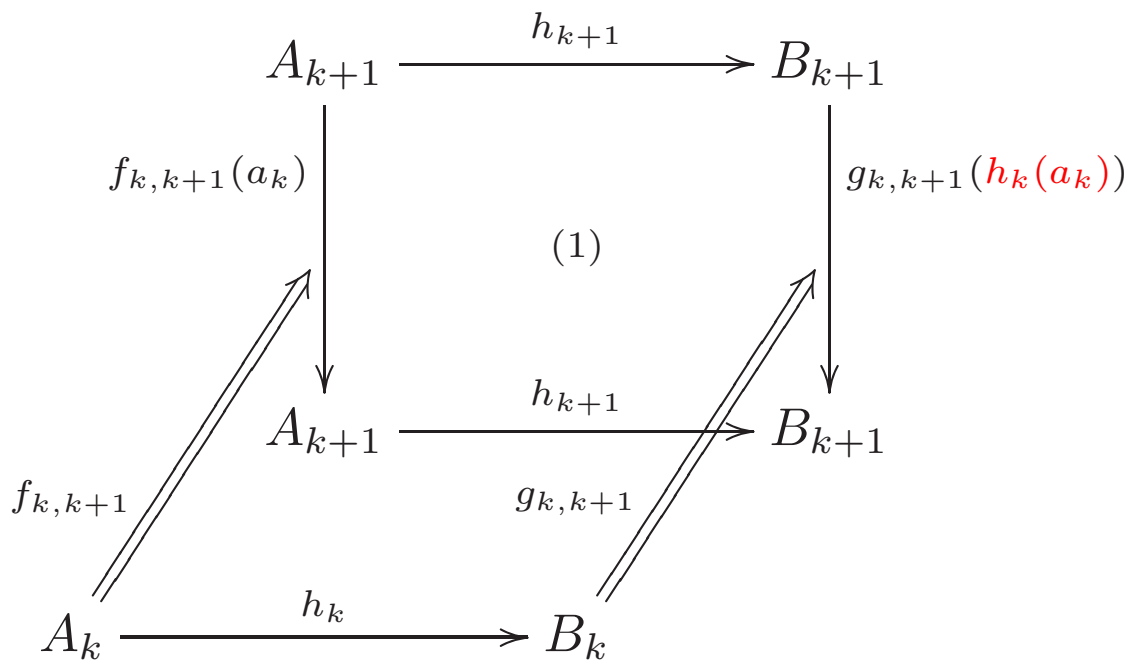

Equations

$$h_{k+1} \circ f_{k,k+1}(a_k) = g_{k,k+1}(h_k(a_k)) \circ h_{k+1} \tag{4.2.2}$$

$$h_{k+1}(f_{k,k+1}(a_k)(a_{k+1})) = g_{k,k+1}(h_k(a_k))(h_{k+1}(a_{k+1})) \tag{4.2.3}$$

express commutativity of diagram (1). However for morphism $(h_k, h_{k+1})$, $k > 1$, diagram (4.2.1) is not complete. Assuming similar diagram for morphism $(h_k, h_{k+1})$, this diagram on the top layer has form
(4.2.4)

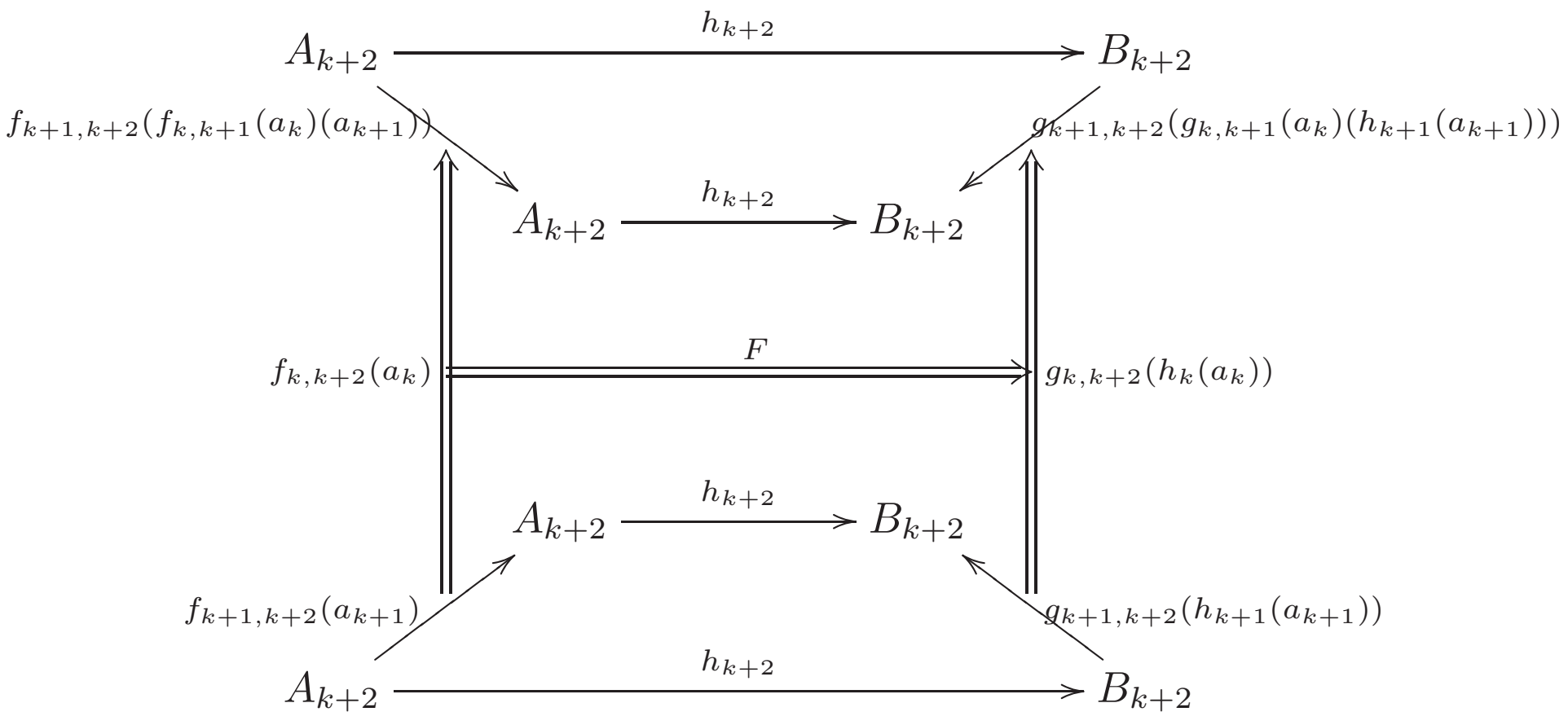

Unfortunately, the diagram (4.2.4) is not too informative. It is evident that there exists morphism from ${}^*A_{k+2}$ into ${}^*B_{k+2}$, which maps $f_{k,k+2}(a_k)$ into $g_{k,k+2}(h_k(a_k))$. However, the structure of this morphism is not clear from the diagram. We need consider map from ${}^*A_{k+2}$ into ${}^*B_{k+2}$, like we have done this in theorem 4.1.3.

Theorem 4.2.2. *Since representation $f_{i+1,i+2}$ is effective, then $(h_i, {}^*h_{i+2})$ is morphism of representations from representation $f_{i,i+2}$ into representation $g_{i,i+2}$ of $\Omega_i$-algebra.*

Proof. Consider the diagram

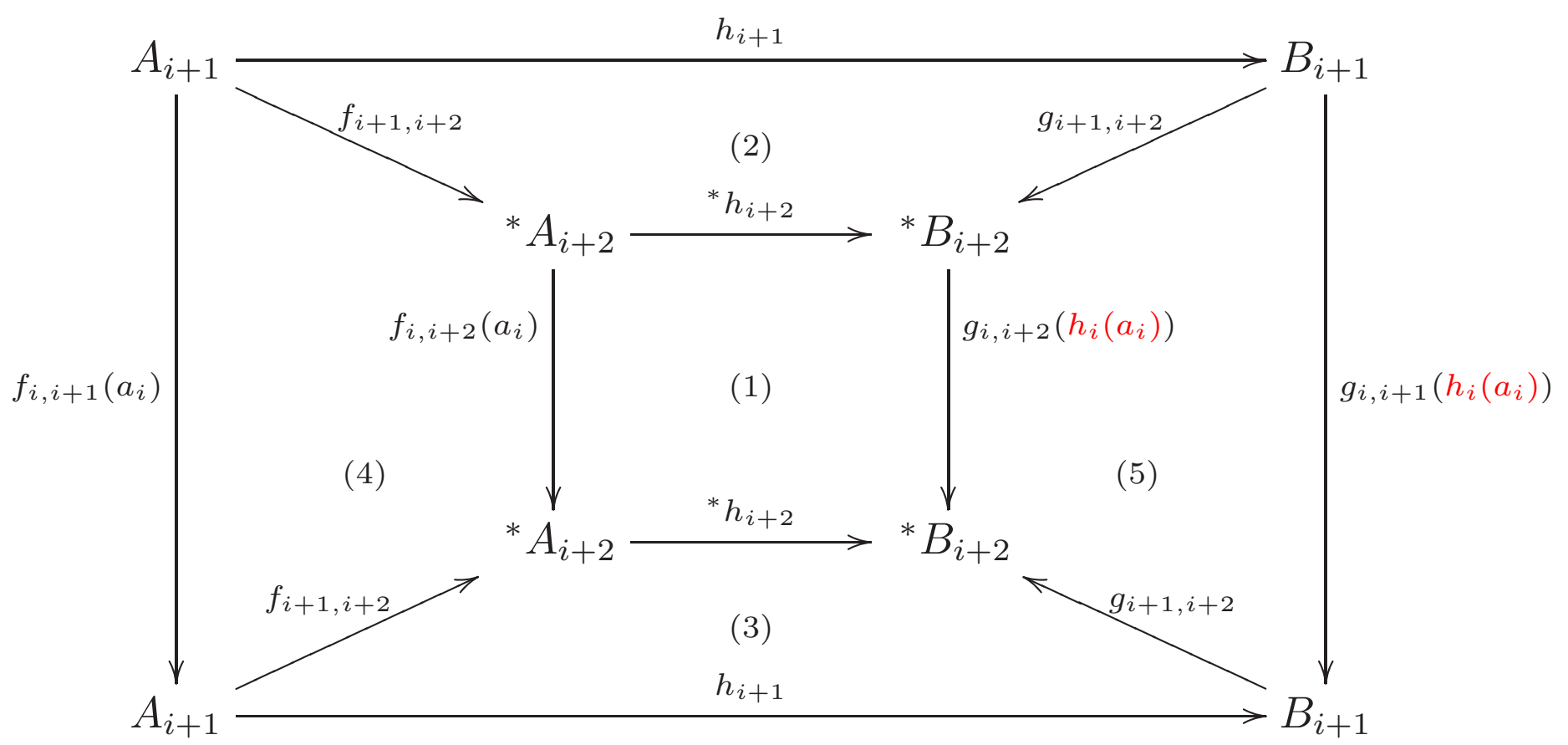

The existence of map ${}^*h_{i+2}$ and commutativity of the diagram (2) and (3) follows from effectiveness of map $f_{i+1,i+2}$ and theorem 2.2.8. Commutativity of diagrams (4) and (5) follows from theorem 4.1.3.

From commutativity of the diagram (4) it follows that

(4.2.5) $$f_{i+1,i+2}\circ f_{i,i+1}(a_i)=f_{i,i+2}(a_i)\circ f_{i+1,i+2}$$

From equation (4.2.5) it follows that

(4.2.6) $${}^*h_{i+2}\circ f_{i+1,i+2}\circ f_{i,i+1}(a_i)={}^*h_{i+2}\circ f_{i,i+2}(a_i)\circ f_{i+1,i+2}$$

From commutativity of diagram (3) it follows that

(4.2.7) $${}^*h_{i+2}\circ f_{i+1,i+2}=g_{i+1,i+2}\circ h_{i+1}$$

From equation (4.2.7) it follows

(4.2.8) $${}^*h_{i+2}\circ f_{i+1,i+2}\circ f_{i,i+1}(a_i)=g_{i+1,i+2}\circ h_{i+1}\circ f_{i,i+1}(a_i)$$

From equations (4.2.6) and (4.2.8) it follows that

(4.2.9) $${}^*h_{i+2}\circ f_{i,i+2}(a_i)\circ f_{i+1,i+2}=g_{i+1,i+2}\circ h_{i+1}\circ f_{i,i+1}(a_i)$$

From commutativity of the diagram (5) it follows that

(4.2.10) $$g_{i+1,i+2}\circ g_{i,i+1}(h_i(a_i))=g_{i,i+2}(h_i(a_i))\circ g_{i+1,i+2}$$

From equation (4.2.10) it follows that

(4.2.11) $$g_{i+1,i+2}\circ g_{i,i+1}(h_i(a_i))\circ h_{i+1}=g_{i,i+2}(h_i(a_i))\circ g_{i+1,i+2}\circ h_{i+1}$$

From commutativity of the diagram (2) it follows that

(4.2.12) $${}^*h_{i+2}\circ f_{i+1,i+2}=g_{i+1,i+2}\circ h_{i+1}$$

From equation (4.2.12) it follows that

(4.2.13) $$g_{i,i+2}(h_i(a_i))\circ {}^*h_{i+2}\circ f_{i+1,i+2}=g_{i,i+2}(h_i(a_i))\circ g_{i+1,i+2}\circ h_{i+1}$$

From equations (4.2.11) and (4.2.13) it follows that

(4.2.14) $$g_{i+1,i+2}\circ g_{i,i+1}(h_i(a_i))\circ h_{i+1}=g_{i,i+2}(h_i(a_i))\circ {}^*h_{i+2}\circ f_{i+1,i+2}$$

External diagram is diagram (4.2.1) when $i = 1$. Therefore, external diagram is commutative

$$h_{i+1} \circ f_{i,i+1}(a_i) = g_{i,i+1}(h_i(a_i)) \circ h_{i+1} \tag{4.2.15}$$

From equation (4.2.15) it follows that

$$g_{i+1,i+2} \circ h_{i+1} \circ f_{i,i+1}(a_i) = g_{i+1,i+2} \circ g_{i,i+1}(h_i(a_i)) \circ h_{i+1}(a_{i+1}) \tag{4.2.16}$$

From equations (4.2.9), (4.2.14) and (4.2.16) it follows that

$$^*h_{i+2} \circ f_{i,i+2}(a_i) \circ f_{i+1,i+2} = g_{i,i+2}(h_i(a_i)) \circ {}^*h_{i+2} \circ f_{i+1,i+2} \tag{4.2.17}$$

Because the map $f_{i+1,i+2}$ is injection, then from equation (4.2.17) it follows that

$$^*h_{i+2} \circ f_{i,i+2}(a_i) = g_{i,i+2}(h_i(a_i)) \circ {}^*h_{i+2} \tag{4.2.18}$$

From equation (4.2.18) commutativity of the diagram (1) follows. This proves the statement of theorem. □

Theorems 4.1.3 and 4.2.2 are true for any layer of tower of representations. In each particular case we need properly show sets and direction of the map. Meaning of given theorems is that all maps in tower of representations act coherently.

The theorem 4.2.2 states that unknown map on the diagram (4.2.4) is the map $^*h_{i+2}$.

Theorem 4.2.3. *Consider the set of $\Omega_i$-algebras $A_i$, $B_i$, $C_i$, $i = 1, ..., n$. Let*

$$\overline{p} :(\overline{A},\overline{f}) \rightarrow (\overline{B},\overline{g})$$
$$\overline{q} :(\overline{B},\overline{g}) \rightarrow (\overline{C},\overline{h})$$

*be morphisms of tower of representations. There exists morphism of representations of $\Omega$-algebra*

$$\overline{r} : (\overline{A},\overline{f}) \rightarrow (\overline{C},\overline{h})$$

*where $r_k = q_k p_k$, $k = 1, ..., n$. We call morphism $\overline{r}$ of tower of representations from $\overline{f}$ into $\overline{h}$* **product of morphisms $\overline{p}$ and $\overline{q}$ of tower of representations**.

Proof. For each $k$, $k = 2, ..., n$, we represent statement of theorem using diagram

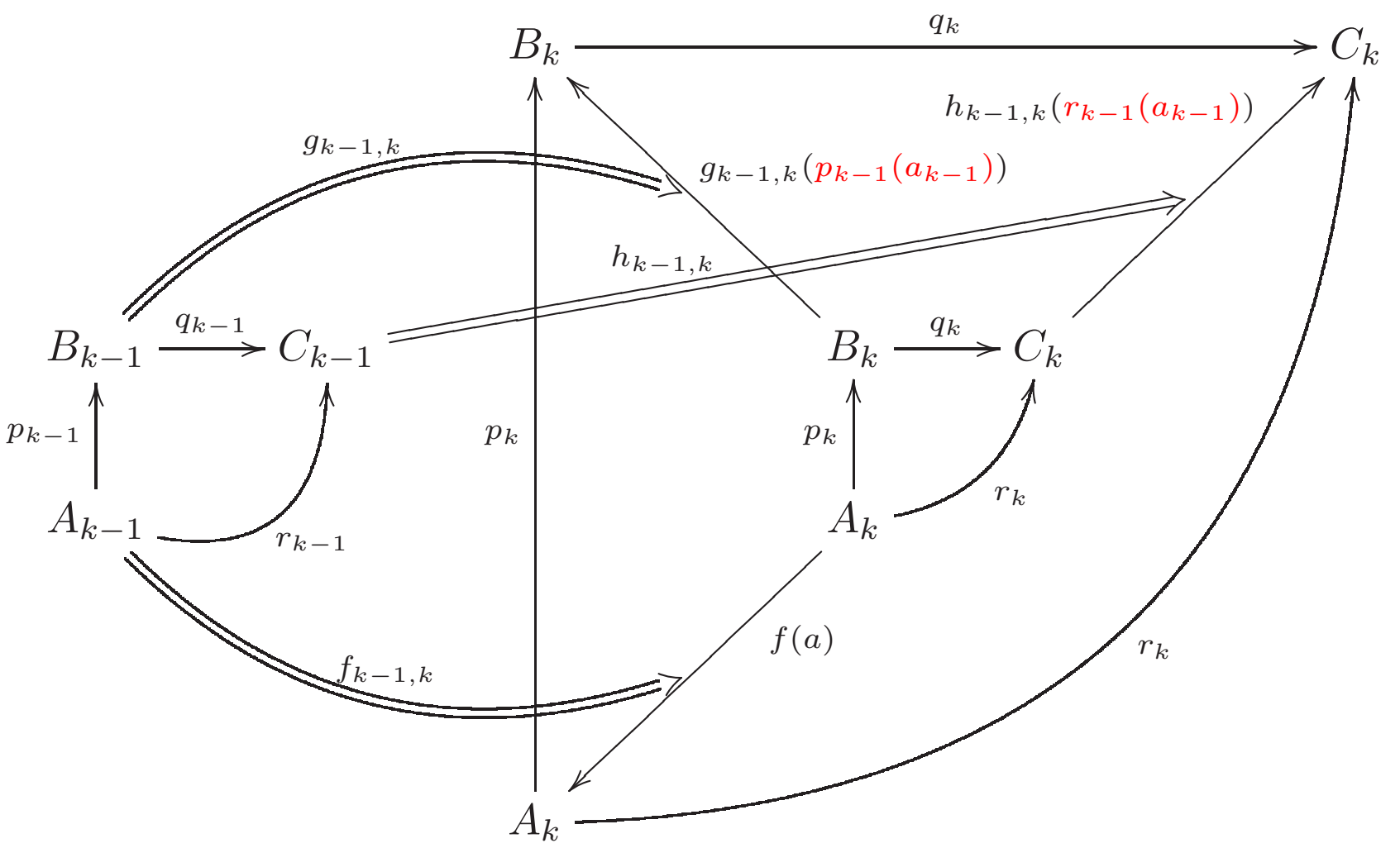

Map $r_{k-1}$ is homomorphism of $\Omega_{k-1}$-algebra $A_{k-1}$ into $\Omega_{k-1}$-algebra $C_{k-1}$. We need to show that tuple of maps $(r_{k-1}, r_k)$ satisfies to (4.2.2):

$$\begin{aligned} r_k(f_{k-1,k}(a_{k-1})a_k) &= q_k p_k(f_{k-1,k}(a_{k-1})a_k) \\ &= q_k(g_{k-1,k}(p_{k-1}(a_{k-1}))p_k(a_k)) \\ &= h_{k-1,k}(q_{k-1}p_{k-1}(a_{k-1}))q_k p_k(a_k)) \\ &= h_{k-1,k}(r(a_{k-1}))r_k(a_k) \end{aligned}$$

□

## 4.3. Automorphism of Tower of Representations

Definition 4.3.1. Let $(\overline{A}, \overline{f})$ be tower of representations of $\overline{\Omega}$-algebras. The morphism of tower of representations $(\mathrm{id}, h_2, ..., h_n)$ such, that for each $k$, $k = 2$, ..., $n$, $h_k$ is endomorphism of $\Omega_k$-algebra $A_k$ is called **endomorphism of tower of representations** $\overline{f}$. □

Definition 4.3.2. Let $(\overline{A}, \overline{f})$ be tower of representations of $\overline{\Omega}$-algebras. The morphism of tower of representations $(\mathrm{id}, h_2, ..., h_n)$ such, that for each $k$, $k = 2$, ..., $n$, $h_k$ is automorphism of $\Omega_k$-algebra $A_k$ is called **automorphism of tower of representations** $f$. □

Theorem 4.3.3. *Let $(\overline{A}, \overline{f})$ be tower of representations of $\overline{\Omega}$-algebras. The set of automorphisms of the representation $\overline{f}$ forms group.*

Proof. Let $\overline{r}_{[1]}$, $\overline{p}_{[1]}$ be automorphisms of the tower of representations $\overline{f}$. According to definition 4.3.2, for each $k$, $k = 2$, ..., $n$, maps $r_k$, $p_k$ are automorphisms of $\Omega_k$-algebra $A_k$. According to theorem II.3.2 ([13], p. 57), for each $k$, $k = 2$, ..., $n$, the map $r_k p_k$ is automorphism of $\Omega_k$-algebra $A_k$. From the theorem 4.2.3 and the definition 4.3.2, it follows that product of automorphisms $\overline{r}_{[1]}\overline{p}_{[1]}$ of the tower of representations $\overline{f}$ is automorphism of the tower of representations $\overline{f}$.

According to proof of the theorem 2.4.5, for any $k$, $k = 2$, ..., $n$, the product of automorphisms of $\Omega_i$-algebra is associative. Therefore, the product of automorphisms of tower of representations is associative.

Let $\overline{r}_{[1]}$ be an automorphism of the tower of representations $\overline{f}$. According to definition 4.3.2 for each $k$, $k = 2$, ..., $n$, the map $r_k$ is automorphism of $\Omega_k$-algebra $A_k$. Therefore, for each $k$, $k = 2$, ..., $n$, the map $r_k^{-1}$ is automorphism of $\Omega_k$-algebra $A_k$. The equation (4.2.3) is true for automorphism $\overline{r}_{[1]}$. Assume $a'_k = r_k(a_k)$, $k = 2$, ..., $n$. Since $r_k$, $k = 2$, ..., $n$, is automorphism then $a_k = r_k^{-1}(a'_k)$ and we can write (4.2.3) in the form

$$h_{k+1}(f_{k,k+1}(h_k^{-1}(a'_k))(h_{k+1}(a'_{k+1}))) = g_{k,k+1}(a'_k)(a'_{k+1}) \tag{4.3.1}$$

Since the map $h_{k+1}$ is automorphism of $\Omega_{k+1}$-algebra $A_{k+1}$, then from the equation (4.3.1) it follows that

$$f_{k,k+1}(h_k^{-1}(a'_k)(h_{k+1}(a'_{k+1}))) = h_{k+1}^{-1}(g_{k,k+1}(a'_k)(a'_{k+1})) \tag{4.3.2}$$

The equation (4.3.2) corresponds to the equation (4.2.3) for the map $\overline{r}_{[1]}^{-1}$. Therefore, map $\overline{r}_{[1]}^{-1}$ is automorphism of the representation $\overline{f}$. □

## 4.4. Examples of Tower of Representations

**4.4.1. Affine Space.** Let $\overset{\circ}{A}$ be a set of points. Let $\overline{A}$ be vector space over field $k$. Single transitive representation of the vector space $\overline{A}$ in the set of points $\overset{\circ}{A}$ is called Affine space. If $A \in \overset{\circ}{A}$ and $\overline{v} \in \overline{A}$, then we use an expression $A + \overline{v}$ to denote the action of vector $\overline{v}$ at the point $A$.

The statement that the representation is single transitive means that for any points $A, B \in \overset{\circ}{A}$ there exists a unique vector $\overline{v} \in \overline{A}$ such that

$$B = A + \overline{v} \tag{4.4.1}$$

If the vector $\overline{v}$ satisfies to the equation (4.4.1), then we denote the vector $\overline{v}$ as $\overrightarrow{AB}$

$$B = A + \overrightarrow{AB} \tag{4.4.2}$$

Since the vector space $\overline{A}$ is the representation of the field $k$ in the Abelian group $\overline{A}$. then the affine space is tower of representations.

**4.4.2. Free Algebra over Ring.** According to the definitions [14], p. 1, [12], p. 4, free algebra $A$ over commutative ring $R$ is free $R$-module in which the product defined as bilinear map of $R$-module $A$. This definition is unsatisfactory because endomorphisms generated by representation are not endomorphisms of algebra. However, we can change the definition of $R$-algebra in such a way that all requirements to representation have been met.

According to definition [10]-3.1.10, the map

$$l(a) : b \in A \to ab \in A \tag{4.4.3}$$

is linear map. Therefore, the map (4.4.3) is endomorphism of $R$-module $A$. Endomorphism (4.4.3) is called **left shift of $R$-module** $A$.

Theorem 4.4.1. *The representation*

$$f_{2,3} : a \in A \to l(a) \in {}^\star A \tag{4.4.4}$$

*of $R$-module $A$ in $R$-module $A$ is equivalent to structure of $R$-algebra $A$.*

Proof. Let the structure of $R$-module $A$ is generated by effective representation

$$f_{1,2} : R \to {}^\star A$$

of ring $R$ in Abelian group $A$.

The map $f_{2,3}(a)$ is left shift of $R$-module $A$. Since left shift is endomorphism of $R$-module $A$, then

$$\begin{aligned} f_{2,3}(a)(b_1 + b_2) &= f_{2,3}(a)(b_1) + f_{2,3}(a)(b_2) \\ f_{2,3}(a)(rb) &= r f_{2,3}(a)(b) \end{aligned} \tag{4.4.5}$$

Since the map (4.4.4) is homomorphism of Abelian groups, then

$$f_{2,3}(a_1 + a_2)(b) = (f_{2,3}(a_1) + f_{2,3}(a_2))(b) = f_{2,3}(a_1)(b) + f_{2,3}(a_2)(b) \tag{4.4.6}$$

According to the theorem 4.1.2 there is representation

$$f_{1,3} : R \to {}^{\star\star} A$$

which according to the equation (4.1.2) is defined by rule

$$f_{1,3}(r)(f_{2,3}(a)) = f_{2,3}(f_{1,2}(r)(a)) \tag{4.4.7}$$

We do not assume that the representation $f_{2,3}$ is effective.[4.1] Therefore, the representation $f_{1,3}$ also can be not effective. However, we can, just as in the case of the map $f_{1,2}$, identify $f_{1,3}(r)$ and $r$. From the equation (4.4.7), it follows that

$$rf_{2,3}(a)(b) = f_{2,3}(ra)(b) \tag{4.4.8}$$

From equations (4.4.5), (4.4.6), (4.4.8) and the definition [10]-3.1.10, it follows that the map $f_{2,3}$ is bilinear map. Therefore, the map $f_{2,3}$ determines the product in $R$-module $A$ according to rule

$$ab = f_{2,3}(a)(b)$$

□

Corollary 4.4.2. *The tower of representations*

$$((f_{1,2}, f_{2,3}), (R, A, A)) \tag{4.4.9}$$

*generates the structure of $R$-algebra $A$.* □

We point out that, in the definition of $R$-algebra $A$, domain of the map $f_{2,3}$ is not all set of endomorphisms of Abelian group $A$, but set of endomorphisms of $R$-module $A$. We assume that in the definition 4.1.1 each element of the tower of representations is an universal algebra. In this case **diagram of representations** for $R$-algebra has form

$$\begin{array}{ccccc} R & \xrightarrow[*]{f_{1,2}} & A & \xrightarrow[*]{f_{2,3}} & A \\ & & & & \uparrow{\scriptstyle *\, g_{1,2}} \\ & & & & R \end{array} \tag{4.4.10}$$

On diagram of representations (4.4.10), we assume that we at first construct the representation $g_{1,2}$ generating the structure of $R$-module $A$. The representation $f_{1,2}$ generates the structure of $R$-module $A$. According to the definition, the representation $f_{2,3}$ is homomorphism of $R$-module $A$ into the set of endomorphisms of $R$-module $A$.

**4.4.3. Group Algebra.** If $R$ is a ring and $G$ is multiplicative group, then free $R$-module $R[G]$ with basis $G$ is called group ring ([5], p. 383). From the theorem [7]-13.1.1, it follows that we can extend the product from group $G$ to module $R[G]$. Therefore, module $R[G]$ is free $R$-algebra.

If $R$ is commutative ring, then module $R[G]$ is called **group algebra** of group $G$ ([6], §8). Since group $G$ is the basis of group algebra, then we can write element of algebra $R[G]$ in form

$$\sum_{g\in G} a_g g \quad a_g \in R \tag{4.4.11}$$

The product of arbitary elements of $R[G]$ has form

$$\left(\sum_{g\in G} a_g g\right)\left(\sum_{g\in G} b_g g\right) = \sum_{g\in G} c_g g \qquad c_g = \sum_{u\in G} a_u b_{u^{-1}g} \tag{4.4.12}$$

[4.1] If the representation $f_{2,3}$ is not effective, then $R$-algebra $A$ has zero divisors.

We can consider coordinates of element $a \in R[g]$ as map

$$a : G \rightarrow R \quad a(g) = a_g$$

**4.4.4. Module over Algebra.** Let $A_2$ be free algebra over field $A_1$. Considering the algebra $A_2$ as a ring, we can determine free module $A_3$ over the algebra $A_2$. We can describe the module $A_3$ using the diagram of representations

$$\begin{array}{ccccc} A_2 & \xrightarrow[*]{f_{2,3}} & A_2 & \xrightarrow[*]{f_{3,4}} & A_3 \\ \uparrow {\scriptstyle f_{1,2}} & & \uparrow {\scriptstyle f_{1,2}} & & \\ A_1 & & A_1 & & \end{array}$$

## 4.5. Generating Set of Tower of Representations

Definition 4.5.1. Tower of representations $(\overline{A}, \overline{f})$. is called **tower of effective representations**, if for any $i$ the representation $f_{i,i+1}$ is effective. □

Theorem 4.5.2. *Consider the tower of representations $(\overline{A}, \overline{f})$. Let representations $f_{i,i+1}$, ..., $f_{i+k-1,i+k}$ be effective. Then representation $f_{i,i+k}$ is effective.*

Proof. We will prove the statement of theorem by induction.

Let representations $f_{i,i+1}$, $f_{i+1,i+2}$ be effective. Assume that transformation $f_{i,i+1}(a_i)$ is not identity transformation. Then there exists $a_{i+1} \in A_{i+1}$ such that $f_{i,i+1}(a_i)(a_{i+1}) \neq a_{i+1}$. Because the representation $f_{i+1,i+2}$ is effective, then transformations $f_{i+1,i+2}(a_{i+1})$ and $f_{i+1,i+2}(f_{i,i+1}(a_i)(a_{i+1}))$ do not coincide. According to construction in theorem 4.1.2, the transformation $f_{i,i+2}(a_i)$ is not identity transformation.

Assume the statement of theorem is true for $k-1$ representations. Let $f_{i,i+1}$, ..., $f_{i+k-1,i+k}$ be effective. According assumption of induction, representations $f_{i,i+k-1}$, $f_{i+k-1,i+k}$ are effective. According to proven above, the representation $f_{i,i+k}$ is effective. □

Theorem 4.5.3. *Consider tower $(\overline{A}, \overline{f})$ of representations of $\overline{\Omega}$-algebras. Let identity transformation*

$$\delta_{i+2} : A_{i+2} \rightarrow A_{i+2}$$

*of $\Omega_{i+2}$-algebra $A_{i+2}$ belong to representation $f_{i+1,i+2}$. Let representations $f_{i,i+1}$, $f_{i+1,i+2}$ be effective. Then representation $f'_{i,i+2}$ defined in the theorem 4.1.4, is effective.*

Proof. Let there exist $a_{i+1} \in A_{i+1}$ such that

$$f_{i+1,i+2}(a_{i+1}) = \delta_{i+2}$$

Assume that the representation $f'_{i,i+2}$ is not effective. Then there exist

(4.5.1) $$a_{i\cdot 1}, a_{i\cdot 2} \in A_i \quad a_{i\cdot 1} \neq a_{i\cdot 2}$$

such that

(4.5.2) $$f'_{i,i+2}(a_{i\cdot 1})(a_{i+2}) = f'_{i,i+2}(a_{i\cdot 2})(a_{i+2})$$

From equations (4.1.9), (4.5.2), it follows that

(4.5.3) $$f_{i+1,i+2}(f_{i,i+1}(a_{i\cdot 1})(a_{i+1}))(a_{i+2}) = f_{i+1,i+2}(f_{i,i+1}(a_{i\cdot 2})(a_{i+1}))(a_{i+2})$$

Since $a_{i+2}$ is arbitrary, then, from (4.5.3), it follows that

$$f_{i+1,i+2}(f_{i,i+1}(a_{i\cdot 1})(a_{i+1})) = f_{i+1,i+2}(f_{i,i+1}(a_{i\cdot 2})(a_{i+1})) \tag{4.5.4}$$

Since the representation $f_{i,i+1}$ is effective, then, from the condition (4.5.1), it follows that

$$f_{i,i+1}(a_{i\cdot 1})(a_{i+1}) \ne f_{i,i+1}(a_{i\cdot 2})(a_{i+1}) \tag{4.5.5}$$

From the condition (4.5.5) and the equation (4.5.4), it follows that the representation $f_{i+1,i+2}$ is not effective. The contradiction completes the proof of the theorem. □

We construct the basis of the tower of representations in a similar way that we constructed the basis of representation in the section 2.7.

We will write elements of tower of representations $(\overline{A},\overline{f})$ as tuple

$$\overline{a} = (a_1, ..., a_3) \quad a_i \in A_i \quad i = 1, ..., n$$

Definition 4.5.4. Let $(\overline{A},\overline{f})$ be tower of representations. The tuple of sets

$$\overline{N}_{[1]} = (N_2 \subset A_2, ..., N_n \subset A_n)$$

is called **tuple of stable sets of tower of representations** $\overline{f}$, if

$$f_{i-1,i}(a_{i-1})(a_i) \in N_i \quad i = 2, ..., n$$

for every $a_1 \in A_1$, $a_2 \in N_2$, ..., $a_n \in N_n$. We also will say that tuple of sets

$$\overline{N}_{[1]} = (N_2 \subset A_2, ..., N_n \subset A_n)$$

is stable relative to tower of representations $\overline{f}$. □

Theorem 4.5.5. *Let $\overline{f}$ be tower of representations. Let set $N_i \subset A_i$ be subalgebra of $\Omega_i$-algebra $A_i$, $i = 2$, ..., $n$. Let tuple of sets*

$$\overline{N}_{[1]} = (N_2 \subset A_2, ..., N_n \subset A_n)$$

*be stable relative to tower of representations $\overline{f}$. Then there exists representation*

$$((A_1, N_2, ..., N_n), (f_{N_2,1,2}, ..., f_{N_n,n-1,n})) \tag{4.5.6}$$

*such that*

$$f_{N_i,i-1,i}(a_{i-1}) = f_{i-1,i}(a_{i-1})|_{N_i} \quad i = 2, ..., n$$

*The tower of representations (4.5.6) is called* **tower of subrepresentations**.

Proof. Let $\omega_{i-1,1}$ be $m$-ary operation of $\Omega_{i-1}$-algebra $A_{i-1}$, $i = 2$, ..., $n$. Then for any $a_{i-1,1}$, ..., $a_{i-1,m} \in N_{i-1}$ [4.2] and any $a_i \in N_i$

$$\begin{aligned}
&(f_{N_i,i-1,i}(a_{i-1,1})...f_{N_i,i-1,i}(a_{i-1,m})\omega_{i-1,1})(a_i)\\
=&(f_{i-1,i}(a_{i-1,1})...f_{i-1,i}(a_{i-1,m})\omega_{i-1,1})(a_i)\\
=&f_{i-1,i}(a_{i-1,1}...a_{i-1,m}\omega_{i-1,1})(a_i)\\
=&f_{N_i,i-1,i}(a_{i-1,1}...a_{i-1,m}\omega_{i-1,1})(a_i)
\end{aligned}$$

[4.2] Assume $N_1 = A_1$.

Let $\omega_{i,2}$ be $m$-ary operation of $\Omega_i$-algebra $A_i$, $i = 2, ..., n$. Then for any $a_{i,1}, ..., a_{i,n} \in N_i$ and any $a_{i-1} \in N_{i-1}$

$$\begin{aligned}&f_{N_i,i-1,i}(a_{i-1})(a_{i,1})...f_{N_i,i-1,i}(a_{i-1})(a_{i,m})\omega_{i,2}\\ =&f_{i-1,i}(a_{i-1})(a_{i,1})...f_{i-1,i}(a_{i-1})(a_{i,m})\omega_{i,2}\\ =&f_{i-1,i}(a_{i-1})(a_{i,1}...a_{i,m}\omega_{i,2})\\ =&f_{N_i,i-1,i}(a_{i-1})(a_{i,1}...a_{i,m}\omega_{i,2})\end{aligned}$$

We proved the statement of theorem. □

From theorem 4.5.5, it follows that if map $(f_{N_2,1,2}, ..., f_{N_n,1,n})$ is tower of subrepresentations of tower of representations $\overline{f}$, then map

$$(id : A_1 \to A_1, id_2 : N_2 \to A_2, ..., id_n : N_n \to A_n)$$

is morphism of towers of representations.

THEOREM 4.5.6. *The set*[4.3] $\mathcal{B}_{\overline{f}}$ *of all towers of subrepresentations of tower of representations* $\overline{f}$ *generates a closure system on tower of representations* $\overline{f}$ *and therefore is a complete lattice.*

PROOF. Let for given $\lambda \in \Lambda$, $K_{\lambda,i}$, $i = 2, ..., n$, be subalgebra of $\Omega_i$-algebra $A_i$ that is stable relative to representation $f_{i-1,i}$. We determine the operation of intersection on the set $\mathcal{B}_{\overline{f}}$ according to rule

$$\bigcap f_{K_{\lambda,i-1,i}} = f_{\cap K_{\lambda,i-1,i}} \quad i = 2, ..., n$$

$$\bigcap \overline{K}_\lambda = \left(K_1 = A_1, K_2 = \bigcap K_{\lambda,2}, ..., K_n = \bigcap K_{\lambda,n}\right)$$

$\cap K_{\lambda,i}$ is subalgebra of $\Omega_i$-algebra $A_i$. Let $a_i \in \cap K_{\lambda,i}$. For any $\lambda \in \Lambda$ and for any $a_{i-1} \in K_{i-1}$,

$$f_{i-1,i}(a_{i-1})(a_i) \in K_{\lambda,i}$$

Therefore,

$$f_{i-1,i}(a_{i-1})(a_i) \in K_i$$

Repeating this construction in the order of increment $i$, $i = 2, ..., n$, we see that $(K_1, ..., K_n)$ is tuple of stable sets of tower of representations $\overline{f}$. Therefore, we determined the operation of intersection of towers of subrepresentations properly. □

We denote the corresponding closure operator by $\overline{J}(\overline{f})$. If we denote $\overline{X}_{[1]}$ the tuple of sets $(X_2 \subset A_2, ..., X_n \subset A_n)$ then $\overline{J}(\overline{f}, \overline{X}_{[1]})$ is the intersection of all tuples $(K_1, ..., K_n)$ stable with respect to representation $\overline{f}$ and such that for $i = 2$, ..., $n$, $K_i$ is subalgebra of $\Omega_i$-algebra $A_i$ containing $X_i$.[4.4]

[4.3] This definition is similar to definition of the lattice of subalgebras ([13], p. 79, 80)

[4.4] For $n = 2$, $J_2(f_{1,2}, X_2) = J_{f_{1,2}}(X_2)$. It would be easier to use common notation in sections 2.7 and 4.6. However I think that using of vector notation in section 2.7 is premature.

THEOREM 4.5.7. *Let*[4.5] $\overline{f}$ *be the tower of representations. Let* $X_i \subset A_i$, $i = 2$, ..., $n$. *Assume* $Y_1 = A_1$. *Step by step increasing the value of* $i$, $i = 2$, ..., $n$, *we define a subsets* $X_{i,m} \subset A_i$ *by induction on* $m$.

$$X_{i,0} = X_i$$
$$x \in X_{i,m} => x \in X_{i,m+1}$$
$$x_1 \in X_{i,m}, ..., x_p \in X_{i,m}, \omega \in \Omega_i(p) => x_1...x_p\omega \in X_{i,m+1}$$
$$x_i \in X_{i,m}, x_{i-1} \in Y_{i-1} => f_{i-1,i}(x_{i-1})(x_i) \in X_{i,m+1}$$

*For each value of* $i$, *we assume*

$$Y_i = \bigcup_{m=0}^{\infty} X_{i,m}$$

*Then*

$$\overline{Y} = (Y_1, ..., Y_n) = \overline{J}(\overline{f}, \overline{X}_{[1]})$$

PROOF. For each value of $i$ the proof of the theorem coincides with the proof of theorem 2.6.4. Because to define stable subset of $\Omega_i$-algebra $A_i$ we need only certain stable subset of $\Omega_{i-1}$-algebra $A_{i-1}$, we have to find stable subset of $\Omega_{i-1}$-algebra $A_{i-1}$ before we do this in $\Omega_i$-algebra $A_i$. □

$\overline{J}(\overline{f}, \overline{X}_{[1]})$ is called **tower of subrepresentations of tower of representations $\overline{f}$ generated by tuple of sets** $\overline{X}_{[1]}$, and $\overline{X}_{[1]}$ is a **tuple of generating sets of tower subrepresentations** $\overline{J}(\overline{f}, \overline{X}_{[1]})$. In particular, a **tuple of generating sets of tower of representations $\overline{f}$** is a tuple $(X_2 \subset A_2, ..., X_n \subset A_n)$ such that $\overline{J}(\overline{f}, \overline{X}_{[1]}) = \overline{A}$.

It is easy to see that definition of the tuple of generating sets of tower of representations does not depend on whether representations of tower are effective or not. For this reason hereinafter we will assume that representations of the tower are effective and we will use convention for effective representation in remark 2.1.7.

We also will use notation

$$\overline{r}_{[1]} \circ \overline{a}_{[1]} = (r_2 \circ a_2, ..., r_n \circ a_n)$$

for image of tuple of elements $a_2 \in A_2$, ..., $a_n \in A_n$ under the endomorphism of tower of effective representations. According to the definition of product of maps, for any endomorphisms $\overline{r}_{[1]}$, $\overline{s}_{[1]}$ the following equation is true

$$(\overline{r}_{[1]} \circ \overline{s}_{[1]}) \circ \overline{a}_{[1]} = \overline{r}_{[1]} \circ (\overline{s}_{[1]} \circ \overline{a}_{[1]}) \tag{4.5.7}$$

The equation (4.5.7) is associative law for $\circ$ and allows us to write expression

$$\overline{r}_{[1]} \circ \overline{s}_{[1]} \circ \overline{a}_{[1]}$$

without brackets.

From theorem 4.5.7, it follows next definition.

DEFINITION 4.5.8. Let $(X_2 \subset A_2, ..., X_n \subset A_n)$ be tuple of sets. For each tuple of elements $\overline{a}$, $\overline{a} \in \overline{J}(\overline{f}, \overline{X}_{[1]})$, there exists tuple of $\overline{\Omega}$-words defined according to following rule.

(1) If $a_1 \in A_1$, then $a_1$ is $\Omega_1$-word.
(2) If $a_i \in X_i$, $i = 2$, ..., $n$, then $a_i$ is $\Omega_i$-word.

[4.5]The statement of theorem is similar to the statement of theorem 5.1, [13], p. 79.

(3) If $a_{i,1}$, ..., $a_{i,p}$ are $\Omega_i$-words, $i=2$, ..., $n$, and $\omega\in\Omega_i(p)$, then $a_{i,1}...a_{i,p}\omega$ is $\Omega_i$-word.
(4) If $a_i$ is $\Omega_i$-word, $i=2$, ..., $n$, and $a_{i-1}$ is $\Omega_{i-1}$-word, then $a_{i-1}a_i$ is $\Omega_i$-word.

**Tuple of $\overline{\Omega}$-words**[4.6]

$$\overline{w}(\overline{f},\overline{X}_{[1]},\overline{a})=(w_1(\overline{f},\overline{X}_{[1]},a_1),...,w_n(\overline{f},\overline{X}_{[1]},a_n))$$

represents given element $\overline{a}\in\overline{J}(\overline{f},\overline{X}_{[1]})$.[4.7] Denote $\overline{w}(\overline{f},\overline{X}_{[1]})$ **the set of tuples of $\overline{\Omega}$-words of tower of representations $\overline{J}(\overline{f},\overline{X}_{[1]})$.** □

Representation of $a_i\in A_i$ as $\Omega_i$-word is ambiguous. If $a_{i,1}$, ..., $a_{i,p}$ are $\Omega_i$-words, $\omega\in\Omega_i(p)$ and $a_{i-1}\in A_{i-1}$, then $\Omega_i$-words $a_{i-1}a_{i,1}...a_{i,p}\omega$ and $a_{i-1}a_{i,1}...a_{i-1}a_{i,p}\omega$ represent the same element of $\Omega_i$-algebra $A_i$. It is possible that there exist equations related with a character of a representation. For instance, if $\omega$ is the operation of $\Omega_{i-1}$-algebra $A_{i-1}$ and the operation of $\Omega_i$-algebra $A_i$, then we can request that $\Omega_i$-words $a_{i-1,1}...a_{i-1,p}\omega a_i$ and $a_{i-1,1}a_i...a_{i-1,p}a_i\omega$ represent the same element of $\Omega_i$-algebra $A_i$. Listed above equations for each value $i$, $i=2$, ..., $n$, determine equivalence $r_i$ on the set of $\Omega_i$-words $W_i(\overline{f},\overline{X}_{[1]})$.According to the construction, equivalence $r_i$ on the set of $\Omega_i$-words $W_i(\overline{f},\overline{X}_{[1]})$ depends not only on the choice of the set $X_i$, but also on the choice of the set $X_{i-1}$.

Theorem 4.5.9. *Endomorphism $\overline{r}_{[1]}$ of tower of representations $\overline{f}$ forms the map of $\overline{\Omega}_{[1]}$-words*

$$\overline{w}_{[1]}[\overline{f},\overline{X}_{[1]},\overline{r}_{[1]}]:\overline{W}_{[1]}(\overline{f},\overline{X}_{[1]})\rightarrow\overline{w}_{[1]}(\overline{f},\overline{X}'_{[1]})\quad\overline{X}_{[1]}\subset\overline{A}_{[1]}\quad\overline{X}'_{[1]}=\overline{r}_{[1]}(\overline{X}_{[1]})$$

*such that for any $i$, $i=2$, ..., $n$,*

(1) *If* $a_i\in X_i$, $a'_i=r_i\circ a_i$, *then*
$$w_i[\overline{f},\overline{X}_{[1]},\overline{r}_{[1]}](a_i)=a'_i$$
(2) *If*
$$a_{i,1},...,a_{i,n}\in w_i(\overline{f},\overline{X}_{[1]})$$
$$a'_{i,1}=w_i[\overline{f},\overline{X}_{[1]},\overline{r}_{[1]}](a_{i,1})\quad...\quad a'_{i,p}=w_i[\overline{f},\overline{X}_{[1]},\overline{r}_{[1]}](a_{i,p})$$
*then for operation $\omega\in\Omega_i(p)$ holds*
$$w_i[\overline{f},\overline{X}_{[1]},\overline{r}_{[1]}](a_{i,1}...a_{i,p}\omega)=a'_{i,1}...a'_{i,p}\omega$$

---

[4.6]According to our design
$$w_2(\overline{f},\overline{X}_{[1]},a_2)=w_2(f_{1,2},X_2,a_2)\quad w_3(\overline{f},\overline{X}_{[1]},a_3)=w_3((f_{1,2},f_{2,3}),(X_2,X_3),a_3)$$

[4.7]These notations have small difference with notations offered earlier in the theory of representation of universal algebra. Namely, we added a word of $\Omega_1$-algebra into tuple of $\overline{\Omega}$-words. However it is evident that this word is always trivial
$$w_1(\overline{f},\overline{X}_{[1]},a_1)=a_1$$

(3) *If*

$$a_i\in w_i(\overline{f},\overline{X}_{[1]})\qquad a'_i=w_i[\overline{f},\overline{X}_{[1]},\overline{r}_{[1]}](a_i)$$
$$a_{i-1}\in w_{i-1}(\overline{f},\overline{X}_{[1]})\quad a'_{i-1}=w_{i-1}[\overline{f},\overline{X}_{[1]},\overline{r}_{[1]}](a_{i-1})$$

*then*

$$w_i[\overline{f},\overline{X}_{[1]},\overline{r}_{[1]}](a_{i-1}a_i)=a'_{i-1}a'_i$$

Proof. Statements (1), (2) of the theorem are true by definition of the endomorpism $r_i$. The statement (3) of the theorem follows from the equation (4.2.3). □

Remark 4.5.10. Let $\overline{r}_{[1]}$ be endomorphism of tower of representations $\overline{f}$. Let

$$\overline{a}_{[1]}\in\overline{J}_{[1]}(\overline{f},\overline{X}_{[1]})\quad \overline{a}'_{[1]}=\overline{r}_{[1]}\circ\overline{a}_{[1]}\quad \overline{X}'_{[1]}=\overline{r}_{[1]}\circ\overline{X}_{[1]}$$

The theorem 4.5.9 states that $\overline{a}'_{[1]}\in\overline{J}_{[1]}(\overline{f},\overline{X}'_{[1]})$ . The theorem 4.5.9 also states that $\overline{\Omega}_{[1]}$-word representing $\overline{a}_{[1]}$ relative $\overline{X}_{[1]}$ and $\overline{\Omega}_{[1]}$-word representing $\overline{a}'_{[1]}$ relative $\overline{X}'_{[1]}$ are generated according to the same algorithm. This allows considering of the set of $\overline{\Omega}_{[1]}$-words $\overline{w}_{[1]}(\overline{f},\overline{X}'_{[1]},\overline{a}'_{[1]})$ as map

$$\overline{W}_{[1]}(\overline{f},\overline{X}_{[1]},\overline{a}_{[1]}):\overline{X}'_{[1]}\to\overline{w}_{[1]}(\overline{f},\overline{X}'_{[1]},\overline{a}'_{[1]})$$
$$\overline{W}_{[1]}(\overline{f},\overline{X}_{[1]},\overline{a}_{[1]})(\overline{X}'_{[1]})=\overline{W}_{[1]}(\overline{f},\overline{X}_{[1]},\overline{a}_{[1]})\circ\overline{X}'_{[1]}$$

such that, if for certain endomorphism $\overline{r}_{[1]}$

$$\overline{X}'_{[1]}=\overline{r}_{[1]}\circ\overline{X}_{[1]}\quad \overline{a}'_{[1]}=\overline{r}_{[1]}\circ\overline{a}_{[1]}$$

then

$$\overline{W}_{[1]}(\overline{f},\overline{X}_{[1]},\overline{a}_{[1]})\circ\overline{X}'_{[1]}=\overline{w}_{[1]}(\overline{f},\overline{X}'_{[1]},\overline{a}'_{[1]})=\overline{a}'_{[1]}$$

Tuple of maps[4.8]

$$(W_2(\overline{f},\overline{X}_{[1]},a_2),...,W_n(\overline{f},\overline{X}_{[1]},a_n))$$

is called **tuple of coordinates of element $\overline{a}$ relative to tuple of sets $\overline{X}_{[1]}$**. Denote $\overline{W}_{[1]}(\overline{f},\overline{X}_{[1]})$ **the set of tuples of coordinates of tower of representations** $\overline{J}(\overline{f},\overline{X}_{[1]})$. Similarly, we consider coordinates of tuple of sets $\overline{B}_{[1]}\subset\overline{J}_{[1]}(f,X)$ relative to the tuple of sets $\overline{X}_{[1]}$

$$\overline{W}_{[1]}(\overline{f},\overline{X}_{[1]},\overline{B}_{[1]})=(W_2(\overline{f},\overline{X}_{[1]},B_2),...,W_n(\overline{f},\overline{X}_{[1]},B_n))$$
$$W_k(\overline{f},\overline{X}_{[1]},B_k)=\{W_k(\overline{f},\overline{X}_{(2...k)},a_k):a_k\in B_k\}=(W_k(\overline{f},\overline{X}_{(2...k)},a_k),a_k\in B_k)$$

□

Theorem 4.5.11. *There is a structure of $\Omega_k$-algebra on the set of coordinates* $W_k(\overline{f},\overline{X}_{[1]})$.

---

[4.8] According to our design

$$W_2(\overline{f},\overline{X}_{[1]},a_2)=W_2(f_{1,2},X_2,a_2)\quad W_3(\overline{f},\overline{X}_{[1]},a_3)=W_3((f_{1,2},f_{2,3}),(X_2,X_3),a_3)$$

We do not care about the map $W_1(\overline{f},\overline{X},a_1)$ because this map is trivial

$$W_1(\overline{f},\overline{X}_{[1]},a_1)\circ\overline{X}'_{[1]}=a_1 \tag{4.5.8}$$

Proof. Let $\omega\in\Omega_k(n)$. Then for any $m_1, ..., m_n\in J_k(\overline{f},\overline{X}_{[1]})$ , we assume

$$W_k(\overline{f},\overline{X}_{[1]},m_1)...W_k(\overline{f},\overline{X}_{[1]},m_n)\omega=W_k(\overline{f},\overline{X}_{[1]},m_1...m_n\omega) \tag{4.5.9}$$

According to the remark 4.5.10,

$$\begin{aligned}(W_k(\overline{f},\overline{X}_{[1]},m_1)...W_k(\overline{f},\overline{X}_{[1]},m_n)\omega)\circ\overline{X}_{[1]}&=W_k(\overline{f},\overline{X}_{[1]},m_1...m_n\omega)\circ\overline{X}_{[1]}\\&=w(\overline{f},\overline{X}_{[1]},m_1...m_n\omega)\end{aligned} \tag{4.5.10}$$

follows from the equation (4.5.9). According to rule (3) of the definition 4.5.8, from the equation (4.5.10), it follows that

$$\begin{aligned}&(W_k(\overline{f},\overline{X}_{[1]},m_1)...W_k(\overline{f},\overline{X}_{[1]},m_n)\omega)\circ\overline{X}_{[1]}\\&=w_k(\overline{f},\overline{X}_{[1]},m_1)...w_k(\overline{f},\overline{X}_{[1]},m_n)\omega\\&=(W_k(\overline{f},\overline{X}_{[1]},m_1)\circ\overline{X}_{[1]})...(W_k(\overline{f},\overline{X}_{[1]},m_n)\circ\overline{X}_{[1]})\omega\end{aligned} \tag{4.5.11}$$

From the equation (4.5.11), it follows that the operation $\omega$ defined by the equation (4.5.9) on the set of coordinates $W_k(\overline{f},\overline{X}_{[1]})$ is defined properly. □

Theorem 4.5.12. *For $k=2, ..., n$, there exists the representation of $\Omega_{k-1}$-algebra $W_{k-1}(\overline{f},\overline{X}_{[1]})$ in $\Omega_k$-algebra $W_k(\overline{f},\overline{X}_{[1]})$.*[4.9]

Proof. Let $a_{k-1}\in J_{k-1}(\overline{f},\overline{X}_{[1]})$. Then for any $a_k\in J_k(\overline{f},\overline{X}_{[1]})$, we assume

$$W_{k-1}(\overline{f},\overline{X}_{[1]},a_{k-1})W_k(\overline{f},\overline{X}_{[1]},a_k)=W_k(\overline{f},\overline{X}_{[1]},a_{k-1}a_k) \tag{4.5.12}$$

According to the remark 4.5.10,

$$\begin{aligned}(W_{k-1}(\overline{f},\overline{X}_{[1]},a_{k-1})W_k(\overline{f},\overline{X}_{[1]},a_k))\circ\overline{X}_{[1]}&=W_k(\overline{f},\overline{X}_{[1]},am)\circ\overline{X}_{[1]}\\&=w_k(\overline{f},\overline{X}_{[1]},a_{k-1}a_k)\end{aligned} \tag{4.5.13}$$

follows from the equation (4.5.12). According to rule (4) of the definition 4.5.8, from the equation (4.5.13), it follows that

$$\begin{aligned}&(W_{k-1}(\overline{f},\overline{X}_{[1]},a_{k-1})W_k(\overline{f},\overline{X}_{[1]},a_k))\circ\overline{X}_{[1]}\\&=w_{k-1}(\overline{f},\overline{X}_{[1]},a_{k-1})w_k(\overline{f},\overline{X}_{[1]},a_k)\\&=(W_{k-1}(\overline{f},\overline{X}_{[1]},a_{k-1})\circ\overline{X}_{[1]})(W_k(\overline{f},\overline{X}_{[1]},a_k)\circ\overline{X}_{[1]})\end{aligned} \tag{4.5.14}$$

From the equation (4.5.14), it follows that the representation (4.5.12) of $\Omega_{k-1}$-algebra $W_{k-1}(\overline{f},\overline{X}_{[1]})$ in $\Omega_k$-algebra $W_k(\overline{f},\overline{X}_{[1]})$ is defined properly. □

Corollary 4.5.13. *Tuple of $\overline{\Omega}$-algebras $\overline{W}(\overline{f},\overline{X}_{[1]})$ forms the tower of representations.* □

Theorem 4.5.14. *Let $\overline{f}$ be tower of representations. For given sets $X_k\subset A_k$, $X'_k\subset A_k$, $k=2, ..., n$, consider tuple of maps*

$$\overline{R}_{[1]}=(R_2,...,R_n)$$

*such that for any $k=2, ..., n$, the map*

$$R_k:X_k\rightarrow X'_k$$

[4.9] According to equation (4.5.8), we identify $\Omega_1$-algebras $W_1(\overline{f},\overline{X}_{[1]})$ and $A_1$.

*agree with the structure of representation $f_{k-1,k}$, i. e.*[4.10]

$$\omega\in\Omega_k(p),\ \ x_{1,k},\ \ ...,\ \ x_{p,k},\ \ x_{1,k}...x_{p,k}\omega\in X_k,\ \ R_1(x_{1,k}...x_{p,k}\omega)\in X'_k$$
$$=> R_1(x_{1,k}...x_{p,k}\omega)=R_1(x_{1,k})...R_1(x_{p,k})\omega$$
$$a_{k-1}\in X_{k-1},\ \ a_k\in X_k,\ \ R_1(a_{k-1}a_k)\in X'_k$$
$$=> R_k(a_{k-1}a_k)=R_{k-1}(a_{k-1})R_k(a_k)$$

*Consider the map of $\overline{\Omega}_{[1]}$-words*

$$\overline{w}_{[1]}[\overline{f},\overline{X}_{[1]},\overline{R}_{[1]}]:\overline{w}_{[1]}(\overline{f},\overline{X}_{[1]})\rightarrow\overline{w}_{[1]}(\overline{f},\overline{X}'_{[1]})$$

*that satisfies conditions (1), (2), (3) of theorem 4.5.9 and such that*

$$\overline{x}_{[1]}\in\overline{X}_{[1]} => w[\overline{f},\overline{X}_{[1]},\overline{R}_{[1]}](\overline{x}_{[1]})=\overline{R}_{[1]}(\overline{x}_{[1]})$$

*For each $k$, $k=2$, ..., $n$, there exists endomorphism of $\Omega_k$-algebra $A_k$*

$$r_k:A_k\rightarrow A_k$$

*defined by rule*

$$r_k\circ a_k=w_k[\overline{f},\overline{X}_{[1]},\overline{R}_{[1]}](w_k(\overline{f},\overline{X}_{[1]},a_k))\tag{4.5.15}$$

*Tuple of endomorphisms $\overline{r}_{[1]}$ is morphism of towers of representations $\overline{J}(\overline{f},\overline{X}_{[1]})$ and $\overline{J}(\overline{f},\overline{X}'_{[1]})$.*

Proof. If $n=2$, then tower of representations $\overline{f}$ is representation of $\Omega_1$-algebra $A_1$ in $\Omega_2$-algebra $A_2$. The statement of theorem is corollary of theorem 2.6.11.

Let the statement of theorem be true for $n-1$. We do not change notation in theorem when we move from one layer to another because a word in $\Omega_{n-1}$-algebra $A_{n-1}$ does not depend on word in $\Omega_n$-algebra $A_n$. We prove the theorem by induction over complexity of $\Omega_n$-word.

If $w_n(\overline{f},\overline{X}_{[1]},a_n)=a_n$, then $a_n\in X_n$. According to condition (1) of theorem 4.5.9,

$$\begin{aligned}r_n\circ a_n&=w_n[\overline{f},\overline{X}_{[1]},\overline{R}_{[1]}](w_n(\overline{f},\overline{X}_{[1]},a_n))\\&=w_n[\overline{f},\overline{X}_{[1]},\overline{R}_{[1]}](a_n)\\&=R_n(a_n)\end{aligned}$$

Therefore, maps $r_n$ and $R_n$ coinside on the set $X_n$, and the map $r_n$ agrees with structure of $\Omega_n$-algebra.

Let $\omega\in\Omega_n(p)$. Let the statement of induction be true for

$$a_{n,1},...,a_{n,p}\in J_n(\overline{f},\overline{X}_{[1]})$$

Let

$$w_{n,1}=w_n(\overline{f},\overline{X}_{[1]},a_{n,1})\quad ...\quad w_{n,p}=w_n(\overline{f},\overline{X}_{[1]},a_{n,p})$$

If

$$a_n=a_{n,1}...a_{n,p}\omega$$

then according to condition (3) of definition 4.5.8,

$$w_n[\overline{f},\overline{X}_{[1]},a_n]=w_{n,1}...w_{n,p}\omega$$

[4.10] We assume $X_1=A_1$, $R_1(a_1)=a_1$.

According to condition (2) of theorem 4.5.9,

$$\begin{aligned} r_n \circ a_n &= w_n[\overline{f}, \overline{X}_{[1]}, \overline{R}_{[1]}](w_n(\overline{f}, \overline{X}_{[1]}, a_n)) \\ &= w_n(\overline{f}, \overline{X}_{[1]}, \overline{R}_{[1]}](w_{n,1}...w_{n,p}\omega) \\ &= w_n(\overline{f}, \overline{X}_{[1]}, \overline{R}_{[1]}](w_{n,1})...w_n(\overline{f}, \overline{X}_{[1]}, \overline{R}_{[1]}](w_{n,p})\omega \\ &= (r_n \circ a_{n,1})...(r_n \circ a_{n,p})\omega \end{aligned}$$

Therefore, the map $r_n$ is endomorphism of $\Omega_n$-algebra $A_n$.

Let the statement of induction be true for

$$\begin{array}{ll} a_n \in J_n(\overline{f}, \overline{X}_{[1]}) & w_n(\overline{f}, \overline{X}_{[1]}, a_n) = m_n \\ a_{n-1} \in J_{n-1}(\overline{f}, \overline{X}_{[1]}) & w_{n-1}(\overline{f}, \overline{X}_{[1]}, a_{n-1}) = m_{n-1} \end{array}$$

According to condition (4) of definition 4.5.8,

$$w_n(\overline{f}, \overline{X}_{[1]}, a_{n-1}a_n) = m_{n-1}m_n$$

According to condition (3) of theorem 4.5.9,

$$\begin{aligned} r_n \circ (a_{n-1}a_n) &= w_n(\overline{f}, \overline{X}_{[1]}, \overline{R}_{[1]}](w_n(\overline{f}, \overline{X}_{[1]}, a_{n-1}a_n)) \\ &= w_n(\overline{f}, \overline{X}_{[1]}, \overline{R}_{[1]}](m_{n-1}m_n) \\ &= w_{n-1}(\overline{f}, \overline{X}_{[1]}, \overline{R}_{[1]}](m_{n-1})w_n(\overline{f}, (r_1, R_2, ..., R_n), \overline{X}_{[1]})(m_n) \\ &= (r_{n-1} \circ a_{n-1})(r_n \circ a_n) \end{aligned}$$

From equation (4.2.3) it follows that the map $\overline{r}$ is morphism of the tower of representations $\overline{f}$. □

Remark 4.5.15. The theorem 4.5.14 is the theorem of extension of map. The only statement we know about the tuple of sets $\overline{X}_{[1]}$ is the statement that $\overline{X}_{[1]}$ is the tuple of generating sets of the tower of representations $\overline{f}$. However, between the elements of the set $X_k$, $k = 2, ..., n$, there may be relationships generated by either operations of $\Omega_k$-algebra $A_k$, or by transformation of representation $f_{k-1,k}$. Therefore, any map of tuple of sets $\overline{X}_{[1]}$, in general, cannot be extended to an endomorphism of tower of representations $\overline{f}$.[4.11] However, if for any $k$, $k = 2, ..., n$, the map $R_k$ is coordinated with the structure of representation $f_{k-1,k}$, then we can construct an extension of this map and this extension is endomorphism of tower of representations $\overline{f}$. □

Definition 4.5.16. Let $\overline{X}_{[1]}$ be the tuple of generating sets of tower of representations $\overline{f}$. Let $\overline{r}_{[1]}$ be the endomorphism of the tower of representations $\overline{f}$. The set of coordinates $\overline{W}_{[1]}(\overline{f}, \overline{X}_{[1]}, \overline{r}_{[1]} \circ \overline{X}_{[1]})$ is called **coordinates of endomorphism of tower of representations**. □

Definition 4.5.17. Let $\overline{X}_{[1]}$ be the tuple of generating sets of tower of representations $\overline{f}$. Let $\overline{r}_{[1]}$ be the endomorphism of the tower of representations $\overline{f}$. Let $\overline{a}_{[1]} \in \overline{A}_{[1]}$. We define **superposition of coordinates** of the tower of representations $\overline{f}$ and the element $\overline{a}_{[1]}$ as coordinates defined according to rule

$$\overline{W}_{[1]}(\overline{f}, \overline{X}_{[1]}, \overline{a}_{[1]}) \circ \overline{W}_{[1]}(\overline{f}, \overline{X}_{[1]}, \overline{r}_{[1]} \circ \overline{X}_{[1]}) = \overline{W}_{[1]}(\overline{f}, \overline{X}_{[1]}, \overline{r}_{[1]} \circ \overline{a}_{[1]}) \tag{4.5.16}$$

[4.11] In the theorem 4.6.7, requirements to tuple of generating sets are more stringent. Therefore, the theorem 4.6.7 says about extension of arbitrary map. A more detailed analysis is given in the remark 4.6.9.

Let $\overline{Y}_{[1]} \subset \overline{A}_{[1]}$. We define superposition of coordinates of the tower of representations $\overline{f}$ and the tuple of sets $\overline{Y}_{[1]}$ according to rule

$$\begin{aligned}&\overline{W}_{[1]}(\overline{f},\overline{X}_{[1]},\overline{Y}_{[1]})\circ\overline{W}_{[1]}(\overline{f},\overline{X}_{[1]},\overline{r}_{[1]}\circ\overline{X}_{[1]})\\ =&(\overline{W}_{[1]}(\overline{f},\overline{X}_{[1]},\overline{a}_{[1]})\circ\overline{W}_{[1]}(\overline{f},\overline{X}_{[1]},\overline{r}_{[1]}\circ\overline{X}_{[1]}),\overline{a}_{[1]}\in\overline{Y}_{[1]})\end{aligned}\tag{4.5.17}$$

□

Theorem 4.5.18. *Endomorphism $\overline{r}_{[1]}$ of representation $\overline{f}$ generates the map of coordinates of tower of representations*

$$\overline{W}_{[1]}[\overline{f},\overline{X}_{[1]},\overline{r}_{[1]}]:\overline{W}_{[1]}(\overline{f},\overline{X}_{[1]})\rightarrow\overline{W}_{[1]}(\overline{f},\overline{X}_{[1]})\tag{4.5.18}$$

*such that*

(4.5.19)

$$\overline{W}_{[1]}(\overline{f},\overline{X}_{[1]},\overline{a}_{[1]})\rightarrow\overline{W}_{[1]}[\overline{f},\overline{X}_{[1]},\overline{r}_{[1]}]\star\overline{W}_{[1]}(\overline{f},\overline{X}_{[1]},\overline{a}_{[1]})=\overline{W}_{[1]}(\overline{f},\overline{X}_{[1]},\overline{r}_{[1]}\circ\overline{a}_{[1]})$$

$$\overline{W}_{[1]}[\overline{f},\overline{X}_{[1]},\overline{r}_{[1]}]\star\overline{W}_{[1]}(\overline{f},\overline{X}_{[1]},\overline{a}_{[1]})=\overline{W}_{[1]}(\overline{f},\overline{X}_{[1]},\overline{a}_{[1]})\circ\overline{W}_{[1]}(\overline{f},\overline{X}_{[1]},\overline{r}_{[1]}\circ\overline{X}_{[1]})$$

Proof. According to the remark 4.5.10, we consider equations (4.5.16), (4.5.18) relative to given tuple of generating sets $\overline{X}_{[1]}$. The tuple of words

$$\overline{W}_{[1]}(\overline{f},\overline{X}_{[1]},\overline{a}_{[1]})\circ\overline{X}_{[1]}=\overline{w}_{[1]}(\overline{f},\overline{X}_{[1]},\overline{a}_{[1]})\tag{4.5.20}$$

corresponds to coordinates $\overline{W}_{[1]}(\overline{f},\overline{X}_{[1]},\overline{a}_{[1]})$; the tuple of words

$$\overline{W}_{[1]}(\overline{f},\overline{X}_{[1]},\overline{r}_{[1]}\circ\overline{a}_{[1]})\circ X=\overline{w}_{[1]}(\overline{f},\overline{X}_{[1]},\overline{r}_{[1]}\circ\overline{a}_{[1]})\tag{4.5.21}$$

corresponds to coordinates $\overline{W}_{[1]}(\overline{f},\overline{X}_{[1]},\overline{r}_{[1]}\circ\overline{a}_{[1]})$. Therefore, in order to prove the theorem, it is sufficient to show that the map $\overline{W}_{[1]}[\overline{f},\overline{X}_{[1]},\overline{r}_{[1]}]$ corresponds to map $\overline{w}_{[1]}[\overline{f},\overline{X}_{[1]},\overline{r}_{[1]}]$. We prove this statement by induction over number $n$ of universal algebras in tower of representations.

If $n=2$, then tower of representations $\overline{f}$ is representation of $\Omega_1$-algebra $A_1$ in $\Omega_2$-algebra $A_2$. The statement of theorem is corollary of theorem 2.6.15. Let the statement of theorem be true for $n-1$. We prove the theorem by induction over complexity of $\Omega_n$-word.

If $a_n\in X_n$, $a'_n=r_n\circ a_n$, then, according to equations (4.5.20), (4.5.21), maps $W_n[\overline{f},\overline{X}_{[1]},\overline{r}_{[1]}]$ and $w_n[\overline{f},\overline{X}_{[1]},\overline{r}_{[1]}]$ are coordinated.

Let for $a_{n,1}$, ..., $a_{n,p}\in X_n$ maps $W_n[\overline{f},\overline{X}_{[1]},\overline{r}_{[1]}]$ and $w_n[\overline{f},\overline{X}_{[1]},\overline{r}_{[1]}]$ are coordinated. Let $\omega\in\Omega_n(p)$. According to the theorem 2.6.9

$$W_n(\overline{f},\overline{X}_{[1]},a_{n,1}...a_{n,p}\omega)=W_n(\overline{f},\overline{X}_{[1]},a_{n,1})...W_n(\overline{f},\overline{X}_{[1]},a_{n,p})\omega\tag{4.5.22}$$

Because $r_n$ is endomorphism of $\Omega_n$-algebra $A_n$, then from the equation (4.5.22) it follows that

(4.5.23)

$$\begin{aligned}W_n(\overline{f},\overline{X}_{[1]},r_n\circ(a_{n,1}...a_{n,p}\omega))&=W_n(\overline{f},\overline{X}_{[1]},(r_n\circ a_{n,1})...(r_n\circ a_{n,p})\omega)\\&=W_n(\overline{f},\overline{X}_{[1]},r_n\circ a_{n,1})...W(f,X,r_n\circ a_{n,p})\omega\end{aligned}$$

From equations (4.5.22), (4.5.23) and the statement of induction, it follows that the maps $W_n[\overline{f},\overline{X}_{[1]},\overline{r}_{[1]}]$ and $w_n[\overline{f},\overline{X}_{[1]},\overline{r}_{[1]}]$ are coordinated for $a_n=a_{n,1}...a_{n,p}\omega$.

Let for $a_{n,1}\in A_n$ maps $W_n[\overline{f},\overline{X}_{[1]},\overline{r}_{[1]}]$ and $w_n[\overline{f},\overline{X}_{[1]},\overline{r}_{[1]}]$ are coordinated. Let $a\in A$. According to the theorem 4.5.12

$$W_n(\overline{f},\overline{X}_{[1]},a_{n-1}a_{n,1})=W_{n-1}(\overline{f},\overline{X}_{[1]},a_{n-1})W_n(\overline{f},\overline{X}_{[1]},a_{n,1}) \tag{4.5.24}$$

Because $\overline{r}_{[1]}$ is endomorphism of tower of representations $\overline{f}$, then from the equation (4.5.24) it follows that

$$\begin{aligned}W_n(\overline{f},\overline{X}_{[1]},r_n\circ(a_{n-1}a_{n,1}))&=W_n(\overline{f},\overline{X}_{[1]},(r_{n-1}\circ a_{n-1})(r_n\circ a_{n,1}))\\&=W_{n-1}(\overline{f},\overline{X}_{[1]},r_{n-1}\circ a_{n-1})W_n(\overline{f},\overline{X}_{[1]},r_n\circ a_{n,1})\end{aligned} \tag{4.5.25}$$

From equations (4.5.24), (4.5.25) and the statement of induction, it follows that maps $W_n[\overline{f},\overline{X}_{[1]},\overline{r}_{[1]}]$ and $w_n[\overline{f},\overline{X}_{[1]},\overline{r}_{[1]}]$ are coordinated for $a_{n,2}=a_{n-1}a_{n,1}$. □

Corollary 4.5.19. *Let $\overline{X}_{[1]}$ be the tuple of generating sets of the tower of representations $\overline{f}$. Let $\overline{r}_{[1]}$ be the endomorphism of the tower of representations $\overline{f}$. The map $\overline{W}_{[1]}[\overline{f},\overline{X}_{[1]},\overline{r}_{[1]}]$ is endomorphism of tower of representations $\overline{W}(\overline{f},\overline{X}_{[1]})$.* □

Hereinafter we will identify map $\overline{W}_{[1]}[\overline{f},\overline{X}_{[1]},\overline{r}_{[1]}]$ and the set of coordinates $\overline{W}_{[1]}(\overline{f},\overline{X}_{[1]},\overline{r}_{[1]}\circ\overline{X}_{[1]})$.

Theorem 4.5.20. *Let $\overline{X}_{[1]}$ be the tuple of generating sets of the tower of representations $\overline{f}$. Let $\overline{r}_{[1]}$ be the endomorphism of the tower of representations $\overline{f}$. Let $\overline{Y}_{[1]}\subset\overline{A}_{[1]}$. Then*

$$\overline{W}_{[1]}(\overline{f},\overline{X}_{[1]},\overline{Y}_{[1]})\circ\overline{W}_{[1]}(\overline{f},\overline{X}_{[1]},\overline{r}_{[1]}\circ\overline{X}_{[1]})=\overline{W}_{[1]}(\overline{f},\overline{X}_{[1]},\overline{r}_{[1]}\circ\overline{Y}_{[1]}) \tag{4.5.26}$$

$$\overline{W}_{[1]}[\overline{f},\overline{X}_{[1]},\overline{r}_{[1]}]\star\overline{W}_{[1]}(\overline{f},\overline{X}_{[1]},\overline{Y}_{[1]})=\overline{W}_{[1]}(\overline{f},\overline{X}_{[1]},\overline{r}_{[1]}\circ\overline{Y}_{[1]}) \tag{4.5.27}$$

Proof. The equation (4.5.26) follows from the equation

$$\overline{r}_{[1]}\circ\overline{Y}_{[1]}=(\overline{r}_{[1]}\circ\overline{a}_{[1]},\overline{a}_{[1]}\in\overline{Y}_{[1]})$$

as well from equations (4.5.16), (4.5.17). The equation (4.5.27) is corollary of equations (4.5.26), (4.5.19). □

Theorem 4.5.21. *Let $\overline{X}_{[1]}$ be the tuple of generating sets of the tower of representations $\overline{f}$. Let $\overline{r}_{[1]}$, $\overline{s}_{[1]}$ be the endomorphisms of the tower of representations $\overline{f}$. Then*

$$\overline{W}_{[1]}(\overline{f},\overline{X}_{[1]},\overline{s}_{[1]}\circ\overline{X}_{[1]})\circ\overline{W}_{[1]}(\overline{f},\overline{X}_{[1]},\overline{r}_{[1]}\circ\overline{X}_{[1]})=\overline{W}_{[1]}(\overline{f},\overline{X}_{[1]},\overline{r}_{[1]}\circ\overline{s}_{[1]}\circ\overline{X}_{[1]}) \tag{4.5.28}$$

$$\overline{W}_{[1]}[\overline{f},\overline{X}_{[1]},\overline{r}_{[1]}]\star\overline{W}_{[1]}[\overline{f},\overline{X}_{[1]},\overline{s}_{[1]}]=\overline{W}_{[1]}[\overline{f},\overline{X}_{[1]},\overline{r}_{[1]}\circ\overline{s}_{[1]}] \tag{4.5.29}$$

Proof. The equation (4.5.28) follows from the equation (4.5.26), if we assume $\overline{Y}_{[1]}=\overline{s}_{[1]}\circ\overline{X}_{[1]}$. The equation (4.5.29) follows from the equation (4.5.28) and

chain of equations

$$
\begin{aligned}
&(\overline{W}_{[1]}[\overline{f},\overline{X}_{[1]},\overline{r}_{[1]}]\star\overline{W}_{[1]}[\overline{f},\overline{X}_{[1]},\overline{s}_{[1]}])\star\overline{W}_{[1]}(\overline{f},\overline{X}_{[1]},\overline{Y}_{[1]})\\
=&\overline{W}_{[1]}[\overline{f},\overline{X}_{[1]},\overline{r}_{[1]}]\star(\overline{W}_{[1]}[\overline{f},\overline{X}_{[1]},\overline{s}_{[1]}]\star\overline{W}_{[1]}(\overline{f},\overline{X}_{[1]},\overline{Y}_{[1]}))\\
=&(\overline{W}_{[1]}(\overline{f},\overline{X}_{[1]},\overline{Y}_{[1]})\circ\overline{W}_{[1]}(\overline{f},\overline{X}_{[1]},\overline{s}_{[1]}\circ\overline{X}_{[1]}))\circ\overline{W}_{[1]}(\overline{f},\overline{X}_{[1]},\overline{r}_{[1]}\circ\overline{X}_{[1]})\\
=&\overline{W}_{[1]}(\overline{f},\overline{X}_{[1]},\overline{Y}_{[1]})\circ(\overline{W}_{[1]}(\overline{f},\overline{X}_{[1]},\overline{s}_{[1]}\circ\overline{X}_{[1]})\circ\overline{W}_{[1]}(\overline{f},\overline{X}_{[1]},\overline{r}_{[1]}\circ\overline{X}_{[1]}))\\
=&\overline{W}_{[1]}(\overline{f},\overline{X}_{[1]},\overline{Y}_{[1]})\circ\overline{W}_{[1]}(\overline{f},\overline{X}_{[1]},\overline{r}_{[1]}\circ\overline{s}_{[1]}\circ\overline{X}_{[1]})\\
=&\overline{W}_{[1]}(\overline{f},\overline{X}_{[1]},\overline{r}_{[1]}\circ\overline{s}_{[1]})\star\overline{W}_{[1]}(\overline{f},\overline{X}_{[1]},\overline{Y}_{[1]})
\end{aligned}
$$

□

We can generalize the definition of the superposition of coordinates and assume that one of the factors is a tuple of sets of $\overline{\Omega}_{[1]}$-words. Accordingly, the definition of the superposition of coordinates has the form

$$
\begin{aligned}
&\overline{w}_{[1]}(\overline{f},\overline{X}_{[1]},\overline{Y}_{[1]})\circ\overline{W}_{[1]}(\overline{f},\overline{X}_{[1]},\overline{r}_{[1]}\circ\overline{X}_{[1]})\\
=&\overline{W}_{[1]}(\overline{f},\overline{X}_{[1]},\overline{Y}_{[1]})\circ\overline{w}_{[1]}(\overline{f},\overline{X}_{[1]},\overline{r}_{[1]}\circ\overline{X}_{[1]})\\
=&\overline{w}_{[1]}(\overline{f},\overline{X}_{[1]},\overline{r}_{[1]}\circ\overline{Y}_{[1]})
\end{aligned}
$$

The following forms of writing an image of a tuple of sets $\overline{Y}_{[1]}$ under endomorphism $\overline{r}_{[1]}$ are equivalent

$$
\begin{aligned}
\overline{r}_{[1]}\circ\overline{Y}_{[1]}&=\overline{W}_{[1]}(\overline{f},\overline{X}_{[1]},Y)\circ(\overline{r}_{[1]}\circ\overline{X}_{[1]})\\
&=\overline{W}_{[1]}(\overline{f},\overline{X}_{[1]},\overline{Y}_{[1]})\circ(\overline{W}_{[1]}(\overline{f},\overline{X}_{[1]},\overline{r}_{[1]}\circ\overline{X}_{[1]})\circ\overline{X}_{[1]})
\end{aligned}
\tag{4.5.30}
$$

From equations (4.5.26), (4.5.30), it follows that

$$
\begin{aligned}
&(\overline{W}_{[1]}(\overline{f},\overline{X}_{[1]},\overline{Y}_{[1]})\circ\overline{W}_{[1]}(\overline{f},\overline{X}_{[1]},\overline{r}_{[1]}\circ\overline{X}_{[1]}))\circ\overline{X}_{[1]}\\
=&\overline{W}_{[1]}(\overline{f},\overline{X}_{[1]},\overline{Y}_{[1]})\circ(\overline{W}_{[1]}(\overline{f},\overline{X}_{[1]},\overline{r}_{[1]}\circ\overline{X}_{[1]})\circ\overline{X}_{[1]})
\end{aligned}
\tag{4.5.31}
$$

The equation (4.5.31) is associative law for composition and allows us to write expression

$$\overline{W}_{[1]}(\overline{f},\overline{X}_{[1]},\overline{Y}_{[1]})\circ\overline{W}_{[1]}(\overline{f},\overline{X}_{[1]},\overline{r}_{[1]}\circ\overline{X}_{[1]})\circ X$$

without brackets.

Definition 4.5.22. Let $(X_2\subset A_2,...,X_n\subset A_n)$ be tuple of generating sets of tower of representations $\overline{f}$. Let map $\overline{r}_{[1]}$ be endomorphism of tower of representations $\overline{f}$. Let the tuple of sets $\overline{X}'_{[1]}=\overline{r}_{[1]}\circ\overline{X}_{[1]}$ be image of tuple of sets $\overline{X}_{[1]}$ under the map $\overline{r}_{[1]}$. Endomorphism $\overline{r}_{[1]}$ of tower of representations $\overline{f}$ is called **regular on tuple of generating sets** $\overline{X}_{[1]}$, if the tuple of sets $\overline{X}'_{[1]}$ is tuple of generating sets of tower of representations $\overline{f}$. Otherwise, endomorphism $\overline{r}_{[1]}$ is called **singular on tuple of generating sets** $\overline{X}_{[1]}$, □

Definition 4.5.23. Endomorphism of tower of representations $\overline{f}$ is called **regular**, if it is regular on any tuple of generating sets. □

Theorem 4.5.24. *Automorphism $\overline{r}$ of tower of representations $\overline{f}$ is regular endomorphism.*

Proof. Let $\overline{X}_{[1]}$ be tuple of generating sets of tower of representations $\overline{f}$. Let $\overline{X}'_{[1]}=\overline{r}_{[1]}\circ\overline{X}_{[1]}$.

According to theorem 4.5.18 endomorphism $\overline{r}_{[1]}$ forms the map of $\overline{\Omega}_{[1]}$-words $\overline{w}_{[1]}[\overline{f},\overline{X}_{[1]},\overline{r}_{[1]}]$.

Let $\overline{a}'_{[1]}\in\overline{A}_{[1]}$. Since $\overline{r}_{[1]}$ is automorphism, then there exists $\overline{a}_{[1]}\in\overline{A}_{[1]}$, $\overline{r}_{[1]}\circ\overline{a}_{[1]}=\overline{a}'_{[1]}$. According to definition 4.5.8, $\overline{w}_{[1]}(\overline{f},\overline{X}_{[1]},\overline{a}_{[1]})$ is tuple of words representing $\overline{a}_{[1]}$ relative to tuple of generating sets $\overline{X}_{[1]}$. According to theorem 4.5.18, $\overline{w}_{[1]}(\overline{f},\overline{X}'_{[1]},\overline{a}'_{[1]})$ is tuple of words representing $\overline{a}'_{[1]}$ relative to tuple of sets $\overline{X}'_{[1]}$

$$\overline{w}_{[1]}(\overline{f},\overline{X}'_{[1]},\overline{a}'_{[1]})=\overline{w}_{[1]}[\overline{f},\overline{X}_{[1]},\overline{r}](\overline{w}_{[1]}(\overline{f},\overline{X}_{[1]},\overline{a}_{[1]}))$$

Therefore, $\overline{X}'_{[1]}$ is generating set of representation $\overline{f}$. According to definition 4.5.23, automorphism $\overline{r}_{[1]}$ is regular. □

## 4.6. Basis of Tower of Representations

Definition 4.6.1. If the tuple of sets $\overline{X}_{[1]}$ is tuple of generating sets of tower of representations $\overline{f}$, then any tuple of sets $\overline{Y}_{[1]}$, $X_i\subset Y_i\subset A_i$, $i=2$, ..., $n$, also is tuple of generating sets of tower of representations $\overline{f}$. If there exists tuple of minimal sets $\overline{X}_{[1]}$ generating the tower of representations $\overline{f}$, then the tuple of sets $\overline{X}_{[1]}$ is called **basis of tower of representations** $\overline{f}$. □

Theorem 4.6.2. *We define a basis of tower of representations by induction over $n$. For $n=2$, the basis of tower of representations is the basis of representation $f_{1,2}$. If tuple of sets $\overline{X}_{[1,n]}$ is basis of tower of representations $\overline{f}_{[1]}$, then the tuple of generating sets $\overline{X}_{[1]}$ of tower of representations $\overline{f}$ is basis iff for any $a_n\in X_n$ the tuple of sets $(X_2,...,X_{n-1},X_n\setminus\{a_n\})$ is not tuple of generating sets of tower of representations $\overline{f}$.*

Proof. For $n=2$, the statement of the theorem is corollary of the theorem 2.7.2.

Let $n>2$. Let $\overline{X}_{[1]}$ be tuple of generating sets of tower of representations $\overline{f}$. Let tuple of sets $\overline{X}_{[1,n]}$ be basis of tower of representations $\overline{f}_{[n]}$. Assume that for some $a_n\in X_n$ there exist word

$$w_n=w_n(a_n,\overline{f},(X_1,...,X_{n-1},X_n\setminus\{a_n\}))\tag{4.6.1}$$

Consider $a'_n\in A_n$ such that it has word

$$w'_n=w_n(a'_n,\overline{f},\overline{X}_{[1]})\tag{4.6.2}$$

that depends on $a_n$. According to the definition 2.6.6, any occurrence $a_n$ into word $w'_n$ can be substituted by the word $w_n$. Therefore, the word $w'_n$ does not depend on $a_n$, and the tuple of sets $(X_2,...,X_{n-1},X_n\setminus\{a_n\})$ is the tuple of generating sets of tower of representations $\overline{f}$. Therefore, $\overline{X}_{[1]}$ is not basis of the tower of representations $\overline{f}$. □

Remark 4.6.3. The proof of the theorem 4.6.2 gives us effective method for constructing the basis of tower of representations $\overline{f}$. We start to build a basis in the lowest layer. When the basis is constructed in layer $i$, $i=2$, ..., $n-1$, we can proceed to the construction of basis in layer $i+1$. □

Remark 4.6.4. We write a basis also in following form

$$X_k=(x_k,x_k\in X_k)\quad k=2,...,n$$

If basis is finite, then we also use notation

$$X_k=(x_{k\cdot i},i\in I_k)=(x_{k\cdot 1},...,x_{k\cdot p_k})\quad k=2,...,n$$

□

Remark 4.6.5. Coordinates $a_k\in A_k$ relative to basis $\overline{X}_{[1]}$ of the tower of representations $\overline{f}$ are defined ambiguously. Just as in the remark 2.7.5, we consider **tuple of equivalence** $\overline{\rho}_{[1]}(\overline{f})$ **generated by tower of representations** $\overline{f}$. If during the construction, we obtain the equality of two $\Omega_k$-word relative to given basis, then we can say without worrying about the equivalence $\rho_k(\overline{f})$ that these $\Omega_k$-words are equal. □

Theorem 4.6.6. *Automorphism of the tower of representations $\overline{f}$ maps a basis of the tower of representations $\overline{f}$ into basis.*

Proof. For $n=2$, the statement of theorem is corollary of the theorem 4.6.6.

Let the map $\overline{r}$ be automorphism of the tower of representations $\overline{f}$. Let the tuple of sets $\overline{X}_{[1]}$ be a basis of the tower of representations $\overline{f}$. Let $\overline{X}'_{[1]}=\overline{r}_{[1]}(\overline{X}_{[1]})$.

Assume that the tuple of sets $\overline{X}'_{[1]}$ is not basis. According to the theorem 4.6.2 there exist $i$, $i=2$, ..., $n$, and $a'_i\in X'_i$ such that the tuple of sets $\overline{X}''_{[1]}$[4.12] is tuple of generating sets of the tower of representations $\overline{f}$. According to the theorem 4.3.3 the map $\overline{r}^{-1}$ is automorphism of the tower of representations $\overline{f}$. According to the theorem 4.5.24 and definition 4.5.23, the tuple of sets $\overline{X}'''_{[1]}$[4.13] is tuple of generating sets of the tower of representations $\overline{f}$. The contradiction completes the proof of the theorem. □

Theorem 4.6.7. *Let $\overline{X}_{[1]}$ be the basis of the representation $\overline{f}$. Let*

$$\overline{R}_{[1]}:\overline{X}_{[1]}\to\overline{X}'_{[1]}$$

*be arbitrary map of the tuple of sets $\overline{X}_{[1]}$. Consider the map of $\overline{\Omega}_{[1]}$-words*

$$\overline{w}_{[1]}[\overline{f},\overline{X}_{[1]},\overline{R}_{[1]}]:\overline{w}_{[1]}(\overline{f},\overline{X}_{[1]})\to\overline{w}_{[1]}(\overline{f},\overline{X}'_{[1]})$$

*that satisfies conditions (1), (2), (3) of the theorem 4.5.9 and such that*

$$\overline{x}_{[1]}\in\overline{X}_{[1]}=>w[\overline{f},\overline{X}_{[1]},\overline{R}_{[1]}](\overline{x}_{[1]})=\overline{R}_{[1]}(\overline{x}_{[1]})$$

*There exists unique endomorphism of representation $\overline{f}$*[4.14]

$$\overline{r}_{[1]}:\overline{A}_{[1]}\to\overline{A}_{[1]}$$

*defined by rule*

$$\overline{r}_{[1]}\circ\overline{a}_{[1]}=w[\overline{f},\overline{X}_{[1]},\overline{R}_{[1]}](\overline{w}_{[1]}(\overline{f},\overline{X}_{[1]},\overline{a}_{[1]}))$$

Proof. The statement of theorem is corollary theorems 2.6.7, 2.6.11. □

[4.12] $X''_j=X'_j$, $j\ne i$, $X''_i=X'_i\setminus\{x'_i\}$
[4.13] $X'''_j=X_j$, $j\ne i$, $X'''_i=X_i\setminus\{x_i\}$
[4.14] This statement is similar to the theorem [1]-4.1, p. 135.

Corollary 4.6.8. *Let* $\overline{X}_{[1]}$, $\overline{X'}_{[1]}$ *be the bases of the representation* $\overline{f}$. *Let* $\overline{r}_{[1]}$ *be the automorphism of the representation* $\overline{f}$ *such that* $\overline{X'}_{[1]} = \overline{r}_{[1]} \circ \overline{X}_{[1]}$ . *Automorphism* $\overline{r}_{[1]}$ *is uniquely defined.* □

Remark 4.6.9. The theorem 4.6.7, as well as the theorem 4.5.14, is the theorem of extension of map. However in this theorem, $\overline{X}_{[1]}$ is not arbitrary tuple of generating sets of the tower of representations, but basis. According to remarks 4.6.3, 2.7.3, we cannot determine coordinates of any element of basis through the remaining elements of the same basis. Therefore, we do not need to coordinate the map of the basis with representation. □

Theorem 4.6.10. *The set of coordinates* $\overline{W}_{[1]}(\overline{f}, \overline{X}_{[1]}, \overline{X}_{[1]})$ *corresponds to identity transformation*

$$\overline{W}_{[1]}[\overline{f}, \overline{X}_{[1]}, \overline{E}_{[1]}] = \overline{W}_{[1]}(\overline{f}, \overline{X}_{[1]}, \overline{X}_{[1]})$$

Proof. The statement of the theorem follows from the equation

$$\overline{a}_{[1]} = \overline{W}_{[1]}(\overline{f}, \overline{X}_{[1]}, \overline{a}_{[1]}) \circ \overline{X}_{[1]} = \overline{W}_{[1]}(\overline{f}, \overline{X}_{[1]}, \overline{a}_{[1]}) \circ \overline{W}_{[1]}(\overline{f}, \overline{X}_{[1]}, \overline{X}_{[1]}) \circ \overline{X}_{[1]}$$

□

Theorem 4.6.11. *Let* $\overline{W}_{[1]}(\overline{f}, \overline{X}_{[1]}, \overline{r}_{[1]} \circ \overline{X}_{[1]})$ *be the set of coordinates of automorphism* $\overline{r}_{[1]}$. *There exists set of coordinates* $\overline{W}_{[1]}(\overline{f}, \overline{r}_{[1]} \circ \overline{X}_{[1]}, \overline{X}_{[1]})$ *, corresponding to automorphism* $\overline{r}_{[1]}^{-1}$. *The set of coordinates* $\overline{W}_{[1]}(\overline{f}, \overline{r}_{[1]} \circ \overline{X}_{[1]}, \overline{X}_{[1]})$ *satisfy to equation*[4.15]

$$\overline{W}_{[1]}(\overline{f}, \overline{X}_{[1]}, \overline{r}_{[1]} \circ \overline{X}_{[1]}) \circ \overline{W}_{[1]}(\overline{f}, \overline{r}_{[1]} \circ \overline{X}_{[1]}, \overline{X}_{[1]}) = \overline{W}_{[1]}(\overline{f}, \overline{X}_{[1]}, \overline{r}_{[1]} \circ \overline{X}_{[1]}) \tag{4.6.3}$$

$$\overline{W}_{[1]}[\overline{f}, \overline{X}_{[1]}, \overline{r}_{[1]}^{-1}] = \overline{W}_{[1]}[\overline{f}, \overline{X}_{[1]}, , R]^{-1} = \overline{W}_{[1]}(\overline{f}, \overline{r}_{[1]} \circ \overline{X}_{[1]}, \overline{X}_{[1]})$$

Proof. Since $\overline{R}_{[1]}$ is automorphism of the tower of representations $\overline{f}$, then, according to the theorem 4.6.6, the set $\overline{r}_{[1]} \circ \overline{X}_{[1]}$ is a basis of the tower of representations $\overline{f}$. Therefore, there exists the set of coordinates $\overline{W}_{[1]}(\overline{f}, \overline{r}_{[1]} \circ \overline{X}_{[1]}, \overline{X}_{[1]})$ . The equation (4.6.3) follows from the chain of equations

$$\begin{aligned}
&\overline{W}_{[1]}(\overline{f}, \overline{X}_{[1]}, \overline{r}_{[1]} \circ \overline{X}_{[1]}) \circ \overline{W}_{[1]}(\overline{f}, \overline{r}_{[1]} \circ \overline{X}_{[1]}, \overline{X}_{[1]})\\
=&\overline{W}_{[1]}(\overline{f}, \overline{X}_{[1]}, \overline{r}_{[1]} \circ \overline{X}_{[1]}) \circ \overline{W}_{[1]}(\overline{f}, \overline{X}_{[1]}, \overline{r}_{[1]}^{-1} \circ \overline{X}_{[1]})\\
=&\overline{W}_{[1]}(\overline{f}, \overline{X}_{[1]}, \overline{r}_{[1]}^{-1} \circ \overline{r}_{[1]} \circ \overline{X}_{[1]}) = \overline{W}_{[1]}(\overline{f}, \overline{X}_{[1]}, \overline{X}_{[1]})
\end{aligned}$$

□

Theorem 4.6.12. *The group of automorphisms* $G(\overline{f})$ *of the tower of effective representations* $\overline{f}$ *generates effective representation in the tower of representations* $\overline{f}$.

Proof. From the corollary 4.6.8, it follows that if automorphism $\overline{r}_{[1]}$ maps a basis $\overline{X}_{[1]}$ into a basis $\overline{X'}_{[1]}$, then the set of coordinates $\overline{W}_{[1]}(\overline{f}, \overline{X}_{[1]}, \overline{X'}_{[1]})$ uniquely determines an automorphism $\overline{r}_{[1]}$. From the theorem 4.5.18, it follows that the set of coordinates $\overline{W}_{[1]}(\overline{f}, \overline{X}_{[1]}, \overline{X'}_{[1]})$ determines the map of coordinates relative to the basis $\overline{X}_{[1]}$ under automorphism of the tower of representations $\overline{f}$. From the equation (4.5.30), it follows that automorphism $\overline{r}_{[1]}$ acts from the right on elements

[4.15] See also remark 2.7.12.

of $\Omega_k$-algebra $A_k$, $k = 2, ..., n$. From the equation (4.5.28), it follows that the representation of group is $*T$-representation. According to the theorem 4.6.10 the set of coordinates $\overline{W}_{[1]}(\overline{f}, \overline{X}_{[1]}, \overline{X}_{[1]})$ corresponds to identity transformation. From the theorem 4.6.11, it follows that the set of coordinates $\overline{W}_{[1]}(\overline{f}, \overline{r}_{[1]} \circ \overline{X}_{[1]}, \overline{X}_{[1]})$ corresponds to transformation, inverse to transformation $\overline{W}_{[1]}(\overline{f}, \overline{X}_{[1]}, \overline{r}_{[1]} \circ \overline{X}_{[1]})$. □

## 4.7. Examples of Basis of Tower of Representations

**4.7.1. Affine Space.** Consider affine space $\overset{\circ}{A}$ (the subsection 4.4.1).

The Abelian group $\overline{A}$ acts single transitive on the set $\overset{\circ}{A}$. From construction in subsection 2.8.2, it follows that the basis of the set $\overset{\circ}{A}$ relative to representation of the Abelian group $\overline{A}$ consists of one point. This point is usually denoted by the letter $O$ and is called **origin of coordinate system of affine space**. Therefore, an arbitrary point $A \in \overset{\circ}{A}$ can be represented using vector $\overrightarrow{OA} \in \overline{A}$

Let $\overline{\overline{e}}$ be the basis of the vector space $\overline{A}$. Then the vector $\overrightarrow{OA}$ has form

$$\overrightarrow{OA} = a^i \overline{e}_i$$

The set $(a_i, i \in I)$ is called **coordinates of point** $A$ **of affine space** $\overset{\circ}{A}$ **relative to basis** $(O, \overline{\overline{e}})$.

**4.7.2. Coordinates of Vector in $R$-Module.** Let a left module $\overline{M}$ over the ring $R$ ([1], p. 117) have a finite basis

$$\overline{\overline{e}} = (\overline{e}_1, ..., \overline{e}_n)$$

According to the definition 2.7.1, we can present any vector $\overline{v} \in \overline{M}$ in the form

$$\overline{v} = v^i \overline{e}_i \tag{4.7.1}$$

According to the theorem 2.7.2, $\overline{e}_i \in \overline{\overline{e}}$ cannot be presented such way using other vectors of basis. However this does not mean that we cannot write equation

$$c^i \overline{e}_i = 0 \tag{4.7.2}$$

where not all $c^i$ are 0. To make equation (4.7.2) possible, it is necessary to require that the numbers $c^1$, ..., $c^n$ are irreversible in the ring $R$. From equations (4.7.1), (4.7.2), it follows that besides expansion (4.7.1) relative to basis $\overline{\overline{e}}$, the vector $\overline{v}$ has also expansion

$$\overline{v} = (v^i + c^i) \overline{e}_i \tag{4.7.3}$$

It is easy to see that the set of tuples $(c^1, ..., c^n)$ generates a left $R$-module. Therefore, the set of coordinates in $R$-module $\overline{M}$ is affinne space, and coordinates of vector $\overline{v}$ form an an affine plane in this space.

It is possible that the set of coordinates of an arbitrary representation of universal algebra can also be regarded as a tower of representations.

**4.7.3. Basis of $R$-algebra.** Usually, when we consider the $R$-algebra $A$, we choose a basis $\overline{\overline{e}}$ of corresponding $R$-module $A$. This choice is convenient because if $R$ is commutative division algebra, then expansion of the vector is unique relative to basis of $R$-vector space. This, in particular, allows us to define product by specifying structural constants of algebra relative to given basis.

In general, the basis of $R$-module $A$ may appear a generating set. For instance, if in vector space $H$, where we consider quaternion algebra over real field, we consider the basis

$$e_{\mathbf{0}} = 1 \quad e_{\mathbf{1}} = i \quad e_{\mathbf{2}} = j \quad e_{\mathbf{3}} = k \tag{4.7.4}$$

then in the algebra $H$ the following equation is true

$$\begin{gathered} e_{\mathbf{0}} = -e_{\mathbf{1}}e_{\mathbf{1}} = -e_{\mathbf{2}}e_{\mathbf{2}} \\ e_{\mathbf{3}} = e_{\mathbf{1}}e_{\mathbf{2}} \end{gathered} \tag{4.7.5}$$

Therefore, the set $(e_{\mathbf{1}}, e_{\mathbf{2}})$ is a basis of algebra $H$. Ambiguity of representation of quaternion relative to the given basis is consequence of the equation (4.7.5). Namely, we can present a quaternion $a \in H$ as

$$a = (a^{\mathbf{0}} - a^{\mathbf{4}})e_{\mathbf{1}}e_{\mathbf{1}} + a^{\mathbf{4}}e_{\mathbf{2}}e_{\mathbf{2}} + a^{\mathbf{1}}e_{\mathbf{1}} + a^{\mathbf{2}}e_{\mathbf{2}} + a^{\mathbf{3}}e_{\mathbf{1}}e_{\mathbf{2}}$$

where $a^{\mathbf{4}}$ is arbitrary.

**4.7.4. Module over Algebra.** Let $A_3$ be free module over the $A_1$-algebra $A_2$ (the subsection 4.4.4). Let $\overline{\overline{e}}_{32}$ be basis of algebra $A_3$ over algebra $A_2$. A vector $a_3 \in A_3$ has representation

$$a_3 = a_3^j \ \overline{e}_{32\cdot j} = \begin{pmatrix} a_3^1 & ... & a_3^n \end{pmatrix} \begin{pmatrix} \overline{e}_{32\cdot 1} \\ ... \\ \overline{e}_{32\cdot n} \end{pmatrix} \tag{4.7.6}$$

Let $\overline{\overline{e}}_{21}$ be basis of algebra $A_2$ over field $A_1$. Because $a_3^j \in A_2$, we can write their coordinates relative to basis $\overline{\overline{e}}_{21}$

$$a_3^j = a_3^{ji} \ \overline{e}_{21\cdot i} = \begin{pmatrix} a_3^{j1} & ... & a_3^{jm} \end{pmatrix} \begin{pmatrix} \overline{e}_{21\cdot n} \\ ... \\ \overline{e}_{21\cdot m} \end{pmatrix} \tag{4.7.7}$$

From equations (4.7.6), (4.7.7) it follows

(4.7.8)

$$\begin{gathered} a_3 = a_3^{ji} \ \overline{e}_{21\cdot i} \ \overline{e}_{32\cdot j} = \\ \begin{pmatrix} \begin{pmatrix} a_3^{11} & ... & a_3^{1m} \end{pmatrix} \begin{pmatrix} \overline{e}_{21\cdot n} \\ ... \\ \overline{e}_{21\cdot m} \end{pmatrix} & ... & \begin{pmatrix} a_3^{n1} & ... & a_3^{nm} \end{pmatrix} \begin{pmatrix} \overline{e}_{21\cdot n} \\ ... \\ \overline{e}_{21\cdot m} \end{pmatrix} \end{pmatrix} \begin{pmatrix} \overline{e}_{32\cdot 1} \\ ... \\ \overline{e}_{32\cdot n} \end{pmatrix} \end{gathered}$$

Equation (4.7.8) shows the structure of coordinates in vector space $A_3$ over field $A_1$. It is easy to see that vectors

$$\overline{e}_{31\cdot ij} = \overline{e}_{21\cdot i} \ \overline{e}_{32\cdot j}$$

are linear independent over field $A_1$. Therefore, we build the basis $\overline{\overline{e}}_{31}$ of vector space $A_3$ over field $A_1$. Therefore, we can rewrite equation (4.7.8) as

$$a_3 = a_3^{ji}\ \overline{e}_{31\cdot ij} = \begin{pmatrix} a_3^{11} & ... & a_3^{1m} & ... & a_3^{n1} & ... & a_3^{nm} \end{pmatrix} \begin{pmatrix} \overline{e}_{31\cdot 11} \\ ... \\ \overline{e}_{31\cdot 1m} \\ ... \\ \overline{e}_{31\cdot n1} \\ ... \\ \overline{e}_{31\cdot nm} \end{pmatrix} \tag{4.7.9}$$

It is easy to see that we can identify vector $\overline{e}_{31\cdot ij}$ with tensor product $\overline{e}_{21\cdot i} \otimes \overline{e}_{32\cdot j}$.

## 4.8. Representations in Category

Definition 4.8.1. Let for any objects $B$ and $C$ of category $\mathcal{B}$ the structure of $\Omega$-algebra is defined on the set of morphisms $Mor(B, C)$. The set of homomorphisms of $\Omega$-algebra

$$f_{B,C} : A \to Mor(B, C)$$

is called **representation of $\Omega$-algebra $A$ in category $\mathcal{B}$**. □

If we assume that the set $Mor(B, C)$ is defined only when $B = C$, then we get definition of representation. The difference in definitions is that we do not restrict ourselves by transformations of the set $B$, but consider $\Omega$-algebra of maps from the set $B$ into set $C$. At first sight there is no fundamental differences between considered theories. However we can see that representation of vector space in category of bundles lead us to connection theory.

CHAPTER 5

# References

[1] Serge Lang, Algebra, Springer, 2002

[2] S. Burris, H.P. Sankappanavar, A Course in Universal Algebra, Springer-Verlag (March, 1982),
eprint http://www.math.uwaterloo.ca/ snburris/htdocs/ualg.html
(The Millennium Edition)

[3] P. K. Rashevsky, Riemann Geometry and Tensor Calculus,
Moscow, Nauka, 1967

[4] A. G. Kurosh, Lectures on General Algebra, Chelsea Pub Co, 1965

[5] Aleksandr Vasilevich Mikhalev; Günter Pilz; The concise handbook of algebra;
Kluwer Academic Publishers, 2002

[6] I. R. Shafarevich, Basic notions of algebra,
Translated from the Russian by M. Reid,
Springer, 2005

[7] Aleks Kleyn, Lectures on Linear Algebra over Division Ring,
eprint arXiv:math.GM/0701238 (2010)

[8] Aleks Kleyn, Fibered Universal Algebra,
eprint arXiv:math.DG/0702561 (2007)

[9] Aleks Kleyn, Introduction into Geometry over Division Ring,
eprint arXiv:0906.0135 (2010)

[10] Aleks Kleyn, Linear Maps of Free Algebra,
eprint arXiv:1003.1544 (2010)

[11] Aleks Kleyn.
Linear Algebra over Division Ring: System of Linear Equations.
CreateSpace Independent Publishing Platform, 2012; ISBN-13: 978-1477631812

[12] John C. Baez, The Octonions,
eprint arXiv:math.RA/0105155 (2002)

[13] Paul M. Cohn, Universal Algebra, Springer, 1981

[14] Richard D. Schafer, An Introduction to Nonassociative Algebras, Dover Publications, Inc., New York, 1995

CHAPTER 6

# Index

CHAPTER 7

# Special Symbols and Notations

$w(f, X)$ set of $\Omega_2$-words of representation $J(f, X)$ 26

$\overline{w}(\overline{f}, \overline{X}_{[1]})$ set of tuples of $\overline{\Omega}_{[1]}$-words of tower of representations $\overline{J}(\overline{f}, \overline{X}_{[1]})$ 60

$\overline{w}(\overline{f}, \overline{X}_{[1]}, \overline{a})$ tuple of words of element $\overline{a}*$ relative to tuple of sets $\overline{X}_{[1]}$ 60

$\delta$ identical transformation 7

# Представление универсальной алгебры

Александр Клейн

Aleks_Kleyn@MailAPS.org
http://AleksKleyn.dyndns-home.com:4080/
http://sites.google.com/site/AleksKleyn/
http://arxiv.org/a/kleyn_a_1
http://AleksKleyn.blogspot.com/

Аннотация. Теория представлений универсальной алгебры является естественным развитием теории универсальной алгебры. Морфизм представления - это отображение, сохраняющее структуру представления. Изучение морфизмов представлений ведёт к понятиям множества образующих и базиса представления. В книге рассмотрено понятие башни левосторонних представлений $\Omega_i$-алгебр, $i = 1, ..., n$, как множество согласованных левосторонних представлений $\Omega_i$-алгебр.

# Оглавление

Глава 1

# Предисловие

## 1.1. Предисловие

В статьях я часто рассматриваю вопросы, связанные с представлением универсальной алгебры. Вначале это были небольшие наброски, которые я многократно исправлял и переписывал. Но постепенно появлялись новые наблюдения. В результате вспомогательный инструмент превратился в стройную теорию.

Я это понял, когда я работал над книгой [9], и решил посвятить отдельную книгу вопросам, связанным с представлением универсальной алгебры. Изучение теории представлений универсальной алгебры показывает, что эта теория имеет много общего с теорией универсальной алгебры.

Основным толчком к более глубокому изучению представлений универсальной алгебры послужило определение векторного пространства как представление поля в абелевой группе. Я обратил внимание, что это определение меняет роль линейного отображения. По сути, линейное отображение - это отображение, которое сохраняет структуру представления. Эту конструкцию легко обобщить на произвольное представление универсальной алгебры. Таким образом появилось понятие морфизма представлений.

Множество невырожденных автоморфизмов векторного пространства порождает группу. Эта группа действует однотранзитивно на множестве базисов векторного пространства. Это утверждение является фундаментом теории инвариантов векторного пространства.

Возникает естественный вопрос. Можно ли обобщить эту конструкцию на произвольное представление? Базис - это не единственное множество, которое порождает векторное пространство. Если мы к множеству векторов базиса добавим произвольный вектор, то новое множество по прежнему порождает тоже самое векторное пространство, но базисом не является. Это утверждение является исходной точкой, от которой я начал изучение множества образующих представления. Множество образующих представления - это ещё одна интересная параллель теории представлений с теорией универсальной алгебры.

Множество автоморфизмов представления является лупой. Неассоциативность произведения порождает многочисленные вопросы, которые требуют дополнительное исследование. Все эти вопросы ведут к необходимости понимания теории инвариантов заданного представления.

Если мы рассматриваем теорию представлений универсальной алгебры как расширение теории универсальной алгебры, то почему не рассмотреть представление одного представления в другом представлении. Так появилась концепция башни представлений. Самый удивительный факт - это то, что все отображения в башне представлений действуют согласовано.

## 1.2. Соглашения

Соглашение 1.2.1. *Функция и отображение - синонимы. Однако существует традиция соответствие между кольцами или векторными пространствами называть отображением, а отображение поля действительных чисел или алгебры кватернионов называть функцией. Я тоже следую этой традиции, хотя встречается текст, в котором неясно, какому термину надо отдать предпочтение.* □

Соглашение 1.2.2. *Так как число универсальных алгебр в башне представлений переменно, то мы будем пользоваться векторными обозначениями для башни представлений. Множество* $(A_1, ..., A_n)$ $\Omega_i$*-алгебр* $A_i$, $i = 1, ..., n$ *мы будем обозначать* $\overline{A}$. *Множество представлений* $(f_{1,2}, ..., f_{n-1,n})$ *этих алгебр мы будем обозначать* $\overline{f}$. *Так как разные алгебры имеют разный тип, мы также будем говорить о множестве* $\overline{\Omega}$*-алгебр. По отношению к множеству* $\overline{A}$ *мы также будем пользоваться матричными обозначениями, предложенными в разделе [7]-2.1. Например, символом* $\overline{A}_{[1]}$ *мы будем обозначать множество* $\overline{\Omega}$*-алгебр* $(A_2, ..., A_n)$. *В соответствующем обозначении* $(\overline{A}_{[1]}, \overline{f})$ *башни представлений подразумевается, что* $\overline{f} = (f_{2,3}, ..., f_{n-1,n})$. □

Соглашение 1.2.3. *Так как мы пользуемся векторными обозначениями для элементов башни представлений, необходимо соглашение о записи операций. Предполагается, что операции выполняются покомпонентно. Например,*

$$\overline{r}(\overline{a}) = (r_1(a_1), ..., r_n(a_n))$$

□

Замечание 1.2.4. Я считаю диаграммы отображений важным инструментом. Однако временами возникает желание увидеть диаграмму трёх мерной, что увеличило бы её выразительную мощность. Кто знает какие сюрпризы готовит будущее. В 1992 на конференции в Казани я рассказывал своим коллегам какие преимущества имеет компьютерная подготовка статей. Спустя 8 лет из письма из Казани я узнал, что теперь можно готовить статьи с помощью LaTeX. □

Глава 2

# Представление универсальной алгебры

## 2.1. Представление универсальной алгебры

Определение 2.1.1. Пусть на множестве $M$ определена структура $\Omega_2$-алгебры ([2, 13]). Эндоморфизм $\Omega_2$-алгебры

$$t : M \rightarrow M$$

называется **преобразованием универсальной алгебры** $M$.[2.1] □

Мы будем обозначать $\delta$ тождественное преобразование.

Определение 2.1.2. Пусть ${}^*M$ - множество **левосторонних преобразований**

$$u' = tu$$

$\Omega_2$-алгебры $M$. Пусть на множестве ${}^*M$ определена структура $\Omega_1$-алгебры. Гомоморфизм

$$f : A \rightarrow {}^*M \tag{2.1.1}$$

$\Omega_1$-алгебры $A$ в $\Omega_1$-алгебру ${}^*M$ называется **левосторонним представлением $\Omega_1$-алгебры** $A$ или $A*$**-представлением** в $\Omega_2$-алгебре $M$. □

Определение 2.1.3. Пусть $M^*$ - множество **правосторонних преобразований**

$$u' = ut$$

$\Omega_2$-алгебры $M$. Пусть на множестве $M^*$ определена структура $\Omega_1$-алгебры. Гомоморфизм

$$f : A \rightarrow M^*$$

$\Omega_1$-алгебры $A$ в $\Omega_1$-алгебру $M^*$ называется **правосторонним представлением $\Omega_1$-алгебры** $A$ или $*A$**-представлением** в $\Omega_2$-алгебре $M$. □

Мы распространим на теорию представлений соглашение, описанное в замечаниях [7]-2.2.15, [11]-2.2.15. Мы можем записать принцип двойственности в следующей форме

Теорема 2.1.4 (принцип двойственности). *Любое утверждение, справедливое для левостороннего представления $\Omega_1$-алгебры $A$, будет справедливо для правостороннего представления $\Omega_1$-алгебры $A$.*

[2.1] Если множество операций $\Omega_2$-алгебры пусто, то

$$t : M \rightarrow M$$

является отображением.

Замечание 2.1.5. Существует две формы записи преобразования $\Omega_2$-алгебры $M$. Если мы пользуемся операторной записью, то преобразование $A$ записывается в виде $Aa$ или $aA$, что соответствует левостороннему преобразованию или правостороннему преобразованию. Если мы пользуемся функциональной записью, то преобразование $A$ записывается в виде $A(a)$ независимо от того, это левостороннее или правостороннее преобразование. Эта запись согласована с принципом двойственности.

Это замечание является основой следующего соглашения. Когда мы пользуемся функциональной записью, мы не различаем левостороннее и правостороннее преобразование. Мы будем обозначать ${}^*M$ множество преобразований $\Omega_2$-алгебры $M$. Пусть на множестве ${}^*M$ определена структура $\Omega_1$-алгебры. Пусть $A$ является $\Omega_1$-алгеброй. Мы будем называть гомоморфизм

$$f : A \to {}^*M \tag{2.1.2}$$

**представлением $\Omega_1$-алгебры $A$ в $\Omega_2$-алгебре** $M$. Мы будем также пользоваться записью

$$f : A \xrightarrow{*} M$$

для обозначения представления $\Omega_1$-алгебры $A$ в $\Omega_2$-алгебре $M$.

Соответствие между операторной записью и функциональной записью однозначно. Мы можем выбирать любую форму записи, которая удобна для изложения конкретной темы. □

Существует несколько способов описать представление. Мы можем указать отображение $f$, имея в виду что область определения - это $\Omega_1$-алгебра $A$ и область значений - это $\Omega_1$-алгебра ${}^*M$. Либо мы можем указать $\Omega_1$-алгебру $A$ и $\Omega_2$-алгебру $M$, имея в виду что нам известна структура отображения $f$.[2.2]

Диаграмма

$$\begin{array}{ccc} M & \xrightarrow{\quad f(a) \quad} & M \\ & f \Big\Uparrow & \\ & A & \end{array}$$

означает, что мы рассматриваем представление $\Omega_1$-алгебры $A$. Отображение $f(a)$ является образом $a \in A$.

Определение 2.1.6. Мы будем называть представление $\Omega_1$-алгебры $A$ **эффективным**, если отображение (2.1.2) - изоморфизм $\Omega_1$-алгебры $A$ в ${}^*M$. □

Замечание 2.1.7. Если левостороннее представление $\Omega_1$-алгебры эффективно, мы можем отождествлять элемент $\Omega_1$-алгебры с его образом и записывать левостороннее преобразование, порождённое элементом $a \in A$, в форме

$$v' = av$$

Если правостороннее представление $\Omega_1$-алгебры эффективно, мы можем отождествлять элемент $\Omega_1$-алгебры с его образом и записывать правостороннее преобразование, порождённое элементом $a \in A$, в форме

$$v' = va$$

□

[2.2] Например, мы рассматриваем векторное пространство $\overline{V}$ над полем $D$ ( определение [7]-5.1.4 ).

Определение 2.1.8. Мы будем называть представление $\Omega_1$-алгебры **транзитивным**, если для любых $a, b \in V$ существует такое $g$, что

$$a = f(g)(b)$$

Мы будем называть представление $\Omega_1$-алгебры **однотранзитивным**, если оно транзитивно и эффективно. □

Теорема 2.1.9. *Представление однотранзитивно тогда и только тогда, когда для любых $a, b \in M$ существует одно и только одно $g \in A$ такое, что $a = f(g)(b)$*

Доказательство. Следствие определений 2.1.6 и 2.1.8. □

## 2.2. Морфизм представлений универсальной алгебры

Теорема 2.2.1. *Пусть $A$ и $B$ - $\Omega_1$-алгебры. Представление $\Omega_1$-алгебры $B$*

$$g : B \to {}^{\star}M$$

*и гомоморфизм $\Omega_1$-алгебры*

$$h : A \to B \tag{2.2.1}$$

*определяют представление $f$ $\Omega_1$-алгебры $A$*

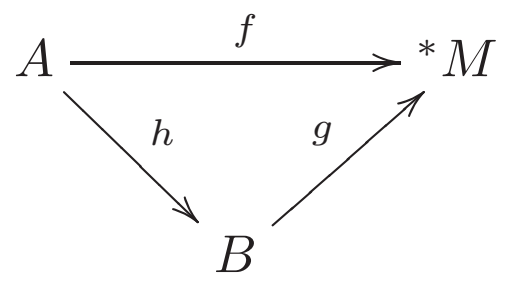

Доказательство. Отображение $f$ является гомоморфизмом $\Omega_1$-алгебры $A$ в $\Omega_1$-алгебру ${}^{*}M$, так как отображение $g$ является гомоморфизмом $\Omega_1$-алгебры $B$ в $\Omega_1$-алгебру ${}^{*}M$. □

Если мы изучаем представление $\Omega_1$-алгебры в $\Omega_2$-алгебрах $M$ и $N$, то нас интересуют отображения из $M$ в $N$, сохраняющие структуру представления.

Определение 2.2.2. Пусть

$$f : A \xrightarrow{*} M$$

представление $\Omega_1$-алгебры $A$ в $\Omega_2$-алгебре $M$ и

$$g : B \xrightarrow{*} N$$

представление $\Omega_1$-алгебры $B$ в $\Omega_2$-алгебре $N$. Пара отображений

$$(r : A \to B, R : M \to N) \tag{2.2.2}$$

таких, что

- $r$ - гомоморфизм $\Omega_1$-алгебры
- $R$ - гомоморфизм $\Omega_2$-алгебры
-
$$R \circ f(a) = g(r(a)) \circ R \tag{2.2.3}$$

называется **морфизмом представлений из $f$ в $g$**. Мы также будем говорить, что определён **морфизм представлений $\Omega_1$-алгебры в $\Omega_2$-алгебре**. □

Замечание 2.2.3. Мы можем рассматривать пару отображений $r$, $R$ как отображение
$$F : A \cup M \to B \cup N$$
такое, что
$$F(A) = B \qquad F(M) = N$$
Поэтому в дальнейшем мы будем говорить, что дано отображение $(r, R)$. □

Определение 2.2.4. Если представления $f$ и $g$ совпадают, то морфизм представлений $(r, R)$ называется **морфизмом представления** $f$. □

Для произвольного $m \in M$ равенство (2.2.3) имеет вид
$$R(f(a)(m)) = g(r(a))(R(m)) \tag{2.2.4}$$

Замечание 2.2.5. Рассмотрим морфизм представлений (2.2.2). Мы можем обозначать элементы множества $B$, пользуясь буквой по образцу $b \in B$. Но если мы хотим показать, что $b$ является образом элемента $a \in A$, мы будем пользоваться обозначением $r(a)$. Таким образом, равенство
$$r(a) = r(a)$$
означает, что $r(a)$ (в левой части равенства) является образом $a \in A$ (в правой части равенства). Пользуясь подобными соображениями, мы будем обозначать элемент множества $N$ в виде $R(m)$. Мы будем следовать этому соглашению, изучая соотношения между гомоморфизмами $\Omega_1$-алгебр и отображениями между множествами, где определены соответствующие представления. □

Замечание 2.2.6. Мы можем интерпретировать (2.2.4) двумя способами

- Пусть преобразование $f(a)$ отображает $m \in M$ в $f(a)(m)$. Тогда преобразование $g(r(a))$ отображает $R(m) \in N$ в $R(f(a)(m))$.
- Мы можем представить морфизм представлений из $f$ в $g$, пользуясь диаграммой

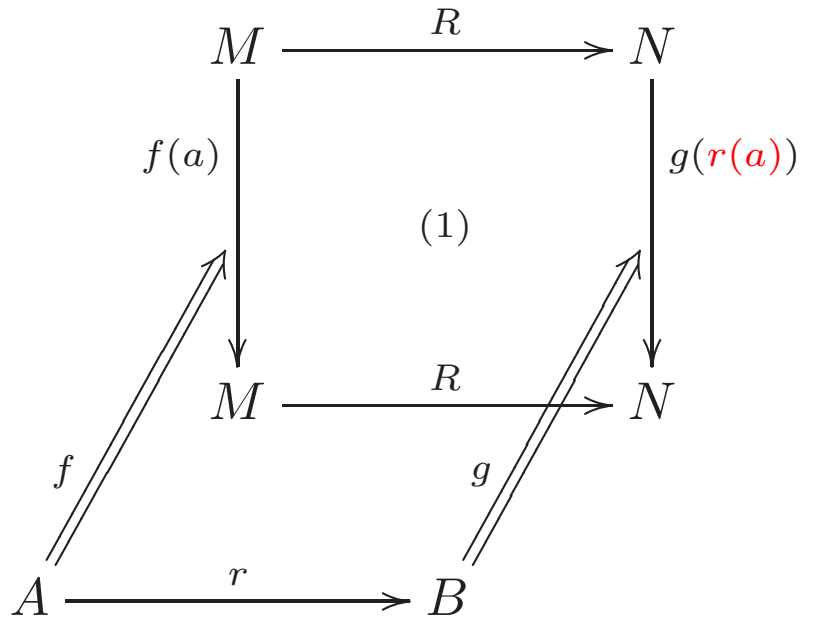

Из (2.2.3) следует, что диаграмма (1) коммутативна.

□

Теорема 2.2.7. *Рассмотрим представление*
$$f : A \twoheadrightarrow\!\!{*}\; M$$
*$\Omega_1$-алгебры $A$ и представление*
$$g : B \twoheadrightarrow\!\!{*}\; N$$

*$\Omega_1$-алгебры $B$. Морфизм*

$$h : A \longrightarrow B \qquad H : M \longrightarrow N$$

*представлений из $f$ в $g$ удовлетворяет соотношению*

$$H \circ \omega(f(a_1), ..., f(a_n)) = \omega(g(h(a_1)), ..., g(h(a_n))) \circ H \tag{2.2.5}$$

*для произвольной $n$-арной операции $\omega$ $\Omega_1$-алгебры.*

Доказательство. Так как $f$ - гомоморфизм, мы имеем

$$H \circ \omega(f(a_1), ..., f(a_n)) = H \circ f(\omega(a_1, ..., a_n)) \tag{2.2.6}$$

Из (2.2.3) и (2.2.6) следует

$$H \circ \omega(f(a_1), ..., f(a_n)) = g(h(\omega(a_1, ..., a_n))) \circ H \tag{2.2.7}$$

Так как $h$ - гомоморфизм, из (2.2.7) следует

$$H \circ \omega(f(a_1), ..., f(a_n)) = g(\omega(h(a_1), ..., h(a_n))) \circ H \tag{2.2.8}$$

Так как $g$ - гомоморфизм, из (2.2.8) следует (2.2.5). □

Теорема 2.2.8. *Пусть отображение*

$$h : A \to B \quad H : M \to N$$

*является морфизмом из представления*

$$f : A \xrightarrow{*} M$$

*$\Omega_1$-алгебры $A$ в представление*

$$g : B \xrightarrow{*} N$$

*$\Omega_1$-алгебры $B$. Если представление $f$ эффективно, то отображение*

$$^*H : {}^*M \to {}^*N$$

*определённое равенством*

$$^*H(f(a)) = g(h(a)) \tag{2.2.9}$$

*является гомоморфизмом $\Omega_1$-алгебры.*

Доказательство. Так как представление $f$ эффективно, то для выбранного преобразования $f(a)$ выбор элемента $a$ определён однозначно. Следовательно, преобразование $g(h(a))$ в равенстве (2.2.9) определено корректно.

Так как $f$ - гомоморфизм, мы имеем

$$^*H(\omega(f(a_1), ..., f(a_n))) = {}^*H(f(\omega(a_1, ..., a_n))) \tag{2.2.10}$$

Из (2.2.9) и (2.2.10) следует

$$^*H(\omega(f(a_1), ..., f(a_n))) = g(h(\omega(a_1, ..., a_n))) \tag{2.2.11}$$

Так как $h$ - гомоморфизм, из (2.2.11) следует

$$^*H(\omega(f(a_1), ..., f(a_n))) = g(\omega(h(a_1), ..., h(a_n))) \tag{2.2.12}$$

Так как $g$ - гомоморфизм,

$$^*H(\omega(f(a_1), ..., f(a_n))) = \omega(g(h(a_1)), ..., g(h(a_n))) = \omega({}^*H(f(a_1)), ..., {}^*H(f(a_n)))$$

следует из (2.2.12). Следовательно, отображение $^*H$ является гомоморфизмом $\Omega_1$-алгебры. □

Теорема 2.2.9. *Пусть*

$$f : A \xrightarrow{*} M$$

*однотранзитивное представление $\Omega_1$-алгебры $A$ и*

$$g : B \xrightarrow{*} N$$

*однотранзитивное представление $\Omega_1$-алгебры $B$. Если отображение*

$$h : A \to B$$

*является гомоморфизмом $\Omega_1$-алгебры, то существует морфизм представлений из $f$ в $g$*

$$h : A \to B \quad H : M \to N$$

Доказательство. Выберем гомоморфизм $h$. Выберем элемент $m \in M$ и элемент $n \in N$. Чтобы построить отображение $H$, рассмотрим следующую диаграмму

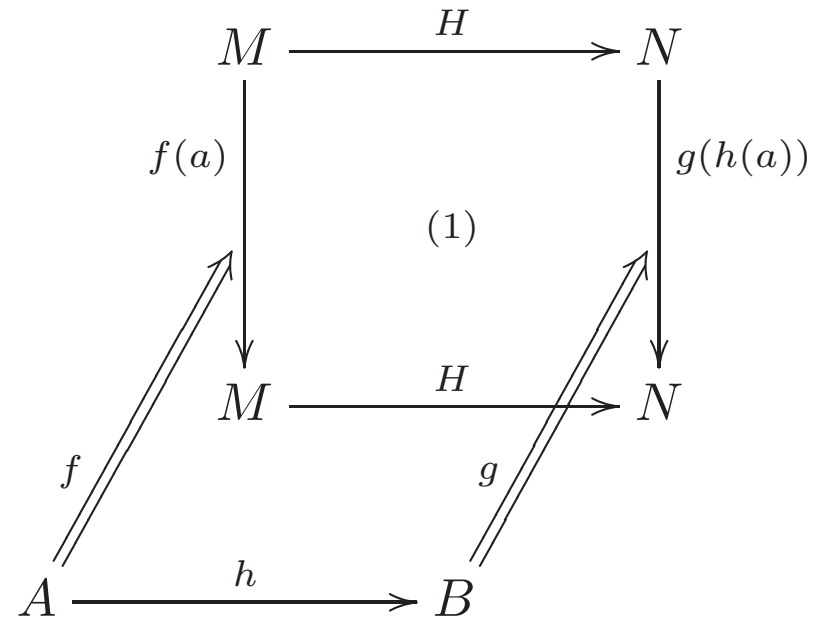

Из коммутативности диаграммы (1) следует

$$H(f(a)m) = h(g(a))H(m)$$

Для произвольного $m' \in M$ однозначно определён $a \in A$ такой, что $m' = f(a)m$. Следовательно, мы построили отображении $H$, которое удовлетворяет равенству (2.2.3). □

Теорема 2.2.10. *Если представление*

$$f : A \xrightarrow{*} M$$

*$\Omega_1$-алгебры $A$ однотранзитивно и представление*

$$g : B \xrightarrow{*} N$$

*$\Omega_1$-алгебры $B$ однотранзитивно, то для заданного гомоморфизма $\Omega_1$-алгебры*

$$h : A \longrightarrow B$$

*гомоморфизм $\Omega_2$-алгебры*

$$H : M \longrightarrow N$$

*такой, что $(h, H)$ является морфизмом представлений из $f$ в $g$, определён однозначно с точностью до выбора образа $n = H(m) \in N$ заданного элемента $m \in M$.*

Доказательство. Из доказательства теоремы 2.2.9 следует, что выбор гомоморфизма $h$ и элементов $m \in M$, $n \in N$ однозначно определяет отображение $H$. □

Теорема 2.2.11. *Если представление*

$$f : A \xrightarrow{*} M$$

*$\Omega_1$-алгебры $A$ однотранзитивно, то для любого эндоморфизма $h$ $\Omega_1$-алгебры $A$ существует морфизм представления $f$*

$$h : A \to A \quad H : M \to M$$

Доказательство. Рассмотрим следующую диаграмму

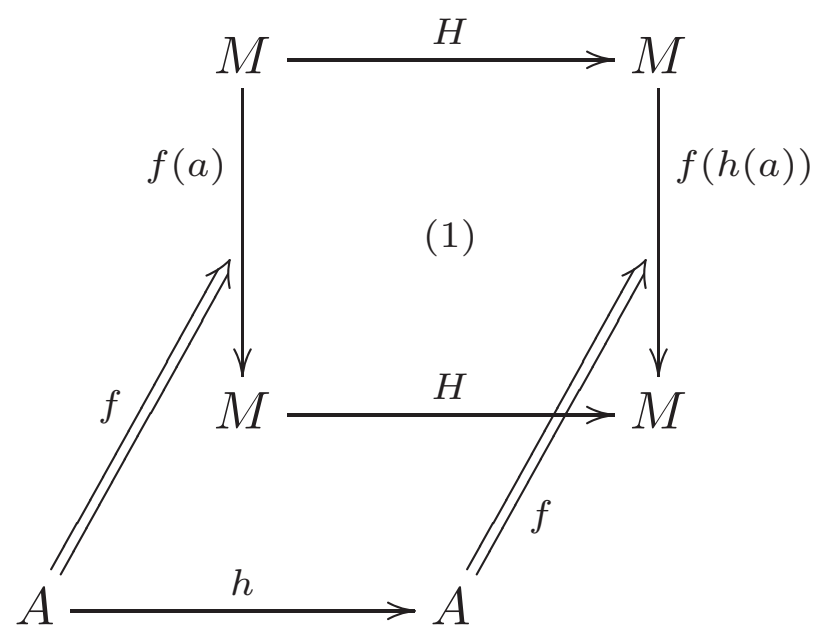

Утверждение теоремы является следствием теоремы 2.2.9. □

Теорема 2.2.12. *Пусть*

$$f : A \xrightarrow{*} M$$

*представление $\Omega_1$-алгебры $A$,*

$$g : B \xrightarrow{*} N$$

*представление $\Omega_1$-алгебры $B$,*

$$h : C \to {}^* L$$

*представление $\Omega_1$-алгебры $C$. Пусть определены морфизмы представлений $\Omega_1$-алгебры*

$$p : A \longrightarrow B \qquad P : M \longrightarrow N$$

$$q : B \longrightarrow C \qquad Q : N \longrightarrow L$$

*Тогда определён морфизм представлений $\Omega_1$-алгебры*

$$r : A \longrightarrow C \qquad R : M \longrightarrow L$$

*где $r = qp$, $R = QP$. Мы будем называть морфизм $(r, R)$ представлений из $f$ в $h$* **произведением морфизмов** $(p, P)$ **и** $(q, Q)$ **представлений универсальной алгебры**.

Доказательство. Мы можем представить утверждение теоремы, пользуясь диаграммой

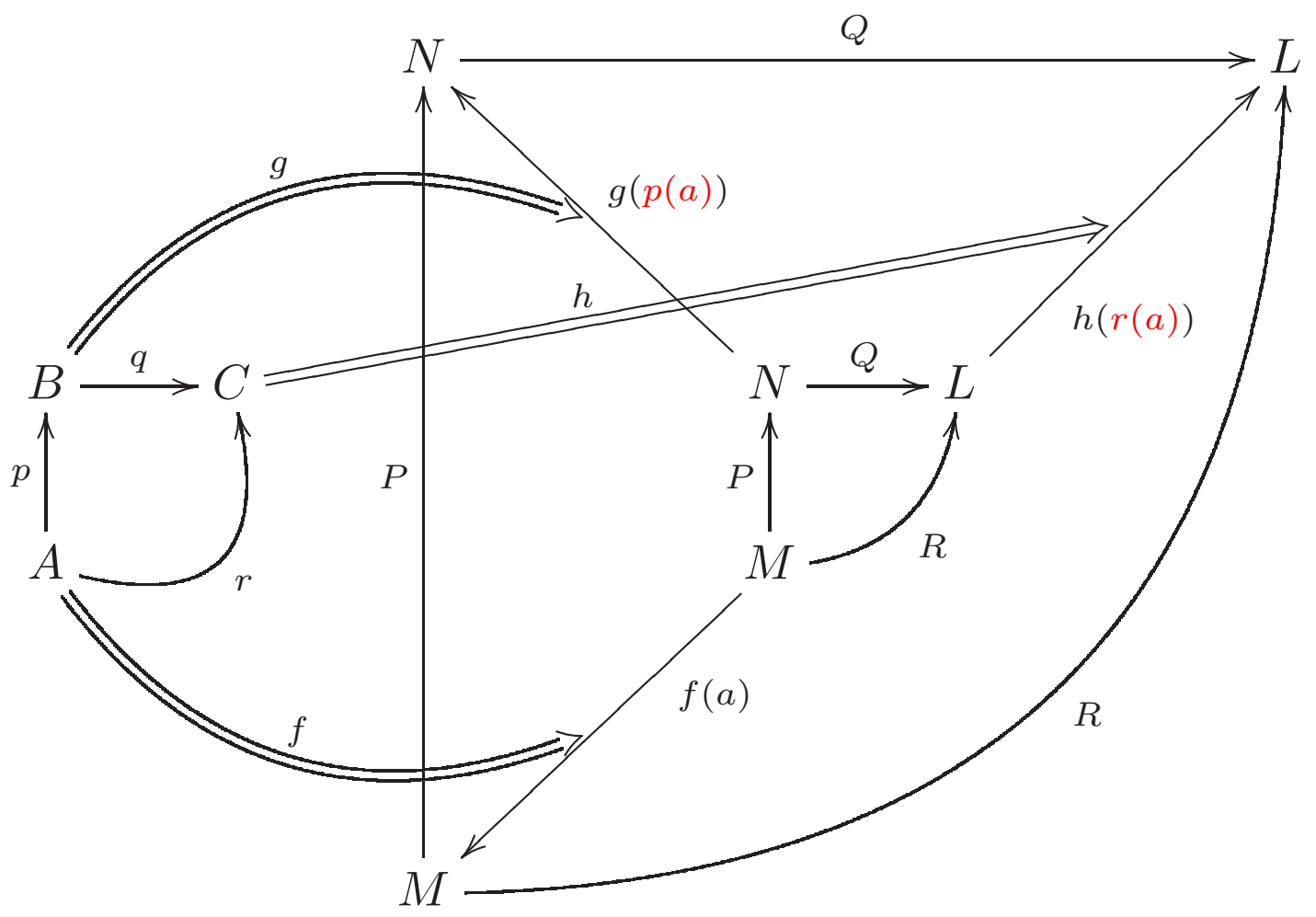

Отображение $r$ является гомоморфизмом $\Omega_1$-алгебры $A$ в $\Omega_1$-алгебру $C$. Нам надо показать, что пара отображений $(r,R)$ удовлетворяет (2.2.3):

$$\begin{aligned}R(f(a)m) &= QP(f(a)m)\\ &= Q(g(p(a))P(m))\\ &= h(qp(a))QP(m))\\ &= h(r(a))R(m)\end{aligned}$$

□

Определение 2.2.13. Допустим $\mathcal{A}$ категория $\Omega_1$-алгебр. Мы определим **категорию** $\mathcal{A}*$ **левосторонних представлений** универсальной алгебры из категории $\mathcal{A}$. Объектами этой категории являются левосторонние представлениями $\Omega_1$-алгебры. Морфизмами этой категории являются морфизмы левосторонних представлений $\Omega_1$-алгебры. □

Определение 2.2.14. Пусть на множестве $M$ определена эквивалентность $S$. Преобразование $f$ называется **согласованным с эквивалентностью** $S$, если из условия $m_1 \equiv m_2(\mathrm{mod} S)$ следует $f(m_1) \equiv f(m_2)(\mathrm{mod}\ S)$. □

Теорема 2.2.15. *Пусть на множестве $M$ определена эквивалентность $S$. Пусть на множестве ${}^*M$ определена $\Omega_1$-алгебра. Если любое преобразование $f \in {}^*M$ согласованно с эквивалентностью $S$, то мы можем определить структуру $\Omega_1$-алгебры на множестве ${}^*(M/S)$.*

Доказательство. Пусть $h = \mathrm{nat}\ S$. Если $m_1 \equiv m_2(\mathrm{mod} S)$ , то $h(m_1) = h(m_2)$. Поскольку $f \in {}^*M$ согласованно с эквивалентностью $S$, то $h(f(m_1)) = h(f(m_2))$. Это позволяет определить преобразование $F$ согласно правилу

$$F([m]) = h(f(m)) \tag{2.2.13}$$

Пусть $\omega$ - n-арная операция $\Omega_1$-алгебры. Пусть $f_1$, ..., $f_n \in {}^\star M$ и

$$F_1([m]) = h(f_1(m)) \quad ... \quad F_n([m]) = h(f_n(m))$$

Согласно условию теоремы, преобразование

$$f = \omega(f_1, ..., f_n) \in {}^*M$$

согласованно с эквивалентностью $S$. Следовательно, из условия $m_1 \equiv m_2(\mathrm{mod} S)$ и определения 2.2.14 следует

$$f(m_1) \equiv f(m_2)(\mathrm{mod} S) \tag{2.2.14}$$
$$\omega(f_1, ..., f_n)(m_1) \equiv \omega(f_1, ..., f_n)(m_2)(\mathrm{mod} S)$$

Следовательно, мы можем определить операцию $\omega$ на множестве ${}^\star(M/S)$ по правилу

$$\omega(F_1, ..., F_n)[m] = h(\omega(f_1, ..., f_n)(m)) \tag{2.2.15}$$

Из определения (2.2.13) и равенства (2.2.14) следует, что мы корректно определили операцию $\omega$ на множестве ${}^\star(M/S)$. □

Теорема 2.2.16. *Пусть*

$$f : A \xrightarrow{*} M$$

*представление $\Omega_1$-алгебры $A$,*

$$g : B \xrightarrow{*} N$$

*представление $\Omega_1$-алгебры $B$. Пусть*

$$r : A \longrightarrow B \qquad R : M \longrightarrow N$$

*морфизм представлений из $f$ в $g$. Положим*

$$s = rr^{-1} \qquad S = RR^{-1}$$

*Тогда для отображений $r$, $R$ существуют разложения, которые можно описать диаграммой*

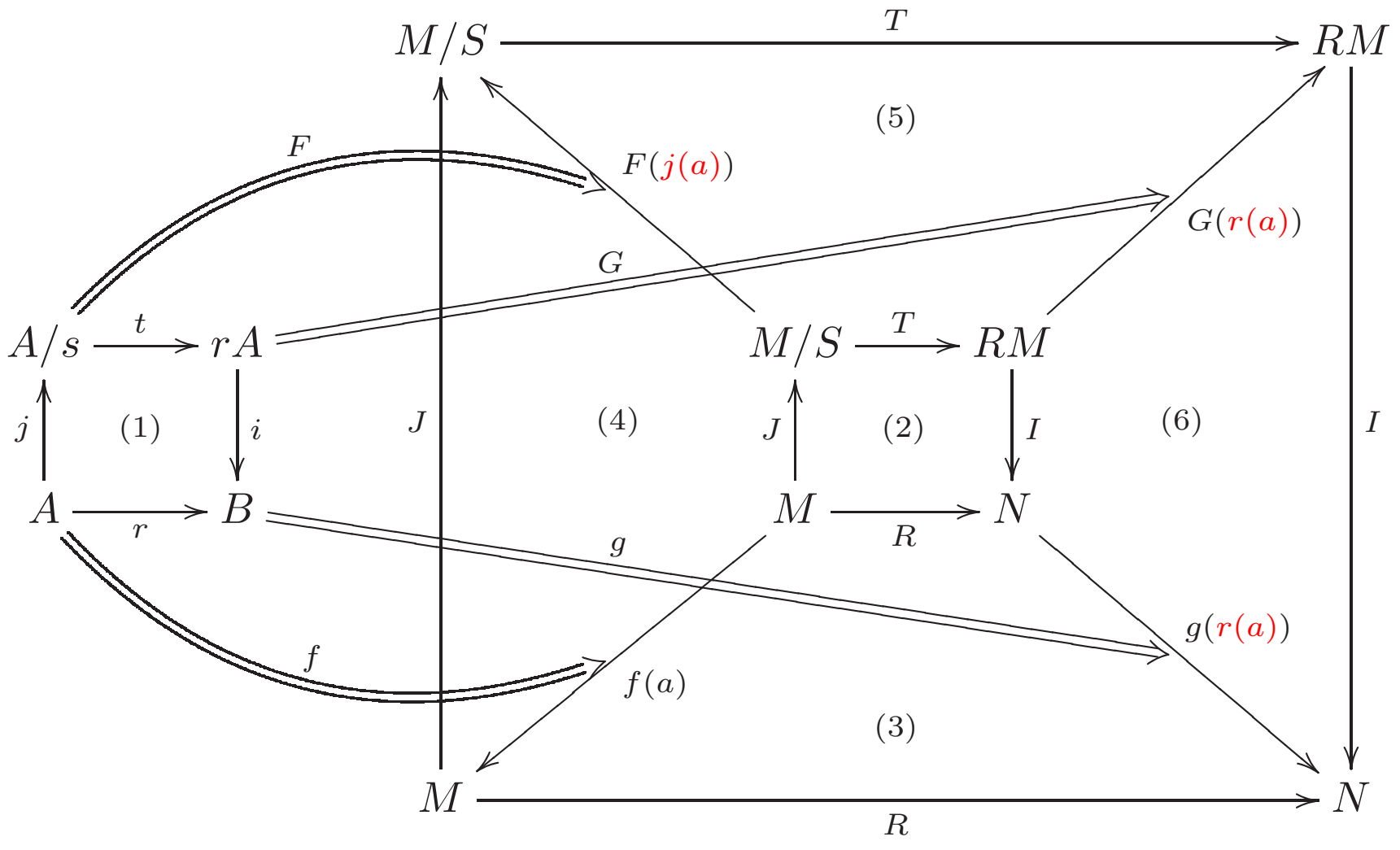

(1) *$s=\ker r$ является конгруэнцией на $A$. Существует разложение гомоморфизма $r$*

$$r=itj \tag{2.2.16}$$

*$j=\mathrm{nat}\ s$ - естественный гомоморфизм*

$$j(a)=j(a)$$

*$t$ - изоморфизм*

$$r(a)=t(j(a)) \tag{2.2.17}$$

*$i$ - вложение*

$$r(a)=i(r(a)) \tag{2.2.18}$$

(2) *$S=\ker R$ является эквивалентностью на $M$. Существует разложение отображения $R$*

$$R=ITJ \tag{2.2.19}$$

*$J=\mathrm{nat}\ S$ - сюръекция*

$$J(m)=J(m)$$

*$T$ - биекция*

$$R(m)=T(J(m)) \tag{2.2.20}$$

*$I$ - вложение*

$$R(m)=I(R(m)) \tag{2.2.21}$$

(3) *$F$ - левостороннее представление $\Omega_1$-алгебры $A/s$ в $M/S$*
(4) *$G$ - левостороннее представление $\Omega_1$-алгебры $rA$ в $RM$*
(5) *$(j,J)$ - морфизм представлений $f$ и $F$*
(6) *$(t,T)$ - морфизм представлений $F$ и $G$*
(7) *$(t^{-1},T^{-1})$ - морфизм представлений $G$ и $F$*
(8) *$(i,I)$ - морфизм представлений $G$ и $g$*
(9) *Существует разложение морфизма представлений*

$$(r,R)=(i,I)(t,T)(j,J) \tag{2.2.22}$$

Доказательство. Существование диаграмм (1) и (2) следует из теоремы II.3.7 ([13], с. 74).

Мы начнём с диаграммы (4).

Пусть $m_1\equiv m_2(\mathrm{mod}\ S)$. Следовательно,

$$R(m_1)=R(m_2) \tag{2.2.23}$$

Если $a_1\equiv a_2(\mathrm{mod}\ s)$, то

$$r(a_1)=r(a_2) \tag{2.2.24}$$

Следовательно, $j(a_1)=j(a_2)$. Так как $(r,R)$ - морфизм представлений, то

$$R(f(a_1)(m_1))=g(r(a_1))(R(m_1)) \tag{2.2.25}$$

$$R(f(a_2)(m_2))=g(r(a_2))(R(m_2)) \tag{2.2.26}$$

Из (2.2.23), (2.2.24), (2.2.25), (2.2.26) следует

$$R(f(a_1)(m_1))=R(f(a_2)(m_2)) \tag{2.2.27}$$

Из (2.2.27) следует

(2.2.28) $$f(a_1)(m_1) \equiv f(a_2)(m_2)(\mathrm{mod}\ S)$$

и, следовательно,

(2.2.29) $$J(f(a_1)(m_1)) = J(f(a_2)(m_2))$$

Из (2.2.29) следует, что отображение

(2.2.30) $$F(j(a))(J(m)) = J(f(a)(m)))$$

определено корректно и является преобразованием множества $M/S$.

Из равенства (2.2.28) (в случае $a_1 = a_2$) следует, что для любого $a$ преобразование согласованно с эквивалентностью $S$. Из теоремы 2.2.15 следует, что на множестве $^*(M/S)$ определена структура $\Omega_1$-алгебры. Рассмотрим $n$-арную операцию $\omega$ и $n$ преобразований

$$F(j(a_i))(J(m)) = J(f(a_i)(m))) \quad i = 1, ..., n$$

пространства $M/S$. Мы положим

$$\omega(F(j(a_1)), ..., F(j(a_n)))(J(m)) = J(\omega(f(a_1), ..., f(a_n)))(m))$$

Следовательно, отображение $F$ является представлением $\Omega_1$-алгебры $A/s$.

Из (2.2.30) следует, что $(j, J)$ является морфизмом представлений $f$ и $F$ (утверждение (5) теоремы).

Рассмотрим диаграмму (5).

Так как $T$ - биекция, то мы можем отождествить элементы множества $M/S$ и множества $MR$, причём это отождествление имеет вид

(2.2.31) $$T(J(m)) = R(m)$$

Мы можем записать преобразование $F(j(a))$ множества $M/S$ в виде

(2.2.32) $$F(j(a)) : J(m) \to F(j(a))(J(m))$$

Так как $T$ - биекция, то мы можем определить преобразование

(2.2.33) $$T(J(m)) \to T(F(j(a))(J(m)))$$

множества $RM$. Преобразование (2.2.33) зависит от $j(a) \in A/s$. Так как $t$ - биекция, то мы можем отождествить элементы множества $A/s$ и множества $rA$, причём это отождествление имеет вид

(2.2.34) $$t(j(a)) = r(a)$$

Следовательно, мы определили отображение

$$G : rA \to {}^\star RM$$

согласно равенству

(2.2.35) $$G(t(j(a)))(T(J(m))) = T(F(j(a))(J(m)))$$

Рассмотрим $n$-арную операцию $\omega$ и $n$ преобразований

$$G(r(a_i))(R(m)) = T(F(j(a_i))(J(m))) \quad i = 1, ..., n$$

пространства $RM$. Мы положим

(2.2.36) $$\omega(G(r(a_1)), ..., G(r(a_n)))(R(m)) = T(\omega(F(j(a_1), ..., F(j(a_n)))(J(m)))$$

Согласно (2.2.35) операция $\omega$ корректно определена на множестве ${}^\star RM$. Следовательно, отображение $G$ является представлением $\Omega_1$-алгебры.

Из (2.2.35) следует, что $(t,T)$ является морфизмом представлений $F$ и $G$ (утверждение (6) теоремы).

Так как $T$ - биекция, то из равенства (2.2.31) следует

$$J(m)=T^{-1}(R(m))\tag{2.2.37}$$

Мы можем записать преобразование $G(r(a))$ множества $RM$ в виде

$$G(r(a)):R(m)\rightarrow G(r(a))(R(m))\tag{2.2.38}$$

Так как $T$ - биекция, то мы можем определить преобразование

$$T^{-1}(R(m))\rightarrow T^{-1}(G(r(a))(R(m)))\tag{2.2.39}$$

множества $M/S$. Преобразование (2.2.39) зависит от $r(a)\in rA$. Так как $t$ - биекция, то из равенства (2.2.34) следует

$$j(a)=t^{-1}(r(a))\tag{2.2.40}$$

Так как по построению диаграмма (5) коммутативна, то преобразование (2.2.39) совпадает с преобразованием (2.2.32). Равенство (2.2.36) можно записать в виде

$$T^{-1}(\omega(G(r(a_1)),...,G(r(a_n)))(R(m)))=\omega(F(j(a_1),...,F(j(a_n)))(J(m))\tag{2.2.41}$$

Следовательно, $(t^{-1},T^{-1})$ является морфизмом представлений $G$ и $F$ (утверждение (7) теоремы).

Диаграмма (6) является самым простым случаем в нашем доказательстве. Поскольку отображение $I$ является вложением и диаграмма (2) коммутативна, мы можем отождествить $n\in N$ и $R(m)$, если $n\in \mathrm{Im}R$. Аналогично, мы можем отождествить соответствующие преобразования.

$$\begin{gathered}g'(i(r(a)))(I(R(m)))=I(G(r(a))(R(m)))\\ \omega(g'(r(a_1)),...,g'(r(a_n)))(R(m))=I(\omega(G(r(a_1),...,G(r(a_n)))(R(m)))\end{gathered}\tag{2.2.42}$$

Следовательно, $(i,I)$ является морфизмом представлений $G$ и $g$ (утверждение (8) теоремы).

Для доказательства утверждения (9) теоремы осталось показать, что определённое в процессе доказательства представление $g'$ совпадает с представлением $g$, а операции над преобразованиями совпадают с соответствующими операциями на ${}^*N$.

$$\begin{aligned}g'(i(r(a)))(I(R(m)))&=I(G(r(a))(R(m))) && \text{by (2.2.42)}\\ &=I(G(t(j(a)))(T(J(m)))) && \text{by (2.2.17), (2.2.20)}\\ &=IT(F(j(a))(J(m))) && \text{by (2.2.35)}\\ &=ITJ(f(a)(m)) && \text{by (2.2.30)}\\ &=R(f(a)(m)) && \text{by (2.2.19)}\\ &=g(r(a))(R(m)) && \text{by (2.2.3)}\end{aligned}$$

$$\begin{aligned}\omega(G(r(a_1)),...,G(r(a_n)))(R(m)) &= T(\omega(F(j(a_1),...,F(j(a_n)))(J(m)))\\ &= T(F(\omega(j(a_1),...,j(a_n)))(J(m)))\\ &= T(F(j(\omega(a_1,...,a_n)))(J(m)))\\ &= T(J(f(\omega(a_1,...,a_n))(m)))\end{aligned}$$

□

Определение 2.2.17. Пусть

$$f : A \xrightarrow{*} M$$

представление $\Omega_1$-алгебры $A$,

$$g : B \xrightarrow{*} N$$

представление $\Omega_1$-алгебры $B$. Пусть

$$r : A \longrightarrow B \qquad R : M \longrightarrow N$$

морфизм представлений из $r$ в $R$ такой, что $f$ - изоморфизм $\Omega_1$-алгебры и $g$ - изоморфизм $\Omega_2$-алгебры. Тогда отображение $(r, R)$ называется **изоморфизмом представлений**. □

Теорема 2.2.18. *В разложении* (2.2.22) *отображение* $(t, T)$ *является изоморфизмом представлений* $F$ *и* $G$.

Доказательство. Следствие определения 2.2.17 и утверждений (6) и (7) теоремы 2.2.16. □

## 2.3. Приведенный морфизм представлений

Из теоремы 2.2.16 следует, что мы можем свести задачу изучения морфизма представлений $\Omega_1$-алгебры к случаю, описываемому диаграммой

(2.3.1)

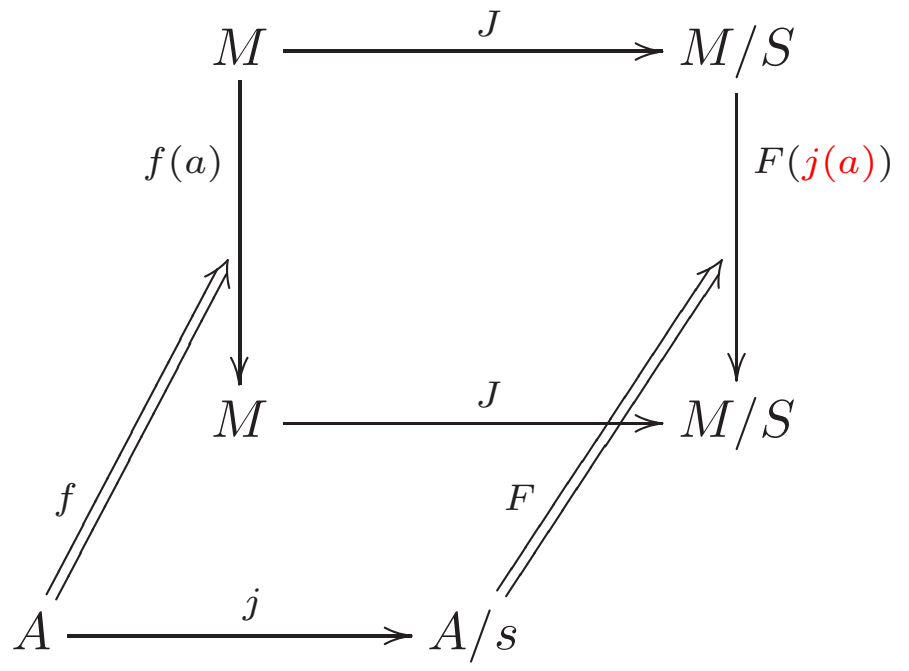

Теорема 2.3.1. *Диаграмма* (2.3.1) *может быть дополнена представлением* $F_1$ $\Omega_1$*-алгебры* $A$ *в множестве* $M/S$ *так, что диаграмма*

(2.3.2)

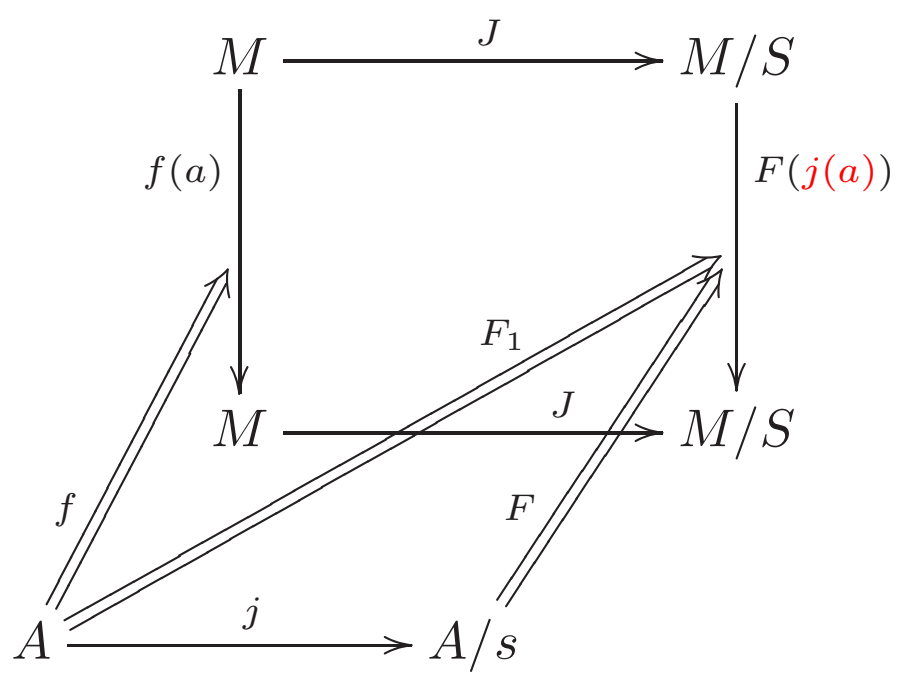

*коммутативна. При этом множество преобразований представления* $F$ *и множество преобразований представления* $F_1$ *совпадают.*

Доказательство. Для доказательства теоремы достаточно положить

$$F_1(a)=F(j(a))$$

Так как отображение $j$ - сюрьекция, то $\mathrm{Im}F_1=\mathrm{Im}F$. Так как $j$ и $F$ - гомоморфизмы $\Omega_1$-алгебры, то $F_1$ - также гомоморфизм $\Omega_1$-алгебры. □

Замечание 2.3.2. Теорема 2.3.1 завершает цикл теорем, посвящённых структуре морфизма представлений $\Omega_1$-алгебры. Из этих теорем следует, что мы можем упростить задачу изучения морфизма представлений $\Omega_1$-алгебры и ограничиться морфизмом представлений вида

$$id:A\longrightarrow A \qquad R:M\longrightarrow N$$

В этом случае мы можем отождествить морфизм $(\mathrm{id},R)$ представлений $\Omega_1$-алгебры и соответствующий гомоморфизм $R$ $\Omega_2$-алгебры и будем называть гомоморфизм $R$ **приведенным морфизмом представлений**. Мы будем пользоваться диаграммой

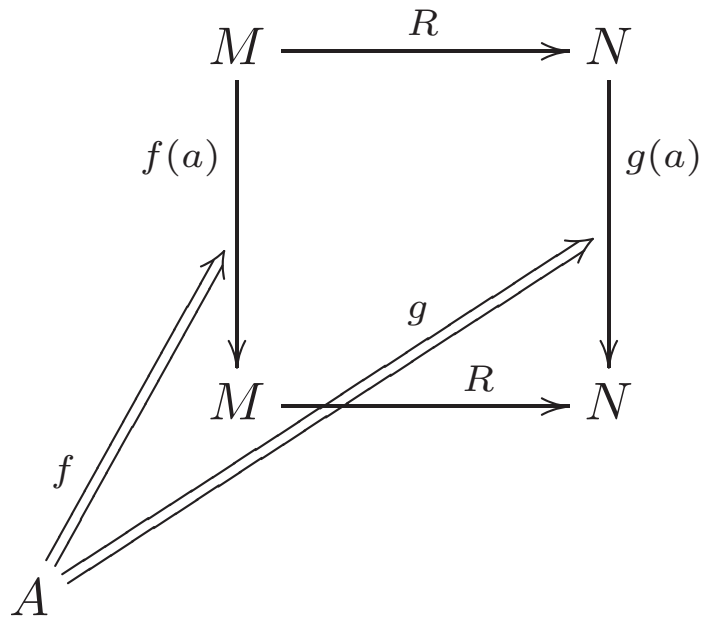

для представления приведенного морфизма $R$ представлений $\Omega_1$-алгебры. Из диаграммы следует

(2.3.3) $$R\circ f(a)=g(a)\circ R$$

□

Теорема 2.3.3. *Пусть отображение*

$$H : M \to N$$

*является приведенным морфизмом из представления*

$$f : A \xrightarrow{*} M$$

*$\Omega_1$-алгебры $A$ в представление*

$$g : A \xrightarrow{*} N$$

*$\Omega_1$-алгебры $A$. Если представление $f$ эффективно, то отображение*

$$^*H : {}^*M \to {}^*N$$

*определённое равенством*

$$^*H(f(a)) = g(a) \tag{2.3.4}$$

*является гомоморфизмом $\Omega_1$-алгебры.*

Доказательство. Теорема является следствием теоремы 2.2.8, если мы положим $h = id$. □

Теорема 2.3.4. *Пусть представления*

$$f : A \xrightarrow{*} M$$

$$g : A \xrightarrow{*} N$$

*$\Omega_1$-алгебры $A$ однотранзитивны. Существует приведенный морфизм представлений из $f$ в $g$*

$$H : M \to N$$

Доказательство. Выберем гомоморфизм $h$. Выберем элемент $m \in M$ и элемент $n \in N$. Чтобы построить отображение $H$, рассмотрим следующую диаграмму

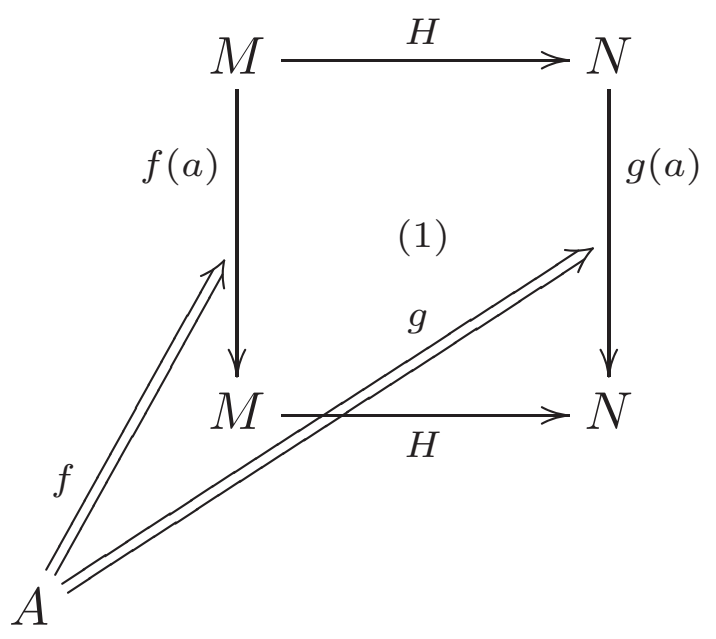

Из коммутативности диаграммы (1) следует

$$H(f(a)m) = g(a)H(m)$$

Для произвольного $m' \in M$ однозначно определён $a \in A$ такой, что $m' = f(a)m$. Следовательно, мы построили отображении $H$, которое удовлетворяет равенству (2.3.3). □

Теорема 2.3.5. *Пусть представления*

$$f : A \xrightarrow{*} M$$

$$g : A \xrightarrow{*} N$$

*$\Omega_1$-алгебры $A$ однотранзитивны. Приведенный морфизм представлений из $f$ в $g$*

$$H : M \longrightarrow N$$

*определён однозначно с точностью до выбора образа $n = H(m) \in N$ заданного элемента $m \in M$.*

Доказательство. Из доказательства теоремы 2.3.4 следует, что выбор элементов $m \in M$, $n \in N$ однозначно определяет отображение $H$. □

Мы дадим следующее определение по аналогии с определением 2.2.13.

Определение 2.3.6. Мы определим **категорию** $A*$ **левосторонних представлений $\Omega_1$-алгебры** $A$. Объектами этой категории являются левосторонними представлениями $\Omega_1$-алгебры $A$. Морфизмами этой категории являются морфизмы $(id, R)$ левосторонних представлений $\Omega_1$-алгебры $A$. □

## 2.4. Автоморфизм представления универсальной алгебры

Определение 2.4.1. Пусть

$$f : A \xrightarrow{*} M$$

представление $\Omega_1$-алгебры $A$ в $\Omega_2$-алгебре $M$. Приведенный морфизм представлений $\Omega_1$-алгебры

$$R : M \to M$$

такой, что $R$ - эндоморфизм $\Omega_2$-алгебры называется **эндоморфизмом представления** $f$. □

Теорема 2.4.2. *Если представление*

$$f : A \xrightarrow{*} M$$

*$\Omega_1$-алгебры $A$ однотранзитивно, то для любых $p$, $q \in M$ существует единственный эндоморфизм*

$$H : M \to M$$

*представления $f$ такой, что $H(p) = q$.*

Доказательство. Рассмотрим следующую диаграмму

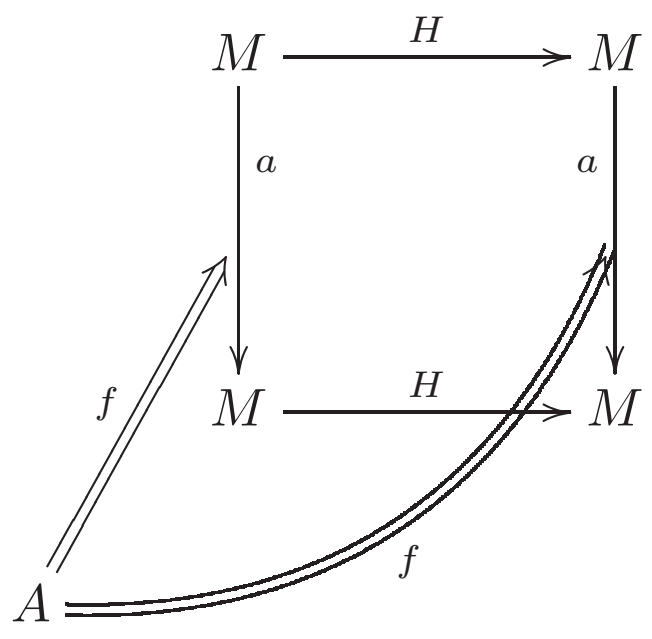

Существование эндоморфизма является следствием теоремы 2.2.9. Единственность эндоморфизма для заданных $p$, $q \in M$ является следствием теоремы 2.2.10, когда $h = \mathrm{id}$. □

Теорема 2.4.3. *Эндоморфизмы представления $f$ порождают полугруппу.*

Доказательство. Из теоремы 2.2.12 следует, что произведение эндоморфизмов $(p, P)$, $(r, R)$ представления $f$ является эндоморфизмом $(pr, PR)$ представления $f$. □

Определение 2.4.4. Пусть

$$f : A \xrightarrow{*} M$$

представление $\Omega_1$-алгебры $A$ в $\Omega_2$-алгебре $M$. Морфизм представлений $\Omega_1$-алгебры

$$R : M \to M$$

такой, что $R$ - автоморфизм $\Omega_2$-алгебры называется **автоморфизмом представления** $f$. □

Теорема 2.4.5. *Пусть*

$$f : A \xrightarrow{*} M$$

*представление $\Omega_1$-алгебры $A$ в $\Omega_2$-алгебре $M$. Множество автоморфизмов представления $f$ порождает* **группу** $GA(f)$ *.*

Доказательство. Пусть $R$, $P$ - автоморфизмы представления $f$. Согласно определению 2.4.4, отображения $R$, $P$ являются автоморфизмами $\Omega_2$-алгебры $M$. Согласно теореме II.3.2, ([13], с. 60), отображение $R \circ P$ является автоморфизмом $\Omega_2$-алгебры $M$. Из теоремы 2.2.12 и определения 2.4.4 следует, что произведение автоморфизмов $R \circ P$ представления $f$ является автоморфизмом представления $f$.

Пусть $R$, $P$, $Q$ - автоморфизмы представления $f$. Из цепочки равенств

$$\begin{aligned}((R \circ P) \circ Q)(a) &= (R \circ P)(Q(a)) = R(P(Q(a))) \\ &= R((P \circ Q)(a)) = (R \circ (P \circ Q))(a)\end{aligned}$$

следует ассоциативность произведения для отображений $R$, $P$, $Q$.[2.3]

Пусть $R$ - автоморфизм представления $f$. Согласно определению 2.4.4 отображение $R$ является автоморфизмом $\Omega_2$-алгебры $M$. Следовательно, отображение $R^{-1}$ является автоморфизмом $\Omega_2$-алгебры $M$. Для автоморфизма $R$ представления справедливо равенство (2.2.4). Положим $m' = R(m)$. Так как $R$ - автоморфизм $\Omega_2$-алгебры, то $m = R^{-1}(m')$ и равенство (2.2.4) можно записать в виде

$$R(f(a')(R^{-1}(m'))) = f(a')(m') \tag{2.4.1}$$

Так как отображение $R$ является автоморфизмом $\Omega_2$-алгебры $M$, то из равенства (2.4.1) следует

$$f(a')(R^{-1}(m')) = R^{-1}(f(a')(m')) \tag{2.4.2}$$

Равенство (2.4.2) соответствует равенству (2.2.4) для отображения $R^{-1}$. Следовательно, отображение $R^{-1}$ является автоморфизмом представления $f$. □

[2.3] При доказательстве ассоциативности произведения я следую примеру полугруппы из [4], с. 20, 21.

## 2.5. Примеры представления универсальной алгебры

**2.5.1. Векторное пространство.** Согласно определениям [7]-5.1.4, [11]-4.1.4, эффективное левостороннее представление тела $D$ в абелевой группе называется **левым $D$-векторным пространством**. Аналогично эффективное представление ассоциативного кольца $R$ в абелевой группе называется **модулем над кольцом $R$ или $R$-модулем**.

**2.5.2. Представление группы на множестве.** Пусть $G$ - абелевая группа, и $M$ - множество. Рассмотрим эфективное представление группы $G$ на множестве $M$. Для заданных $a \in G$, $A \in M$ положим $A \to A + a$. Мы будем также пользоваться записью $a = \vec{AB}$, если

$$B = A + a \tag{2.5.1}$$

Тогда действие группы можно представить в виде

$$B = A + \vec{AB} \tag{2.5.2}$$

Поскольку представление эффективно, то из равенств (2.5.1), (2.5.2) и равенства

$$D = C + a$$

следует, что

$$\overrightarrow{AB} = \overrightarrow{CD}$$

## 2.6. Множество образующих представления

Определение 2.6.1. Пусть

$$f : A \xrightarrow{*} M$$

представление $\Omega_1$-алгебры $A$ в $\Omega_2$-алгебре $M$. Множество $N \subset M$ называется **стабильным множеством представления** $f$, если $f(a)(m) \in N$ для любых $a \in A$, $m \in N$. □

Мы также будем говорить, что множество $M$ стабильно относительно представления $f$.

Теорема 2.6.2. *Пусть*

$$f : A \xrightarrow{*} M$$

*представление $\Omega_1$-алгебры $A$ в $\Omega_2$-алгебре $M$. Пусть множество $N \subset M$ является подалгеброй $\Omega_2$-алгебры $M$ и стабильным множеством представления $f$. Тогда существует представление*

$$f_N : A \to {}^*N$$

*такое, что $f_N(a) = f(a)|_N$. Представление $f_N$ называется* **подпредставлением представления $f$**.

Доказательство. Пусть $\omega_1$ - $n$-арная операция $\Omega_1$-алгебры $A$. Тогда для любых $a_1$, ..., $a_n \in A$ и любого $b \in N$

$$\begin{aligned}(f_N(a_1)...f_N(a_n)\omega_1)(b) &= (f(a_1)...f(a_n)\omega_1)(b)\\ &= f(a_1...a_n\omega_1)(b)\\ &= f_N(a_1...a_n\omega_1)(b)\end{aligned}$$

Пусть $\omega_2$ - $n$-арная операция $\Omega_2$-алгебры $M$. Тогда для любых $b_1$, ..., $b_n \in N$ и любого $a \in A$

$$\begin{aligned}f_N(a)(b_1)...f_N(a)(b_n)\omega_2 &= f(a)(b_1)...f(a)(b_n)\omega_2\\ &= f(a)(b_1...b_n\omega_2)\\ &= f_N(a)(b_1...b_n\omega_2)\end{aligned}$$

Утверждение теоремы доказано. □

Из теоремы 2.6.2 следует, что если $f_N$ - подпредставление представления $f$, то отображение $(id : A \rightarrow A, id_N : N \rightarrow M)$ является морфизмом представлений.

Теорема 2.6.3. *Множество*[2.4] $\mathcal{B}_f$ *всех подпредставлений представления* $f$ *порождает систему замыканий на* $\Omega_2$*-алгебре* $M$ *и, следовательно, является полной структурой.*

Доказательство. Пусть $(K_\lambda)_{\lambda\in\Lambda}$ - семейство подалгебр $\Omega_2$-алгебры $M$, стабильных относительно представления $f$. Операцию пересечения на множестве $\mathcal{B}_f$ мы определим согласно правилу

$$\bigcap f_{K_\lambda} = f_{\cap K_\lambda}$$

Операция пересечения подпредставлений определена корректно. $\cap K_\lambda$ - подалгебра $\Omega_2$-алгебры $M$. Пусть $m \in \cap K_\lambda$. Для любого $\lambda \in \Lambda$ и для любого $a \in A$, $f(a)(m) \in K_\lambda$. Следовательно, $f(a)(m) \in \cap K_\lambda$. Следовательно, $\cap K_\lambda$ - стабильное множество представления $f$. □

Обозначим соответствующий оператор замыкания через $J(f)$. Таким образом, $J(f, X)$ является пересечением всех подалгебр $\Omega_2$-алгебры $M$, содержащих $X$ и стабильных относительно представления $f$.

Теорема 2.6.4. *Пусть*[2.5]

$$f : A \twoheadrightarrow M$$

*представление* $\Omega_1$*-алгебры* $A$ *в* $\Omega_2$*-алгебре* $M$. *Пусть* $X \subset M$. *Определим подмножество* $X_k \subset M$ *индукцией по* $k$.

2.6.4.1: $X_0 = X$
2.6.4.2: $x \in X_k \Rightarrow x \in X_{k+1}$
2.6.4.3: $x_1 \in X_k$, ..., $x_n \in X_k$, $\omega \in \Omega_2(n) \Rightarrow x_1...x_n\omega \in X_{k+1}$
2.6.4.4: $x \in X_k$, $a \in A \Rightarrow f(a)(x) \in X_{k+1}$

*Тогда*

$$\bigcup_{k=0}^{\infty} X_k = J(f, X) \tag{2.6.1}$$

Доказательство. Если положим $U = \cup X_k$, то по определению $X_k$ имеем $X_0 \subset J(f, X)$, и если $X_k \subset J(f, X)$, то $X_{k+1} \subset J(f, X)$. По индукции следует, что $X_k \subset J(f, X)$ для всех $k$. Следовательно,

$$U \subset J(f, X) \tag{2.6.2}$$

[2.4] Это определение аналогично определению структуры подалгебр ([13], стр. 93, 94)
[2.5] Утверждение теоремы аналогично утверждению теоремы 5.1, [13], стр. 94.

Если $a \in U^n$, $a = (a_1, ..., a_n)$, где $a_i \in X_{k_i}$, и если $k = \max\{k_1, ..., k_n\}$, то $a_1...a_n\omega \in X_{k+1} \subset U$. Следовательно, $U$ является подалгеброй $\Omega_2$-алгебры $M$.

Если $m \in U$, то $m \in X_k$ для некоторого $k$. Следовательно, $f(a)(m) \in X_{k+1} \subset U$ для любого $a \in A$. Следовательно, $U$ - стабильное множество представления $f$.

Так как $U$ - подалгеброй $\Omega_2$-алгебры $M$ и стабильное множество представления $f$, то определено подпредставление $f_U$. Следовательно,

$$J(f, X) \subset U \tag{2.6.3}$$

Из (2.6.2), (2.6.3), следует $J(f, X) = U$. □

Определение 2.6.5. $J(f, X)$ называется **подпредставлением, порождённым множеством** $X$, а $X$ - **множеством образующих подпредставления** $J(f, X)$. В частности, **множеством образующих представления** $f$ будет такое подмножество $X \subset M$, что $J(f, X) = M$. □

Нетрудно видеть, что определение множества образующих представления не зависит от того, эффективно представление или нет. Поэтому в дальнейшем мы будем предполагать, что представление эффективно и будем опираться на соглашение для эффективного левостороннего представления в замечании 2.1.7. Мы также будем пользоваться записью

$$R \circ m = R(m)$$

для образа $m \in M$ при эндоморфизме эффективного представления. Согласно определению произведения отображений, для любых эндоморфизмов $R$, $S$ верно равенство

$$(R \circ S) \circ m = R \circ (S \circ m) \tag{2.6.4}$$

Равенство (2.6.4) является законом ассоциативности для $\circ$ и позволяет записать выражение

$$R \circ S \circ m$$

без использования скобок.

Из теоремы 2.6.4 следует следующее определение.

Определение 2.6.6. Пусть $X \subset M$. Для любого $x \in J(f, X)$ существует $\Omega_2$-слово, определённое согласно следующему правилам.

(1) Если $m \in X$, то $m$ - $\Omega_2$-слово.
(2) Если $m_1$, ..., $m_n$ - $\Omega_2$-слова и $\omega \in \Omega_2(n)$, то $m_1...m_n\omega$ - $\Omega_2$-слово.
(3) Если $m$ - $\Omega_2$-слово и $a \in A$, то $am$ - $\Omega_2$-слово.

**$\Omega_2$-слово** $w(f, X, m)$ представляет данный элемент $m \in J(f, X)$. Мы будем отождествлять элемент $m \in J(f, X)$ и соответствующее ему $\Omega_2$-слово, выражая это равенством

$$m = w(f, X, m)$$

Аналогично, для произвольного множества $B \subset J(f, X)$ рассмотрим множество $\Omega_2$-слов[2.6]

$$w(f, X, B) = \{w(f, X, m) : m \in B\}$$

[2.6] Выражение $w(f, X, m)$ является частным случаем выражения $w(f, X, B)$, а именно

$$w(f, X, \{m\}) = \{w(f, X, m)\}$$

Мы будем также пользоваться записью

$$w(f,X,B)=(w(f,X,m),m\in B)$$

Обозначим $w(f,X)$ **множество $\Omega_2$-слов представления** $J(f,X)$. □

Теорема 2.6.7. *Эндоморфизм $R$ представления*

$$f:A\longrightarrow\!\!\!\!\!*\;\; M$$

*порождает отображение $\Omega_2$-слов*

$$w[f,X,R]:w(f,X)\rightarrow w(f,X')\quad X\subset M\quad X'=R\circ X$$

*такое, что*

2.6.7.1: *Если* $m\in X$, $m'=R\circ m$, *то*

$$w[f,X,R](m)=m'$$

2.6.7.2: *Если*

$$m_1,...,m_n\in w(f,X)$$

$$m'_1=w[f,X,R](m_1)\quad ...\quad m'_n=w[f,X,R](m_n)$$

*то для операции $\omega\in\Omega_2(n)$ справедливо*

$$w[f,X,R](m_1...m_n\omega)=m'_1...m'_n\omega$$

2.6.7.3: *Если*

$$m\in w(f,X)\quad m'=w[f,X,R](m)\quad a\in A$$

*то*

$$w[f,X,R](am)=am'$$

Доказательство. Утверждения 2.6.7.1, 2.6.7.2 справедливы в силу определения эндоморфизма $R$. Утверждение
RefItemmap of words of representation, am следует из равенства (2.2.4), так как $r=\mathrm{id}$. □

Замечание 2.6.8. Пусть $R$ - эндоморфизм представления $f$. Пусть

$$m\in J(f,X)\quad m'=R\circ m\quad X'=R\circ X$$

Теорема 2.6.7 утверждает, что $m'\in J_f(X')$ . Теорема 2.6.7 также утверждает, что $\Omega_2$-слово, представляющее $m$, относительно $X$ и $\Omega_2$-слово, представляющее $m'$, относительно $X'$ формируются согласно одному и тому же алгоритму. Это позволяет рассматривать множество $\Omega_2$-слов $w(f,X',m')$ как отображение

$$W(f,X,m):X'\rightarrow w(f,X',m')$$

$$W(f,X,m)(X')=W(f,X,m)\circ X'$$

такое, что, если для некоторого эндоморфизма $R$

$$X'=R\circ X\quad m'=R\circ m$$

то

$$W(f,X,m)\circ X'=w(f,X',m')=m'$$

Отображение $W(f,X,m)$ называется **координатами элемента** $m$ **относительно множества** $X$. Аналогично, мы можем рассмотреть координаты множества $B\subset J(f,X)$ относительно множества $X$

$$W(f,X,B)=\{W(f,X,m):m\in B\}=(W(f,X,m),m\in B)$$

Обозначим $W(f,X)$ **множество координат представления** $J(f,X)$. □

Теорема 2.6.9. *На множестве координат $W(f,X)$ определена структура $\Omega_2$-алгебры.*

Доказательство. Пусть $\omega\in\Omega_2(n)$. Тогда для любых $m_1$, ..., $m_n\in J(f,X)$ положим

$$W(f,X,m_1)...W(f,X,m_n)\omega=W(f,X,m_1...m_n\omega) \tag{2.6.5}$$

Согласно замечанию 2.6.8, из равенства (2.6.5) следует

$$\begin{aligned}(W(f,X,m_1)...W(f,X,m_n)\omega)\circ X&=W(f,X,m_1...m_n\omega)\circ X\\&=w(f,X,m_1...m_n\omega)\end{aligned} \tag{2.6.6}$$

Согласно правилу (2) определения 2.6.6, из равенства (2.6.6) следует

$$\begin{aligned}(W(f,X,m_1)...W(f,X,m_n)\omega)\circ X&=w(f,X,m_1)...w(f,X,m_n)\omega\\&=(W(f,X,m_1)\circ X)...(W(f,X,m_n)\circ X)\omega\end{aligned} \tag{2.6.7}$$

Из равенства (2.6.7) следует корректность определения (2.6.5) операции $\omega$ на множестве координат. □

Теорема 2.6.10. *Определено представление $\Omega_1$-алгебры $A$ в $\Omega_2$-алгебре $W(f,X)$.*

Доказательство. Пусть $a\in A$. Тогда для любого $m\in J(f,X)$ положим

$$aW(f,X,m)=W(f,X,am) \tag{2.6.8}$$

Согласно замечанию 2.6.8, из равенства (2.6.8) следует

$$(aW(f,X,m))\circ X=W(f,X,am)\circ X=w(f,X,am) \tag{2.6.9}$$

Согласно правилу (3) определения 2.6.6, из равенства (2.6.9) следует

$$(aW(f,X,m))\circ X=aw(f,X,m)=a(W(f,X,m)\circ X) \tag{2.6.10}$$

Из равенства (2.6.10) следует корректность определения (2.6.8) представления $\Omega_1$-алгебры $A$ в $\Omega_2$-алгебре $W(f,X)$. □

Теорема 2.6.11. *Пусть*

$$f:A\xrightarrow{*} M$$

*представление $\Omega_1$-алгебры $A$ в $\Omega_2$-алгебре $M$. Для заданных множеств $X\subset M$, $X'\subset M$, пусть отображение*

$$R_1:X\to X'$$

*согласовано со структурой представления $f$, т. е.*

$$\omega \in \Omega_2(n)\ x_1,\ ...,\ x_n,\ x_1...x_n\omega \in X,\ R_1(x_1...x_n\omega) \in X'$$
$$=> R_1(x_1...x_n\omega) = R_1(x_1)...R_1(x_n)\omega$$
$$x \in X,\ a \in A,\ R_1(ax) \in X'$$
$$=> R_1(ax) = aR_1(x)$$

*Рассмотрим отображение $\Omega_2$-слов*

$$w[f, X, R_1] : w(f, X) \to w(f, X')$$

*удовлетворяющее условиям 2.6.7.1, 2.6.7.2, 2.6.7.3, и такое, что*

$$x \in X => w[f, X, R_1](x) = R_1(x)$$

*Существует единственный эндоморфизм $\Omega_2$-алгебры $M$*

$$R : M \to M$$

*определённый правилом*

$$R \circ m = w[f, X, R_1](w(f, X, m))$$

*который является морфизмом представлений $J(f, X)$ и $J(f, X')$.*

Доказательство. Мы будем доказывать теорему индукцией по сложности $\Omega_2$-слова.

Если $w(f, X, m) = m$, то $m \in X$. Согласно условию 2.6.7.1,

$$R \circ m = w[f, X, R_1](w(f, X, m)) = w[f, X, R_1](m) = R_1(m)$$

Следовательно, на множестве $X$ отображения $R$ и $R_1$ совпадают, и отображение $R$ согласовано со структурой представления $f$.

Пусть $\omega \in \Omega_2(n)$. Пусть отображение $R$ определено для $m_1$, ..., $m_n \in J(f, X)$. Пусть

$$w_1 = w(f, X, m_1) \quad ... \quad w_n = w(f, X, m_n)$$

Если $m = m_1...m_n\omega$, то согласно правилу (2) определения 2.6.6,

$$w(f, X, m) = w_1...w_n\omega$$

Согласно условию 2.6.7.2,

$$\begin{aligned} R \circ m &= w[f, X, R_1](w(f, X, m)) = w[f, X, R_1](w_1...w_n\omega) \\ &= w[f, X, R_1](w_1)...w[f, X, R_1](w_n)\omega \\ &= (R \circ m_1)...(R \circ m_n)\omega \end{aligned}$$

Следовательно, отображение $R$ является эндоморфизмом $\Omega_2$-алгебры $M$.

Пусть отображение $R$ определено для $m_1 \in J(f, X)$, $w_1 = w(f, X, m_1)$. Пусть $a \in A$. Если $m = am_1$, то согласно правилу (3) определения 2.6.6,

$$w(f, X, am_1) = aw_1$$

Согласно условию 2.6.7.3,

$$\begin{aligned} R \circ m &= w[f, X, R_1](w(f, X, m)) = w[f, X, R_1](aw_1) \\ &= aw[f, X, R_1](w_1) = aR \circ m_1 \end{aligned}$$

Из равенства (2.2.4) следует, что отображение $R$ является морфизмом представления $f$.

Единственность эндоморфизма $R$, а следовательно, корректность его определения, следует из следующего рассуждения. Допустим, $m \in M$ имеет различные $\Omega_2$-слова относительно множества $X$, например

$$m = x_1 ... x_n \omega = ax \tag{2.6.11}$$

Так как $R$ - эндомоморфизм представления, то из равенства (2.6.11) следует

$$R \circ m = R \circ (x_1 ... x_n \omega) = (R \circ x_1) ... (R \circ x_n) \omega = R \circ (ax) = aR \circ x \tag{2.6.12}$$

Из равенства (2.6.12) следует

$$R \circ m = (R \circ x_1) ... (R \circ x_n) \omega = aR \circ x \tag{2.6.13}$$

Из равенств (2.6.11), (2.6.13) следует, что равенство (2.6.11) сохраняется при отображении. Следовательно, образ $m$ не зависит от выбора координат. □

Замечание 2.6.12. Теорема 2.6.11 - это теорема о продолжении отображения. Единственное, что нам известно о множестве $X$ - это то, что $X$ - множество образующих представления $f$. Однако, между элементами множества $X$ могут существовать соотношения, порождённые либо операциями $\Omega_2$-алгебры $M$, либо преобразованиями представления $f$. Поэтому произвольное отображение множества $X$, вообще говоря, не может быть продолжено до эндоморфизма представления $f$.[2.7] Однако, если отображение $R_1$ согласованно со структурой представления на множестве $X$, то мы можем построить продолжение этого отображения, которое является эндоморфизмом представления $f$. □

Определение 2.6.13. Пусть $X$ - множество образующих представления $f$. Пусть $R$ - эндоморфизм представления $f$. Множество координат $W(f, X, R \circ X)$ называется **координатами эндоморфизма представления**. □

Определение 2.6.14. Пусть $X$ - множество образующих представления $f$. Пусть $R$ - эндоморфизм представления $f$. Пусть $m \in M$. Мы определим **суперпозицию координат** эндоморфизма представления $f$ и элемента $m$ как координаты, определённые согласно правилу

$$W(f, X, m) \circ W(f, X, R \circ X) = W(f, X, R \circ m) \tag{2.6.14}$$

Пусть $Y \subset M$. Мы определим суперпозицию координат эндоморфизма представления $f$ и множества $Y$ согласно правилу

$$W(f, X, Y) \circ W(f, X, R \circ X) = (W(f, X, m) \circ W(f, X, R \circ X), m \in Y) \tag{2.6.15}$$

□

Теорема 2.6.15. *Эндоморфизм $R$ представления*

$$f : A \xrightarrow{*} M$$

*порождает отображение координат представления*

$$W[f, X, R] : W(f, X) \to W(f, X) \tag{2.6.16}$$

[2.7] В теореме 2.7.7, требования к множеству образующих более жёсткие. Поэтому теорема 2.7.7 говорит о продолжении произвольного отображения. Более подробный анализ дан в замечании 2.7.9.

*такое, что*

$$W(f,X,m)\to W[f,X,R]\star W(f,X,m)=W(f,X,R\circ m)$$
$$W[f,X,R]\star W(f,X,m)=W(f,X,m)\circ W(f,X,R\circ X) \tag{2.6.17}$$

Доказательство. Согласно замечанию 2.6.8, мы можем рассматривать равенства (2.6.14), (2.6.16) относительно заданного множества образующих $X$. При этом координатам $W(f,X,m)$ соответствует слово

$$W(f,X,m)\circ X=w(f,X,m) \tag{2.6.18}$$

а координатам $W(f,X,R\circ m)$ соответствует слово

$$W(f,X,R\circ m)\circ X=w(f,X,R\circ m) \tag{2.6.19}$$

Поэтому для того, чтобы доказать теорему, нам достаточно показать, что отображению $W[f,X,R]$ соответствует отображение $w[f,X,R]$. Мы будем доказывать это утверждение индукцией по сложности $\Omega_2$-слова.

Если $m\in X$, $m'=R\circ m$, то, согласно равенствам (2.6.18), (2.6.19), отображения $W[f,X,R]$ и $w[f,X,R]$ согласованы.

Пусть для $m_1$, ..., $m_n\in X$ отображения $W[f,X,R]$ и $w[f,X,R]$ согласованы. Пусть $\omega\in\Omega_2(n)$. Согласно теореме 2.6.9

$$W(f,X,m_1...m_n\omega)=W(f,X,m_1)...W(f,X,m_n)\omega \tag{2.6.20}$$

Так как $R$ - эндоморфизм $\Omega_2$-алгебры $M$, то из равенства (2.6.20) следует

$$\begin{aligned}W(f,X,R\circ(m_1...m_n\omega))&=W(f,X,(R\circ m_1)...(R\circ m_n)\omega)\\&=W(f,X,R\circ m_1)...W(f,X,R\circ m_n)\omega\end{aligned} \tag{2.6.21}$$

Из равенств (2.6.20), (2.6.21) и предположения индукции следует, что отображения $W[f,X,R]$ и $w[f,X,R]$ согласованы для $m=m_1...m_n\omega$.

Пусть для $m_1\in M$ отображения $W[f,X,R]$ и $w[f,X,R]$ согласованы. Пусть $a\in A$. Согласно теореме 2.6.10

$$W(f,X,am_1)=aW(f,X,m_1) \tag{2.6.22}$$

Так как $R$ - эндоморфизм представления $f$, то из равенства (2.6.22) следует

$$W(f,X,R\circ(am_1))=W(f,X,a\,R\circ m_1)=aW(f,X,R\circ m_1) \tag{2.6.23}$$

Из равенств (2.6.22), (2.6.23) и предположения индукции следует, что отображения $W[f,X,R]$ и $w[f,X,R]$ согласованы для $m=am_1$. □

Следствие 2.6.16. *Пусть $X$ - множество образующих представления $f$. Пусть $R$ - эндоморфизм представления $f$. Отображение $W[f,X,R]$ является эндоморфизмом представления $\Omega_1$-алгебры $A$ в $\Omega_2$-алгебре $W(f,X)$.* □

В дальнейшем мы будем отождествлять отображение $W[f,X,R]$ и множество координат $W(f,X,R\circ X)$.

Теорема 2.6.17. *Пусть $X$ - множество образующих представления $f$. Пусть $R$ - эндоморфизм представления $f$. Пусть $Y\subset M$. Тогда*

$$W(f,X,Y)\circ W(f,X,R\circ X)=W(f,X,R\circ Y) \tag{2.6.24}$$

$$W[f,X,R]\star W(f,X,Y)=W(f,X,R\circ Y) \tag{2.6.25}$$

Доказательство. Равенство (2.6.24) является следствием равенства

$$R\circ Y=(R\circ m,m\in Y)$$

а также равенств (2.6.14), (2.6.15). Равенство (2.6.25) является следствием равенств (2.6.24), (2.6.17). □

Теорема 2.6.18. *Пусть $X$ - множество образующих представления $f$. Пусть $R$, $S$ - эндоморфизмы представления $f$. Тогда*

$$W(f,X,S\circ X)\circ W(f,X,R\circ X)=W(f,X,R\circ S\circ X) \tag{2.6.26}$$

$$W[f,X,R]\star W[f,X,S]=W[f,X,R\circ S] \tag{2.6.27}$$

Доказательство. Равенство (2.6.26) следует из равенства (2.6.24), если положить $Y=S\circ X$. Равенство (2.6.27) следует из равенства (2.6.26) и цепочки равенств

$$\begin{aligned}
&(W[f,X,R]\star W[f,X,S])\star W(f,X,Y)\\
=&W[f,X,R]\star(W[f,X,S]\star W(f,X,Y))\\
=&(W(f,X,Y)\circ W(f,X,S\circ X))\circ W(f,X,R\circ X)\\
=&W(f,X,Y)\circ(W(f,X,S\circ X)\circ W(f,X,R\circ X))\\
=&W(f,X,Y)\circ W(f,X,R\circ S\circ X)\\
=&W(f,X,R\circ S)\star W(f,X,Y)
\end{aligned}$$

□

Концепция суперпозиции координат очень проста и напоминает своеобразную машину Тюринга. Если элемент $m\in M$ имеет вид

$$m=m_1...m_n\omega$$

или

$$m=am_1$$

то мы ищем координаты элементов $m_i$, для того чтобы подставить их в соответствующее выражение. Как только элемент $m\in M$ принадлежит множеству образующих $\Omega_2$-алгебры $M$, мы выбираем координаты соответствующего элемента из второго множителя. Поэтому мы требуем, чтобы второй множитель в суперпозиции был множеством координат образа множества образующих $X$.

Мы можем обобщить определение суперпозиции координат и предположить, что один из множителей является множеством $\Omega_2$-слов. Соответственно, определение суперпозиции координат имеет вид

$$w(f,X,Y)\circ W(f,X,R\circ X)=W(f,X,Y)\circ w(f,X,R\circ X)=w(f,X,R\circ Y)$$

Следующие формы записи образа множества $Y$ при эндоморфизме $R$ эквивалентны.

$$R\circ Y=W(f,X,Y)\circ(R\circ X)=W(f,X,Y)\circ(W(f,X,R\circ X)\circ X) \tag{2.6.28}$$

Из равенств (2.6.24), (2.6.28) следует, что

$$(W(f,X,Y)\circ W(f,X,R\circ X))\circ X=W(f,X,Y)\circ(W(f,X,R\circ X)\circ X) \tag{2.6.29}$$

Равенство (2.6.29) является законом ассоциативности для операции композиции и позволяет записать выражение

$$W(f,X,Y)\circ W(f,X,R\circ X)\circ X$$

без использования скобок.

Рассмотрим равенство (2.6.26) где мы видим изменение порядка эндоморфизмов в суперпозиции координат. Это равенство также следует из цепочки равенств, где мы можем непосредственно видеть, когда изменяется порядок эндоморфизмов

$$\begin{aligned} &W(f,X,m)\circ W(f,X,R\circ S\circ X)\circ X\\ =&(R\circ S)\circ m=R\circ(S\circ m)\\ =&R\circ((W(f,X,m)\circ W(f,X,S\circ X))\circ X)\\ =&W(f,X,m)\circ W(f,X,S\circ X)\circ W(f,X,R\circ X)\circ X \end{aligned} \tag{2.6.30}$$

Из равенства (2.6.30) следует, что координаты эндоморфизма действуют на координаты элементы $\Omega_2$-алгебры $M$ справа.

Определение 2.6.19. Пусть $X\subset M$ - множество образующих представления

$$f:A \xrightarrow{*} M$$

Пусть отображение

$$H:M\to M$$

является эндоморфизмом представления $f$. Пусть множество $X'=H\circ X$ является образом множества $X$ при отображении $H$. Эндоморфизм $H$ представления $f$ называется **невырожденным на множестве образующих** $X$, если множество $X'$ является множеством образующих представления $f$. В противном случае, эндоморфизм $H$ представления $f$ называется **вырожденным на множестве образующих** $X$, □

Определение 2.6.20. Эндоморфизм представления $f$ называется **невырожденным**, если он невырожден на любом множестве образующих. □

Теорема 2.6.21. *Автоморфизм $R$ представления*

$$f:A \xrightarrow{*} M$$

*является невырожденным эндоморфизмом.*

Доказательство. Пусть $X$ - множество образующих представления $f$. Пусть $X'=R\circ X$.

Согласно теореме 2.6.7 эндоморфизм $R$ порождает отображение $\Omega_2$-слов $w[f,X,R]$.

Пусть $m'\in M$. Так как $R$ - автоморфизм, то существует $m\in M$, $R\circ m=m'$. Согласно определению 2.6.6, $w(f,X,m)$ - $\Omega_2$-слово, представляющее $m$ относительно множества образующих $X$. Согласно теореме 2.6.7, $w(f,X',m')$ - $\Omega_2$-слово, представляющее $m'$ относительно множества образующих $X'$

$$w(f,X',m')=w[f,X,R](w(f,X,m))$$

Следовательно, $X'$ - множество образующих представления $f$. Согласно определению 2.6.20, автоморфизм $R$ - невырожден. □

## 2.7. Базис представления

Определение 2.7.1. Если множество $X \subset M$ является множеством образующих представления $f$, то любое множество $Y$, $X \subset Y \subset M$ также является множеством образующих представления $f$. Если существует минимальное множество $X$, порождающее представление $f$, то такое множество $X$ называется **базисом представления** $f$. □

Теорема 2.7.2. *Множество образующих $X$ представления $f$ является базисом тогда и только тогда, когда для любого $m \in X$ множество $X \setminus \{m\}$ не является множеством образующих представления $f$.*

Доказательство. Пусть $X$ - множество образующих расслоения $f$. Допустим для некоторого $m \in X$ существует $\Omega_2$-слово

(2.7.1) $$w = w(f, X \setminus \{m\}, m)$$

Рассмотрим элемент $m'$, для которого $\Omega_2$-слово

(2.7.2) $$w' = w(f, X, m')$$

зависит от $m$. Согласно определению 2.6.6, любое вхождение $m$ в $\Omega_2$-слово $w'$ может быть заменено $\Omega_2$-словом $w$. Следовательно, $\Omega_2$-слово $w'$ не зависит от $m$, а множество $X \setminus \{m\}$ является множеством образующих представления $f$. Следовательно, $X$ не является базисом расслоения $f$. □

Замечание 2.7.3. Доказательство теоремы 2.7.2 даёт нам эффективный метод построения базиса представления $f$. Выбрав произвольное множество образующих, мы шаг за шагом исключаем те элементы множества, которые имеют координаты относительно остальных элементов множества. Если множество образующих представления бесконечно, то рассмотренная операция может не иметь последнего шага. Если представление имеет конечное множество образующих, то за конечное число шагов мы можем построить базис этого представления.

Как отметил Кон в [13], стр. 96, 97, представление может иметь неэквивалентные базисы. Например, циклическая группа шестого порядка имеет базисы $\{a\}$ и $\{a^2, a^3\}$, которые нельзя отобразить один в другой эндоморфизмом представления. □

Замечание 2.7.4. Мы будем записывать базис также в виде

$$X = (x, x \in X)$$

Если базис - конечный, то мы будем также пользоваться записью

$$X = (x_i, i \in I) = (x_1, ..., x_n)$$

□

Замечание 2.7.5. Мы ввели $\Omega_2$-слово элемента $x \in M$ относительно множества образующих $X$ в определении 2.6.6. Из теоремы 2.7.2 следует, что если множество образующих $X$ не является базисом, то выбор $\Omega_2$-слова относительно множества образующих $X$ неоднозначен. Но даже если множество образующих $X$ является базисом, то представление $m \in M$ в виде $\Omega_2$-слова неоднозначно. Если $m_1$, ..., $m_n$ - $\Omega_2$-слова, $\omega\in\Omega_2(n)$ и $a\in A$, то[2.8]

$$a(m_1...m_n\omega)=(am_1)...(am_n)\omega \tag{2.7.3}$$

Возможны равенства, связанные со спецификой представления. Например, если $\omega$ является операцией $\Omega_1$-алгебры $A$ и операцией $\Omega_2$-алгебры $M$, то мы можем потребовать, что $\Omega_2$-слова $a_1...a_n\omega x$ и $a_1x...a_nx\omega$ описывают один и тот же элемент $\Omega_2$-алгебры $M$.[2.9]

Помимо перечисленных равенств в $\Omega_2$-алгебре могут существовать соотношения вида

$$w_1(f,X,m)=w_2(f,X,m)\quad m\notin X \tag{2.7.4}$$

Особенностью равенства (2.7.4) является то, что это равенство не может быть приведенно к более простому виду.[2.10]

На множестве $\Omega_2$-слов $w(f,X)$, перечисленные выше равенства определяют **отношение эквивалентности** $\rho(f)$**, порождённое представлением** $f$. Очевидно, что для произвольного $m\in M$ выбор соответствующего $\Omega_2$-слова однозначен с точностью до отношения эквивалентности $\rho(f)$. Тем не менее, если в процессе построений мы получим равенство двух $\Omega_2$-слов относительно заданного базиса, мы можем утверждать, что эти $\Omega_2$-слова равны, не заботясь об отношении эквивалентности $\rho(f)$.

Аналогичное замечание касается отображения $W(f,X,m)$, определённого в замечании 2.6.8.[2.11] □

Теорема 2.7.6. *Автоморфизм представления $f$ отображает базис представления $f$ в базис.*

---

[2.8]Например, пусть $\{\overline{e}_1,\overline{e}_2\}$ - базис векторного пространства над полем $k$. Равенство (2.7.3) принимает форму закона дистрибутивности

$$a(b^1\overline{e}_1+b^2\overline{e}_2)=(ab^1)\overline{e}_1+(ab^2)\overline{e}_2$$

[2.9]Для векторного пространства это требование принимает форму закона дистрибутивности

$$(a+b)\overline{e}_1=a\overline{e}_1+b\overline{e}_1$$

[2.10]Смотри, например, подразделы 4.7.2, 4.7.3.

[2.11]Если базис векторного пространства - конечен, то мы можем представить базис в виде матрицы строки

$$\overline{\overline{e}}=\begin{pmatrix}\overline{e}_1 & ... & \overline{e}_2\end{pmatrix}$$

Мы можем представить отображение $W(f,\overline{\overline{e}},\overline{v})$ в виде матрицы столбца

$$W(f,\overline{\overline{e}},\overline{v})=\begin{pmatrix}v^1\\ ...\\ v^n\end{pmatrix}$$

Тогда

$$W(f,\overline{\overline{e}},\overline{v})(\overline{\overline{e}}')=W(f,\overline{\overline{e}},\overline{v})\begin{pmatrix}\overline{e}'_1 & ... & \overline{e}'_n\end{pmatrix}=\begin{pmatrix}v^1\\ ...\\ v^n\end{pmatrix}\begin{pmatrix}\overline{e}'_1 & ... & \overline{e}'_n\end{pmatrix}$$

имеет вид произведения матриц.

Доказательство. Пусть отображение $R$ - автоморфизм представления $f$. Пусть множество $X$ - базис представления $f$. Пусть $X' = R \circ X$.

Допустим множество $X'$ не является базисом. Согласно теореме 2.7.2 существует $m' \in X'$ такое, что $X' \setminus \{x'\}$ является множеством образующих представления $f$. Согласно теореме 2.4.5 отображение $R^{-1}$ является автоморфизмом представления $f$. Согласно теореме 2.6.21 и определению 2.6.20, множество $X \setminus \{m\}$ является множеством образующих представления $f$. Полученное противоречие доказывает теорему. □

Теорема 2.7.7. *Пусть $X$ - базис представления $f$. Пусть*

$$R_1 : X \to X'$$

*произвольное отображение множества $X$. Рассмотрим отображение $\Omega_2$-слов*

$$w[f, X, R_1] : w(f, X) \to w(f, X')$$

*удовлетворяющее условиям 2.6.7.1, 2.6.7.2, 2.6.7.3, и такое, что*

$$x \in X => w[f, X, R_1](x) = R_1(x)$$

*Существует единственный эндоморфизм представлнения $f$*[2.12]

$$R : M \to M$$

*определённый правилом*

$$R \circ m = w[f, X, R_1](w(f, X, m))$$

Доказательство. Утверждение теоремы является следствием теорем 2.6.7, 2.6.11. □

Следствие 2.7.8. *Пусть $X$, $X'$ - базисы представления $f$. Пусть $R$ - автоморфизм представления $f$ такой, что $X' = R \circ X$. Автоморфизм $R$ определён однозначно.* □

Замечание 2.7.9. Теорема 2.7.7, так же как и теорема 2.6.11, является теоремой о продолжении отображения. Однако здесь $X$ - не произвольное множество образующих представления, а базис. Согласно замечанию 2.7.3, мы не можем определить координаты любого элемента базиса через остальные элементы этого же базиса. Поэтому отпадает необходимость в согласованности отображения базиса с представлением. □

Теорема 2.7.10. *Набор координат $W(f, X, X)$ соответствует тождественному преобразованию*

$$W[f, X, E] = W(f, X, X)$$

Доказательство. Утверждение теоремы следует из равенства

$$m = W(f, X, m) \circ X = W(f, X, m) \circ W(f, X, X) \circ X$$

□

[2.12] Это утверждение похоже на теорему [1]-1, с. 104.

Теорема 2.7.11. *Пусть $W(f,X,R\circ X)$ - множество координат автоморфизма $R$. Определено множество координат $W(f,R\circ X,X)$, соответствующее автоморфизму $R^{-1}$. Множество координаты $W(f,R\circ X,X)$ удовлетворяют равенству*

$$W(f,X,R\circ X)\circ W(f,R\circ X,X)=W(f,X,X) \tag{2.7.5}$$
$$W[f,X,R^{-1}]=W[f,X,R]^{-1}=W(f,R\circ X,X)$$

Доказательство. Поскольку $R$ - автоморфизм представления $f$, то, согласно теореме 2.7.6, множество $R\circ X$ - базис представления $f$. Следовательно, определено множество координат $W(f,R\circ X,X)$. Равенство (2.7.5) следует из цепочки равенств

$$\begin{aligned}W(f,X,R\circ X)\circ W(f,R\circ X,X)&=W(f,X,R\circ X)\circ W(f,X,R^{-1}\circ X)\\&=W(f,X,R^{-1}\circ R\circ X)=W(f,X,X)\end{aligned}$$

□

Замечание 2.7.12. В $\Omega_2$-алгебре $M$ не существует универсального алгоритма определения множества координат $W(f,R\circ X,X)$ для заданного множества $W(f,X,R\circ X)$.[2.13] Мы полагаем, что в теореме 2.7.11 этот алгоритм задан неявно. Очевидно также, что множество $\Omega_2$-слов

$$W(f,X,R\circ X)\circ W(f,R\circ X,X)\circ X \tag{2.7.6}$$

вообще говоря, не совпадает с множеством $\Omega_2$-слов

$$W(f,X,X)\circ X \tag{2.7.7}$$

Теорема 2.7.11 утверждает, что множества $\Omega_2$-слов (2.7.6) и (2.7.7) совпадают с точностью до отношения эквивалентности, порождённой представлением $f$.
□

Теорема 2.7.13. *Группа автоморфизмов $G(f)$ эффективного представления $f$ в $\Omega_2$-алгебре $M$ порождает эффективное представление в $\Omega_2$-алгебре $M$.*

Доказательство. Из следствия 2.7.8 следует, что если автоморфизм $R$ отображает базис $X$ в базис $X'$, то множество координат $W(f,X,X')$ однозначно определяет автоморфизм $R$. Из теоремы 2.6.15 следует, что множество координат $W(f,X,X')$ определяет правило отображения координат относительно базиса $X$ при автоморфизме представления $f$. Из равенства (2.6.28) следует, что автоморфизм $R$ действует справа на элементы $\Omega_2$-алгебры $M$. Из равенства (2.6.26) следует, что представление группы является правосторонним представлением. Согласно теореме 2.7.10 набор координат $W(f,X,X)$ соответствует тождественному преобразованию. Из теоремы 2.7.11 следует, что набор координат $W(f,R\circ X,X)$ соответствует преобразованию, обратному преобразованию $W(f,X,R\circ X)$. □

---

[2.13] В векторном пространстве линейному преобразованию соответствует матрица чисел. Соответственно, обратному преобразованию соответствует обратная матрица.

## 2.8. Примеры базиса представления универсальной алгебры

**2.8.1. Векторное пространство.** Рассмотрим векторное пространство $\overline{V}$ над полем $F$. Если дано множество векторов $\overline{e}_1$, ..., $\overline{e}_n$, то, согласно алгоритму построения координат над векторным пространством, координаты включают такие элементы как $\overline{e}_1+\overline{e}_2$ и $a\overline{e}_1$. Рекурсивно применяя правила, приведенные в определении 2.6.6, мы придём к выводу, что множество векторов $\overline{e}_1$, ..., $\overline{e}_n$, порождает множество линейных комбинаций

$$a^1\overline{e}_1 + ... + a^n\overline{e}_n$$

Согласно теореме 2.7.2, множество векторов $\overline{e}_1$, ..., $\overline{e}_n$, является базисом при условии, если для любого $i$, $i = 1, ..., n$, вектор $\overline{e}_i$ не является линейной комбинацией остальных векторов. Это требование равносильно требованию линейной независимости векторов.

**2.8.2. Представление группы на множестве.** Рассмотрим представление из примера 2.5.2. Мы можем рассматривать множество $M$ как объединение орбит представления группы $G$. В качестве базиса представления можно выбрать множество точек таким образом, что одна и только одна точка принадлежит каждой орбите представления. Если $X$ - базис представления, $A \in X$, $g \in G$, то $\Omega_2$-слово имеет вид $A+g$. Поскольку на множестве $M$ не определены операции, то не существует $\Omega_2$-слово, содержащее различные элементы базиса. Если представление группы $G$ однотранзитивно, то базис представления состоит из одной точки. Этой точкой может быть любая точка множества $M$.

Глава 3

# Представление группы

## 3.1. Представление группы

Группа - одна из немногих алгебр, которая позволяет рассматривать произведение преобразований $\Omega$-алгебры $M$ таким образом, что если преобразования принадлежат представлению, то их произведение также принадлежит представлению. При этом следует помнить, что порядок отображений при суперпозиции зависит от порядка отображений на диаграмме и с какой стороны отображения действуют на элементы множества.

Определение 3.1.1. Пусть ${}^*M$ - группа с произведением

$$(f\circ g)x=f(gx)$$

и $\delta$ - единица группы ${}^\star M$. Пусть $G$ - группа. Мы будем называть гомоморфизм групп

$$f:G\rightarrow{}^\star M \tag{3.1.1}$$

**левосторонним представлением группы** $G$ или $G*$**-представлением группы** в $\Omega$-алгебре $M$, если отображение $f$ удовлетворяет условиям

$$f(ab)u=f(a)(f(b)u) \tag{3.1.2}$$

□

Замечание 3.1.2. Поскольку отображение (3.1.1) - гомоморфизм, то

$$f(ab)u=(f(a)f(b))u \tag{3.1.3}$$

Мы здесь пользуемся соглашением

$$f(a)f(b)=f(a)\circ f(b)$$

Таким образом, концепция представления групп состоит в том, что в каком порядке мы перемножаем элементы группы, в том же порядке перемножаются соответствующие преобразования представления. Из равенств (3.1.2) и (3.1.3) следует

$$(f(a)f(b))u=f(a)(f(b)u) \tag{3.1.4}$$

Равенство (3.1.4) совместно с ассоциативностью произведения преобразований представляет собой **закон ассоциативности** для $G*$-представления. Это позволяет записывать равенство (3.1.4) без использования скобок

$$f(ab)u=f(a)f(b)u$$

□

Определение 3.1.3. Пусть $M^\star$ - группа с произведением

$$x(f\circ g)=(xf)g$$

и $\delta$ - единица группы $M^\star$. Пусть $G$ - группа. Мы будем называть гомоморфизм групп

$$f:G\rightarrow M^* \tag{3.1.5}$$

**правосторонним представлением группы** $G$ или $*G$**-представлением** в $\Omega$-алгебре $M$, если отображение $f$ удовлетворяет условиям

$$uf(ab)=(uf(a))f(b) \tag{3.1.6}$$

□

Замечание 3.1.4. Поскольку отображение (3.1.5) - гомоморфизм, то

$$uf(ab)=u(f(a)f(b)) \tag{3.1.7}$$

Из равенств (3.1.6) и (3.1.7) следует

$$u(f(a)f(b))=(uf(a))f(b) \tag{3.1.8}$$

Равенство (3.1.8) совместно с ассоциативностью произведения преобразований представляет собой **закон ассоциативности** для $*G$-представления. Это позволяет записывать равенство (3.1.8) без использования скобок

$$uf(ab)=uf(a)f(b)$$

□

Определение 3.1.5. Мы будем называть преобразование

$$t:M\rightarrow M$$

**невырожденным преобразованием**, если существует обратное отображение. □

Теорема 3.1.6. *Для любого $g\in G$ преобразование $f(g)$ является невырожденным и удовлетворяет равенству*

$$f(g^{-1})=f(g)^{-1} \tag{3.1.9}$$

Доказательство. На основании (3.1.2) и

$$f(e)=\delta$$

мы можем записать

$$u=\delta(u)=f(gg^{-1})(u)=f(g)(f(g^{-1})(u))$$

Это завершает доказательство. □

Теорема 3.1.7. *Групповая операция определяет два различных представления на группе:*

- **Левый сдвиг** $t_\star$

$$\begin{aligned}b'&=t_\star(a)b=ab\\ b'&=t_\star(a)(b)=ab\end{aligned} \tag{3.1.10}$$

*является $G*$-представлением на множестве*[3.1] $G$

$$t_\star(ab) = t_\star(a) \circ t_\star(b) \tag{3.1.11}$$

• **Правый сдвиг** $_\star t$

$$\begin{aligned} b' &= b\ _\star t(a) = ba \\ b' &= {}_\star t(a)(b) = ba \end{aligned} \tag{3.1.12}$$

*является $*G$-представлением на множестве $G$*

$${}_\star t(ab) = {}_\star t(a) \circ {}_\star t(b) \tag{3.1.13}$$

Доказательство. Равенство (3.1.11) следует из ассоциативности произведения

$$t_\star(ab)c = (ab)c = a(bc) = t_\star(a)(t_\star(b)c) = (t_\star(a) \circ t_\star(b))c$$

Аналогично доказывается равенство (3.1.13). □

Определение 3.1.8. Пусть $G$ - группа. Пусть $f$ - $G*$-представление в $\Omega$-алгебре $M$. Для любого $v \in M$ мы определим **орбиту представления группы** $G$ как множество

$$f(G)v = \{w = f(g)v : g \in G\}$$

□

Так как $f(e) = \delta$, то $v \in f(G)v$.

Теорема 3.1.9. *Если*

$$v \in f(G)u \tag{3.1.14}$$

*то*

$$f(G)u = f(G)v$$

Доказательство. Из (3.1.14) следует существование $a \in G$ такого, что

$$v = f(a)u \tag{3.1.15}$$

Если $w \in f(G)v$, то существует $b \in G$ такой, что

$$w = f(b)v \tag{3.1.16}$$

Подставив (3.1.15) в (3.1.16), мы получим

$$w = f(b)(f(a)u) \tag{3.1.17}$$

На основании (3.1.2) из (3.1.17) следует, что $w \in f(G)u$. Таким образом,

$$f(G)v \subseteq f(G)u$$

На основании (3.1.9) из (3.1.15) следует, что

$$u = f(a)^{-1}v = f(a^{-1})v \tag{3.1.18}$$

Равенство (3.1.18) означает, что $u \in f(G)v$ и, следовательно,

$$f(G)u \subseteq f(G)v$$

Это завершает доказательство. □

[3.1] Левый сдвиг не является представлением группы в группе, так как преобразование $t_\star$ не является гомоморфизмом группы. Аналогичное замечание верно для правого сдвига.

Таким образом, $G*$-представление $f$ в $\Omega$-алгебре $M$ порождает отношение эквивалентности $S$ и орбита $f(G)u$ является классом эквивалентности. Мы будем пользоваться обозначением $M/f(G)$ для фактор множества $M/S$ и мы будем называть это множество **пространством орбит $G*$-представления** $f$.

Теорема 3.1.10. *Если определены $G*$-представление $f_1$ в $\Omega$-алгебре $M_1$ и $G*$-представление $f_2$ в $\Omega$-алгебре $M_2$, то мы можем определить* **прямое произведение $G*$-представлений** *$f_1$ и $f_2$ группы*

$$f = f_1 \times f_2 : G \to M_1 \otimes M_2$$
$$f(g) = (f_1(g), f_2(g))$$

Доказательство. Чтобы показать, что $f$ является представлением, достаточно показать, что $f$ удовлетворяет определению 3.1.1.

$$f(e) = (f_1(e), f_2(e)) = (\delta_1, \delta_2) = \delta$$

$$\begin{aligned} f(ab)u &= (f_1(ab)u_1, f_2(ab)u_2) \\ &= (f_1(a)(f_1(b)u_1), f_2(a)(f_2(b)u_2)) \\ &= f(a)(f_1(b)u_1, f_2(b)u_2) \\ &= f(a)(f(b)u) \end{aligned}$$

□

## 3.2. Однотранзитивное правостороннее представление группы

Определение 3.2.1. Мы будем называть **ядром неэффективности $G*$-представления** множество

$$K_f = \{g \in G : f(g) = \delta\}$$

□

Теорема 3.2.2. *Ядро неэффективности $G*$-представления - это подгруппа группы $G$.*

Доказательство. Допустим $f(a_1) = \delta$ и $f(a_2) = \delta$. Тогда

$$f(a_1a_2)u = f(a_1)(f(a_2)u) = u$$
$$f(a^{-1}) = f^{-1}(a) = \delta$$

□

Теорема 3.2.3. *$G*$-представление* **эффективно** *тогда и только тогда, когда ядро неэффективности $K_f = \{e\}$.*

Доказательство. Утверждение является следствием определений 2.1.6 и 3.2.1 и теоремы 3.2.2. □

Если действие не эффективно, мы можем перейти к эффективному заменив группой $G_1 = G|K_f$, пользуясь факторизацией по ядру неэффективности. Это означает, что мы можем изучать только эффективное действие.

Определение 3.2.4. Рассмотрим $G*$-представление $f$ в $\Omega$-алгебре $M$. **Малая группа** или **группа стабилизации** элемента $x \in M$ - это множество

$$= \{g \in G : f(g)x = x\}$$

Мы будем называть $G*$-представление $f$ **свободным**, если для любого $x \in M$ группа стабилизации $G_x = \{e\}$. □

Теорема 3.2.5. *Если определено свободное $G*$-представление $f$ группы $G$ на $\Omega$-алгебре $A$, то определено взаимно однозначное соответствие между орбитами представления, а также между орбитой представления и группой $G$.*

Доказательство. Допустим для точки $a \in A$ существуют $g_1$, $g_2 \in G$

$$f(g_1)a = f(g_2)a \tag{3.2.1}$$

Умножим обе части равенства (3.2.1) на $f(g_1^{-1})$

$$a = f(g_1^{-1})f(g_2)a$$

Поскольку представление свободное, $g_1 = g_2$. Теорема доказана, так как мы установили взаимно однозначное соответствие между орбитой и группой $G$. □

Определение 3.2.6. Мы будем называть пространство $V$ **однородным пространством группы** $G$, если мы имеем однотранзитивное $G*$-представление на $V$. □

Теорема 3.2.7. *Если мы определим однотранзитивное представление $f$ группы $G$ на $\Omega$-алгебре $A$, то мы можем однозначно определить координаты на $A$, пользуясь координатами на группе $G$.*

*Если $f$ - левостороннее представление, то $f(a)$ эквивалентно левому сдвигу $t_\star(a)$ на группе $G$. Если $f$ - правостороннее представление, то $f(a)$ эквивалентно правому сдвигу ${}_\star t(a)$ на группе $G$.*

Доказательство. Мы выберем точку $v \in A$ и определим координаты точки $w \in A$ как координаты $a \in G$ такого, что $w = f(a)v$. Координаты, определённые таким образом, однозначны с точностью до выбора начальной точки $v \in A$, так как действие эффективно.

Если $f$ - левостороннее представление, мы будем пользоваться записью

$$f(a)v = av$$

Так как запись

$$f(a)(f(b)v) = a(bv) = (ab)v = f(ab)v$$

совместима с групповой структурой, мы видим, что левостороннее представление $f$ эквивалентно левому сдвигу.

Если $f$ - правостороннее представление, мы будем пользоваться записью

$$vf(a) = va$$

Так как запись

$$(vf(b))f(a) = (vb)a = v(ba) = vf(ba)$$

совместима с групповой структурой, мы видим, что правостороннее представление $f$ эквивалентно правому сдвигу. □

Замечание 3.2.8. Мы будем записывать эффективное $G*$-представление в форме

$$v' = t_\star(a)v = av$$

Орбита этого представления имеет вид

$$Gv = t_\star(G)v$$

Мы будем пользоваться обозначением $M/t_\star(G)$ для пространства орбит эффективного $G*$-представления. □

Замечание 3.2.9. Мы будем записывать эффективное $*G$-представление в форме

$$v' = v\ {}_\star t(a) = va$$

Орбита этого представления имеет вид

$$vG = v\ {}_\star t(G)$$

Мы будем пользоваться обозначением $M/{}_\star t(G)$ для пространства орбит эффективного $*G$-представления. □

Теорема 3.2.10. *Свободное $G*$-представление эффективно. Свободное $G*$-представление $f$ в $\Omega$-алгебре $M$ однотранзитивно на орбите.*

Доказательство. Следствие определения 3.2.4. □

Теорема 3.2.11. *Правый и левый сдвиги на группе $G$ перестановочны.*

Доказательство. Это следствие ассоциативности группы $G$

$$(t_\star(a) \circ\ {}_\star t(b))c = a(cb) = (ac)b = ({}_\star t(b) \circ t_\star(a))c$$

□

Теорема 3.2.11 может быть сформулирована следующим образом.

Теорема 3.2.12. *Пусть $G$ - группа. Для любого $a \in G$ отображение $t_\star(a)$ является автоморфизмом представления ${}_\star t$.*

Доказательство. Согласно теореме 3.2.11

$$t_\star(a) \circ {}_\star t(b) = {}_\star t(b) \circ t_\star(a) \tag{3.2.2}$$

Равенство (3.2.2) совпадает с равенством (2.2.3) из определения 2.2.2 при условии $r = id$, $R = t_\star(a)$. □

Теорема 3.2.13. *Пусть $G*$-представление $f$ на $\Omega$-алгебре $M$ однотранзитивно. Тогда мы можем однозначно определить однотранзитивное $*G$-представление $h$ на $\Omega$-алгебре $M$ такое, что диаграмма*

$$\begin{array}{ccc}
M & \xrightarrow{h(a)} & M \\
\downarrow{\scriptstyle f(b)} & & \downarrow{\scriptstyle f(b)} \\
M & \xrightarrow[h(a)]{} & M
\end{array}$$

*коммутативна для любых $a, b \in G$.*[3.2]

[3.2] Это утверждение можно также найти в [3].

Доказательство. Мы будем пользоваться групповыми координатами для точек $v \in M$. Тогда согласно теореме 3.2.7 мы можем записать левый сдвиг $t_\star(a)$ вместо преобразования $f(a)$.

Пусть $v_0, v \in M$. Тогда мы можем найти одно и только одно $a \in G$ такое, что

$$v = v_0 a = v_0\ {}_\star t(a)$$

Мы предположим

$$h(a) = {}_\star t(a)$$

Существует $b \in G$ такое, что

$$w_0 = f(b)v_0 = t_\star(b)v_0 \qquad w = f(b)v = t_\star(b)v$$

Согласно теореме 3.2.11 диаграмма

(3.2.3)
$$\begin{array}{ccc} v_0 & \xrightarrow{h(a)=\ {}_\star t(a)} & v \\ \Big\downarrow{\scriptstyle f(b)=t_\star(b)} & & \Big\downarrow{\scriptstyle f(b)=t_\star(b)} \\ w_0 & \xrightarrow[h(a)=\ {}_\star t(a)]{} & w \end{array}$$

коммутативна.

Изменяя $b$ мы получим, что $w_0$ - это произвольная точка, принадлежащая $M$.

Мы видим из диаграммы, что, если $v_0 = v$, то $w_0 = w$ и следовательно $h(e) = \delta$. С другой стороны, если $v_0 \neq v$, то $w_0 \neq w$ потому, что $G*$-представление $f$ однотранзитивно. Следовательно $*G$-представление $h$ эффективно.

Таким же образам мы можем показать, что для данного $w_0$ мы можем найти $a$ такое, что $w = h(a)w_0$. Следовательно $*G$-представление $h$ однотранзитивно.

В общем случае, произведение преобразований $G*$-представления $f$ не коммутативно и следовательно $*G$-представление $h$ отлично от $G*$-представления $f$. Таким же образом мы можем создать $G*$-представление $f$, пользуясь $*G$-представлением $h$. □

Мы будем называть представления $f$ и $h$ **парными представлениями группы** $G$.

Замечание 3.2.14. Очевидно, что преобразования $t_\star(a)$ и ${}_\star t(a)$ отличаются, если группа $G$ неабелева. Тем не менее, они являются отображениями на. Теорема 3.2.13 утверждает, что, если оба представления правого и левого сдвига существуют на множестве $M$, то мы можем определить два перестановочных представления на множестве $M$. Только правый или левый сдвиг не может представлять оба типа представления. Чтобы понять почему это так, мы можем изменить диаграмму (3.2.3) и предположить $h(a)v_0 = t_\star(a)v_0 = v$ вместо $h(a)v_0 = v_{0\star}t(a) = v$ и проанализировать, какое выражение $h(a)$ имеет в точке $w_0$. Диаграмма

$$\begin{array}{ccc} v_0 & \xrightarrow{h(a)=t_\star(a)} & v \\ \Big\downarrow{\scriptstyle f(b)=t_\star(b)} & & \Big\downarrow{\scriptstyle f(b)=t_\star(b)} \\ w_0 & \xrightarrow[h(a)]{} & w \end{array}$$

эквивалентна диаграмме

$$\begin{array}{ccc} v_0 & \xrightarrow{h(a)=t_\star(a)} & v \\ \uparrow{\scriptstyle f^{-1}(b)=t_\star(b^{-1})} & & \downarrow{\scriptstyle f(b)=t_\star(b)} \\ w_0 & \xrightarrow[h(a)]{} & w \end{array}$$

и мы имеем $w = bv = bav_0 = bab^{-1}w_0$. Следовательно

$$h(a)w_0 = (bab^{-1})w_0$$

Мы видим, что представление $h$ зависит от его аргумента. □

Теорема 3.2.15. *Пусть $f$ и $h$ - парные преставления группы $G$. Для любого $a \in G$ отображение $h(a)$ является автоморфизмом представления $f$.*

Доказательство. Следствие теорем 3.2.12 и 3.2.13. □

Замечание 3.2.16. Существует ли морфизм представлений из $t_\star$ в $t_\star$, отличный от автоморфизма $(\mathrm{id}, {}_\star t(a))$? Если мы положим

$$r(g) = cgc^{-1}$$

$$R(a)(m) = cmac^{-1}$$

то нетрудно убедиться, что отображение $(r, R(a))$ является морфизмом представлений из $t_\star$ в $t_\star$. Но это отображение не является автоморфизмом представления $t_\star$, так как $r \neq \mathrm{id}$. □

Глава 4

# Башня представлений универсальных алгебр

## 4.1. Башня представлений универсальных алгебр

Определение 4.1.1. Рассмотрим множество $\Omega_k$-алгебр $A_k$, $k = 1, ..., n$. Положим $\overline{A} = (A_1, ..., A_n)$. Положим $\overline{f} = (f_{1,2}, ..., f_{n-1,n})$. Множество представлений $f_{k,k+1}$, $k = 1, ..., n$, $\Omega_k$-алгебры $A_k$ в $\Omega_{k+1}$-алгебре $A_{k+1}$ называется **башней $(\overline{A}, \overline{f})$ представлений $\overline{\Omega}$-алгебр**. □

Мы можем проиллюстрировать определение 4.1.1 с помощью диаграммы

(4.1.1)

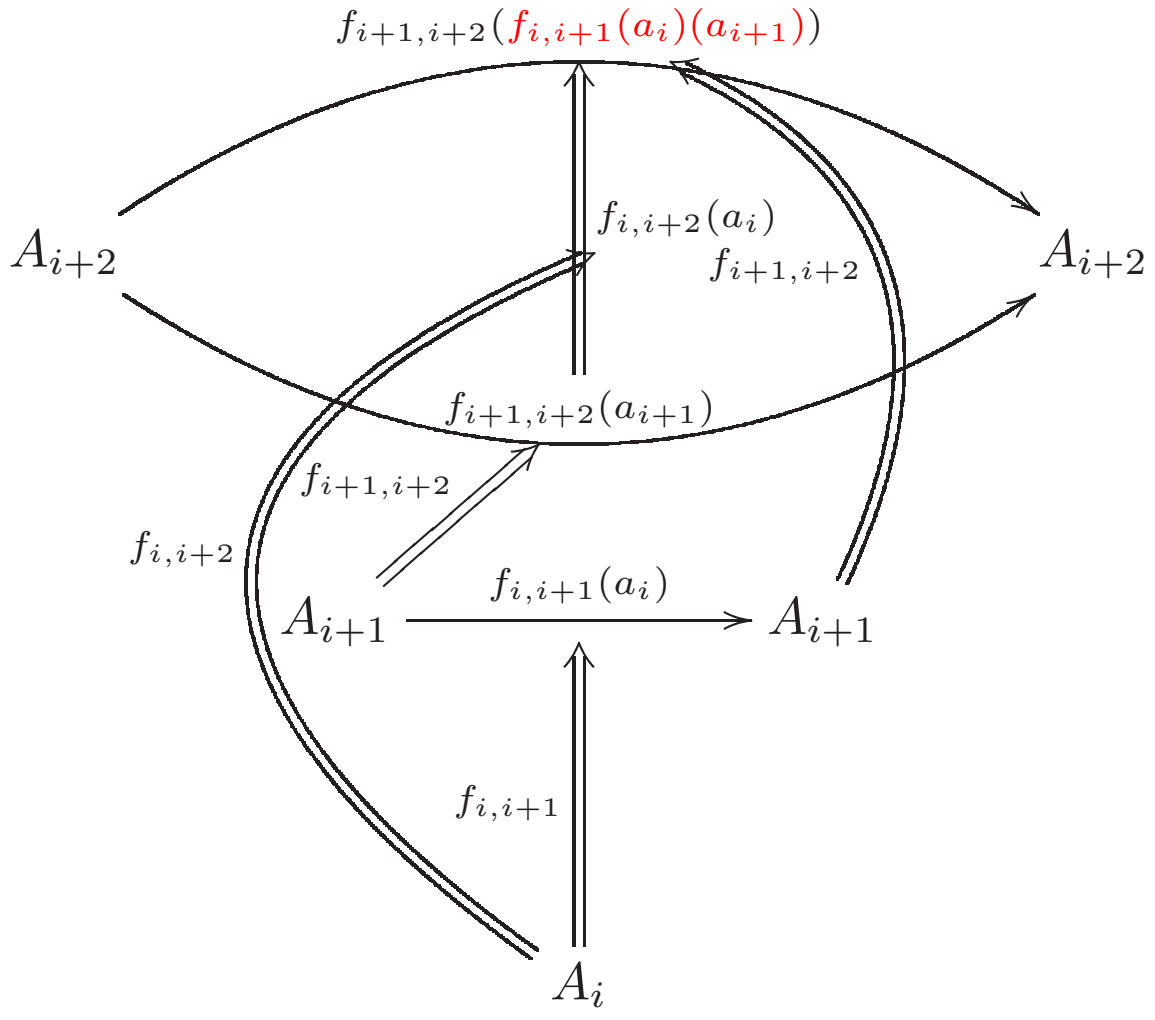

$f_{i,i+1}$ - представление $\Omega_i$-алгебры $A_i$ в $\Omega_{i+1}$-алгебре $A_{i+1}$. $f_{i+1,i+2}$ - представление $\Omega_{i+1}$-алгебры $A_{i+1}$ в $\Omega_{i+2}$-алгебре $A_{i+2}$.

Теорема 4.1.2. *Отображение*

$$f_{i,i+2} : A_i \to {}^{**}A_{i+2}$$

*является представлением* $\Omega_i$*-алгебры* $A_i$ *в* $\Omega_{i+1}$*-алгебре* ${}^*A_{i+2}$.

Доказательство. Произвольному $a_{i+1} \in A_{i+1}$ соответствует автоморфизм $f_{i+1,i+2}(a_{i+1}) \in {}^*A_{i+2}$. Произвольному $a_i \in A_i$ соответствует автоморфизм $f_{i,i+1}(a_i) \in {}^*A_{i+1}$. Так как образом элемента $a_{i+1} \in A_{i+1}$ является элемент $f_{i,i+1}(a_i)(a_{i+1}) \in A_{i+1}$, то тем самым элемент $a_i \in A_i$ порождает преобразование $\Omega_{i+1}$-алгебры ${}^*A_{i+2}$, определённое равенством

$$f_{i,i+2}(a_i)(f_{i+1,i+2}(a_{i+1})) = f_{i+1,i+2}(f_{i,i+1}(a_i)(a_{i+1})) \tag{4.1.2}$$

Пусть $\omega$ - $n$-арная операция $\Omega_i$-алгебры. Так как $f_{i,i+1}$ - гомоморфизм $\Omega_i$-алгебры, то

$$f_{i,i+1}(a_{i,1}...a_{i,n}\omega) = f_{i,i+1}(a_{i,1})...f_{i,i+1}(a_{i,n})\omega \tag{4.1.3}$$

Для $a_1$, ..., $a_n \in A_i$ мы определим операцию $\omega$ на множестве ${}^{**}A_{i+2}$ с помощью равенства

$$\begin{aligned} &f_{i,i+2}(a_{i,1})(f_{i+1,i+2}(a_{i+1}))...f_{i,i+2}(a_{i,n})(f_{i+1,i+2}(a_{i+1}))\omega \\ =&f_{i+1,i+2}(f_{i,i+1}(a_{i,1})(a_{i+1}))...f_{i+1,i+2}(f_{i,i+1}(a_{i,n})(a_{i+1}))\omega \\ =&f_{i+1,i+2}(f_{i,i+1}(a_{i,1})...f_{i,i+1}(a_{i,n})\omega)(a_{i+1}) \end{aligned} \tag{4.1.4}$$

Первое равенство следует из равенства (4.1.2). Второе равенство постулировано на основе требования, что отображение $f_{i+1,i+2}$ является гомоморфизмом $\Omega_i$-алгебры. Следовательно, мы определили структуру $\Omega_i$-алгебры на множестве ${}^{**}A_{i+2}$.

Из (4.1.4) и (4.1.3) следует

$$\begin{aligned} &f_{i,i+2}(a_{i,1})(f_{i+1,i+2}(a_{i+1}))...f_{i,i+2}(a_{i,n})(f_{i+1,i+2}(a_{i+1}))\omega \\ =&f_{i+1,i+2}(f_{i,i+1}(a_{i,1}...a_{i,n}\omega)(a_{i+1})) \\ =&f_{i,i+2}(a_{i,1}...a_{i,n}\omega)(f_{i+1,i+2}(a_{i+1})) \end{aligned} \tag{4.1.5}$$

Из равенства (4.1.5) следует

$$f_{i,i+2}(a_{i,1})...f_{i,i+2}(a_{i,n})\omega = f_{i,i+2}(a_{i,1}...a_{i,n}\omega)$$

Следовательно, отображение $f_{i,i+2}$ является гомоморфизмом $\Omega_i$-алгебры. Следовательно, отображение $f_{i,i+2}$ является представлением $\Omega_i$-алгебры $A_i$ в $\Omega_{i+1}$-алгебре ${}^{*}A_{i+2}$. □

Теорема 4.1.3. *$(id, f_{i+1,i+2})$ является морфизмом представлений $\Omega_i$-алгебры из $f_{i,i+1}$ в $f_{i,i+2}$.*

Доказательство. Рассмотрим более детально диаграмму (4.1.1).

(4.1.6)

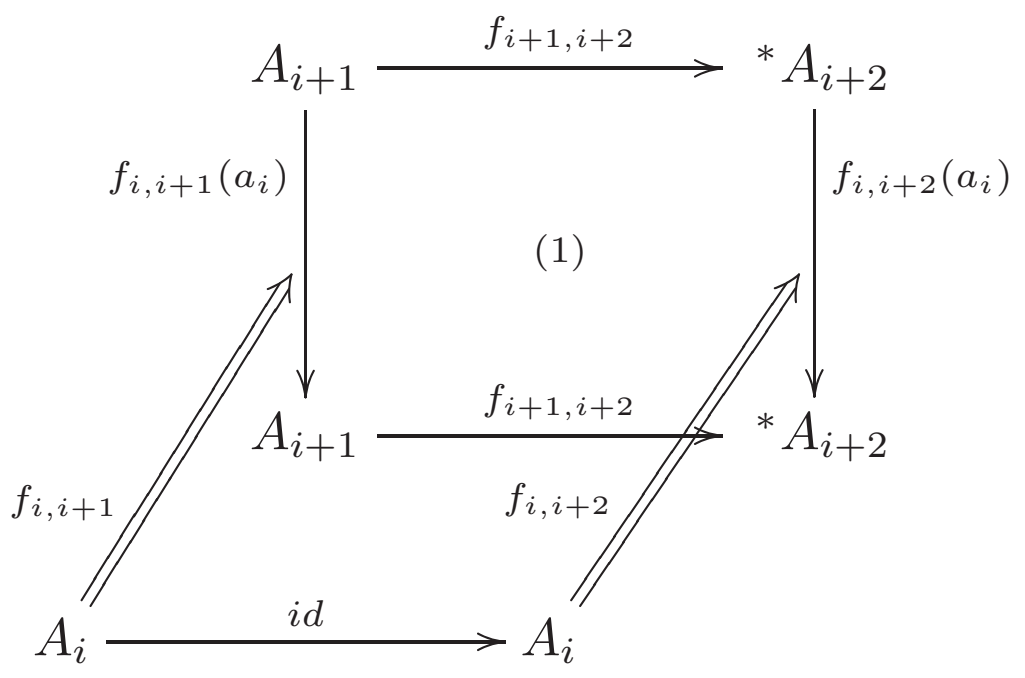

Утверждение теоремы следует из равенства (4.1.2) и определения 2.2.2. □

Теорема 4.1.4. *Рассмотрим башню $(\overline{A}, \overline{f})$ представлений $\overline{\Omega}$-алгебр. Если тождественное преобразование*

$$\delta_{i+2} : A_{i+2} \to A_{i+2}$$

*$\Omega_{i+2}$-алгебры $A_{i+2}$ принадлежит представлению $f_{i+1,i+2}$, то представление $f_{i,i+2}$ $\Omega_i$-алгебры $A_i$ в $\Omega_{i+1}$-алгебре ${}^*A_{i+2}$ может быть продолжено до представления*

$$f'_{i,i+2}:A_i\to {}^*A_{i+2}$$

*$\Omega_i$-алгебры $A_i$ в $\Omega_{i+2}$-алгебре $A_{i+2}$.*

Доказательство. Согласно равенству (4.1.2), для заданного $a_{i+1}\in A_{i+1}$ мы можем определить отображение

(4.1.7)
$$\begin{gathered}f'_{i,i+2}:A_i\to {}^*A_{i+2}\\ f'_{i,i+2}(a_i)(f_{i+1,i+2}(a_{i+1})(a_{i+2}))=f_{i+1,i+2}(f_{i,i+1}(a_i)(a_{i+1}))(a_{i+2})\end{gathered}$$

Пусть существует $a_{i+1}\in A_{i+1}$ такой, что

(4.1.8)
$$f_{i+1,i+2}(a_{i+1})=\delta_{i+2}$$

Тогда из (4.1.7),(4.1.8) следует

(4.1.9)
$$f'_{i,i+2}(a_i)(a_{i+2})=f_{i+1,i+2}(f_{i,i+1}(a_i)(a_{i+1}))(a_{i+2})$$

Пусть $\omega$ - $n$-арная операция $\Omega_i$-алгебры $A_i$. Поскольку $f_{i,i+2}$ - гомоморфизм $\Omega_i$ -алгебры, то из равенства (4.1.2) следует

(4.1.10)
$$\begin{aligned}&f_{i,i+2}(a_{i\cdot 1}...a_{i\cdot n}\omega)(f_{i+1,i+2}(a_{i+1}))\\ =&f_{i,i+2}(a_{i\cdot 1})(f_{i+1,i+2}(a_{i+1}))...f_{i,i+2}(a_{i\cdot n})(f_{i+1,i+2}(a_{i+1}))\omega\\ =&f_{i+1,i+2}(f_{i,i+1}(a_{i\cdot 1})(a_{i+1})...f_{i,i+1}(a_{i\cdot n})(a_{i+1})\omega)\end{aligned}$$

Из равенства (4.1.10) следует
(4.1.11)
$$\begin{aligned}&f_{i,i+2}(a_{i\cdot 1}...a_{i\cdot n}\omega)(f_{i+1,i+2}(a_{i+1})(a_{i+2}))\\ =&f_{i,i+2}(a_{i\cdot 1})(f_{i+1,i+2}(a_{i+1})(a_{i+2}))...f_{i,i+2}(a_{i\cdot n})(f_{i+1,i+2}(a_{i+1})(a_{i+2}))\omega\\ =&f_{i+1,i+2}(f_{i,i+1}(a_{i\cdot 1})(a_{i+1})...f_{i,i+1}(a_{i\cdot n})(a_{i+1})\omega)(a_{i+2})\end{aligned}$$

Из равенств (4.1.8), (4.1.9), (4.1.11) следует

(4.1.12)
$$\begin{aligned}&f'_{i,i+2}(a_{i\cdot 1}...a_{i\cdot n}\omega)(a_{i+2})\\ =&f'_{i,i+2}(a_{i\cdot 1})(a_{i+2})...f'_{i,i+2}(a_{i\cdot n})(a_{i+2})\omega\\ =&f_{i+1,i+2}(f_{i,i+1}(a_{i\cdot 1})(a_{i+1})...f_{i,i+1}(a_{i\cdot n})(a_{i+1})\omega)(a_{i+2})\end{aligned}$$

Равенство (4.1.12) определяет операцию $\omega$ на множестве ${}^*A_{i+2}$

(4.1.13)
$$f'_{i,i+2}(a_{i\cdot 1})...f'_{i,i+2}(a_{i\cdot n})\omega=f'_{i,i+2}(a_{i\cdot 1}...a_{i\cdot n}\omega)$$

При этом отображение $f'_{i,i+2}$ является гомоморфизмом $\Omega_i$-алгебры. Следовательно, мы построили представление $\Omega_i$-алгебры $A_i$ в $\Omega_{i+2}$-алгебре $A_{i+2}$. □

Определение 4.1.5. Рассмотрим башню представлений

$$((A_1,A_2,A_3),(f_{1,2},f_{2,3}))$$

Отображение $f_*=(f_{1,2},f_{1,3})$ называется **представлением $\Omega_1$-алгебры $A_1$ в представлении** $f_{2,3}$. □

Определение 4.1.6. Рассмотрим башню представлений $(\overline{A}, \overline{f})$. Отображение

$$f_* = (f_{1,2}, ..., f_{1,n})$$

называется **представлением $\Omega_1$-алгебры $A_1$ в башне представлений**

$$(\overline{A}_{[1]}, \overline{f}) = ((A_2, ..., A_n), (f_{2,3}, ..., f_{n-1,n}))$$

□

## 4.2. Морфизм башни $T*$-представлений

Определение 4.2.1. Рассмотрим множество $\Omega_k$-алгебр $A_k$, $B_k$, $k = 1$, ..., $n$. Множество отображений $(h_1, ..., h_n)$ называется **морфизмом из башни представлений** $(\overline{A}, \overline{f})$ **в башню представлений** $(\overline{B}, \overline{g})$, если для каждого $k$, $k = 1$, ..., $n-1$, пара отображений $(h_k, h_{k+1})$ является морфизмом представлений из $f_{k,k+1}$ в $g_{k,k+1}$. □

Для любого $k$, $k = 1$, ..., $n-1$, мы имеем диаграмму

(4.2.1)

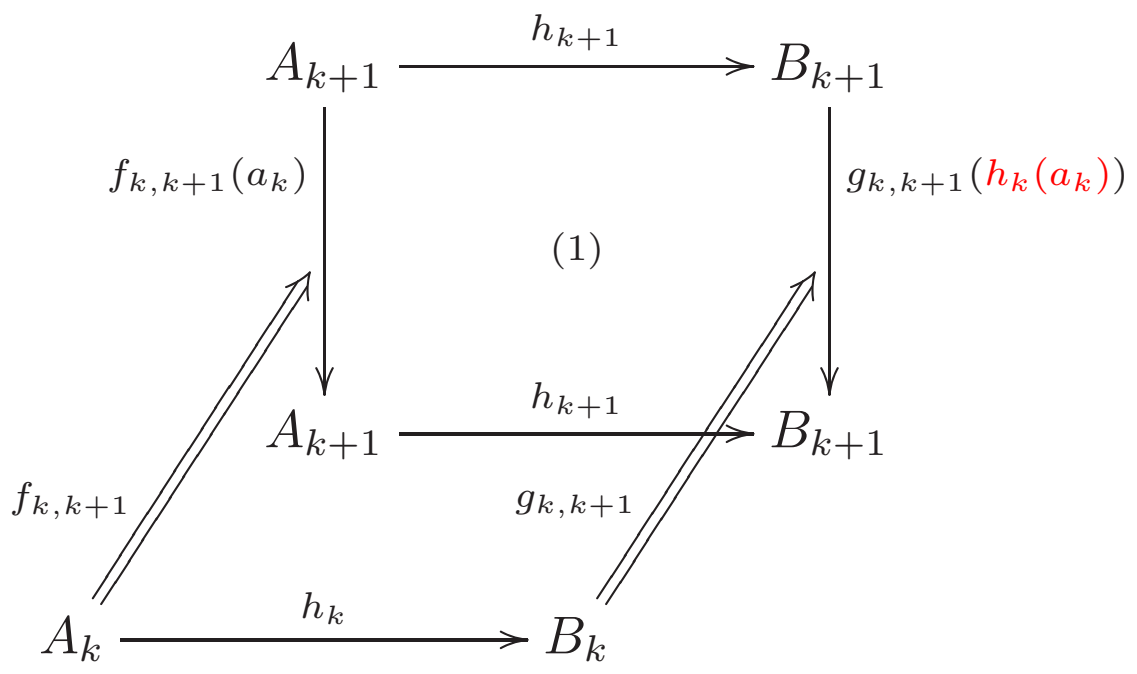

Равенства

$$h_{k+1} \circ f_{k,k+1}(a_k) = g_{k,k+1}(h_k(a_k)) \circ h_{k+1} \tag{4.2.2}$$

$$h_{k+1}(f_{k,k+1}(a_k)(a_{k+1})) = g_{k,k+1}(h_k(a_k))(h_{k+1}(a_{k+1})) \tag{4.2.3}$$

выражают коммутативность диаграммы (1). Однако уже для морфизма $(h_k, h_{k+1})$, $k > 1$, диаграмма (4.2.1) неполна. Учитывая аналогичную диаграмму для морфизма $(h_k, h_{k+1})$ эта диаграмма на верхнем уровне приобретает вид
(4.2.4)

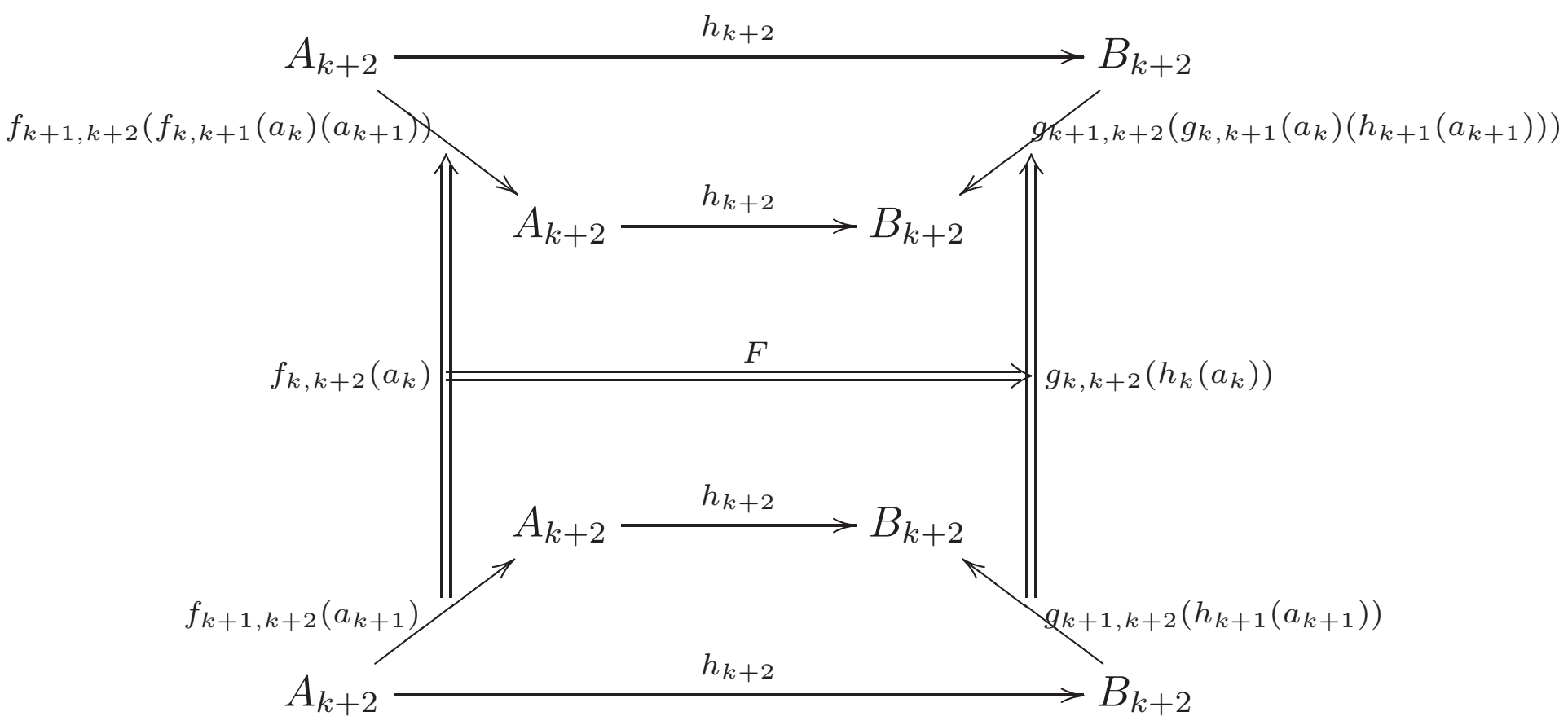

К сожалению, диаграмма (4.2.4) малоинформативна. Очевидно, что существует морфизм из ${}^*A_{k+2}$ в ${}^*B_{k+2}$, отображающий $f_{k,k+2}(a_k)$ в $g_{k,k+2}(h_k(a_k))$. Однако структура этого морфизма из диаграммы неясна. Мы должны рассмотреть отображение из ${}^*A_{k+2}$ в ${}^*B_{k+2}$, так же мы это сделали в теореме 4.1.3.

Теорема 4.2.2. *Если представления $f_{i+1,i+2}$ эффективно, то $(h_i, {}^*h_{i+2})$ является морфизмом представлений из представления $f_{i,i+2}$ в представление $g_{i,i+2}$ $\Omega_i$-алгебры.*

Доказательство. Рассмотрим диаграмму

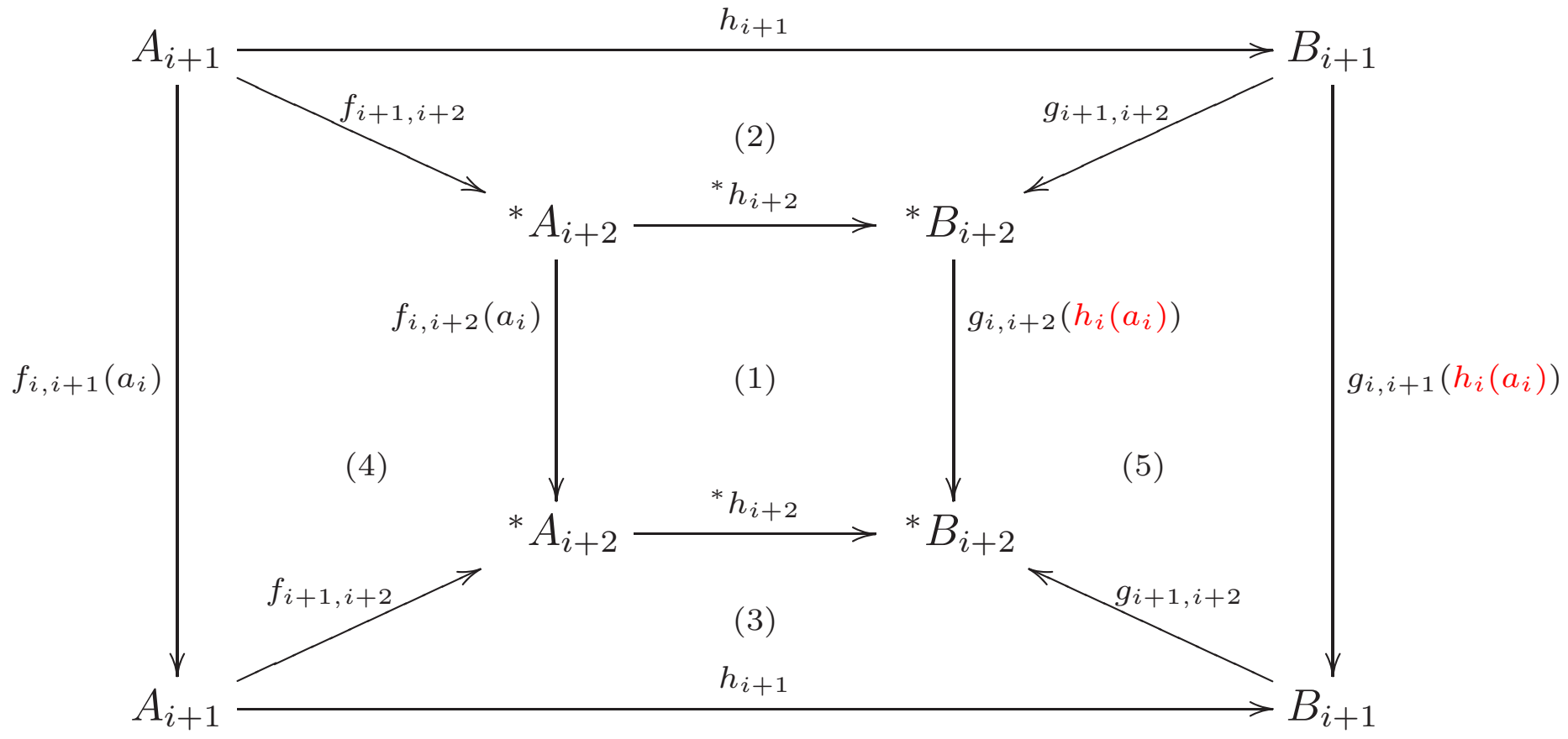

Существование отображения ${}^*h_{i+2}$ и коммутативность диаграмм (2) и (3) следует из эффективности отображения $f_{i+1,i+2}$ и теоремы 2.2.8. Коммутативность диаграмм (4) и (5) следует из теоремы 4.1.3.

Из коммутативности диаграммы (4) следует

$$f_{i+1,i+2} \circ f_{i,i+1}(a_i) = f_{i,i+2}(a_i) \circ f_{i+1,i+2} \tag{4.2.5}$$

Из равенства (4.2.5) следует

$$(4.2.6)\qquad {}^*h_{i+2}\circ f_{i+1,i+2}\circ f_{i,i+1}(a_i)={}^*h_{i+2}\circ f_{i,i+2}(a_i)\circ f_{i+1,i+2}$$

Из коммутативности диаграммы (3) следует

$$(4.2.7)\qquad {}^*h_{i+2}\circ f_{i+1,i+2}=g_{i+1,i+2}\circ h_{i+1}$$

Из равенства (4.2.7) следует

$$(4.2.8)\qquad {}^*h_{i+2}\circ f_{i+1,i+2}\circ f_{i,i+1}(a_i)=g_{i+1,i+2}\circ h_{i+1}\circ f_{i,i+1}(a_i)$$

Из равенств (4.2.6) и (4.2.8) следует

$$(4.2.9)\qquad {}^*h_{i+2}\circ f_{i,i+2}(a_i)\circ f_{i+1,i+2}=g_{i+1,i+2}\circ h_{i+1}\circ f_{i,i+1}(a_i)$$

Из коммутативности диаграммы (5) следует

$$(4.2.10)\qquad g_{i+1,i+2}\circ g_{i,i+1}(h_i(a_i))=g_{i,i+2}(h_i(a_i))\circ g_{i+1,i+2}$$

Из равенства (4.2.10) следует

$$(4.2.11)\qquad g_{i+1,i+2}\circ g_{i,i+1}(h_i(a_i))\circ h_{i+1}=g_{i,i+2}(h_i(a_i))\circ g_{i+1,i+2}\circ h_{i+1}$$

Из коммутативности диаграммы (2) следует

$$(4.2.12)\qquad {}^*h_{i+2}\circ f_{i+1,i+2}=g_{i+1,i+2}\circ h_{i+1}$$

Из равенства (4.2.12) следует

$$(4.2.13)\qquad g_{i,i+2}(h_i(a_i))\circ {}^*h_{i+2}\circ f_{i+1,i+2}=g_{i,i+2}(h_i(a_i))\circ g_{i+1,i+2}\circ h_{i+1}$$

Из равенств (4.2.11) и (4.2.13) следует

$$(4.2.14)\qquad g_{i+1,i+2}\circ g_{i,i+1}(h_i(a_i))\circ h_{i+1}=g_{i,i+2}(h_i(a_i))\circ {}^*h_{i+2}\circ f_{i+1,i+2}$$

Внешняя диаграмма является диаграммой (4.2.1) при $i=1$. Следовательно, внешняя диаграмма коммутативна

$$(4.2.15)\qquad h_{i+1}\circ f_{i,i+1}(a_i)=g_{i,i+1}(h_i(a_i))\circ h_{i+1}$$

Из равенства (4.2.15) следует

$$(4.2.16)\qquad g_{i+1,i+2}\circ h_{i+1}\circ f_{i,i+1}(a_i)=g_{i+1,i+2}\circ g_{i,i+1}(h_i(a_i))\circ h_{i+1}(a_{i+1})$$

Из равенств (4.2.9), (4.2.14) и (4.2.16) следует

$$(4.2.17)\qquad {}^*h_{i+2}\circ f_{i,i+2}(a_i)\circ f_{i+1,i+2}=g_{i,i+2}(h_i(a_i))\circ {}^*h_{i+2}\circ f_{i+1,i+2}$$

Так как отображение $f_{i+1,i+2}$ - инъекция, то из равенства (4.2.17) следует

$$(4.2.18)\qquad {}^*h_{i+2}\circ f_{i,i+2}(a_i)=g_{i,i+2}(h_i(a_i))\circ {}^*h_{i+2}$$

Из равенства (4.2.18) следует коммутативность диаграммы (1), откуда следует утверждение теоремы. □

Теоремы 4.1.3 и 4.2.2 справедливы для любых уровней башни представлений. В каждом конкретном случае надо правильно указывать множества, откуда и куда направленно отображение. Смысл приведенных теорем состоит в том, что все отображения в башне представлений действуют согласовано.

Теорема 4.2.2 утверждает, что неизвестное отображение на диаграмме (4.2.4) является отображением ${}^*h_{i+2}$.

Теорема 4.2.3. *Рассмотрим множество $\Omega_i$-алгебр $A_i$, $B_i$, $C_i$, $i=1, ..., n$. Пусть определены морфизмы башни представлений*

$$\overline{p} :(\overline{A},\overline{f}) \to (\overline{B},\overline{g})$$
$$\overline{q} :(\overline{B},\overline{g}) \to (\overline{C},\overline{h})$$

*Тогда определён морфизм представлений $\Omega$-алгебры*

$$\overline{r} : (\overline{A},\overline{f}) \to (\overline{C},\overline{h})$$

*где $r_k = q_k p_k$, $k=1, ..., n$. Мы будем называть морфизм $\overline{r}$ башни представлений из $\overline{f}$ в $\overline{h}$* **произведением морфизмов $\overline{p}$ и $\overline{q}$ башни представлений**.

Доказательство. Для каждого $k$, $k=2, ..., n$, мы можем представить утверждение теоремы, пользуясь диаграммой

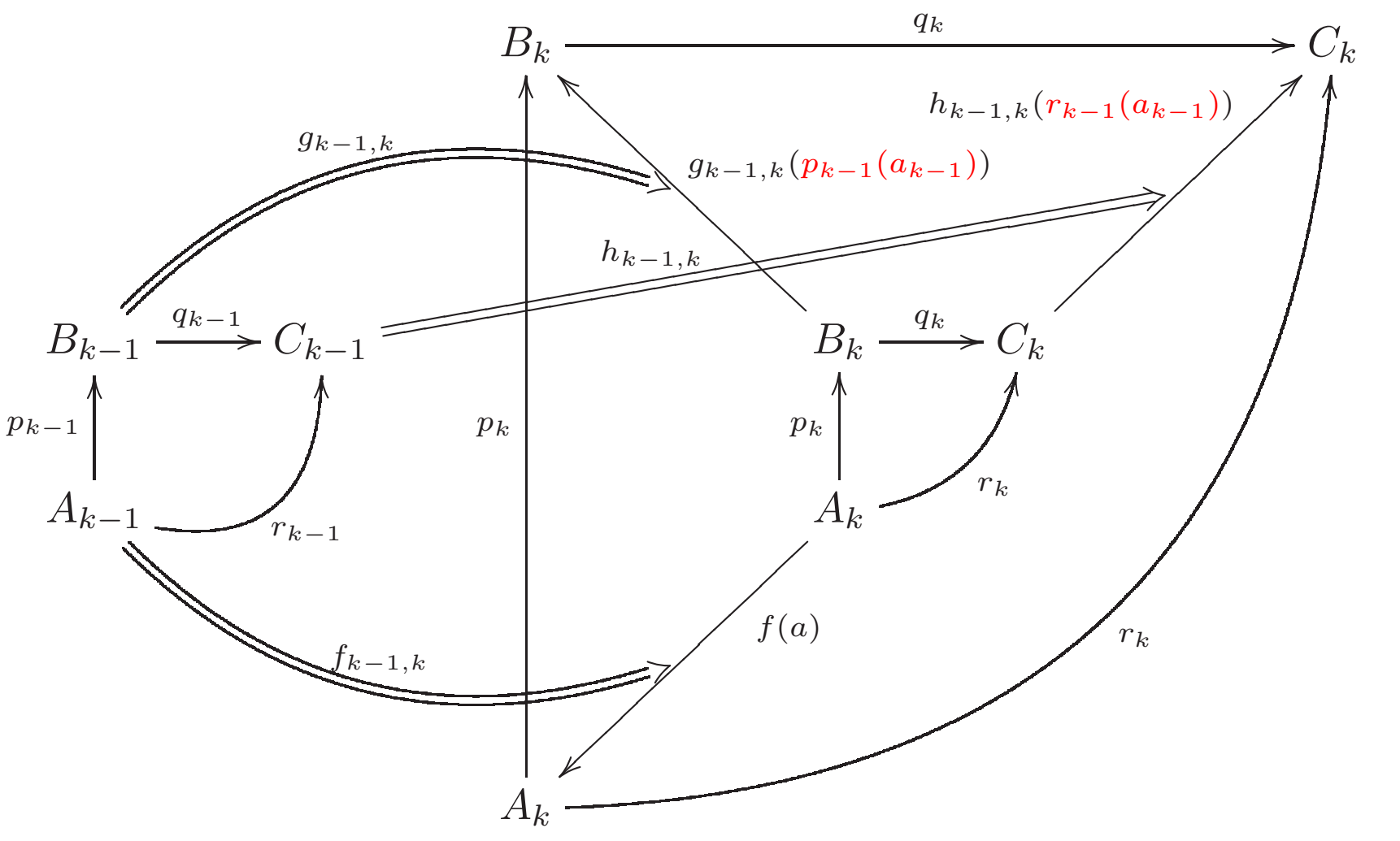

Отображение $r_{k-1}$ является гомоморфизмом $\Omega_{k-1}$-алгебры $A_{k-1}$ в $\Omega_{k-1}$-алгебру $C_{k-1}$. Нам надо показать, что пара отображений $(r_{k-1}, r_k)$ удовлетворяет (4.2.2):

$$\begin{aligned}
r_k(f_{k-1,k}(a_{k-1})a_k) &= q_k p_k(f_{k-1,k}(a_{k-1})a_k)\\
&= q_k(g_{k-1,k}(p_{k-1}(a_{k-1}))p_k(a_k))\\
&= h_{k-1,k}(q_{k-1}p_{k-1}(a_{k-1}))q_k p_k(a_k))\\
&= h_{k-1,k}(r(a_{k-1}))r_k(a_k)
\end{aligned}$$

□

## 4.3. Автоморфизм башни представлений

Определение 4.3.1. Пусть $(\overline{A},\overline{f})$ - башня представлений $\overline{\Omega}$-алгебр. Морфизм башни представлений $(\mathrm{id}, h_2, ..., h_n)$ такой, что для любого $k$, $k=2, ..., n$, $h_k$ - эндоморфизм $\Omega_k$-алгебры $A_k$, называется **эндоморфизмом башни представлений** $\overline{f}$. □

Определение 4.3.2. Пусть $(\overline{A},\overline{f})$ - башня представлений $\overline{\Omega}$-алгебр. Морфизм башни представлений $(\mathrm{id}, h_2, ..., h_n)$ такой, что для любого $k$, $k=2, ..., n$,

$n$, $h_k$ - автоморфизм $\Omega_k$-алгебры $A_k$, называется **автоморфизмом башни представлений** $f$. □

Теорема 4.3.3. *Пусть* $(\overline{A},\overline{f})$ *- башня представлений* $\overline{\Omega}$*-алгебр. Множество автоморфизмов представления* $\overline{f}$ *порождает группу.*

Доказательство. Пусть $\overline{r}_{[1]}$, $\overline{p}_{[1]}$ - автоморфизмы башни представлений $\overline{f}$. Согласно определению 4.3.2 для любого $k$, $k=2$, ..., $n$, отображения $r_k$, $p_k$ являются автоморфизмами $\Omega_k$-алгебры $A_k$. Согласно теореме II.3.2 ([13], с. 60) для любого $k$, $k=2$, ..., $n$, отображение $r_kp_k$ является автоморфизмом $\Omega_k$-алгебры $A_k$. Из теоремы 4.2.3 и определения 4.3.2 следует, что произведение автоморфизмов $\overline{r}_{[1]}\overline{p}_{[1]}$ башни представлений $\overline{f}$ является автоморфизмом башни представлений $\overline{f}$.

Согласно доказательству теоремы 2.4.5, для любого $k$, $k=2$, ..., $n$, произведение автоморфизмов $\Omega_i$-алгебры ассоциативно. Следовательно, ассоциативно произведение автоморфизмов башни представлений.

Пусть $\overline{r}_{[1]}$ - автоморфизм башни представлений $\overline{f}$. Согласно определению 4.3.2 для любого $k$, $k=2$, ..., $n$, отображение $r_k$ является автоморфизмом $\Omega_k$-алгебры $A_k$. Следовательно, для любого $k$, $k=2$, ..., $n$, отображение $r_k^{-1}$ является автоморфизмом $\Omega_k$-алгебры $A_k$. Для автоморфизма $\overline{r}_{[1]}$ справедливо равенство (4.2.3). Положим $a'_k=r_k(a_k)$, $k=2$, ..., $n$. Так как $r_k$, $k=2$, ..., $n$, - автоморфизм, то $a_k=r_k^{-1}(a'_k)$ и равенство (4.2.3) можно записать в виде

$$h_{k+1}(f_{k,k+1}(h_k^{-1}(a'_k))(h_{k+1}(a'_{k+1})))=g_{k,k+1}(a'_k)(a'_{k+1}) \tag{4.3.1}$$

Так как отображение $h_{k+1}$ является автоморфизмом $\Omega_{k+1}$-алгебры $A_{k+1}$, то из равенства (4.3.1) следует

$$f_{k,k+1}(h_k^{-1}(a'_k)(h_{k+1}(a'_{k+1})))=h_{k+1}^{-1}(g_{k,k+1}(a'_k)(a'_{k+1})) \tag{4.3.2}$$

Равенство (4.3.2) соответствует равенству (4.2.3) для отображения $\overline{r}_{[1]}^{-1}$. Следовательно, отображение $\overline{r}_{[1]}^{-1}$ является автоморфизмом представления $\overline{f}$. □

## 4.4. Примеры башни представлений

**4.4.1. Аффинное пространство.** Пусть $\overset{\circ}{A}$ - множество точек. Пусть $\overline{A}$ - векторное пространство над полем $k$. Однотранзитивное представление векторного пространства $\overline{A}$ в множестве точек $\overset{\circ}{A}$ называется аффинным пространством. Если $A\in\overset{\circ}{A}$ и $\overline{v}\in\overline{A}$, то мы будем обозначать действие вектора $\overline{v}$ на точку $A$ выражением $A+\overline{v}$.

Однотранзитивность представления означает, что для любых точек $A,B\in\overset{\circ}{A}$ существует единственный вектор $\overline{v}\in\overline{A}$ такой, что

$$B=A+\overline{v} \tag{4.4.1}$$

Если вектор $\overline{v}$ удовлетворяет равенству (4.4.1), то мы будем также пользоваться обозначением $\overrightarrow{AB}$ для вектора $\overline{v}$

$$B=A+\overrightarrow{AB} \tag{4.4.2}$$

Аффинное пространство является башней представлений, так как векторное пространство $\overline{A}$ является представлением поля $k$ в абелевой группе $\overline{A}$.

**4.4.2. Свободная алгебра над кольцом.** Согласно определениям [14], с. 1, [12], с. 4, свободная алгебра $A$ над коммутативным кольцом $R$ - это свободный $R$-модуль, в котором определено произведение, являющееся билинейным отображением $R$-модуля $A$. Это определение неудовлетворительно, так как эндоморфизмы, порождённые представлением, не являются эндоморфизмами алгебры. Однако мы можем изменить определение $R$-алгебры таким образом, чтобы все требования к представлению были удовлетворены.

Согласно определению [10]-3.1.10, отображение

$$l(a) : b \in A \to ab \in A \tag{4.4.3}$$

является линейным отображением. Следовательно, отображение (4.4.3) является эндоморфизмом $R$-модуля $A$. Мы будем называть эндоморфизм (4.4.3) **левым сдвигом $R$-модуля $A$**.

Теорема 4.4.1. *Представление*

$$f_{2,3} : a \in A \to l(a) \in {}^\star A \tag{4.4.4}$$

*$R$-модуля $A$ в $R$-модуле $A$ эквивалентно структуре $R$-алгебры $A$.*

Доказательство. Пусть структура $R$-модуля $A$ порождена эффективным представлением

$$f_{1,2} : R \to {}^\star A$$

кольца $R$ в абелевой группе $A$.

Отображение $f_{2,3}(a)$ является левым сдвигом $R$-модуля $A$. Поскольку левый сдвиг является эндоморфизмом $R$-модуля $A$, то

$$\begin{aligned} f_{2,3}(a)(b_1 + b_2) &= f_{2,3}(a)(b_1) + f_{2,3}(a)(b_2) \\ f_{2,3}(a)(rb) &= r f_{2,3}(a)(b) \end{aligned} \tag{4.4.5}$$

Поскольку отображение (4.4.4) - гомоморфизм абелевых групп, то

$$f_{2,3}(a_1 + a_2)(b) = (f_{2,3}(a_1) + f_{2,3}(a_2))(b) = f_{2,3}(a_1)(b) + f_{2,3}(a_2)(b) \tag{4.4.6}$$

Согласно теореме 4.1.2 существует представление

$$f_{1,3} : R \to {}^{\star\star} A$$

которое согласно равенству (4.1.2) определено правилом

$$f_{1,3}(r)(f_{2,3}(a)) = f_{2,3}(f_{1,2}(r)(a)) \tag{4.4.7}$$

Мы не предполагаем представление $f_{2,3}$ эффективным.[4.1] Поэтому представление $f_{1,3}$ так же может быть не эффективным. Однако мы можем, также как и в случае отображения $f_{1,2}$, отождествлять $f_{1,3}(r)$ и $r$. Из равенства (4.4.7) следует

$$r f_{2,3}(a)(b) = f_{2,3}(ra)(b) \tag{4.4.8}$$

Из равенств (4.4.5), (4.4.6), (4.4.8) и определения [10]-3.1.10, следует, что отображение $f_{2,3}$ является билинейным отображением. Следовательно, отображение $f_{2,3}$ определяет произведение в $R$-модуле $A$ согласно правилу

$$ab = f_{2,3}(a)(b)$$

□

[4.1] Если представление $f_{2,3}$ не эффективно, то $R$-алгебра $A$ имеет делители нуля.

Следствие 4.4.2. *Башня представлений*

$$((f_{1,2}, f_{2,3}), (R, A, A)) \tag{4.4.9}$$

*порождает структуру $R$-алгебры $A$.* □

Обратим внимание, что при определении $R$-алгебры $A$ область определения отображения $f_{2,3}$ - не всё множество эндоморфизмов абелевой группы $A$, а множество эндоморфизмов $R$-модуля $A$. Мы предполагаем, что в определении 4.1.1 каждый элемент башни представлений является башней представлений. В этом случае, **диаграмма представлений** для $R$-алгебры имеет вид

$$\begin{array}{ccccc} R & \xrightarrow[*]{f_{1,2}} & A & \xrightarrow[*]{f_{2,3}} & A \\ & & & & \uparrow {\scriptstyle *\, g_{1,2}} \\ & & & & R \end{array} \tag{4.4.10}$$

На диаграмме представлений (4.4.10) предполагается, что мы сперва строим представление $g_{1,2}$, порождающее структуру $R$-модуля $A$. Представление $f_{1,2}$ порождает структуру $R$-модуля $A$. Согласно определению, представление $f_{2,3}$ является гомоморфизмом $R$-модуля $A$ в множество эндоморфизмов $R$-модуля $A$.

**4.4.3. Групповая алгебра.** Если $R$ - кольцо и $G$ - мультипликативная группа, то свободный $R$-модуль $R[G]$ с базисом $G$ называется групповым кольцом ([5], с. 383). Из теоремы [7]-13.1.1 следует, что мы можем продолжить умножение с группы $G$ на модуль $R[G]$. Следовательно, модуль $R[G]$ является свободной $R$-алгеброй.

Если $R$ - коммутативное кольцо, то модуль $R[G]$ называется **групповой алгеброй** группы $G$ ([6], с. 74). Поскольку группа $G$ является базисом групповой алгебры, мы можем записать элемент алгебры $R[G]$ в виде

$$\textstyle\sum_{g\in G} a_g g \quad a_g \in R \tag{4.4.11}$$

Произведение прроизвольных элементов $R[G]$ имеет вид

$$\left(\sum_{g\in G} a_g g\right)\left(\sum_{g\in G} b_g g\right) = \sum_{g\in G} c_g g \qquad c_g = \sum_{u\in G} a_u b_{u^{-1}g} \tag{4.4.12}$$

Координаты элемента $a \in R[g]$ можно интерпретировать как отображение

$$a : G \to R \quad a(g) = a_g$$

**4.4.4. Модуль над алгеброй.** Пусть $A_2$ - свободная алгебра над полем $A_1$. Рассматривая алгебру $A_2$ как кольцо, мы можем определить свободный модуль $A_3$ над алгеброй $A_2$. Мы можем описать модуль $A_3$ с помощью диаграммы представлений

$$\begin{array}{ccccc} A_2 & \xrightarrow[*]{f_{2,3}} & A_2 & \xrightarrow[*]{f_{3,4}} & A_3 \\ {\scriptstyle f_{1,2}\,*}\uparrow & & \uparrow{\scriptstyle *\, f_{1,2}} & & \\ A_1 & & A_1 & & \end{array}$$

## 4.5. Множество образующих башни представлений

Определение 4.5.1. Башня представлений $(\overline{A},\overline{f})$. называется **башней эффективных представлений**, если для любого $i$ представление $f_{i,i+1}$ эффективно. □

Теорема 4.5.2. *Рассмотрим башню представлений $(\overline{A},\overline{f})$. Пусть представления $f_{i,i+1}$, ..., $f_{i+k-1,i+k}$ эффективны. Тогда представление $f_{i,i+k}$ эффективно.*

Доказательство. Мы докажем утверждение теоремы по индукции.

Пусть представления $f_{i,i+1}$, $f_{i+1,i+2}$ эффективны. Предположим, что преобразование $f_{i,i+1}(a_i)$ не является тождественным преобразованием. Тогда существует $a_{i+1}\in A_{i+1}$ такой, что $f_{i,i+1}(a_i)(a_{i+1})\ne a_{i+1}$. Так как представление $f_{i+1,i+2}$ эффективно, то преобразования $f_{i+1,i+2}(a_{i+1})$ и $f_{i+1,i+2}(f_{i,i+1}(a_i)(a_{i+1}))$ не совпадают. Согласно построению, выполненному в теореме 4.1.2, преобразование $f_{i,i+2}(a_i)$ не является тождественным преобразованием. Следовательно, представление $f_{i,i+2}$ эффективно.

Предположим утверждение теоремы верно для $k-1$ представлений и пусть $f_{i,i+1}$, ..., $f_{i+k-1,i+k}$ эффективны. Согласно предположению индукции, представления $f_{i,i+k-1}$, $f_{i+k-1,i+k}$ эффективны. Согласно доказанному выше, представление $f_{i,i+k}$ эффективно. □

Теорема 4.5.3. *Рассмотрим башню $(\overline{A},\overline{f})$ представлений $\overline{\Omega}$-алгебр. Пусть тождественное преобразование*

$$\delta_{i+2}:A_{i+2}\to A_{i+2}$$

*$\Omega_{i+2}$-алгебры $A_{i+2}$ принадлежит представлению $f_{i+1,i+2}$. Пусть представления $f_{i,i+1}$, $f_{i+1,i+2}$ эффективны. Тогда представление $f'_{i,i+2}$, определённое в теореме 4.1.4, эффективно.*

Доказательство. Пусть существует $a_{i+1}\in A_{i+1}$ такой, что

$$f_{i+1,i+2}(a_{i+1})=\delta_{i+2}$$

Допустим, что представление $f'_{i,i+2}$ не является эффективным. Тогда существуют

$$a_{i\cdot 1},a_{i\cdot 2}\in A_i\quad a_{i\cdot 1}\ne a_{i\cdot 2}\tag{4.5.1}$$

такие, что

$$f'_{i,i+2}(a_{i\cdot 1})(a_{i+2})=f'_{i,i+2}(a_{i\cdot 2})(a_{i+2})\tag{4.5.2}$$

Из равенств (4.1.9), (4.5.2), следует

$$f_{i+1,i+2}(f_{i,i+1}(a_{i\cdot 1})(a_{i+1}))(a_{i+2})=f_{i+1,i+2}(f_{i,i+1}(a_{i\cdot 2})(a_{i+1}))(a_{i+2})\tag{4.5.3}$$

Так как $a_{i+2}$ произвольно, то из (4.5.3) следует

$$f_{i+1,i+2}(f_{i,i+1}(a_{i\cdot 1})(a_{i+1}))=f_{i+1,i+2}(f_{i,i+1}(a_{i\cdot 2})(a_{i+1}))\tag{4.5.4}$$

Поскольку представление $f_{i,i+1}$ эффективно, то из условия (4.5.1) следует

$$f_{i,i+1}(a_{i\cdot 1})(a_{i+1})\ne f_{i,i+1}(a_{i\cdot 2})(a_{i+1})\tag{4.5.5}$$

Из условия (4.5.5) и равенства (4.5.4) следует, что представление $f_{i+1,i+2}$ не является эффективным. Полученное противоречие доказывает утверждение теоремы. □

Мы строим базис представления башни представлений по той же схеме, что мы построили базис представления в секции 2.7.

Мы будем записывать элементы башни представлений $(\overline{A},\overline{f})$ в виде кортежа

$$\overline{a}=(a_1,...,a_3)\quad a_i\in A_i\quad i=1,...,n$$

Определение 4.5.4. Пусть $(\overline{A},\overline{f})$ - башня представлений. Кортеж множеств

$$\overline{N}_{[1]}=(N_2\subset A_2,...,N_n\subset A_n)$$

называется **кортежем стабильных множеств башни представлений** $\overline{f}$, если

$$f_{i-1,i}(a_{i-1})(a_i)\in N_i\quad i=2,...,n$$

для любых $a_1\in A_1$, $a_2\in N_2$, ..., $a_n\in N_n$. Мы также будем говорить, что кортеж множеств

$$\overline{N}_{[1]}=(N_2\subset A_2,...,N_n\subset A_n)$$

стабилен относительно башни представлений $\overline{f}$. □

Теорема 4.5.5. *Пусть* $\overline{f}$ *- башня представлений. Пусть множество* $N_i\subset A_i$ *является подалгеброй* $\Omega_i$*-алгебры* $A_i$, $i=2$, ..., $n$. *Пусть кортеж множеств*

$$\overline{N}_{[1]}=(N_2\subset A_2,...,N_n\subset A_n)$$

*стабилен относительно башни представлений* $\overline{f}$. *Тогда существует башня представлений*

$$((A_1,N_2,...,N_n),(f_{N_2,1,2},...,f_{N_n,n-1,n}))\tag{4.5.6}$$

*такая, что*

$$f_{N_i,i-1,i}(a_{i-1})=f_{i-1,i}(a_{i-1})|_{N_i}\quad i=2,...,n$$

*Башня представлений* (4.5.6) *называется* **башней подпредставлений**.

Доказательство. Пусть $\omega_{i-1,1}$ - $m$-арная операция $\Omega_{i-1}$-алгебры $A_{i-1}$, $i=2$, ..., $n$. Тогда для любых $a_{i-1,1}$, ..., $a_{i-1,m}\in N_{i-1}$ [4.2] и любого $a_i\in N_i$

$$\begin{aligned}&(f_{N_i,i-1,i}(a_{i-1,1})...f_{N_i,i-1,i}(a_{i-1,m})\omega_{i-1,1})(a_i)\\=&(f_{i-1,i}(a_{i-1,1})...f_{i-1,i}(a_{i-1,m})\omega_{i-1,1})(a_i)\\=&f_{i-1,i}(a_{i-1,1}...a_{i-1,m}\omega_{i-1,1})(a_i)\\=&f_{N_i,i-1,i}(a_{i-1,1}...a_{i-1,m}\omega_{i-1,1})(a_i)\end{aligned}$$

Пусть $\omega_{i,2}$ - $m$-арная операция $\Omega_i$-алгебры $A_i$, $i=2$, ..., $n$. Тогда для любых $a_{i,1}$, ..., $a_{i,n}\in N_i$ и любого $a_{i-1}\in N_{i-1}$

$$\begin{aligned}&f_{N_i,i-1,i}(a_{i-1})(a_{i,1})...f_{N_i,i-1,i}(a_{i-1})(a_{i,m})\omega_{i,2}\\=&f_{i-1,i}(a_{i-1})(a_{i,1})...f_{i-1,i}(a_{i-1})(a_{i,m})\omega_{i,2}\\=&f_{i-1,i}(a_{i-1})(a_{i,1}...a_{i,m}\omega_{i,2})\\=&f_{N_i,i-1,i}(a_{i-1})(a_{i,1}...a_{i,m}\omega_{i,2})\end{aligned}$$

Утверждение теоремы доказано. □

[4.2] Положим $N_1=A_1$.

Из теоремы 4.5.5 следует, что если отображение $(f_{N_2,1,2}, ..., f_{N_n,1,n})$ - башня подпредставлений башни представлений $\overline{f}$, то отображение

$$(id : A_1 \to A_1, id_2 : N_2 \to A_2, ..., id_n : N_n \to A_n)$$

является морфизмом башен представлений.

Теорема 4.5.6. *Множество*[4.3] $\mathcal{B}_{\overline{f}}$ *всех башен подпредставлений башни представлений* $\overline{f}$ *порождает систему замыканий на башне представлений* $\overline{f}$ *и, следовательно, является полной структурой.*

Доказательство. Пусть для данного $\lambda \in \Lambda$, $K_{\lambda,i}$, $i = 2, ..., n$, - подалгебра $\Omega_i$-алгебры $A_i$, стабильная относительно представления $f_{i-1,i}$. Операцию пересечения на множестве $\mathcal{B}_{\overline{f}}$ мы определим согласно правилу

$$\bigcap f_{K_{\lambda,i-1,i}} = f_{\cap K_{\lambda,i-1,i}} \quad i = 2, ..., n$$

$$\bigcap \overline{K}_\lambda = \left(K_1 = A_1, K_2 = \bigcap K_{\lambda,2}, ..., K_n = \bigcap K_{\lambda,n}\right)$$

$\cap K_{\lambda,i}$ - подалгебра $\Omega_i$-алгебры $A_i$. Пусть $a_i \in \cap K_{\lambda,i}$. Для любого $\lambda \in \Lambda$ и для любого $a_{i-1} \in K_{i-1}$,

$$f_{i-1,i}(a_{i-1})(a_i) \in K_{\lambda,i}$$

Следовательно,

$$f_{i-1,i}(a_{i-1})(a_i) \in K_i$$

Повторяя приведенное построение в порядке возрастания $i$, $i = 2, ..., n$, мы видим, что $(K_1, ..., K_n)$ - кортеж стабильных множеств башни представления $\overline{f}$. Следовательно, операция пересечения башен подпредставлений определена корректно. □

Обозначим соответствующий оператор замыкания через $\overline{J}(\overline{f})$. Если мы обозначим через $\overline{X}_{[1]}$ кортеж множеств $(X_2 \subset A_2, ..., X_n \subset A_n)$ то $\overline{J}(\overline{f}, \overline{X}_{[1]})$ является пересечением всех кортежей $(K_1, ..., K_n)$, стабильных относительно представления $\overline{f}$ и таких, что для $i = 2, ..., n$, $K_i$ - подалгебра $\Omega_i$-алгебры $A_i$, содержащая $X_i$.[4.4]

Теорема 4.5.7. *Пусть*[4.5] $\overline{f}$ *- башня представлений. Пусть* $X_i \subset A_i$, $i = 2, ..., n$. *Положим* $Y_1 = A_1$. *Последовательно увеличивая значение* $i$, $i = 2, ..., n$, *определим подмножества* $X_{i,m} \subset A_i$ *индукцией по* $m$.

$$X_{i,0} = X_i$$

$$x \in X_{i,m} => x \in X_{i,m+1}$$

$$x_1 \in X_{i,m}, ..., x_p \in X_{i,m}, \omega \in \Omega_i(p) => x_1...x_p\omega \in X_{i,m+1}$$

$$x_i \in X_{i,m}, x_{i-1} \in Y_{i-1} => f_{i-1,i}(x_{i-1})(x_i) \in X_{i,m+1}$$

*Для каждого значения* $i$ *положим*

$$Y_i = \bigcup_{m=0}^{\infty} X_{i,m}$$

[4.3] Это определение аналогично определению структуры подалгебр ([13], стр. 93, 94)

[4.4] При $n = 2$, $J_2(f_{1,2}, X_2) = J_{f_{1,2}}(X_2)$. Было бы проще использовать единные обозначения в разделах 2.7 и 4.6. Однако использование векторных обозначений в разделе 2.7 мне кажется несвоевременным.

[4.5] Утверждение теоремы аналогично утверждению теоремы 5.1, [13], стр. 94.

*Тогда*

$$\overline{Y}=(Y_1,...,Y_n)=\overline{J}(\overline{f},\overline{X}_{[1]})$$

Доказательство. Для каждого значения $i$ доказательство теоремы совпадает с доказательством теоремы 2.6.4. Так как для определения устойчивого подмножества $\Omega_i$-алгебры $A_i$ нас интересует только некоторое устойчивое подмножество $\Omega_{i-1}$-алгебры $A_{i-1}$, мы должны сперва найти устойчивое подмножество $\Omega_{i-1}$-алгебры $A_{i-1}$. □

$\overline{J}(\overline{f},\overline{X}_{[1]})$ называется **башней подпредставлений башни представлений $\overline{f}$, порождённой кортежем множеств $\overline{X}_{[1]}$**, а $\overline{X}_{[1]}$ - **кортежем множеств образующих башни подпредставлений** $\overline{J}(\overline{f},\overline{X}_{[1]})$. В частности, **кортеж множеств образующих башни представлений** $\overline{f}$ будет такой кортеж $(X_2\subset A_2,...,X_n\subset A_n)$, что $\overline{J}(\overline{f},\overline{X}_{[1]})=\overline{A}$.

Нетрудно видеть, что определение кортежа множеств образующих башни представлений не зависит от того, эффективны представления башни или нет. Поэтому в дальнейшем мы будем предполагать, что представления башни эффективны и будем опираться на соглашение для эффективного представления в замечании 2.1.7.

Мы также будем пользоваться записью

$$\overline{r}_{[1]}\circ\overline{a}_{[1]}=(r_2\circ a_2,...,r_n\circ a_n)$$

для образа кортежа элементов $a_2\in A_2$, ..., $a_n\in A_n$ при эндоморфизме башни эффективных представлений. Согласно определению произведения отображений, для любых эндоморфизмов $\overline{r}_{[1]}$, $\overline{s}_{[1]}$ верно равенство

$$(\overline{r}_{[1]}\circ\overline{s}_{[1]})\circ\overline{a}_{[1]}=\overline{r}_{[1]}\circ(\overline{s}_{[1]}\circ\overline{a}_{[1]}) \tag{4.5.7}$$

Равенство (4.5.7) является законом ассоциативности для $\circ$ и позволяет записать выражение

$$\overline{r}_{[1]}\circ\overline{s}_{[1]}\circ\overline{a}_{[1]}$$

без использования скобок.

Из теоремы 4.5.7 следует следующее определение.

Определение 4.5.8. Пусть $(X_2\subset A_2,...,X_n\subset A_n)$ - кортеж множеств. Для любого кортежа элементов $\overline{a}$, $\overline{a}\in\overline{J}(\overline{f},\overline{X}_{[1]})$, существует кортеж $\overline{\Omega}$-слов, определённых согласно следующему правилу.

(1) Если $a_1\in A_1$, то $a_1$ - $\Omega_1$-слово.
(2) Если $a_i\in X_i$, $i=2$, ..., $n$, то $a_i$ - $\Omega_i$-слово.
(3) Если $a_{i,1}$, ..., $a_{i,p}$ - $\Omega_i$-слова, $i=2$, ..., $n$, и $\omega\in\Omega_i(p)$, то $a_{i,1}...a_{i,p}\omega$ - $\Omega_i$-слово.
(4) Если $a_i$ - $\Omega_i$-слово, $i=2$, ..., 3, и $a_{i-1}$ - $\Omega_{i-1}$-слово, то $a_{i-1}a_i$ - $\Omega_i$-слово.

**Кортеж $\overline{\Omega}$-слов**[4.6]

$$\overline{w}(\overline{f},\overline{X}_{[1]},\overline{a})=(w_1(\overline{f},\overline{X}_{[1]},a_1),...,w_n(\overline{f},\overline{X}_{[1]},a_n))$$

[4.6] Согласно построению

$$w_2(\overline{f},\overline{X}_{[1]},a_2)=w_2(f_{1,2},X_2,a_2)\quad w_3(\overline{f},\overline{X}_{[1]},a_3)=w_3((f_{1,2},f_{2,3}),(X_2,X_3),a_3)$$

представляет данный элемент $\overline{a} \in \overline{J}(\overline{f}, \overline{X}_{[1]})$.[4.7] Обозначим $\overline{w}(\overline{f}, \overline{X}_{[1]})$ **множество кортежей $\overline{\Omega}$-слов башни представлений** $\overline{J}(\overline{f}, \overline{X}_{[1]})$. □

Представление $a_i \in A_i$ в виде $\Omega_i$-слова неоднозначно. Если $a_{i,1}$, ..., $a_{i,p}$ - $\Omega_i$-слова, $\omega \in \Omega_i(p)$ и $a_{i-1} \in A_{i-1}$, то $\Omega_i$-слова $a_{i-1}a_{i,1}...a_{i,p}\omega$ и $a_{i-1}a_{i,1}...a_{i-1}a_{i,p}\omega$ описывают один и тот же элемент $\Omega_i$-алгебры $A_i$. Возможны равенства, связанные со спецификой представления. Например, если $\omega$ является операцией $\Omega_{i-1}$-алгебры $A_{i-1}$ и операцией $\Omega_i$-алгебры $A_i$, то мы можем потребовать, что $\Omega_i$-слова $a_{i-1,1}...a_{i-1,p}\omega a_i$ и $a_{i-1,1}a_i...a_{i-1,p}a_i\omega$ описывают один и тот же элемент $\Omega_i$-алгебры $A_i$. Перечисленные выше равенства для каждого значения $i$, $i = 2$, ..., $n$, определяют отношение эквивалентности $r_i$ на множестве $\Omega_i$-слов $W_i(\overline{f}, \overline{X}_{[1]})$.Согласно построению, отношение эквивалентности $r_i$ на множестве $\Omega_i$-слов $W_i(\overline{f}, \overline{X}_{[1]})$ зависит не только от выбора множества $X_i$, но и от выбора множества $X_{i-1}$.

Теорема 4.5.9. *Эндоморфизм $\overline{r}_{[1]}$ башни представлений $\overline{f}$ порождает отображение $\overline{\Omega}_{[1]}$-слов*

$$\overline{w}_{[1]}[\overline{f}, \overline{X}_{[1]}, \overline{r}_{[1]}] : \overline{W}_{[1]}(\overline{f}, \overline{X}_{[1]}) \to \overline{w}_{[1]}(\overline{f}, \overline{X}'_{[1]}) \quad \overline{X}_{[1]} \subset \overline{A}_{[1]} \quad \overline{X}'_{[1]} = \overline{r}_{[1]}(\overline{X}_{[1]})$$

*такое, что для любого $i$, $i = 2$, ..., $n$,*

(1) *Если* $a_i \in X_i$, $a'_i = r_i \circ a_i$, *то*
$$w_i[\overline{f}, \overline{X}_{[1]}, \overline{r}_{[1]}](a_i) = a'_i$$

(2) *Если*
$$a_{i,1}, ..., a_{i,n} \in w_i(\overline{f}, \overline{X}_{[1]})$$
$$a'_{i,1} = w_i[\overline{f}, \overline{X}_{[1]}, \overline{r}_{[1]}](a_{i,1}) \quad ... \quad a'_{i,p} = w_i[\overline{f}, \overline{X}_{[1]}, \overline{r}_{[1]}](a_{i,p})$$
*то для операции $\omega \in \Omega_i(p)$ справедливо*
$$w_i[\overline{f}, \overline{X}_{[1]}, \overline{r}_{[1]}](a_{i,1}...a_{i,p}\omega) = a'_{i,1}...a'_{i,p}\omega$$

(3) *Если*
$$a_i \in w_i(\overline{f}, \overline{X}_{[1]}) \qquad a'_i = w_i[\overline{f}, \overline{X}_{[1]}, \overline{r}_{[1]}](a_i)$$
$$a_{i-1} \in w_{i-1}(\overline{f}, \overline{X}_{[1]}) \quad a'_{i-1} = w_{i-1}[\overline{f}, \overline{X}_{[1]}, \overline{r}_{[1]}](a_{i-1})$$
*то*
$$w_i[\overline{f}, \overline{X}_{[1]}, \overline{r}_{[1]}](a_{i-1}a_i) = a'_{i-1}a'_i$$

Доказательство. Утверждения (1), (2) теоремы справедливы в силу определения эндоморфизма $r_i$. Утверждение (3) теоремы следует из равенства (4.2.3). □

[4.7]Эти обозначения имеют небольшое различие с обозначениями, предложенными ранее в теории представлении универсальной алгебры. А именно, в кортеж $\overline{\Omega}$-слов мы добавили слово $\Omega_1$-алгебры. Однако очевидно, что это слово всегда тривиально
$$w_1(\overline{f}, \overline{X}_{[1]}, a_1) = a_1$$

Замечание 4.5.10. Пусть $\overline{r}_{[1]}$ - эндоморфизм башни представлений $\overline{f}$. Пусть

$$\overline{a}_{[1]} \in \overline{J}_{[1]}(\overline{f}, \overline{X}_{[1]}) \quad \overline{a}'_{[1]} = \overline{r}_{[1]} \circ \overline{a}_{[1]} \quad \overline{X}'_{[1]} = \overline{r}_{[1]} \circ \overline{X}_{[1]}$$

Теорема 4.5.9 утверждает, что $\overline{a}'_{[1]} \in \overline{J}_{[1]}(\overline{f}, \overline{X}'_{[1]})$ . Теорема 4.5.9 также утверждает, что $\overline{\Omega}_{[1]}$-слово, представляющее $\overline{a}_{[1]}$, относительно $\overline{X}_{[1]}$ и $\overline{\Omega}_{[1]}$-слово, представляющее $\overline{a}'_{[1]}$, относительно $\overline{X}'_{[1]}$ формируются согласно одному и тому же алгоритму. Это позволяет рассматривать множество $\overline{\Omega}_{[1]}$-слов $\overline{w}_{[1]}(\overline{f}, \overline{X}'_{[1]}, \overline{a}'_{[1]})$ как отображение

$$\overline{W}_{[1]}(\overline{f}, \overline{X}_{[1]}, \overline{a}_{[1]}) : \overline{X}'_{[1]} \to \overline{w}_{[1]}(\overline{f}, \overline{X}'_{[1]}, \overline{a}'_{[1]})$$

$$\overline{W}_{[1]}(\overline{f}, \overline{X}_{[1]}, \overline{a}_{[1]})(\overline{X}'_{[1]}) = \overline{W}_{[1]}(\overline{f}, \overline{X}_{[1]}, \overline{a}_{[1]}) \circ \overline{X}'_{[1]}$$

такое, что, если для некоторого эндоморфизма $\overline{r}_{[1]}$

$$\overline{X}'_{[1]} = \overline{r}_{[1]} \circ \overline{X}_{[1]} \quad \overline{a}'_{[1]} = \overline{r}_{[1]} \circ \overline{a}_{[1]}$$

то

$$\overline{W}_{[1]}(\overline{f}, \overline{X}_{[1]}, \overline{a}_{[1]}) \circ \overline{X}'_{[1]} = \overline{w}_{[1]}(\overline{f}, \overline{X}'_{[1]}, \overline{a}'_{[1]}) = \overline{a}'_{[1]}$$

Кортеж отображений[4.8]

$$(W_2(\overline{f}, \overline{X}_{[1]}, a_2), ..., W_n(\overline{f}, \overline{X}_{[1]}, a_n))$$

представляющий данный элемент $\overline{a} \in \overline{J}(\overline{f}, \overline{X}_{[1]})$, называется **кортежем координат элемента $\overline{a}$ относительно кортежа множеств $\overline{X}_{[1]}$**. Обозначим $\overline{W}_{[1]}(\overline{f}, \overline{X}_{[1]})$ **множество кортежей координат башни представлений** $\overline{J}(\overline{f}, \overline{X}_{[1]})$. Аналогично, мы можем рассмотреть координаты кортежа множеств $\overline{B}_{[1]} \subset \overline{J}_{[1]}(f, X)$ относительно кортежа множеств $\overline{X}_{[1]}$

$$\overline{W}_{[1]}(\overline{f}, \overline{X}_{[1]}, \overline{B}_{[1]}) = (W_2(\overline{f}, \overline{X}_{[1]}, B_2), ..., W_n(\overline{f}, \overline{X}_{[1]}, B_n))$$

$$W_k(\overline{f}, \overline{X}_{[1]}, B_k) = \{W_k(\overline{f}, \overline{X}_{(2...k)}, a_k) : a_k \in B_k\} = (W_k(\overline{f}, \overline{X}_{(2...k)}, a_k), a_k \in B_k)$$

□

Теорема 4.5.11. *На множестве координат $W_k(\overline{f}, \overline{X}_{[1]})$ определена структура $\Omega_k$-алгебры.*

Доказательство. Пусть $\omega \in \Omega_k(n)$. Тогда для любых $m_1, ..., m_n \in J_k(\overline{f}, \overline{X}_{[1]})$ положим

$$W_k(\overline{f}, \overline{X}_{[1]}, m_1)...W_k(\overline{f}, \overline{X}_{[1]}, m_n)\omega = W_k(\overline{f}, \overline{X}_{[1]}, m_1...m_n\omega) \tag{4.5.9}$$

Согласно замечанию 4.5.10, из равенства (4.5.9) следует

(4.5.10)

$$\begin{aligned}(W_k(\overline{f}, \overline{X}_{[1]}, m_1)...W_k(\overline{f}, \overline{X}_{[1]}, m_n)\omega) \circ \overline{X}_{[1]} &= W_k(\overline{f}, \overline{X}_{[1]}, m_1...m_n\omega) \circ \overline{X}_{[1]} \\ &= w(\overline{f}, \overline{X}_{[1]}, m_1...m_n\omega)\end{aligned}$$

---

[4.8] Согласно построению

$$W_2(\overline{f}, \overline{X}_{[1]}, a_2) = W_2(f_{1,2}, X_2, a_2) \quad W_3(\overline{f}, \overline{X}_{[1]}, a_3) = W_3((f_{1,2}, f_{2,3}), (X_2, X_3), a_3)$$

Нас не интересует отображение $W_1(\overline{f}, \overline{X}, a_1)$ поскольку оно тривиально

$$W_1(\overline{f}, \overline{X}_{[1]}, a_1) \circ \overline{X}'_{[1]} = a_1 \tag{4.5.8}$$

Согласно правилу (3) определения 4.5.8, из равенства (4.5.10) следует

$$\begin{aligned}&(W_k(\overline{f},\overline{X}_{[1]},m_1)...W_k(\overline{f},\overline{X}_{[1]},m_n)\omega)\circ\overline{X}_{[1]}\\ =&w_k(\overline{f},\overline{X}_{[1]},m_1)...w_k(\overline{f},\overline{X}_{[1]},m_n)\omega\\ =&(W_k(\overline{f},\overline{X}_{[1]},m_1)\circ\overline{X}_{[1]})...(W_k(\overline{f},\overline{X}_{[1]},m_n)\circ\overline{X}_{[1]})\omega\end{aligned}\tag{4.5.11}$$

Из равенства (4.5.11) следует корректность определения (4.5.9) операции $\omega$ на множестве координат $W_k(\overline{f},\overline{X}_{[1]})$. □

Теорема 4.5.12. *Для $k=2,\ ...,\ n$, определено представление $\Omega_{k-1}$-алгебры $W_{k-1}(\overline{f},\overline{X}_{[1]})$ в $\Omega_k$-алгебре $W_k(\overline{f},\overline{X}_{[1]})$.*[4.9]

Доказательство. Пусть $a_{k-1}\in J_{k-1}(\overline{f},\overline{X}_{[1]})$. Тогда для любого $a_k\in J_k(\overline{f},\overline{X}_{[1]})$, положим

$$W_{k-1}(\overline{f},\overline{X}_{[1]},a_{k-1})W_k(\overline{f},\overline{X}_{[1]},a_k)=W_k(\overline{f},\overline{X}_{[1]},a_{k-1}a_k)\tag{4.5.12}$$

Согласно замечанию 4.5.10, из равенства (4.5.12) следует

$$\begin{aligned}(W_{k-1}(\overline{f},\overline{X}_{[1]},a_{k-1})W_k(\overline{f},\overline{X}_{[1]},a_k))\circ\overline{X}_{[1]}&=W_k(\overline{f},\overline{X}_{[1]},am)\circ\overline{X}_{[1]}\\ &=w_k(\overline{f},\overline{X}_{[1]},a_{k-1}a_k)\end{aligned}\tag{4.5.13}$$

Согласно правилу (4) определения 4.5.8, из равенства (4.5.13) следует

$$\begin{aligned}&(W_{k-1}(\overline{f},\overline{X}_{[1]},a_{k-1})W_k(\overline{f},\overline{X}_{[1]},a_k))\circ\overline{X}_{[1]}\\ =&w_{k-1}(\overline{f},\overline{X}_{[1]},a_{k-1})w_k(\overline{f},\overline{X}_{[1]},a_k)\\ =&(W_{k-1}(\overline{f},\overline{X}_{[1]},a_{k-1})\circ\overline{X}_{[1]})(W_k(\overline{f},\overline{X}_{[1]},a_k)\circ\overline{X}_{[1]})\end{aligned}\tag{4.5.14}$$

Из равенства (4.5.14) следует корректность определения (4.5.12) представления $\Omega_{k-1}$-алгебры $W_{k-1}(\overline{f},\overline{X}_{[1]})$ в $\Omega_k$-алгебре $W_k(\overline{f},\overline{X}_{[1]})$. □

Следствие 4.5.13. *Кортеж $\overline{\Omega}$-алгебр $\overline{W}(\overline{f},\overline{X}_{[1]})$ порождает башню представлений.* □

Теорема 4.5.14. *Пусть $\overline{f}$ - башня представлений. Для заданных множеств $X_k\subset A_k$, $X'_k\subset A_k$, $k=2,\ ...,\ n$, рассмотрим кортеж отображений*

$$\overline{R}_{[1]}=(R_2,...,R_n)$$

*таких, что для любого $k=2,\ ...,\ n$, отображение*

$$R_k:X_k\to X'_k$$

*согласовано со структурой представления $f_{k-1,k}$, т. е.*[4.10]

$$\begin{aligned}&\omega\in\Omega_k(p),\ x_{1,k},\ ...,\ x_{p,k},\ x_{1,k}...x_{p,k}\omega\in X_k,\ R_1(x_{1,k}...x_{p,k}\omega)\in X'_k\\ =>&R_1(x_{1,k}...x_{p,k}\omega)=R_1(x_{1,k})...R_1(x_{p,k})\omega\\ &a_{k-1}\in X_{k-1},\ a_k\in X_k,\ R_1(a_{k-1}a_k)\in X'_k\\ =>&R_k(a_{k-1}a_k)=R_{k-1}(a_{k-1})R_k(a_k)\end{aligned}$$

4.9 Согласно равенству (4.5.8), мы можем отождествить $\Omega_1$-алгебры $W_1(\overline{f},\overline{X}_{[1]})$ и $A_1$.

4.10 Мы полагаем $X_1=A_1$, $R_1(a_1)=a_1$.

*Рассмотрим отображение $\overline{\Omega}_{[1]}$-слов*

$$\overline{w}_{[1]}[\overline{f},\overline{X}_{[1]},\overline{R}_{[1]}]:\overline{w}_{[1]}(\overline{f},\overline{X}_{[1]})\rightarrow\overline{w}_{[1]}(\overline{f},\overline{X}'_{[1]})$$

*удовлетворяющее условиям* (1), (2), (3) *теоремы 4.5.9, и такое, что*

$$\overline{x}_{[1]}\in\overline{X}_{[1]}=>w[\overline{f},\overline{X}_{[1]},\overline{R}_{[1]}](\overline{x}_{[1]})=\overline{R}_{[1]}(\overline{x}_{[1]})$$

*Для каждого $k$, $k=2$, ..., $n$, существует эндоморфизм $\Omega_k$-алгебры $A_k$*

$$r_k:A_k\rightarrow A_k$$

*определённый правилом*

$$r_k\circ a_k=w_k[\overline{f},\overline{X}_{[1]},\overline{R}_{[1]}](w_k(\overline{f},\overline{X}_{[1]},a_k))\tag{4.5.15}$$

*Кортеж эндоморфизмов $\overline{r}_{[1]}$ является морфизмом башен представлнений $\overline{J}(\overline{f},\overline{X}_{[1]})$ и $\overline{J}(\overline{f},\overline{X}'_{[1]})$.*

Доказательство. При $n=2$ башня представлений $\overline{f}$ является представлением $\Omega_1$-алгебры $A_1$ в $\Omega_2$-алгебре $A_2$. Утверждение теоремы является следствием теоремы 2.6.11.

Пусть утверждение теоремы верно для $n-1$. Обозначения в теореме не меняются при переходе от одного уровня к другому, так как слово в $\Omega_{n-1}$-алгебре $A_{n-1}$ не зависит от слова в $\Omega_n$-алгебре $A_n$. Мы будем доказывать теорему индукцией по сложности $\Omega_n$-слова.

Если $w_n(\overline{f},\overline{X}_{[1]},a_n)=a_n$, то $a_n\in X_n$. Согласно условию (1) теоремы 4.5.9,

$$\begin{aligned}r_n\circ a_n&=w_n[\overline{f},\overline{X}_{[1]},\overline{R}_{[1]}](w_n(\overline{f},\overline{X}_{[1]},a_n))\\&=w_n[\overline{f},\overline{X}_{[1]},\overline{R}_{[1]}](a_n)\\&=R_n(a_n)\end{aligned}$$

Следовательно, на множестве $X_n$ отображения $r_n$ и $R_n$ совпадают, и отображение $r_n$ согласовано со структурой $\Omega_n$-алгебры.

Пусть $\omega\in\Omega_n(p)$. Пусть предположение индукции верно для

$$a_{n,1},...,a_{n,p}\in J_n(\overline{f},\overline{X}_{[1]})$$

Пусть

$$w_{n,1}=w_n(\overline{f},\overline{X}_{[1]},a_{n,1})\quad...\quad w_{n,p}=w_n(\overline{f},\overline{X}_{[1]},a_{n,p})$$

Если

$$a_n=a_{n,1}...a_{n,p}\omega$$

то согласно условию (3) определения 4.5.8,

$$w_n[\overline{f},\overline{X}_{[1]},a_n]=w_{n,1}...w_{n,p}\omega$$

Согласно условию (2) теоремы 4.5.9,

$$\begin{aligned}r_n\circ a_n&=w_n[\overline{f},\overline{X}_{[1]},\overline{R}_{[1]}](w_n(\overline{f},\overline{X}_{[1]},a_n))\\&=w_n(\overline{f},\overline{X}_{[1]},\overline{R}_{[1]}](w_{n,1}...w_{n,p}\omega)\\&=w_n(\overline{f},\overline{X}_{[1]},\overline{R}_{[1]}](w_{n,1})...w_n(\overline{f},\overline{X}_{[1]},\overline{R}_{[1]}](w_{n,p})\omega\\&=(r_n\circ a_{n,1})...(r_n\circ a_{n,p})\omega\end{aligned}$$

Следовательно, отображение $r_n$ является эндоморфизмом $\Omega_n$-алгебры $A_n$.

Пусть предположение индукции верно для

$$a_n \in J_n(\overline{f},\overline{X}_{[1]}) \qquad w_n(\overline{f},\overline{X}_{[1]},a_n)=m_n$$
$$a_{n-1} \in J_{n-1}(\overline{f},\overline{X}_{[1]}) \qquad w_{n-1}(\overline{f},\overline{X}_{[1]},a_{n-1})=m_{n-1}$$

Согласно условию (4) определения 4.5.8,

$$w_n(\overline{f},\overline{X}_{[1]},a_{n-1}a_n)=m_{n-1}m_n$$

Согласно условию (3) теоремы 4.5.9,

$$\begin{aligned}r_n\circ(a_{n-1}a_n)&=w_n(\overline{f},\overline{X}_{[1]},\overline{R}_{[1]}](w_n(\overline{f},\overline{X}_{[1]},a_{n-1}a_n))\\&=w_n(\overline{f},\overline{X}_{[1]},\overline{R}_{[1]}](m_{n-1}m_n)\\&=w_{n-1}(\overline{f},\overline{X}_{[1]},\overline{R}_{[1]}](m_{n-1})w_n(\overline{f},(r_1,R_2,...,R_n),\overline{X}_{[1]})(m_n)\\&=(r_{n-1}\circ a_{n-1})(r_n\circ a_n)\end{aligned}$$

Из равенства (4.2.3) следует, что отображение $\overline{r}$ является морфизмом башни представлений $\overline{f}$. □

Замечание 4.5.15. Теорема 4.5.14 - это теорема о продолжении отображения. Единственное, что нам известно о кортеже множеств $\overline{X}_{[1]}$ - это то, что $\overline{X}_{[1]}$ - кортеж множеств образующих башни представлений $\overline{f}$. Однако, между элементами множества $X_k$, $k=2$, ..., $n$, могут существовать соотношения, порождённые либо операциями $\Omega_k$-алгебры $A_k$, либо преобразованиями представления $f_{k-1,k}$. Поэтому произвольное отображение кортежа множеств $\overline{X}_{[1]}$, вообще говоря, не может быть продолжено до эндоморфизма башни представлений $\overline{f}$.[4.11] Однако, если для каждого $k$, $k=2$, ..., $n$, отображение $R_k$ согласованно со структурой представления $f_{k-1,k}$, то мы можем построить продолжение этого отображения, которое является эндоморфизмом башни представлений $\overline{f}$. □

Определение 4.5.16. Пусть $\overline{X}_{[1]}$ - кортеж множеств образующих башни представлений $\overline{f}$. Пусть $\overline{r}_{[1]}$ - эндоморфизм башни представлений $\overline{f}$. Множество координат $\overline{W}_{[1]}(\overline{f},\overline{X}_{[1]},\overline{r}_{[1]}\circ\overline{X}_{[1]})$ называется **координатами эндоморфизма башни представлений**. □

Определение 4.5.17. Пусть $\overline{X}_{[1]}$ - кортеж множеств образующих башни представлений $\overline{f}$. Пусть $\overline{r}_{[1]}$ - эндоморфизм башни представлений $\overline{f}$. Пусть $\overline{a}_{[1]}\in\overline{A}_{[1]}$. Мы определим **суперпозицию координат** эндоморфизма башни представлений $\overline{f}$ и элемента $\overline{a}_{[1]}$ как координаты, определённые согласно правилу

$$\overline{W}_{[1]}(\overline{f},\overline{X}_{[1]},\overline{a}_{[1]})\circ\overline{W}_{[1]}(\overline{f},\overline{X}_{[1]},\overline{r}_{[1]}\circ\overline{X}_{[1]})=\overline{W}_{[1]}(\overline{f},\overline{X}_{[1]},\overline{r}_{[1]}\circ\overline{a}_{[1]}) \tag{4.5.16}$$

Пусть $\overline{Y}_{[1]}\subset\overline{A}_{[1]}$. Мы определим суперпозицию координат эндоморфизма башни представлений $\overline{f}$ и кортежа множеств $\overline{Y}_{[1]}$ согласно правилу

$$\begin{aligned}&\overline{W}_{[1]}(\overline{f},\overline{X}_{[1]},\overline{Y}_{[1]})\circ\overline{W}_{[1]}(\overline{f},\overline{X}_{[1]},\overline{r}_{[1]}\circ\overline{X}_{[1]})\\&=(\overline{W}_{[1]}(\overline{f},\overline{X}_{[1]},\overline{a}_{[1]})\circ\overline{W}_{[1]}(\overline{f},\overline{X}_{[1]},\overline{r}_{[1]}\circ\overline{X}_{[1]}),\overline{a}_{[1]}\in\overline{Y}_{[1]})\end{aligned} \tag{4.5.17}$$

[4.11] В теореме 4.6.7, требования к кортежу множеств образующих более жёсткие. Поэтому теорема 4.6.7 говорит о продолжении произвольного отображения. Более подробный анализ дан в замечании 4.6.9.

□

Теорема 4.5.18. *Эндоморфизм $\overline{r}_{[1]}$ представления $\overline{f}$ порождает отображение координат башни представлений*

$$\overline{W}_{[1]}[\overline{f},\overline{X}_{[1]},\overline{r}_{[1]}]:\overline{W}_{[1]}(\overline{f},\overline{X}_{[1]})\to\overline{W}_{[1]}(\overline{f},\overline{X}_{[1]}) \tag{4.5.18}$$

*такое, что*

(4.5.19)

$$\overline{W}_{[1]}(\overline{f},\overline{X}_{[1]},\overline{a}_{[1]})\to\overline{W}_{[1]}[\overline{f},\overline{X}_{[1]},\overline{r}_{[1]}]\star\overline{W}_{[1]}(\overline{f},\overline{X}_{[1]},\overline{a}_{[1]})=\overline{W}_{[1]}(\overline{f},\overline{X}_{[1]},\overline{r}_{[1]}\circ\overline{a}_{[1]})$$

$$\overline{W}_{[1]}[\overline{f},\overline{X}_{[1]},\overline{r}_{[1]}]\star\overline{W}_{[1]}(\overline{f},\overline{X}_{[1]},\overline{a}_{[1]})=\overline{W}_{[1]}(\overline{f},\overline{X}_{[1]},\overline{a}_{[1]})\circ\overline{W}_{[1]}(\overline{f},\overline{X}_{[1]},\overline{r}_{[1]}\circ\overline{X}_{[1]})$$

Доказательство. Согласно замечанию 4.5.10, мы можем рассматривать равенства (4.5.16), (4.5.18) относительно заданного кортежа множеств образующих $\overline{X}_{[1]}$. При этом координатам $\overline{W}_{[1]}(\overline{f},\overline{X}_{[1]},\overline{a}_{[1]})$ соответствует кортеж слов

$$\overline{W}_{[1]}(\overline{f},\overline{X}_{[1]},\overline{a}_{[1]})\circ\overline{X}_{[1]}=\overline{w}_{[1]}(\overline{f},\overline{X}_{[1]},\overline{a}_{[1]}) \tag{4.5.20}$$

а координатам $\overline{W}_{[1]}(\overline{f},\overline{X}_{[1]},\overline{r}_{[1]}\circ\overline{a}_{[1]})$ соответствует кортеж слов

$$\overline{W}_{[1]}(\overline{f},\overline{X}_{[1]},\overline{r}_{[1]}\circ\overline{a}_{[1]})\circ X=\overline{w}_{[1]}(\overline{f},\overline{X}_{[1]},\overline{r}_{[1]}\circ\overline{a}_{[1]}) \tag{4.5.21}$$

Поэтому для того, чтобы доказать теорему, нам достаточно показать, что отображению $\overline{W}_{[1]}[\overline{f},\overline{X}_{[1]},\overline{r}_{[1]}]$ соответствует отображение $\overline{w}_{[1]}[\overline{f},\overline{X}_{[1]},\overline{r}_{[1]}]$. Мы докажем утверждение индукцией по числу $n$ универсальных алгебр в башне представлений.

При $n=2$ башня представлений $\overline{f}$ является представлением $\Omega_1$-алгебры $A_1$ в $\Omega_2$-алгебре $A_2$. Утверждение теоремы является следствием теоремы 2.6.15. Пусть утверждение теоремы верно для $n-1$. Мы будем доказывать теорему индукцией по сложности $\Omega_n$-слова.

Если $a_n\in X_n$, $a'_n=r_n\circ a_n$, то, согласно равенствам (4.5.20), (4.5.21), отображения $W_n[\overline{f},\overline{X}_{[1]},\overline{r}_{[1]}]$ и $w_n[\overline{f},\overline{X}_{[1]},\overline{r}_{[1]}]$ согласованы.

Пусть для $a_{n,1}$, ..., $a_{n,p}\in X_n$ отображения $W_n[\overline{f},\overline{X}_{[1]},\overline{r}_{[1]}]$ и $w_n[\overline{f},\overline{X}_{[1]},\overline{r}_{[1]}]$ согласованы. Пусть $\omega\in\Omega_n(p)$. Согласно теореме 2.6.9

$$W_n(\overline{f},\overline{X}_{[1]},a_{n,1}...a_{n,p}\omega)=W_n(\overline{f},\overline{X}_{[1]},a_{n,1})...W_n(\overline{f},\overline{X}_{[1]},a_{n,p})\omega \tag{4.5.22}$$

Так как $r_n$ - эндоморфизм $\Omega_n$-алгебры $A_n$, то из равенства (4.5.22) следует

(4.5.23)

$$\begin{aligned}W_n(\overline{f},\overline{X}_{[1]},r_n\circ(a_{n,1}...a_{n,p}\omega))&=W_n(\overline{f},\overline{X}_{[1]},(r_n\circ a_{n,1})...(r_n\circ a_{n,p})\omega)\\&=W_n(\overline{f},\overline{X}_{[1]},r_n\circ a_{n,1})...W(f,X,r_n\circ a_{n,p})\omega\end{aligned}$$

Из равенств (4.5.22), (4.5.23) и предположения индукции следует, что отображения $W_n[\overline{f},\overline{X}_{[1]},\overline{r}_{[1]}]$ и $w_n[\overline{f},\overline{X}_{[1]},\overline{r}_{[1]}]$ согласованы для $a_n=a_{n,1}...a_{n,p}\omega$.

Пусть для $a_{n,1}\in A_n$ отображения $W_n[\overline{f},\overline{X}_{[1]},\overline{r}_{[1]}]$ и $w_n[\overline{f},\overline{X}_{[1]},\overline{r}_{[1]}]$ согласованы. Пусть $a_{n-1}\in A_{n-1}$. Согласно теореме 4.5.12

$$W_n(\overline{f},\overline{X}_{[1]},a_{n-1}a_{n,1})=W_{n-1}(\overline{f},\overline{X}_{[1]},a_{n-1})W_n(\overline{f},\overline{X}_{[1]},a_{n,1}) \tag{4.5.24}$$

Так как $\overline{r}_{[1]}$ - эндоморфизм башни представлений $\overline{f}$, то из равенства (4.5.24) следует

(4.5.25)

$$\begin{aligned}W_n(\overline{f},\overline{X}_{[1]},r_n\circ(a_{n-1}a_{n,1}))=&W_n(\overline{f},\overline{X}_{[1]},(r_{n-1}\circ a_{n-1})(r_n\circ a_{n,1}))\\=&W_{n-1}(\overline{f},\overline{X}_{[1]},r_{n-1}\circ a_{n-1})W_n(\overline{f},\overline{X}_{[1]},r_n\circ a_{n,1})\end{aligned}$$

Из равенств (4.5.24), (4.5.25) и предположения индукции следует, что отображения $W_n[\overline{f},\overline{X}_{[1]},\overline{r}_{[1]}]$ и $w_n[\overline{f},\overline{X}_{[1]},\overline{r}_{[1]}]$ согласованы для $a_{n,2}=a_{n-1}a_{n,1}$. □

Следствие 4.5.19. *Пусть $\overline{X}_{[1]}$ - кортеж множеств образующих башни представлений $\overline{f}$. Пусть $\overline{r}_{[1]}$ - эндоморфизм башни представлений $\overline{f}$. Отображение $\overline{W}_{[1]}[\overline{f},\overline{X}_{[1]},\overline{r}_{[1]}]$ является эндоморфизмом башни представлений $\overline{W}(\overline{f},\overline{X}_{[1]})$.* □

В дальнейшем мы будем отождествлять отображение $\overline{W}_{[1]}[\overline{f},\overline{X}_{[1]},\overline{r}_{[1]}]$ и множество координат $\overline{W}_{[1]}(\overline{f},\overline{X}_{[1]},\overline{r}_{[1]}\circ\overline{X}_{[1]})$.

Теорема 4.5.20. *Пусть $\overline{X}_{[1]}$ - кортеж множеств образующих башни представлений $\overline{f}$. Пусть $\overline{r}_{[1]}$ - эндоморфизм башни представлений $\overline{f}$. Пусть $\overline{Y}_{[1]}\subset\overline{M}_{[1]}$. Тогда*

(4.5.26) $$\overline{W}_{[1]}(\overline{f},\overline{X}_{[1]},\overline{Y}_{[1]})\circ\overline{W}_{[1]}(\overline{f},\overline{X}_{[1]},\overline{r}_{[1]}\circ\overline{X}_{[1]})=\overline{W}_{[1]}(\overline{f},\overline{X}_{[1]},\overline{r}_{[1]}\circ\overline{Y}_{[1]})$$

(4.5.27) $$\overline{W}_{[1]}[\overline{f},\overline{X}_{[1]},\overline{r}_{[1]}]\star\overline{W}_{[1]}(\overline{f},\overline{X}_{[1]},\overline{Y}_{[1]})=\overline{W}_{[1]}(\overline{f},\overline{X}_{[1]},\overline{r}_{[1]}\circ\overline{Y}_{[1]})$$

Доказательство. Равенство (4.5.26) является следствием равенства

$$\overline{r}_{[1]}\circ\overline{Y}_{[1]}=(\overline{r}_{[1]}\circ\overline{a}_{[1]},\overline{a}_{[1]}\in\overline{Y}_{[1]})$$

а также равенств (4.5.16), (4.5.17). Равенство (4.5.27) является следствием равенств (4.5.26), (4.5.19). □

Теорема 4.5.21. *Пусть $\overline{X}_{[1]}$ - кортеж множеств образующих башни представлений $\overline{f}$. Пусть $\overline{r}_{[1]}$, $\overline{s}_{[1]}$ - эндоморфизмы башни представлений $\overline{f}$. Тогда*

(4.5.28)

$$\overline{W}_{[1]}(\overline{f},\overline{X}_{[1]},\overline{s}_{[1]}\circ\overline{X}_{[1]})\circ\overline{W}_{[1]}(\overline{f},\overline{X}_{[1]},\overline{r}_{[1]}\circ\overline{X}_{[1]})=\overline{W}_{[1]}(\overline{f},\overline{X}_{[1]},\overline{r}_{[1]}\circ\overline{s}_{[1]}\circ\overline{X}_{[1]})$$

(4.5.29) $$\overline{W}_{[1]}[\overline{f},\overline{X}_{[1]},\overline{r}_{[1]}]\star\overline{W}_{[1]}[\overline{f},\overline{X}_{[1]},\overline{s}_{[1]}]=\overline{W}_{[1]}[\overline{f},\overline{X}_{[1]},\overline{r}_{[1]}\circ\overline{s}_{[1]}]$$

Доказательство. Равенство (4.5.28) следует из равенства (4.5.26), если положить $\overline{Y}_{[1]}=\overline{s}_{[1]}\circ\overline{X}_{[1]}$. Равенство (4.5.29) следует из равенства (4.5.28) и цепочки равенств

$$\begin{aligned}&(\overline{W}_{[1]}[\overline{f},\overline{X}_{[1]},\overline{r}_{[1]}]\star\overline{W}_{[1]}[\overline{f},\overline{X}_{[1]},\overline{s}_{[1]}])\star\overline{W}_{[1]}(\overline{f},\overline{X}_{[1]},\overline{Y}_{[1]})\\=&\overline{W}_{[1]}[\overline{f},\overline{X}_{[1]},\overline{r}_{[1]}]\star(\overline{W}_{[1]}[\overline{f},\overline{X}_{[1]},\overline{s}_{[1]}]\star\overline{W}_{[1]}(\overline{f},\overline{X}_{[1]},\overline{Y}_{[1]}))\\=&(\overline{W}_{[1]}(\overline{f},\overline{X}_{[1]},\overline{Y}_{[1]})\circ\overline{W}_{[1]}(\overline{f},\overline{X}_{[1]},\overline{s}_{[1]}\circ\overline{X}_{[1]}))\circ\overline{W}_{[1]}(\overline{f},\overline{X}_{[1]},\overline{r}_{[1]}\circ\overline{X}_{[1]})\\=&\overline{W}_{[1]}(\overline{f},\overline{X}_{[1]},\overline{Y}_{[1]})\circ(\overline{W}_{[1]}(\overline{f},\overline{X}_{[1]},\overline{s}_{[1]}\circ\overline{X}_{[1]})\circ\overline{W}_{[1]}(\overline{f},\overline{X}_{[1]},\overline{r}_{[1]}\circ\overline{X}_{[1]}))\\=&\overline{W}_{[1]}(\overline{f},\overline{X}_{[1]},\overline{Y}_{[1]})\circ\overline{W}_{[1]}(\overline{f},\overline{X}_{[1]},\overline{r}_{[1]}\circ\overline{s}_{[1]}\circ\overline{X}_{[1]})\\=&\overline{W}_{[1]}(\overline{f},\overline{X}_{[1]},\overline{r}_{[1]}\circ\overline{s}_{[1]})\star\overline{W}_{[1]}(\overline{f},\overline{X}_{[1]},\overline{Y}_{[1]})\end{aligned}$$

□

Мы можем обобщить определение суперпозиции координат и предположить, что один из множителей является кортежем множеств $\overline{\Omega}_{[1]}$-слов. Соответственно, определение суперпозиции координат имеет вид

$$\begin{aligned}&\overline{w}_{[1]}(\overline{f},\overline{X}_{[1]},\overline{Y}_{[1]})\circ\overline{W}_{[1]}(\overline{f},\overline{X}_{[1]},\overline{r}_{[1]}\circ\overline{X}_{[1]})\\=&\overline{W}_{[1]}(\overline{f},\overline{X}_{[1]},\overline{Y}_{[1]})\circ\overline{w}_{[1]}(\overline{f},\overline{X}_{[1]},\overline{r}_{[1]}\circ\overline{X}_{[1]})\\=&\overline{w}_{[1]}(\overline{f},\overline{X}_{[1]},\overline{r}_{[1]}\circ\overline{Y}_{[1]})\end{aligned}$$

Следующие формы записи образа кортежа множеств $\overline{Y}_{[1]}$ при эндоморфизме $\overline{r}_{[1]}$ эквивалентны.

$$\begin{aligned}\overline{r}_{[1]}\circ\overline{Y}_{[1]}&=\overline{W}_{[1]}(\overline{f},\overline{X}_{[1]},Y)\circ(\overline{r}_{[1]}\circ\overline{X}_{[1]})\\&=\overline{W}_{[1]}(\overline{f},\overline{X}_{[1]},\overline{Y}_{[1]})\circ(\overline{W}_{[1]}(\overline{f},\overline{X}_{[1]},\overline{r}_{[1]}\circ\overline{X}_{[1]})\circ\overline{X}_{[1]})\end{aligned}\tag{4.5.30}$$

Из равенств (4.5.26), (4.5.30) следует, что

$$\begin{aligned}&(\overline{W}_{[1]}(\overline{f},\overline{X}_{[1]},\overline{Y}_{[1]})\circ\overline{W}_{[1]}(\overline{f},\overline{X}_{[1]},\overline{r}_{[1]}\circ\overline{X}_{[1]}))\circ\overline{X}_{[1]}\\=&\overline{W}_{[1]}(\overline{f},\overline{X}_{[1]},\overline{Y}_{[1]})\circ(\overline{W}_{[1]}(\overline{f},\overline{X}_{[1]},\overline{r}_{[1]}\circ\overline{X}_{[1]})\circ\overline{X}_{[1]})\end{aligned}\tag{4.5.31}$$

Равенство (4.5.31) является законом ассоциативности для операции композиции и позволяет записать выражение

$$\overline{W}_{[1]}(\overline{f},\overline{X}_{[1]},\overline{Y}_{[1]})\circ\overline{W}_{[1]}(\overline{f},\overline{X}_{[1]},\overline{r}_{[1]}\circ\overline{X}_{[1]})\circ X$$

без использования скобок.

Определение 4.5.22. Пусть $(X_2\subset A_2,...,X_n\subset A_n)$ - кортеж множеств образующих башни представлений $\overline{f}$. Пусть отображение $\overline{r}_{[1]}$ является эндоморфизмом башни представления $\overline{f}$. Пусть кортеж множеств $\overline{X}'_{[1]}=\overline{r}_{[1]}\circ\overline{X}_{[1]}$ является образом кортежа множеств $\overline{X}_{[1]}$ при отображении $\overline{r}_{[1]}$. Эндоморфизм $\overline{r}_{[1]}$ башни представлений $\overline{f}$ называется **невырожденным на кортеже множеств образующих** $\overline{X}_{[1]}$, если кортеж множеств $\overline{X}'_{[1]}$ является кортежем множеств образующих башни представлений $\overline{f}$. В противном случае, эндоморфизм $\overline{r}_{[1]}$ называется **вырожденным на кортеже множеств образующих** $\overline{X}_{[1]}$, □

Определение 4.5.23. Эндоморфизм башни представлений $\overline{f}$ называется **невырожденным**, если он невырожден на любом кортеже множеств образующих. □

Теорема 4.5.24. *Автоморфизм $\overline{r}$ башни представлений $\overline{f}$ является невырожденным эндоморфизмом.*

Доказательство. Пусть $\overline{X}_{[1]}$ - кортеж множеств образующих башни представлений $\overline{f}$. Пусть $\overline{X}'_{[1]}=\overline{r}_{[1]}\circ\overline{X}_{[1]}$.

Согласно теореме 4.5.18 эндоморфизм $\overline{r}_{[1]}$ порождает отображение $\overline{\Omega}_{[1]}$-слов $\overline{w}_{[1]}[\overline{f},\overline{X}_{[1]},\overline{r}_{[1]}]$.

Пусть $\overline{a}'_{[1]}\in\overline{A}_{[1]}$. Так как $\overline{r}_{[1]}$ - автоморфизм, то существует $\overline{a}_{[1]}\in\overline{A}_{[1]}$, $\overline{r}_{[1]}\circ\overline{a}_{[1]}=\overline{a}'_{[1]}$. Согласно определению 4.5.8, $\overline{w}_{[1]}(\overline{f},\overline{X}_{[1]},\overline{a}_{[1]})$ - кортеж $\overline{\Omega}_{[1]}$-слов, представляющих $\overline{a}_{[1]}$ относительно кортежа множеств образующих $\overline{X}_{[1]}$.

Согласно теореме 4.5.18, $\overline{w}_{[1]}(\overline{f}, \overline{X}'_{[1]}, \overline{a}'_{[1]})$ - кортеж $\overline{\Omega}_{[1]}$-слов, представляющих $\overline{a}'_{[1]}$ относительно кортежа множеств $\overline{X}'_{[1]}$

$$\overline{w}_{[1]}(\overline{f}, \overline{X}'_{[1]}, \overline{a}'_{[1]}) = \overline{w}_{[1]}[\overline{f}, \overline{X}_{[1]}, \overline{r}](\overline{w}_{[1]}(\overline{f}, \overline{X}_{[1]}, \overline{a}_{[1]}))$$

Следовательно, $\overline{X}'_{[1]}$ - множество образующих представления $\overline{f}$. Согласно определению 4.5.23, автоморфизм $\overline{r}_{[1]}$ - невырожден. □

## 4.6. Базис башни представлений

Определение 4.6.1. Если кортеж множеств $\overline{X}_{[1]}$ является кортежем множеств образующих башни представлений $\overline{f}$, то любой кортеж множеств $\overline{Y}_{[1]}$, $X_i \subset Y_i \subset A_i$, $i = 2, ..., n$, также является кортежем множеств образующих башни представлений $\overline{f}$. Если существует кортеж минимальных множеств $\overline{X}_{[1]}$, порождающих башню представлений $\overline{f}$, то такой кортеж множеств $\overline{X}_{[1]}$ называется **базисом башни представлений** $\overline{f}$. □

Теорема 4.6.2. *Базис башни представлений определён индукцией по $n$. При $n = 2$ базис башни представлений является базисом представления $f_{1,2}$. Если кортеж множеств $\overline{X}_{[1,n]}$ является базисом башни представлений $\overline{f}_{[1]}$, то кортеж множеств образующих $\overline{X}_{[1]}$ башни представлений $\overline{f}$ является базисом тогда и только тогда, когда для любого $a_n \in X_n$ кортеж множеств $(X_2, ..., X_{n-1}, X_n \backslash \{a_n\})$ не является кортежем множеств образующих башни представлений $\overline{f}$.*

Доказательство. При $n = 2$ утверждение теоремы является следствием теоремы 2.7.2.

Пусть $n > 2$. Пусть $\overline{X}_{[1]}$ - кортеж множеств образующих башни расслоений $\overline{f}$. Пусть кортеж множеств $\overline{X}_{[1,n]}$ является базисом башни расслоений $\overline{f}_{[n]}$. Допустим для некоторого $a_n \in X_n$ существует слово

$$w_n = w_n(a_n, \overline{f}, (X_1, ..., X_{n-1}, X_n \setminus \{a_n\})) \tag{4.6.1}$$

Рассмотрим $a'_n \in A_n$, для которого слово

$$w'_n = w_n(a'_n, \overline{f}, \overline{X}_{[1]}) \tag{4.6.2}$$

зависит от $a_n$. Согласно определению 2.6.6, любое вхождение $a_n$ в слово $w'_n$ может быть заменено словом $w_n$. Следовательно, слово $w'_n$ не зависит от $a_n$, а кортеж множеств $(X_2, ..., X_{n-1}, X_n \setminus \{a_n\})$ является кортежем множеств образующих башни представлений $\overline{f}$. Следовательно, $\overline{X}_{[1]}$ не является базисом башни представлений $\overline{f}$. □

Замечание 4.6.3. Доказательство теоремы 4.6.2 даёт нам эффективный метод построения базиса башни представлений $\overline{f}$. Мы начинаем строить базис в самом нижнем слое. Когда базис построен в слое $i$, $i = 2, ..., n-1$, мы можем перейти к построению базиса в слое $i+1$. □

Замечание 4.6.4. Мы будем записывать базис также в виде

$$X_k = (x_k, x_k \in X_k) \quad k = 2, ..., n$$

Если базис - конечный, то мы будем также пользоваться записью

$$X_k = (x_{k\cdot i}, i \in I_k) = (x_{k\cdot 1}, ..., x_{k\cdot p_k}) \quad k = 2, ..., n$$

□

Замечание 4.6.5. Координаты $a_k \in A_k$ относительно базиса $\overline{X}_{[1]}$ башни представлений $\overline{f}$ определены неоднозначно. Так же как в замечании 2.7.5, мы рассмотрим **кортеж отношений эквивалентности $\overline{\rho}_{[1]}(\overline{f})$, порождённых башней представлений $\overline{f}$**. Если в процессе построений мы получим равенство двух $\Omega_k$-слов относительно заданного базиса, мы можем утверждать, что эти $\Omega_k$-слова равны, не заботясь об отношении эквивалентности $\rho_k(\overline{f})$. □

Теорема 4.6.6. *Автоморфизм башни представлений $\overline{f}$ отображает базис башни представлений $\overline{f}$ в базис.*

Доказательство. При $n = 2$ утверждение теоремы является следствием теоремы 4.6.6.

Пусть отображение $\overline{r}$ - автоморфизм башни представлений $\overline{f}$. Пусть кортеж множеств $\overline{X}_{[1]}$ - базис башни представлений $\overline{f}$. Пусть $\overline{X}'_{[1]} = \overline{r}_{[1]}(\overline{X}_{[1]})$.

Допустим кортеж множеств $\overline{X}'_{[1]}$ не является базисом. Согласно теореме 4.6.2 существуют $i$, $i = 2, ..., n$, и $a'_i \in X'_i$ такие, что кортеж множеств $\overline{X}''_{[1]}$ [4.12] является кортежем множеств образующих башни представлений $\overline{f}$. Согласно теореме 4.3.3 отображение $\overline{r}^{-1}$ является автоморфизмом башни представлений $\overline{f}$. Согласно теореме 4.5.24 и определению 4.5.23, кортеж множеств $\overline{X}'''_{[1]}$ [4.13] является кортежем множеств образующих представления $\overline{f}$. Полученное противоречие доказывает теорему. □

Теорема 4.6.7. *Пусть $\overline{X}_{[1]}$ - базис представления $\overline{f}$. Пусть*

$$\overline{R}_{[1]} : \overline{X}_{[1]} \to \overline{X}'_{[1]}$$

*произвольное отображение кортежа множеств $\overline{X}_{[1]}$. Рассмотрим отображение $\overline{\Omega}_{[1]}$-слов*

$$\overline{w}_{[1]}[\overline{f}, \overline{X}_{[1]}, \overline{R}_{[1]}] : \overline{w}_{[1]}(\overline{f}, \overline{X}_{[1]}) \to \overline{w}_{[1]}(\overline{f}, \overline{X}'_{[1]})$$

*удовлетворяющее условиям (1), (2), (3) теоремы 4.5.9, и такое, что*

$$\overline{x}_{[1]} \in \overline{X}_{[1]} => w[\overline{f}, \overline{X}_{[1]}, \overline{R}_{[1]}](\overline{x}_{[1]}) = \overline{R}_{[1]}(\overline{x}_{[1]})$$

*Существует единственный эндоморфизм представления $\overline{f}$* [4.14]

$$\overline{r}_{[1]} : \overline{A}_{[1]} \to \overline{A}_{[1]}$$

*определённый правилом*

$$\overline{r}_{[1]} \circ \overline{a}_{[1]} = w[\overline{f}, \overline{X}_{[1]}, \overline{R}_{[1]}](\overline{w}_{[1]}(\overline{f}, \overline{X}_{[1]}, \overline{a}_{[1]}))$$

Доказательство. Утверждение теоремы является следствием теорем 2.6.7, 2.6.11. □

Следствие 4.6.8. *Пусть $\overline{X}_{[1]}$, $\overline{X}'_{[1]}$ - базисы представления $\overline{f}$. Пусть $\overline{r}_{[1]}$ - автоморфизм представления $\overline{f}$ такой, что $\overline{X}'_{[1]} = \overline{r}_{[1]} \circ \overline{X}_{[1]}$ . Автоморфизм $\overline{r}_{[1]}$ определён однозначно.* □

---

[4.12] $X''_j = X'_j$, $j \neq i$, $X''_i = X'_i \setminus \{x'_i\}$

[4.13] $X'''_j = X_j$, $j \neq i$, $X'''_i = X_i \setminus \{x_i\}$

[4.14] Это утверждение похоже на теорему [1]-1, с. 104.

Замечание 4.6.9. Теорема 4.6.7, так же как и теорема 4.5.14, является теоремой о продолжении отображения. Одако здесь $\overline{X}_{[1]}$ - не произвольный кортеж множеств образующих башни представлений, а базис. Согласно замечаниям 4.6.3, 2.7.3, мы не можем определить координаты любого элемента базиса через остальные элементы этого же базиса. Поэтому отпадает необходимость в согласованности отображения базиса с представлением. □

Теорема 4.6.10. *Набор координат $\overline{W}_{[1]}(\overline{f},\overline{X}_{[1]},\overline{X}_{[1]})$ соответствует тождественному преобразованию*

$$\overline{W}_{[1]}[\overline{f},\overline{X}_{[1]},\overline{E}_{[1]}]=\overline{W}_{[1]}(\overline{f},\overline{X}_{[1]},\overline{X}_{[1]})$$

Доказательство. Утверждение теоремы следует из равенства

$$\overline{a}_{[1]}=\overline{W}_{[1]}(\overline{f},\overline{X}_{[1]},\overline{a}_{[1]})\circ\overline{X}_{[1]}=\overline{W}_{[1]}(\overline{f},\overline{X}_{[1]},\overline{a}_{[1]})\circ\overline{W}_{[1]}(\overline{f},\overline{X}_{[1]},\overline{X}_{[1]})\circ\overline{X}_{[1]}$$

□

Теорема 4.6.11. *Пусть $\overline{W}_{[1]}(\overline{f},\overline{X}_{[1]},\overline{r}_{[1]}\circ\overline{X}_{[1]})$ - множество координат автоморфизма $\overline{r}_{[1]}$. Определено множество координат $\overline{W}_{[1]}(\overline{f},\overline{r}_{[1]}\circ\overline{X}_{[1]},\overline{X}_{[1]})$, соответствующее автоморфизму $\overline{r}_{[1]}^{-1}$. Множество координат $\overline{W}_{[1]}(\overline{f},\overline{r}_{[1]}\circ\overline{X}_{[1]},\overline{X}_{[1]})$ удовлетворяет равенству*[4.15]

$$\overline{W}_{[1]}(\overline{f},\overline{X}_{[1]},\overline{r}_{[1]}\circ\overline{X}_{[1]})\circ\overline{W}_{[1]}(\overline{f},\overline{r}_{[1]}\circ\overline{X}_{[1]},\overline{X}_{[1]})=\overline{W}_{[1]}(\overline{f},\overline{X}_{[1]},\overline{r}_{[1]}\circ\overline{X}_{[1]}) \tag{4.6.3}$$

$$\overline{W}_{[1]}[\overline{f},\overline{X}_{[1]},\overline{r}_{[1]}^{-1}]=\overline{W}_{[1]}[\overline{f},\overline{X}_{[1]},,R]^{-1}=\overline{W}_{[1]}(\overline{f},\overline{r}_{[1]}\circ\overline{X}_{[1]},\overline{X}_{[1]})$$

Доказательство. Поскольку $\overline{r}_{[1]}$ - автоморфизм башни представлений $\overline{f}$, то, согласно теореме 4.6.6, множество $\overline{r}_{[1]}\circ\overline{X}_{[1]}$ - базис башни представлений $\overline{f}$. Следовательно, определено множество координат $\overline{W}_{[1]}(\overline{f},\overline{r}_{[1]}\circ\overline{X}_{[1]},\overline{X}_{[1]})$. Равенство (4.6.3) следует из цепочки равенств

$$\begin{aligned}&\overline{W}_{[1]}(\overline{f},\overline{X}_{[1]},\overline{r}_{[1]}\circ\overline{X}_{[1]})\circ\overline{W}_{[1]}(\overline{f},\overline{r}_{[1]}\circ\overline{X}_{[1]},\overline{X}_{[1]})\\=&\overline{W}_{[1]}(\overline{f},\overline{X}_{[1]},\overline{r}_{[1]}\circ\overline{X}_{[1]})\circ\overline{W}_{[1]}(\overline{f},\overline{X}_{[1]},\overline{r}_{[1]}^{-1}\circ\overline{X}_{[1]})\\=&\overline{W}_{[1]}(\overline{f},\overline{X}_{[1]},\overline{r}_{[1]}^{-1}\circ\overline{r}_{[1]}\circ\overline{X}_{[1]})=\overline{W}_{[1]}(\overline{f},\overline{X}_{[1]},\overline{X}_{[1]})\end{aligned}$$

□

Теорема 4.6.12. *Группа автоморфизмов $G(\overline{f})$ башни эффективных представлений $\overline{f}$ порождает эффективное представление в башне представлений $\overline{f}$.*

Доказательство. Из следствия 4.6.8 следует, что если автоморфизм $\overline{r}_{[1]}$ отображает базис $\overline{X}_{[1]}$ в базис $\overline{X}'_{[1]}$, то множество координат $\overline{W}_{[1]}(\overline{f},\overline{X}_{[1]},\overline{X}'_{[1]})$ однозначно определяет автоморфизм $\overline{r}_{[1]}$. Из теоремы 4.5.18 следует, что множество координат $\overline{W}_{[1]}(\overline{f},\overline{X}_{[1]},\overline{X}'_{[1]})$ определяет правило отображения координат относительно базиса $\overline{X}_{[1]}$ при автоморфизме башни представлений $\overline{f}$. Из равенства (4.5.30) следует, что автоморфизм $\overline{r}_{[1]}$ действует справа на элементы $\Omega_k$-алгебры $A_k$, $k=2$, ..., $n$. Из равенства (4.5.28) следует, что представление группы является $*T$-представлением. Согласно теореме 4.6.10 набор координат $\overline{W}_{[1]}(\overline{f},\overline{X}_{[1]},\overline{X}_{[1]})$ соответствует тождественному преобразованию. Из теоремы

[4.15] Смотри также замечание 2.7.12.

4.6.11 следует, что набор координат $\overline{W}_{[1]}(\overline{f}, \overline{r}_{[1]} \circ \overline{X}_{[1]}, \overline{X}_{[1]})$ соответствует преобразованию, обратному преобразованию $\overline{W}_{[1]}(\overline{f}, \overline{X}_{[1]}, \overline{r}_{[1]} \circ \overline{X}_{[1]})$. □

## 4.7. Примеры базиса башни представлений

**4.7.1. Аффинное пространство.** Рассмотрим аффинное пространство $\overset{\circ}{A}$ (раздел 4.4.1).

Абелевая группа $\overline{A}$ действует однотранзитивно на множестве $\overset{\circ}{A}$. Из построений в подразделе 2.8.2 следует, что базис множества $\overset{\circ}{A}$ относительно представления абелевой группы $\overline{A}$ состоит из одной точки. Эту точку обычно обозначают буквой $O$ и называют **началом системы координат аффинного пространства**. Следовательно, произвольную точку $A \in \overset{\circ}{A}$ можно представить с помощью вектора $\overrightarrow{OA} \in \overline{A}$

Пусть $\overline{\overline{e}}$ - базис векторного пространства $\overline{A}$. Тогда вектор $\overrightarrow{OA}$ имеет вид

$$\overrightarrow{OA} = a^i \overline{e}_i$$

Множество $(a_i, i \in I)$ называется **координатами точки** $A$ **аффинного пространства** $\overset{\circ}{A}$ **относительно базиса** $(O, \overline{\overline{e}})$.

**4.7.2. Координаты вектора в $R$-модуле.** Пусть левый модуль $\overline{M}$ над кольцом $R$ ([1], с. 93) имеет конечный базис

$$\overline{\overline{e}} = (\overline{e}_1, ..., \overline{e}_n)$$

Согласно определению 2.7.1 любой вектор $\overline{v} \in \overline{M}$ можно представить в виде

$$\overline{v} = v^i \overline{e}_i \tag{4.7.1}$$

Согласно теореме 2.7.2 $\overline{e}_i \in \overline{\overline{e}}$ не может быть представлен подобным образом через остальные векторы базиса. Однако это не означает, что мы не можем записать равенство

$$c^i \overline{e}_i = 0 \tag{4.7.2}$$

где не все $c^i$ равны 0. Для того, чтобы равенство (4.7.2) было возможным, необходимо потребовать, чтобы числа $c^1$, ..., $c^n$ были необратимы в кольце $R$. Из равенств (4.7.1), (4.7.2) следует, что вектор $\overline{v}$ кроме разложения (4.7.1) относительно базиса $\overline{\overline{e}}$, имеет также разложение

$$\overline{v} = (v^i + c^i) \overline{e}_i \tag{4.7.3}$$

Нетрудно заметить, что множество кортежей $(c^1, ..., c^n)$ порождает левый $R$-модуль. Следовательно, множество координат в $R$-модуле $\overline{M}$ является аффинным пространством, а координаты вектора $\overline{v}$ являются аффинной плоскостью в этом пространстве.

Не исключено, что множество координат произвольного представления универсальной алгебры также можно рассматривать как некоторую башню представлений.

**4.7.3. Базис $R$-алгебры.** Обычно, когда мы рассматриваем $R$-алгебру $A$, мы выбираем базис $\overline{\overline{e}}$ соответствующего $R$-модуля $A$. Этот выбор удобен, так как если $R$ - коммутативная алгебра с делением, то разложение вектора однозначно относительно базиса $R$-векторного пространства. Это, в частности, позволяет описать операции произведения, указав структурные константы алгебры относительно заданого базиса.

В общем случае, базис $R$-модуля $A$ может оказаться множеством образующих. Например, если в векторном пространстве $H$, в котором задана алгебра кватернионов над полем действительных чисел, рассмотреть базис

$$e_{\mathbf{0}} = 1 \quad e_{\mathbf{1}} = i \quad e_{\mathbf{2}} = j \quad e_{\mathbf{3}} = k \tag{4.7.4}$$

то в алгебре $H$ верно равенство

$$\begin{gathered} e_{\mathbf{0}} = -e_{\mathbf{1}}e_{\mathbf{1}} = -e_{\mathbf{2}}e_{\mathbf{2}} \\ e_{\mathbf{3}} = e_{\mathbf{1}}e_{\mathbf{2}} \end{gathered} \tag{4.7.5}$$

Следовательно, множество $(e_{\mathbf{1}}, e_{\mathbf{2}})$ является базисом алгебры $H$. Следствием равенства (4.7.5) является неоднозначность представления кватерниона относительно заданного базиса. А именно, кватернион $a \in H$ можно записать в виде

$$a = (a^{\mathbf{0}} - a^{\mathbf{4}})e_{\mathbf{1}}e_{\mathbf{1}} + a^{\mathbf{4}}e_{\mathbf{2}}e_{\mathbf{2}} + a^{\mathbf{1}}e_{\mathbf{1}} + a^{\mathbf{2}}e_{\mathbf{2}} + a^{\mathbf{3}}e_{\mathbf{1}}e_{\mathbf{2}}$$

где $a^{\mathbf{4}}$ - произвольно.

**4.7.4. Модуль над алгеброй.** Пусть $A_3$ - свободный модуль над $A_1$-алгеброй $A_2$ (раздел 4.4.4). Пусть $\overline{\overline{e}}_{32}$ - базис алгебры $A_3$ над алгеброй $A_2$. Вектор $a_3 \in A_3$ имеет представление

$$a_3 = a_3^j\ \overline{e}_{32\cdot j} = \begin{pmatrix} a_3^1 & ... & a_3^n \end{pmatrix} \begin{pmatrix} \overline{e}_{32\cdot 1} \\ ... \\ \overline{e}_{32\cdot n} \end{pmatrix} \tag{4.7.6}$$

Пусть $\overline{\overline{e}}_{21}$ - базис алгебры $A_2$ над полем $A_1$. Так как $a_3^j \in A_2$, то мы можем записать их координаты относительно базиса $\overline{\overline{e}}_{21}$

$$a_3^j = a_3^{ji}\ \overline{e}_{21\cdot i} = \begin{pmatrix} a_3^{j1} & ... & a_3^{jm} \end{pmatrix} \begin{pmatrix} \overline{e}_{21\cdot n} \\ ... \\ \overline{e}_{21\cdot m} \end{pmatrix} \tag{4.7.7}$$

Из равенств (4.7.6), (4.7.7) следует

(4.7.8)

$$a_3 = a_3^{ji}\ \overline{e}_{21\cdot i}\ \overline{e}_{32\cdot j} = \begin{pmatrix} \begin{pmatrix} a_3^{11} & ... & a_3^{1m} \end{pmatrix} \begin{pmatrix} \overline{e}_{21\cdot n} \\ ... \\ \overline{e}_{21\cdot m} \end{pmatrix} & ... & \begin{pmatrix} a_3^{n1} & ... & a_3^{nm} \end{pmatrix} \begin{pmatrix} \overline{e}_{21\cdot n} \\ ... \\ \overline{e}_{21\cdot m} \end{pmatrix} \end{pmatrix} \begin{pmatrix} \overline{e}_{32\cdot 1} \\ ... \\ \overline{e}_{32\cdot n} \end{pmatrix}$$

Равенство (4.7.8) показывает структуру координат в векторном пространстве $A_3$ над полем $A_1$. Нетрудно убедиться, что векторы

$$\overline{e}_{31\cdot ij}=\overline{e}_{21\cdot i}\ \overline{e}_{32\cdot j}$$

линейно независимы над полем $A_1$. Следовательно, мы построили базис $\overline{\overline{e}}_{31}$ векторного пространстве $A_3$ над полем $A_1$. Следовательно, мы можем переписать равенство (4.7.8) в виде

$$a_3=a_3^{ji}\ \overline{e}_{31\cdot ij}=\begin{pmatrix}a_3^{11} & ... & a_3^{1m} & ... & a_3^{n1} & ... & a_3^{nm}\end{pmatrix}\begin{pmatrix}\overline{e}_{31\cdot 11}\\ ...\\ \overline{e}_{31\cdot 1m}\\ ...\\ \overline{e}_{31\cdot n1}\\ ...\\ \overline{e}_{31\cdot nm}\end{pmatrix}\tag{4.7.9}$$

Нетрудно убедиться, что вектор $\overline{e}_{31\cdot ij}$ можно отожлествить с тензорным произведением $\overline{e}_{21\cdot i}\otimes\overline{e}_{32\cdot j}$.

## 4.8. Представление в категории

Определение 4.8.1. Пусть для любых объектов $B$ и $C$ категории $\mathcal{B}$ на множестве морфизмов $Mor(B,C)$ определена структура $\Omega$-алгебры. Множество гомоморфизмов $\Omega$-алгебры

$$f_{B,C}:A\to Mor(B,C)$$

называется **представлением $\Omega$-алгебры $A$ в категории $\mathcal{B}$**. □

Если предположить, что множество $Mor(B,C)$ определенно только, когда $B=C$, то мы получим определение представления. Различие в определениях состоит в том, что мы не ограничиваем себя преобразованиями множества $B$, а рассматриваем $\Omega$-алгебру отображений из множества $B$ в множество $C$. На первый взгляд нет принципиальных различий между рассматриваемыми теориями. Однако нетрудно заметить, что представление векторного пространства в категории расслоений приводит к теории связностей.

Глава 5

# Список литературы

[1] Серж Ленг, Алгебра, М. Мир, 1968

[2] S. Burris, H.P. Sankappanavar, A Course in Universal Algebra, Springer-Verlag (March, 1982),
eprint http://www.math.uwaterloo.ca/ snburris/htdocs/ualg.html
(The Millennium Edition)

[3] П. К. Рашевский, Риманова геометрия и тензорный анализ,
М., Наука, 1967

[4] А. Г. Курош, Общая алгебра, (лекции 1969 - 70 учебного года), М., МГУ, 1970

[5] Aleksandr Vasilevich Mikhalev; Günter Pilz; The concise handbook of algebra;
Kluwer Academic Publishers, 2002

[6] И. Р. Шафаревич, Основные понятия алгебры,
В сборнике "Современные проблемы математики. Фундаментальные направления (Итоги науки и техники. ВИНИТИ АН СССР", М., 1986,
**11** 3 - 289 М., Мир, 1980

[7] Александр Клейн, Лекции по линейной алгебре над телом,
eprint arXiv:math.GM/0701238 (2010)

[8] Александр Клейн, Расслоенная универсальная алгебра,
eprint arXiv:math.DG/0702561 (2007)

[9] Александр Клейн, Введение в геометрию над телом,
eprint arXiv:0906.0135 (2010)

[10] Александр Клейн, Линейные отображения свободной алгебры,
eprint arXiv:1003.1544 (2010)

[11] Александр Клейн.
Линейная алгебра над телом: Система линейных уравнений.
CreateSpace Independent Publishing Platform, 2014; ISBN-13: 978-1502982476

[12] John C. Baez, The Octonions,
eprint arXiv:math.RA/0105155 (2002)

[13] П. Кон, Универсальная алгебра, М., Мир, 1968

[14] Richard D. Schafer, An Introduction to Nonassociative Algebras, Dover Publications, Inc., New York, 1995

Глава 6

# Предметный указатель

Глава 7

# Специальные символы и обозначения